\documentclass[11pt]{smfbook}
\usepackage{amsmath,amssymb,smfthm,stmaryrd,latexsym}
\usepackage[frenchb]{babel}
\usepackage[latin1]{inputenc}
\input diagrams.tex

\def\Y{\mathfrak{Y}}
\def\R{\mathbb{R}}
\def\S{\mathfrak{S}}
\def\D{\boldsymbol{D}}

\def\C{\mathbb{C}}
\def\CC{\mathfrak{C}}
\def\Cp{\mathbb{C}_p}
\def\Fq{\mathbb{F}_q}
\def\Fqb{\overline{\mathbb{F}}_q}
\def\L{\mathcal{L}}
\def\La{\Lambda}
\def\A{\text{{\bf A}}} 
\def\Aff{\mathbb{A}}
\def\O{\mathcal{O}}
\def\F{\mathcal{F}}
\def\I{\mathcal{I}}

\def\*{^\times }
\def\dpt{\displaystyle}
\def\GL{\text{\bf GL}}
\def\G{\mathcal{G}}
 
\def\Gm{\mathbb{G}_m}         
\def\l{\lambda}
\def\g{\gamma}
\def\BB{\mathbb{B}}

\def\a{\alpha}
\def\b{\beta}
\def\d{\delta}
\def\s{\sigma}
\def\ph{\varphi}

\def\e{\epsilon}
\def\ssi{\Leftrightarrow}
\def\lssi{\Longleftrightarrow}
\def\impl{\Rightarrow}
\def\limpl{\Longrightarrow}
\def\drt{\rightarrow}
\def\ldrt{\longrightarrow}
\def\Q{\mathbb{Q}}
\def\Qb{\overline{\mathbb{Q}}}
\def\Ql{\mathbb{Q}_\ell}
\def\Qlb{\overline{\mathbb{Q}_\ell}}
\def\Qp{\mathbb{Q}_p}
\def\Qpb{\overline{\mathbb{Q}_p}}
\def\Zp{\mathbb{Z}_p}
\def\Z{\mathbb{Z}}
\def\N{\mathbb{N}}
\def\mm{\mathfrak{m}}
\def\Zl{\mathbb{Z}_\ell}
\def\Hom{\text{Hom}}
\def\Homf{\underline{\text{Hom}}}

\def\Rpsi{\text{R}\Psi_{\bar{\eta}}}
\def\rpsi{\text{R}\Psi}
\def\Rth{\text{R}\Theta_{\bar{\eta}}}
\def\rth{\text{R}\Theta}
\def\Coh{\text{Coh}}
\def\Cohf{\underline{\text{Coh}}}

\def\Gal{\boldsymbol{\text{\bf Gal}}}
\def\={\! = \!}

\newcommand{\citeg}[1]{{\bf\cite{#1}}}

\def\spec{\boldsymbol{Spec}}
\def\spf{\boldsymbol{Spf}}

\def\limp{\underset{\longleftarrow}{\text{ lim }}\;}
\def\limi{\underset{\longrightarrow}{\text{ lim }}\;}

\def\Ab{\boldsymbol{Ab}}

\def\iso{\xrightarrow{\;\sim\;}}

\def\End{\text{End}}
\def\Aut{\text{Aut}}
\def\Gl{\text{GL}}
\def\im{\text{im}}
\def\Sh{\text{Sh}}
\def\m1{^{-1}}
\def\xrig{\xrightarrow}
\def\pK{K^p}

\def\Res{\text{Res}}
\def\Eb{\breve{E}}

\def\M{\breve{\mathcal{M}}}
\def\Mb{\overline{\mathcal{M}}}
\def\XX{\mathbb{X}}
\def\X{\mathfrak{X}}
\def\Sb{\overline{S}}
\def\kb{\overline{k}}
\def\St{\widetilde{S}(\phi)}
\def\GG{\Gamma}
\def\Ext{\text{Ext}}
\def\bc{\backslash}
\def\AA{\mathcal{A}}
\def\Spa{\text{Spa}}
\def\Lpf{\Lambda_\bullet -\mathcal{F}\!\text{sc}}
\def\LpGf{\Lambda_\bullet-G -\mathcal{F}\!\text{sc}}
\def\Lf{\Lambda -\mathcal{F}\!\text{sc}}
\def\LGf{\Lambda -G-\mathcal{F}\!\text{sc}}
\def\Lmf{\Lambda/\mm^n -\mathcal{F}\!\text{sc}}
\def\Xet{_{/X_{\text{{\'e}t}}}}
\def\et{\text{{\'e}t}}
\def\U{\mathbb{U}}
\def\Fix{\text{Fix}(KgK)}

\def\Irr{\text{Irr}}
\def\Groth{\text{Groth}}
\def\L{\,^{L}}
\def\Lie{\text{Lie}}
\def\LL{\mathcal{L}}

\newarrow{Equalto} =====

\title{Correspondances de Langlands locales dans la 
cohomologie des espaces de Rapoport Zink} 
\author{Laurent Fargues}
\email{lfargues@math.jussieu.fr}

\begin{document}

\frontmatter
\renewcommand{\labelitemi}{$\bullet$}
\newtheorem{Fait}{Fait}[section]

\maketitle

\tableofcontents

\chapter*{Introduction}

\section*{Motivation générale}

  Dans l'article \cite{Langlands1} Langlands a prédit l'existence d'un
  lien entre les motifs de Grothendieck définis sur des corps de
  nombres et certaines représentations automorphes des groupes
  linéaires. Pour certaines formes modulaires un cas particulier bien
  connu de cette correspondance est réalisé dans la cohomologie des
  courbes modulaires. Plus généralement, soit $G$ un groupe réductif
  défini sur $\Q$ possédant une donnée de Shimura $(G,X)$. Langlands a
  donné une description conjecturale du motif découpé par une
  représentation automorphe $\Pi$ de $G$ dans la variété de Shimura
  $\Sh(G,X)$ en fonction du L-paramètre (conjectural) de $\Pi$
  (\cite{Langlands1}). 

 Dans cette thèse, nous démontrons des analogues locaux en des places
 non-archimédiennes de cette conjecture. Qu'entendons nous par
 analogue local ? En la place archimédienne de $\Q$ 
 l'objet local équivalent naturel de la variété de Shimura $\Sh (G,X)$
 est l'espace symétrique hermitien $X$ uniformisant la variété
 analytique $\Sh (G,X) (\C)$, tandis que les objets associés aux 
 représentations automorphes de $G$ sont les représentations
 irréductibles du groupe de Lie $G(\R)$ au sens d'Harisch-Chandra. Les travaux de
 Schmid étudient ainsi la contribution des séries discrète de $G(\R)$
 dans la cohomologie $L^2$ de $X$ (l'analogue local de l'identification de
 l'action du groupe de Galois motivique consisterait dans ce cas à
 déterminer une structure de Hodge sur ces espaces de cohomologie).
Les analogues non-archimédiens en un premier $p$ de $\Q$ de l'espace
 symétrique $X$  sont les espaces rigides $p$-adiques de Rapoport-Zink
 (\cite{RZ}) uniformisant certains ouverts rigides des variétés de
 Shimura $\Sh (G,X)$ (ils ne sont pas définis pour tous les couples
 $(G,X)$). 

\section*{Espaces de Rapoport-Zink}

Pour comprendre la  définition des espaces de Rapoport-Zink
rappelons que 
 la plupart des variétés de Shimura devraient se décrire
 comme des espaces de modules de motifs. L'espace $X$ s'interprète ainsi
 comme la contrepartie locale archimédienne de ce point de vue puisque
 c'est un espace de modules de structures de Hodge. Supposons que
 les variétés $\Sh(G,X)$ s'interprètent comme espaces de modules de variétés
 abéliennes. Les espaces de Rapoport-Zink sont alors la contrepartie
 $p$-adique de ce point de vue : ce sont des espaces de modules de
 groupes p-divisibles munis de structures de niveaux. 

 Contrairement à $X$ les espaces de Rapoport-Zink ne sont pas simplement
 connexes; il s'agit d'une tour d'espaces rigides $(\M_K)_K$ indexée
 par les sous groupes compacts ouverts de $G(\Qp)$ telle que si $K'$ est
 un  sous groupe compact ouvert de $K$, le morphisme de transition
 $\M_{K'}\drt \M_K$ est un revêtement étale. Cette tour est munie
 d'une action de $G(\Qp)$. Un deuxième groupe algébrique $p$-adique
 entre en jeux,  une forme intérieure d'un sous groupe
 de Levi de la forme intérieure quasi-déployée de $G_{\Qp}$. On le note
 $J_b$. Celui-ci agit sur chaque espace $\M_K$ de manière compatible
 aux morphismes de transition. Le
 troisième groupe entrant en jeu est le groupe de Weil $W_E$ d'une
 extension de degré fini de $\Qp$. 

\section*{Esquisse des résultats}

Considérons la cohomologie étale
 $\ell$-adique à support compact au sens de Berkovich des $\M_K$ pour
 un $\ell \neq p$, 
$$
H^\bullet_c (\M_K,\Qlb)
$$
Lorsque $K$ varie celle ci est munie d'une action lisse de
$G(\Qp)\times J_b\times W_E$.

 Les résultats principaux de cette thèse
concernent des espaces de Rapoport-Zink pour lesquels le groupe $J_b$ est une
forme intérieure de $G_{\Qp}$. Ce sont ceux qui uniformisent les
strates supersingulières des variétés de Shimura de type P.E.L. non
ramifiées définies dans \cite{Ko4}. Le premier résultat concerne le
cas où $G(\Qp)=\GL_n (F)$ 
pour une extension $F|\Qp$ non ramifiée, et $J_b=D^\times$ où $D$ est
une algèbre à division sur $F$. Soit $\text{JL}$ la correspondance de
Jacquet-Langlands qui associe à une représentation irréductible de
$D^\times $ une représentation de la série discrète  de
$\GL_n (F)$. Soit $\rho$ une représentation  de $J_b$ telle
que $\text{JL} (\rho)$ soit supercuspidale. Soit $\pi$ une représentation
supercuspidale de $G(\Qp)$. D'après la correspondance de Langlands
locale (\cite{Har4}), à $\pi$ est associé un certain L-paramètre
$\psi$. Nous démontrons que  $\rho\otimes\pi$ intervient
virtuellement dans la somme alternée de la cohomologie des espaces de
Rapoport-Zink ssi $\check{\pi}=\text{JL}(\rho)$. Nous démontrons  qu'alors la
représentation galoisienne découpée par $\rho\otimes \pi$ s'exprime 
simplement en fonction du L-paramètre $\psi$. 

Nous démontrons le même type de résultat lorsque $G(\Qp)=J_b= \GL_n
(F)$, la correspondance de Jacquet-Langlands étant alors l'identité. 

Nous démontrons également le premier cas de loi de réciprocité
non abélienne locale construite géométriquement et 
associée à des groupes autres que des formes
intérieures des groupes linéaires. Plus précisément, il s'agit du cas
où $J_b=G(\Qp)=GU(3)$ est le groupe de similitudes unitaires non ramifié en trois
variables. Contrairement au cas précédents, les L-paquets de
représentations de $GU(3)$ ne sont pas tous de cardinal $1$. Le
résultat est alors du même type que précédemment après sommation sur
un tel L-paquet.

\section*{Conjectures de Kottwitz}

Se basant sur la description conjecturale de Langlands du motif, et
en particulier de la représentation galoisienne dans la cohomologie
$\ell$-adique, associé à une représentation automorphe de $G$ dans la
variété de Shimura $\Sh (G,X)$, ainsi que sur la formule conjecturale 
d'Arthur pour la multiplicité d'une représentation automorphe dans un
A-paquet discret de $G$, Kottwitz a formulé une conjecture concernant
la contribution d'une représentation ``discrète'' de $G(\Qp)\times
J_b$ dans les espaces de cohomologie $H^\bullet_c (\M_K,\Qlb)$ ainsi que de
la représentation galoisienne associée (\cite{Rapoport2} et
\cite{ECM}). Les résultats de cette thèse sont des cas particuliers de
ces conjectures. 

\section*{Énoncés précis}

Les résultats principaux de cette thèse sont les suivants :

\begin{enonce*}{Théorème}[théorèmes III.\ref{theo_princ1} et III.\ref{theo_princ2}]
  Soit $(V,F,b,\mu)$ une donnée locale  de Rapoport-Zink de type
  E.L. non ramifiée simple basique (I.\ref{EL_RZ}). Soit $E\subset \Qpb$ le corps
  réflex associé.  Soient $\M_K=\M_K (b,\mu)$ les espaces de
  Rapoport-Zink associés (I.\ref{esp_RZ_ass}).
\begin{enumerate}
\item
Supposons le groupe réductif $J_b$ anisotrope modulo
  son centre, c'est à dire
 $J_b=D^\times$ pour une algèbre à division
  $D$ sur $F$ d'invariant calculé en fonction de $\mu$. 
Notons $\text{JL}$ la correspondance de Jacquet-Langlands. 
Soit $\pi$ une
  représentation irréductible de $J_b$
telle que $\text{JL}(\pi)$ soit supercuspidale. L'égalité suivante est
vérifiée dans le
  groupe de Grothendieck $\Groth ( G(\Qp)\times W_E)$ 
\begin{eqnarray*}
 \sum_i (-1)^i 
[\underset{K}{\limi} \Hom_{J_b} ( H^i_c ( \M_K\hat{\otimes}\Cp,\Qlb),
\pi)]_{cusp} 
 = 
[\text{\text{JL}} (\pi)] \otimes [r_\mu\circ \tilde{\s}_\ell (\text{\text{JL}}
(\pi))_{|E}]\; |.|^{\frac{-\sum_\tau {p_\tau q_\tau}}{2}} 
\end{eqnarray*}
où
$\tilde{\s}_\ell$ la correspondance de Langlands locale ``absolue''
$\ell$-adique pour le
groupe linéaire (\ref{Lang_Rev}) et $r_\mu$ la représentation du
L-groupe associé (\ref{le_r_mu}) sur $E$. 
\item Supposons le groupe réductif $J_b$ égal à $G=\GL_n$. Soit $\pi$ une
  représentation irréductible supercuspidale de $J_b$. Il y a une
  égalité dans $\Groth ( G(\Qp)\times W_E)$ 
$$
\sum_i (-1)^i 
[\underset{K}{\limi} \Hom_{J_b} ( H^i_c ( \M_K\hat{\otimes}\Cp,\Qlb),
\pi)]_{cusp} = [\pi]\otimes  [r_\mu\circ \tilde{\s}_\ell (
\pi)_{|E}]\; |.|^{-\sum_\tau \frac{p_\tau q_\tau}{2}}
$$
\end{enumerate}
\end{enonce*}

Dans cet énoncé, la donnée de Rapoport-Zink locale est l'analogue
non-archimédien local d'une donnée de Shimura. 

\begin{enonce*}{Théorème}[théorème III.\ref{theo_princ3}]
  Soit $(F,*,V,<.,.>,b,\mu)$ une donnée locale de type P.E.L. non
  ramifiée simple basique (I.\ref{PEL_RZ}). Supposons $p\neq 2$, $[F:\Qp]/2$
  impair et $\dim_F (V)=3$. Soit $E\subset \Qpb$ le corps réflex
  associé.
Soient $\M_K=\M_k(b,\mu)$ les espaces de Rapoport-Zink associés
  (I.\ref{esp_RZ_ass_PEL}). 
\begin{enumerate}
\item 
 Soit
  $\psi:W_{\Qp}\ldrt \L G$ une classe de conjugaison de L-paramètres
  vérifiant : $\im (\psi )$ n'est pas contenue dans $\L B$ pour un
 sous groupe de  Borel $B$ de $G$. Soit $\Pi (\psi )$ le L-paquet
  supercuspidal de représentations de $G(\Qp)$ associé
  (\ref{paquets_sup}). Rappelons que $J_b= G(\Qp)=GU(3)$. 
  L'égalité suivante est vérifiée dans 
 $\Groth (G(\Qp)\times W_{E})$ : 
$$ 
\sum_{\pi\in \Pi (\psi)}\left [\underset{K}{\limi} \Hom_{J_b} \left ( H^\bullet_c
( \M_K\hat{\otimes}\Cp,\Qlb),\pi\right ) \right ]_{cusp} = \sum_{\pi\in \Pi
  (\psi)} [\pi]\otimes [ r_\mu \circ \psi_{|W_E}]\;
|.|^{-\frac{\sum_\tau p_\tau q_\tau}{4}}
$$ 
En particulier, si $\Pi (\psi)$ est stable la conjecture de Kottwitz
(\cite{Rapoport2},\cite{ECM}). 
est vérifiée. 
\item
Soit $\xi$ un caractère du groupe endoscopique $H$ de $U(3)$
(appendice \ref{Resultats_Rog}) et $St_H (\xi)$ la représentation de Steinberg
associée. Soit $\Pi ( St_H (\xi)) = \{ \pi^2 (\xi), \pi^s (\xi) \}$ le
L-paquet transfert endoscopique à $U(3)$ (\ref{Class_sup}). Soit
$\Pi=\{\pi^2,\pi^s\}$ un L-paquet de représentations de $G(\Qp)$ tel
  que $\pi^2_{|U(3)}=\pi^2 (\xi),\pi^s_{|U(3)}=\pi^s (\xi)$. A $\Pi$
  est associé un L-paramètre $\psi:W_{\Qp}\times \text{SL}_2 (\C)
  \ldrt \L G$. L'égalité suivante est vérifiée 
dans $\Groth (G(\Qp)\times W_{E})$
  : 
  \begin{eqnarray*}
  &\,&  \left [\underset{K}{\limi} \Hom_{J_b} \left ( H^\bullet_c
    (\M_{K}\hat{\otimes} \Cp ,\Qlb),\pi^s\right )\right ]_{cusp} + \left[
    \underset{K}{\limi} \text{Ext}^{\, \bullet}_{J_b-\text{lisse}} (
    H^\bullet_c (\M_{K}\hat{\otimes} \Cp ,\Qlb),\pi^2)\right ]_{cusp}
    \\
 &=& [\pi^s]\otimes [r_\mu\circ \psi_{|W_E}] \;
    |.|^{-\frac{\sum_\tau p_\tau q_\tau }{4}}
  \end{eqnarray*}
\item 
Soit $\chi$ un caractère de $G(\Qp)$ et $St_G(\chi)$ la représentation
de Steinberg associée. L'égalité suivante est vérifiée 
 dans $\Groth ( G(\Q)\times
W_{E})$ : 
$$
\left [ \text{Ext}^\bullet_{J_b-\text{lisse}} \left (  H^\bullet_c
( \M_K\hat{\otimes}\Cp,\Qlb), St_G (\chi)\right )\right ]_{cusp} = 0 
$$
\end{enumerate}
\end{enonce*}
\vspace{3mm}

\section*{Quelques rappels sur la méthode de Harris,Taylor et Boyer}
\subsection*{Rapide historique}

Lubin et Tate avaient déjà considéré (mais formulé différemment bien sûr) le cas
particulier $G=\GL_1$ pour 
lequel les espaces $\M$ sont de dimension $0$, et donné ainsi une
construction géométrique de la loi de réciprocité locale
d'Artin. Carayol a ensuite formulé dans \cite{Carayol1} des
conjectures concernant les espaces de déformation de Lubin-Tate de
dimension supérieure. Ces conjectures ont été démontrées dans le cas
des corps locaux de caractéristique $p$ par Boyer
(\cite{Boyer}). Harris et Taylor on alors démontré la conjecture de
Langlands locale pour les groupes linéaires sur des extensions de
degré fini de $\Qp$ en utilisant cette approche (\cite{Har4}), puisque
dans leur cas la représentation $r_\mu$ est la représentation standard
du groupe linéaire. 

\subsection*{Quelques points de \cite{Har4}}

Pour comprendre la tactique utilisée pour démontrer nos théorèmes il
nous faut 
rappeler brièvement certains points de  
(\cite{Har4}). Harris et Taylor 
considèrent des espaces de déformation de $\O_F$-modules
divisibles de dimension $1$ et de hauteur $g$, munis d'une structure
 de Drinfeld de niveau $m$, et notés $\spf (R_{F,g,m})$ (ce sont des
 schémas formels). Avec les notations du théorème
 III.\ref{theo_princ1} énoncé précédemment, cela correspond à un choix
 particulier du cocaractère $\mu$. Les espaces $\M_K$ associés sont
 alors une union disjointe infinie des fibres génériques au sens des
 espaces analytiques des $\spf (R_{F,g,m})$. Leur résultat est énoncé en
 termes de cycles évanescents $\ell$-adiques au sens de Berkovich,
 cependant les résultats de la section \ref{tors_ell} de cette thèse
 permettent de faire le lien avec la cohomologie à support compact des
 espaces analytiques $\M_K$. Avec nos notations, ils démontrent que si
 $\rho$ est une représentation de $J_b$ telle que $\text{JL} (\rho)$
 est supercuspidale alors, 
$$
[\Hom_{J_b} ( H^\bullet_c (
\M_K\hat{\otimes}\Cp,\Qlb),\rho)]=[\text{JL} (\rho)]\otimes [r_\ell
(\text{JL} (\rho ))]
$$
où $r_\ell
(\text{JL} (\rho ))$ est une représentation $\ell$-adique de
$W_F$. Ils montrent alors que la correspondance définie sur les
représentations supercuspidales 
$$
\pi\longmapsto r_\ell (\pi)
$$
est une correspondance de Langlands locale tordue par un caractère non
ramifié. Leur démonstration utilise certaines variétés de Shimura, dites
``simples'', possédant de bon modèles entiers en tout niveau. Elles
paramètrent des variétés abéliennes $A$ munies de structures
additionnelles (polarisation, action d'une algèbre) et de structures de
niveau définies en $p$ en utilisant la théorie des structures de
niveau de Drinfeld.  

L'un des points clef est que 
celles-ci sont 
stratifiées par la classe d'isogénie du groupe p-divisible
universel $A[p^\infty]$. La strate minimale, ou strate supersingulière, est de
dimension $0$, et le complété formel de la variété de Shimura le long
de cette strate est uniformisé par une union disjointe de $\spf
(R_{F,g,m})$. Un point clef de \cite{Har4} remarqué par Boyer dans
\cite{Boyer} dans le cadre des $\mathcal{D}$-modules elliptiques et
transposé par Harris et Taylor dans le cadre des groupes p-divisibles,
est que la cohomologie des autres strates est induite parabolique en
$p$, c'est à dire comme représentation de $G(\Qp)$. Cela est une
conséquence du fait qu'en restriction aux autres strates le groupe
p-divisible universel (pas $A[p^\infty]$, mais un déduit de celui-ci)
 est une extension d'un groupe étale par un
groupe du même type que celui associé à la strate supersingulière. 
Utilisant que 
les
déformations de tels groupes p-divisibles s'expriment simplement en
termes des déformations de leur partie connexe, l'astuce de Boyer s'en
déduit facilement. Une fois l'astuce de Boyer démontrée on peut facilement relier 
 la partie supercuspidale en
$p$ de la 
cohomologie  de la variété de Shimura  à celle des espaces $\spf
(R_{F,g,m})^{an}$, ce qui est la première étape de la récurrence de
\cite{Har4}. 

\section*{L'approche utilisée}

Dans \cite{ECM}
M. Harris  esquisse un programme  visant à généraliser
l'approche de \cite{Har4}. 
 Certains points de cette thèse sont inspirés de ce
programme. Nous utilisons une approche similaire à celle de \cite{Har4}
basée sur les travaux de Kottwitz sur les variétés de Shimura de type
P.E.L. (\cite{Ko1}) et ceux de Rapport et Zink sur l'uniformisation
$p$-adique de ces variétés (\cite{RZ}). Les
principaux points de \cite{Har4} ne se transposant  pas dans ce cadre
général  et les solutions que nous y avons apportées  sont les suivants :
\begin{itemize}
\item Contrairement aux variétés ``simples'' de \cite{Har4},
  on ne sait construire en général de bon modèles entiers des
 variétés de Shimura de type P.E.L. non ramifiées de
  \cite{Ko1}  qu'en niveau maximal
  (compact hyperspécial, ou parahorique) en $p$. Cela  est lié au fait
  qu'on ne sait pas définir de notion de structure de niveau de
  Drinfeld sur les groupes p-divisbles de dimension plus grande que
  un. L'alternative consiste à travailler au niveau des espaces
  rigides dès que l'on n'est plus en niveau maximal, et par exemple à
  définir le tube au dessus d'une strate comme étant l'image
  réciproque par l'application de changement de niveau vers un niveau
  maximal du tube rigide au dessus de la strate considérée. 
\item Du point de vue des cycles évanescents, 
  le lien entre la cohomologie des espaces de Rapoport-Zink et
  celle de la strate supersingulière via l'uniformisation $p$-adique
  est quasi-immédiat dans \cite{Har4}. En effet, la fibre spéciale des
  espaces $\M_K$ de \cite{Har4} est de dimension $0$. Le lien entre
  les deux est donc un ``problème de combinatoire''. Le problème dans
  notre cadre est que non seulement nous travaillons d'un point de vue
  rigide, mais en général l'uniformisation de \cite{RZ} n'est pas une
  simple réécriture du théorème de Serre-Tate; il y a des familles de  déformations
  par quasi-isogénie ``continues horizontales'' (c'est à dire non
  discrètes sur un espace 
  de paramètres de dimension
  $0$)  des groupes $p$-divisibles que nous considérons. 
Pour y remédier, 
nous
  démontrons dans la partie II une suite spectrale de
  Hochschild-Serre liée à l'uniformisation $p$-adique rigide qui relie
  la cohomologie des $\M_K$ à celle de la strate
  supersingulière. 
\item  L'astuce de Boyer n'admet pas de généralisation immédiate. Soit
  par exemple $X$ un groupe p-divisible sur la clôture algébrique
  d'un corps fini possédant deux pentes $\l_1,\l_2$ avec $0<\l_2<\l_2
  < 1$, et une filtration 
$$
0\ldrt X_{\l_2} \ldrt X \ldrt X_{\l_1} \ldrt 0
$$
où $X_{\l_1}$, resp. $X_{\l_2}$, est isoclin de pente $\l_1$, resp.
$\l_2$. Les déformations de $X$ ne se décrivent pas simplement à
partir de celles de $X_{\l_2}$ et de $X_{\l_1}$ comme c'était le cas
lorsque $X_{\l_1}$ était étale. 

La résolution de ce problème est à l'origine de la partie la plus
originale de cette thèse, la partie III. Dans \cite{ECM} M.Harris
propose de démontrer que la cohomologie des strates non basiques est
induite parabolique en $p$ en utilisant une formule des traces de
Lefschetz en géométrie rigide. Malheureusement, à part en dimension
$1$ (\cite{Hu4},\cite{St}) une telle formule des traces est loin
d'exister (il faudrait comprendre les termes associés aux points fixes
dans le bord des compactifications de Huber au sens des espaces
adiques). La solution que nous apportons à ce problème est la suivante
: utilisant les travaux de Fujiwara (\cite{Fuj}) nous démontrons une
formule des traces de Lefschetz pour la cohomologie des tubes rigides
au dessus des strates des variétés de Shimura que nous considérons. La
méthode consiste alors à comparer cette formule des traces
``$p$-adique'' à la formule des traces de Lefschetz archimédienne
topologique associée à l'uniformisation complexe des variétés de
Shimura. Utilisant 
 la conjugaison stable en une place auxiliaire différente de
$p$ pour séparer les classes de conjugaison elliptiques régulières 
des non-elliptiques en $p$,  
on montre le résultat. Les techniques utilisées au passage sont
très diverses : théorie des types, cycles évanescents rigides...
\item Reste la partie purement automorphe du programme : comparaison
  de la formule des traces d'Arthur pour les groupes unitaires de
  signature quelconque à l'infini. Nous utilisons pour cela des
  résultats non publiés de Harris et Labesse généralisant ceux de
  Clozel en signature $(1,n-1)$ (\cite{Clozel1}). Ces résultats sont
  assujettis à certaines conditions sur les groupes unitaires. Ce sont
  ces conditions qui imposent les restrictions sur $J_b$ dans le
  théorèmes \ref{theo_princ1} et \ref{theo_princ2}. Pour les résultats
  concernant $U(3)$ nous utilisons les résultats de Rogawski
  (\cite{Rogawski1}). 

Une meilleur connaissance des la stabilisation de la formule des traces
pour les groupes unitaires permettrait sûrement d'obtenir des
résultats complets sur $U(n)$ pour $n$ quelconque. 
\item La partie représentations galoisienne, c'est à dire la
  détermination du fameux $r_\mu\circ \psi$ est achevée en
  déterminant  la restriction à un groupe de décomposition des
  représentations galoisiennes globales découpées par certaines
  représentations automorphes dans les variétés de Shimura associées
  aux groupes unitaires. Nous utilisons pour cela l'existence
  (résultats non publiés de Harris et Labesse, résultats de Rogawski
  sur $U(3)$) d'un changement de base quadratique stable pour de
  telles représentations automorphes, la détermination par Kottwitz de
  cette représentation en presque toutes les places non ramifiées
  ainsi que l'un des résultats principaux de \cite{Har4} généralisant
  le théorème principal de Clozel dans \cite{Clozel1}. 

Bien qu'en annexe, ce résultat peut être considéré en soi même
comme un résultat indépendant  de cette thèse. 
\end{itemize}

\section*{Description  des différentes parties}
\subsection*{Partie I}
\begin{itemize}
\item Dans le premier chapitre nous rappelons et explicitons les
  définitions de base concernant les variétés de Shimura de type
  P.E.L. non ramifiées de \cite{Ko1}. 
\item La section \ref{Generalites_RZ} du chapitre 2 contient les définitions et
  propriétés de base
  concernant les espaces de Rapoport-Zink non ramifiés sans structure
  de niveau (en tant que schémas formels). 

La section \ref{sect_theo_fin}  contient le premier résultat de cette thèse,
le théorème \ref{TheoFini}. Il s'agit
d'une amélioration du théorème de finitude de 2.18 de \cite{RZ} basée
sur des résultats récents de Zink (\cite{Zink1}). 
Nous montrons que 
`` $J_b\bc \M (b,\mu)$''  est quasi-compact (pour tout $(b,\mu)$) 
, ou encore que sous
l'action de $J_b$ il n'y a qu'un nombre fini d'orbites de composantes
irréductibles dans la fibre spéciale de $\M (b,\mu)$. Le demi-plan
supérieur de Drinfeld possède une décomposition cellulaire indexée par
l'immeuble de Bruhat-Tits de $PGL_n$. Comme dans \cite{Mi1} nous
donnons une description des points sur $\bar{\mathbb{F}}_p$ de $\M
(b,\mu)$  comme un sous ensemble de l'immeuble de $G$ sur
$\widehat{\Qp^{nr}}$ 
(section \ref{Description_en_termes_cristaux}  et \ref{Lien_imm_G}).
 Le théorème \ref{TheoFini} s'interprète
alors comme un théorème de finitude sur cet immeuble. Il nous
permettra plus loin d'obtenir des théorèmes de finitude concernant la 
cohomologie à support compact des $\M (b,\mu)$ comme $J_b$-modules. 

La section 2.5 contient les définitions et propriétés de base
concernant les espaces $\M_K (b,\mu)$ (structures de niveau, action sur
ces structures). 

Dans cette thèse nous travaillons avec deux théories différentes pour
les variétés rigides et leur cohomologie étale $\ell$-adique à support
compact : la théorie des espaces analytiques de Berkovich et celle des
espaces adiques de Huber. La première est la version surconvergente de
la seconde: le topos étale analytique coïncide avec le topos étale
adique surconvergent. Les raisons pour lesquelles nous n'avons pas
fixé une des deux théories sont multiples :
\begin{itemize}
\item Le théorème de lissité de l'action sur la cohomologie
  $\ell$-adique  à support
  compact est formulé dans le cadre de la théorie de
  Berkovich. La démonstration de ce théorème n'étant pas publiée (\cite{Berk4}), nous en avons mis une
  démonstration en annexe \ref{coho_l}.
\item Les propriétés des cycles évanescents analytiques de Berkovich
  sont clairement exposées dans les deux articles \cite{Berk2} et
  \cite{Berk3}. Les propriétés des cycles évanescents adiques de Huber
  sont plus dispersées. 
\item Cependant, les théorèmes de finitude de Huber sont plus
  généraux (il a par exemple une définition des faisceaux
  constructibles). Nous utilisons notamment les théorèmes de
  finitude de 
  \cite{Hu2}.  
\end{itemize}
Nous avons donc mis en annexe \ref{an_ad} quelques propriétés de bases
qui nous l'espérons permettrons au lecteur de passer aisément d'une théorie à
l'autre. 

La définition et les propriétés de base (lissité des actions par exemple) de la
cohomologie $\ell$-adique  des espaces de Rapoport-Zink sont données
dans la section \ref{cohomo_de_MK}. Le lemme
\ref{lemme_annulation_Ext} la proposition \ref{Jfin} et le corollaire
\ref{Extf} constituent les trois théorèmes de finitude concernant
l'action de $J_b$ sur ces espaces de cohomologie. 
\item Le chapitre 3 contient des rappels et des précisions sur les
  résultats principaux de \cite{RZ} concernant l'uniformisation des
  variétés de Shimura de type P.E.L.. Pour paramétrer les strates de
  ces variétés nous utilisons l'ensemble $B(G,\mu)$ de Kottwitz
  (\cite{Ko3}) des classes d'isomorphismes d'isocristaux munis de
  structures additionnelles vérifiant le théorème de Mazur généralisé
  (\cite{RR}) sur la position relative des polygones de Hodge et
  Newton (le polygone de Hodge étant fixé par $\mu$ via la condition de
  Kottwitz). Nous donnons la description explicite de ces ensembles
  (ainsi que des groupes $J_b$ associés) 
  dans le cas de $\GL_n$ dans la section \ref{Isoc_GL}, et dans le cas
  unitaire dans l'appendice \ref{app_cl_sig}. 

On retiendra 
une propriété importante de 
  l'application d'uniformisation rigide énoncée
  dans ce chapitre : c'est un isomorphisme local
  (\ref{iso_loc_unif}). Cependant (remarque \ref{rem_u_pas_tout}),  
 à part dans le cas de la strate
  basique, les ouverts rigides uniformisés sont beaucoup plus petit
  que les tubes au dessus des strates 
 (c'est déjà clair dans le cas des
  variétés de \cite{Har4}, où l'uniformisation de Rapoport-Zink
  uniformise des complétés formels de la strate $\mu$-ordinaire,
  ouverte dans la fibre spéciale, le long de sous schémas de dimension $0$ de cette fibre
  spéciale).  
  C'est l'un des problèmes
  principaux qu'il faudrait  surmonter pour comprendre la contribution des strates
  non basiques. 
\end{itemize}

\subsection*{Partie II} 
 
  D'après la formule de Matsushima,
  la cohomologie de nos variétés de Shimura s'exprime en termes de 
représentations automorphes du groupe de similitudes unitaires $G$. 
  Une des idées qui est au coeur de cette thèse est que 
le facteur en $p$ 
(après restriction de la représentation galoisienne à un groupe 
de décomposition en $p$, c'est à dire comme 
comme représentation de $G(\Qp)\times W_{E_\nu}$ où
  $E_\nu$ est le complété en $p$ du corps réflex) de  certains
  morceaux de la cohomologie des variétés de Shimura que nous
  considérons ne devraient dépendre que de la cohomologie des objets
  locaux en $p$ que sont les espaces rigides de Rapoport-Zink. 

Dans \cite{Har4} cela est une conséquence du théorème de comparaison
de Berkovich entre cycles évanescents algébriques et analytiques
rigides (comparaison entre le site étale et étale formel d'un schéma
hensélien). Comme expliqué auparavant, nous utilisons une
approche différente. Nous démontrons l'existence d'une suite spectrale
de type Hochschild-Serre  reliant la cohomologie des espaces de
Rapoport-Zink à celle des tubes qu'ils uniformisent
$p$-adiquement. M.Harris avait déjà démontré l'existence d'une telle
suite spectrale dans le cas du demi-plan supérieur de Drinfeld
(\cite{Har1}). Il utilisait pour cela la décomposition cellulaire 
du demi-plan et des résultats de Berkovich sur les ``$G$-modules
discrets'' (non publiés). Une telle décomposition n'existe pas en
général. 
Notre méthode est différente et 
s'applique à tous les cas d'uniformisation connus (et conjecturaux,
confère \cite{ECM} conjecture 4.2.). 

La principale difficulté dans l'établissement d'une telle suite
spectrale provient des coefficients $\ell$-adiques. En effet, les
espaces $\M_K$ ne sont pas quasicompacts et en général
$$
H^\bullet_c (\M_K,\Z_\ell)\neq \underset{n}{\limp}
H^\bullet_c (\M_K,\Z/\ell^n \Z)
$$
à cause des problèmes d'interversion de limite inductive (sur le
support compact des sections) et de limite projective (sur les
coefficients). Nous contournons ce problème en nous inspirant de
l'article (\cite{Hu2}). Nous utilisons le foncteur dérivé $\R\pi_*$ de
l'annexe \ref{coho_l} (ainsi que ses versions équivariantes) 
introduit par  Ekedeahl dans \cite{Eke} qui envoi des faisceaux
$\ell$-adiques sur des complexes de faisceaux de $\Zl$-modules. 
Grâce au lemme \ref{Rpp}
nous
pouvons alors travailler avec les foncteurs dérivés des sections à
support compact (noté $\GG_!$) 
 dans la catégorie dérivée de tels faisceaux
de $\Zl$-modules, éliminant ainsi les problèmes de faisceaux
$\ell$-adiques. Cela permet de démontrer le théorème général
\ref{Hosch}. Le théorème principal \ref{sui_spe_hhhh} s'en déduit. 

L'application fondamentale est la formule de Matsushima $p$-adique
calculant la cohomologie de la strate basique en fonction de la
cohomologie des espaces de Rapoport-Zink basiques (ceux qui
interviennent dans les résultats principaux de la thèse) et des
 représentations automorphes d'une forme intérieure $I^\phi$
de $G$ (corollaire \ref{coro_Matsushima}) : 
il y a une égalité dans le groupe de Grothendieck
$\text{Groth}(G(\A_f)\times W_{E_\nu})$  
\begin{eqnarray*}
&\;& 
\sum_{i,j \atop { \Pi\in\mathcal{T}(I^\phi) \atop \Pi_\infty =\check{\rho}}}
(-1)^{i+j} \left [ \underset{K_p}{\limi}
  \text{Ext}^{\;i}_{J_b\text{-lisse}}
 \left ( H^j_c(\M_{K_p}\otimes \C_p,\Qlb)(N),\Pi_p\right ) \right ] 
\otimes 
\left [ \Pi^p \right ] \\
&=& \sum_i (-1)^i \left [ \underset{K}{\limi} H^i(
  \Sh_K(G,X)^{an}_{\text{basique}},\LL_{\rho}^{an}) \right ] 
\end{eqnarray*}
où $\mathcal{T}(I^\phi)$ désigne l'ensemble des représentations
automorphes de $I^\phi$, et $I^\phi$  vérifie $I^\phi
(\A_f^p)=G(\A_f^p)$, $I^\phi (\Qp)=J_b$. 

C'est la comparaison géométrique 
entre cette formule de Matsushima $p$-adique et
celle archimédienne
couplée à la comparaison analytique 
entre la formule des traces d'Arthur de $G$ et sa
forme intérieure $I^\phi$ qui permettra de démontrer 
dans la  partie IV les résultats principaux. 

\subsection*{Partie III} 

Le but de cette partie est de démontrer l'analogue de l'astuce de
Boyer en comparant deux formules des traces : l'une $p$-adique, c'est
à dire utilisant l'uniformisation $p$-adique, l'autre topologique
utilisant l'uniformisation archimédienne par l'espace symétrique
hermitien $X$ (un cas particulier de la formule des traces d'Arthur
telle qu'elle est reformulée dans \cite{Gor_Mac1}). 

\subsubsection*{Chapitre 6} 

Dans ce chapitre nous démontrons la formule $p$-adique
(théorème \ref{Lefschetz_speciale}).
 Plusieurs 
points sont au coeur de cette formule : le théorème de Fujiwara, le
fait  que la fibre des cycles évanescents algébriques ne dépend que du
complété formel en ce point... nous renvoyons à l'introduction de ce
chapitre pour plus de détails.

Pour démontrer une telle formule nous utilisons la formule des traces
de Lefschetz 
 modulo $p$, et plus particulièrement le théorème de Fujiwara
 (\cite{Fuj}). Pour cela nous utilisons les cycles évanescents
 analytiques rigides de Berkovich.  Cependant nos variétés de Shimura 
et les correspondances de Hecke associées 
ne possèdent pas de bon modèles entiers en tout niveau en $p$. L'idée
consiste alors à spécialiser (au sens de \cite{Fuj})
l'image directe sur la fibre générique
vers un niveau maximal des correspondances de Hecke en tant que
correspondances cohomologiques. Les correspondances cohomologiques
agissent sur la cohomologie,
 cependant nous avons besoin de
vérifier que la suite spectrale des cycles évanescents est
équivariante vis à vis de ces deux opérations : image directe
propre et spécialisation. C'est une des raisons pour laquelle nous
développons un formalisme de spécialisation des correspondances
cohomologiques rigides dans l'esprit de \cite{Fuj} et montrons que
celui-ci est compatible à l'image directe et à l'isomorphisme de
Berkovich. 

Le théorème  \ref{Lefschetz_speciale} est alors une application de ce
formalisme couplé au théorème de Fujiwara (\cite{Fuj}). 

Cependant un point technique doit être soulevé : même en niveau
compact hyperspécial $C_0$ seules les correspondances de Hecke à support
dans $C_0$ possèdent un bon modèle
entier (pour appliquer le théorème de Fujiwara une des deux applications 
de projections de la correspondance doit être quasi-finie). 
 Nous voudrions à priori appliquer notre formule des traces en mettant en
$p$ un pseudo-coefficient d'une représentation supercuspidale. Cela
n'est donc pas possible. L'alternative consiste à utiliser la théorie
des types de Bushnell et Kutzko.
 En effet, si $(J,\l)$ est un type supercuspidal 
l'idempotent associé
$e_\l$ est à support dans $J$ et donc dans un sous groupe compact
hyperspécial. Le problème est presque résolu à la remarque suivante
près : $e_\l$ ne permet pas de séparer les représentations
supercuspidales dans une même classe d'équivalence
inertielle. A tous les types supercuspidaux connus est associé un
couple $(\tilde{J},\tilde{\l})$ où $\tilde{J}$ est un sous groupe
compact modulo le centre de $G(\Qp)$ contenant $J$ et $\tilde{\l}$ une
extension de $\l$ à $\tilde{J}$. On remarque alors que l'on peut
séparer toutes les supercuspidales par les fonctions $e_\l*\delta_g$
où $g\in \tilde{J}$ (annexe \ref{Annexe_types}). Il suffit donc de faire agir
l'automorphisme associé à $g\in \tilde{J}$ sur nos modèles
entiers. Remarquant que pour les groupes $G(\Qp)$ avec lesquels nous
travaillons tous les groupes $\tilde{J}$ sont contenus dans le
normalisateur d'un sous groupe parahorique de $G(\Qp)$ on vérifie
que l'on peut s'en sortir (\ref{ext_para_lu})
en utilisant les modèles entiers définis en
niveau parahorique de \cite{RZ}. 

\subsubsection*{Chapitre 7} 

Dans ce chapitre nous redémontrons un cas très 
 particulier de la formule
des traces topologique d'Arthur telle qu'elle est reformulée dans
 \cite{Gor_Mac1} (théorème \ref{Lefschetz_generique}). 
En effet, en supposant la variété de Shimura compacte et la correspondance de
 Hecke à support dans les éléments réguliers en une place finie la formule prend une
 forme très simple. 
 
Bien que cette
 formule puisse se déduire de \cite{Gor_Mac1}, la lecture de
 \cite{Gor_Mac1} (plutôt technique)
 ne la laisse pas du tout transparaître (les auteurs
 de \cite{Gor_Mac1} s'intéressent surtout aux points fixes
 dans le bord). C'est pourquoi nous avons inclus sa démonstration.  
 
 Nous vérifions de plus en utilisant la théorie de Hodge $p$-adique 
(théorème \ref{lehteoint}) que si l'on
 disposait d'une formule des traces en géométrie rigide de la forme
une somme sur des points fixes locaux naïfs + des termes de Lefschetz
 au bord (au sens des compactifications de Huber des espaces adiques,
 confère \cite{Hu4} et \cite{St}), la distribution somme sur les points fixes naïfs
 associée à des strates non basiques est à support dans des éléments
 non-elliptiques réguliers ce qui ramènerait la démonstration de l'astuce de
 Boyer à montrer que la distribution associée aux termes au bord est
 ``induite'' en $p$. 

\subsubsection*{Chapitre 8} 
Celui-ci contient quelques propositions destinées à montrer le
théorème principal du chapitre 9.

\subsubsection*{Chapitre 9}

Il contient le théorème principal de la partie III (et l'un des
principaux de cette thèse), l'équivalent de l'astuce de Boyer : le
théorème \ref{coho_ind_en_p}. 

La démonstration consiste à comparer les deux formules des traces :
topologique et $p$-adique. La formule des traces topologique
\ref{Lefschetz_generique} 
pour la trace d'une fonction $f$ dans l'algèbre de Hecke globale 
 est une combinaison d'intégrales orbitales 
$\mathcal{O}_\g (f_p) \times \mathcal{O}_\g (f^p)$ où $\g\in G(\Q)$. 
La formule $p$-adique \ref{Lefschetz_speciale}
est une combinaison linéaire de termes de la forme $a_{\g}(f_p) \times
\mathcal{O}_\g (f^p)$ où $\g$ est une classe de conjugaison dans un
groupe $I^\phi (\Q)$ qui se transfert (via l'action sur la cohomologie
$\ell$-adique hors $p$ d'une variété abélienne) en un élément de
$G(\A_f^p)$ et permet de définir $\mathcal{O}_\g (f^p)$. Quant à $
a_{\g}(f_p)$ il s'agit d'une certaine trace sur 
la cohomologie étale d'un espace rigide de
déformation (par isomorphismes) de groupes p-divisibles. 
On suppose qu'en une place hors $p\;$ $f$ est à support dans les
éléments réguliers. Toutes les classes de conjugaison 
$\g$ précédentes sont alors régulières.

Si $f_p$ est dans l'algèbre de Hecke d'un type supercuspidal, dans la
formule des traces archimédienne seules des classes de conjugaison 
elliptiques en $p$ contribuent de façon non nulle.
Choisissons une place auxiliaire $w$ de $\Q$ différente de $p$.
 Égalant les deux
formules des traces et faisant varier les fonctions test en $w$ on en
déduit que si la classe d'isogénie $\phi$ ainsi que $\g\in I^\phi
(\Q)$ contribuent de façon non nulle alors $\g$ est conjugué dans
$G(\Q_w)$ à un élément de $G(\Q)$ elliptique dans $G(\Qp)$. 
Regardons alors $\g$ comme un élément de $I^\phi (\Qp)\subset J_b$
(où $b$ est déduit de $\phi$). La
forme intérieure quasidéployée de $J_b$ est le centralisateur du
morphisme des pentes $M(b)$
au sens de \cite{Ko2}, une forme intérieure d'un
sous groupe de Levi de $G_{\Qp}$. Le transfert de $\g$ vers cette
forme intérieure est donc un élément de $M(b)(\Qp)$ stablement
conjugué à un élément elliptique de $G(\Qp)$. Cela n'est donc possible
que si $M(b)=G_{\Qp}$, c'est à dire si $\phi$ intervient dans 
la strate basique. On en déduit donc que seuls ces termes
interviennent. Appliquant de nouveau la formule des traces $p$-adique
(cette fois ci à la strate basique) on obtient le résultat. 

\subsection*{Partie IV} 
Elle contient la démonstration des résultats annoncés au début de cette
introduction en utilisant la formule de Matsushima $p$-adique de la
partie II, le théorème \ref{coho_ind_en_p} du III et les résultats de
l'annexe A.

%%% Local Variables: 
%%% mode: latex
%%% TeX-master: "thesef"
%%% End: 

\mainmatter

\part{Vari{\'e}t{\'e}s de Shimura de type P.E.L. et espaces
de Rapoport-Zink} 
\chapter{Variétés de Shimura de type P.E.L. non ramifiées} 
\section{Donnée de Shimura de type P.E.L.}

Soit $\mathcal{D}=(B,*,V,<.,.>,h)$ une donnée de Shimura de type
P.E.L.(\citeg{Ko1}). Plus précisément : 
\begin{itemize}
\item $B$ est une algèbre simple sur un corps de nombres $F$
\item $*$ est une involution positive sur $B$ : $\forall b\in B \;
  tr(bb^*)>0$ 
\item $(V,<.,.>)$ est un $B$-module hermitien où $<.,.>:V\times V \drt
  \Q$ est symplectique
\item $h:\C\drt \End_B(V)_\R$ est un morphisme d'algèbres {\`a} involutions
  où $\End_B(V)$ est muni de l'adjonction par rapport {\`a} la forme
  symplectique  $<.,.>$ et $\C$ de la conjugaison complexe notée $c$. De plus
  $<.,h(i).>$ est un produit scalaire.
\end{itemize}
Notons 
$$
G=\{ g\in\Gl_B(V)\; |\; gg^*=c(g)\in \Q^* \}
$$
le groupe des similitudes unitaires défini sur $\Q$ associé et
$G_1=\ker (c)$ le 
groupe unitaire.  
     
Le morphisme $h$ définit une $\Q$ structure de Hodge polarisée munie d'une action 
de $B$. Cette structure de Hodge munie de structures additionnelles est
un point dans la variation de structures de Hodge associée aux
variétés de Shimura que nous allons utiliser.

Nous noterons $X$ la classe de $G(\R)$ conjugaison de $h:\C^\times 
\drt G_{\R}$. Au morphisme  $h$ on peut associer 
$$
\mu_h:\C\*\drt G_\C
$$
(ou plutôt sa classe de conjugaison) définissant la filtration de Hodge
sur $V_\C=V_0\oplus V_1$ où $\mu_h (z)=z$ sur $V_1$ et $\mu_h=1$ sur
$V_0$. Dans toute la suite nous suivrons la convention homologique de
Kottwitz, c'est {\`a} dire l'inverse de \cite{Del_Var_Sh}. 

Nous noterons $\Qb$ la clôture algébrique de $\Q$ dans $\C$. 
Soit $E$ le corps reflex associé {\`a} $\mathcal{D}$. Le corps $E$
est le corps de 
définition de la classe de conjugaison de $\mu$, 
$$
E=\Q[tr_\C(b;V_0)\;|\; b\in B] \subset \Qb
$$

\subsection{Le cas (A)}\label{Le_cas_A}

Nécessairement, $F$ est soit un corps
C.M. soit un corps totalement réel. 
Supposons que $F$ est un corps C.M., ce que nous appellerons le cas
$(A)$. 
 
Fixons $\Phi\subset \Hom (F,\C)$ un type C.M. de $F$. Il existe alors 
des signatures $(p_\tau,q_\tau)_{\tau\in\Phi}$ telles que $p_\tau+q_\tau=n$ soit 
indépendant de $\tau$ dans $\Phi$ et 
$$
G_{1/\R}\simeq \prod_{\tau\in\Phi} U(p_\tau,q_\tau)
$$ 
De plus
$$
G_{\C}\simeq \prod_{\tau\in\Phi} \Gl_n(\C) \times \C\*
$$
le dernier facteur représentant le facteur de similitude d'un élément 
de $G(\C)$. 

Avec ces notations, 

$$
h(z)=\prod_{\tau\in\Phi} \left ( 
\begin{array}{cccccc}
z \\
& \ddots \\
& & z \\
& & & \bar{z} \\
& & & & \ddots \\
& & & & & \bar{z}
\end{array}
\right )
$$
$$\;\;\;\;\;\;\;\;\;\;\;\hspace{5mm}
\;\;\underbrace{\hspace{1.8cm}}_{p_\tau} \underbrace{\hspace{1.8cm}}_{q_\tau}
$$
et 
$$
\mu(z)=\prod_{\tau \in\Phi} \left ( 
\begin{array}{cccccc}
z \\
& \ddots \\
& & z \\
& & & 1 \\
& & & & \ddots \\
& & & & & 1 
\end{array}
\right )
\times (z) 
$$

On a la relation
$$
n = [B :F]^{1/2} \; \text{rg}_B V 
$$

Le corps $E$ peut alors se décrire {\`a} partir de ces données. Faisons
opérer $\Gal(\Qb | \Q)$ sur les $\Phi$-uplets de couples d'entiers
$(a_\tau,b_\tau)$ de la fa{\c c}on suivante : 
$$\forall \s\in \Gal ( \Qb |
\Q)\;\;\; \s. (a_\tau,b_\tau)_{\tau\in \Phi}= (a'_\tau,b'_\tau)_{\tau\in
\Phi}$$ où 
$$\forall \tau\in \Phi \;\;
(a'_\tau,b'_\tau) = \left \{ (a_{\s^{-1}\tau},b_{\s^{-1}\tau}) \text{
si } \s^{-1}\tau\in \Phi  \atop 
(b_{c\s^{-1}\tau},a_{c\s^{-1}\tau}) \text{
si } \s^{-1}\tau\notin \Phi \right.
$$
$\Gal(\Qb |E)$ est alors le stabilisateur de $(p_\tau,q_\tau)_{\tau\in
\Phi}$ pour cette action. 

\begin{exem} Notons $F^\Phi$ le corps réflex de $(F,\Phi)$ défini
  par l'égalité
$$
\Gal( \Qb | F^\Phi)= \{\tau\in \Gal (\Qb|\Q)\; |\; \tau \Phi=\Phi \;\}
$$ 
Supposons qu'il existe $\tau_0\in \Phi$ tel que $q_{\tau_0}\neq 0$ et
$\forall \tau\neq \tau_0\; (p_\tau,q_\tau)=(n,0)$. Supposons de plus que
$p_{\tau_0}\neq n/2$. Le corps
 $E$ est alors le composé 
$$
E= F^\Phi. \tau_0 (F) 
$$ 
C'est par exemple le cas de l'article \cite{Clozel1} dans
lequel $(p_{\tau_0},q_{\tau_0})=(1,n-1)$.  

Par contre, dans  \cite{Har4}, $(p_{\tau_0},q_{\tau_0})=(1,n-1)$ tandis que $\forall \tau\neq \tau_0 \; (p_\tau,q_\tau)=(0,n)$. On a encore dans dans ce cas l{\`a} 
$ E= F^\Phi. \tau_0 (F)$.
\end{exem}
\begin{exem} 
Notons $F^+$ le sous-corps totalement réel maximal de $F$. Supposons
que $F$ est de la forme $F^+.\mathcal{K}$ où $\mathcal{K}|\Q$ est une
extension quadratique imaginaire et que le C.M. type $\Phi$ est induit
à partir d'un C.M. type de $\mathcal{K}$. Le corps réflex de $(F,\Phi)$ est
alors $\mathcal{K}$. Nous ferons parfois cette hypothèse.

Par exemple,sous cette hypothèse, dans l'exemple précédent $E=\tau_0 (F)$ ce qui est le cas de \cite{Har4}.  
\end{exem}

\subsection{Le cas (C)}

Considérons maintenant le cas où $F$ est un corps totalement réel. Dans ce cas
l{\`a},  $G_1$ est, soit une forme d'un produit de groupes
symplectiques, 
soit 
une forme de produits de groupes orthogonaux suivant que $B_\R$ est
un produit d'algèbres $M_n(\R)$ ou $M_n(\mathbb{H})$ (cas (C) et (D) de 
\citeg{Ko1} respectivement). Nous ne considérerons que le premier
cas. Alors, 
$$
G_{1/\R}\simeq \prod_{\tau: F\hookrightarrow \R} Sp_n(\R)
$$
où $Sp_n=\{g\in \Gl_n \;|\; {}^t g J g= J\},\; J=\left (
\begin{array}{cc}
0 & -I_n \\
I_n & 0 \\
\end{array} \right )
$.
Avec ces notations, 
$$
h(a+ib)=\prod_{\tau:F\hookrightarrow  \R}
\left (
\begin{array}{cc}
a I_n & -b I_n\\
b I_n  & a I_n  \\ 
\end{array}
\right )
$$
et  
$$
\mu(z)=\prod_{\tau:F:\hookrightarrow \R} 
\left (
\begin{array}{cc}
z I_n  & 0 \\
0 & I_n \\ 
\end{array}
\right )
$$

\section{Variétés de Shimura sur $E$}

A la donnée $\mathcal{D}$ est associée une tour $(\Sh_K)_K$, pour des 
sous-groupes compacts ouverts $K$ de 
$G(\A_f)$ suffisamment petits, de variétés quasiprojectives
lisses sur $E$ classifiant les quadruplets $(A,\l,\iota,\overline{\eta})$ où
:
\begin{itemize}
\item $A$ est un schéma abélien {\`a} isogénie près 
\item $\l$ est une polarisation $\Q\*$ homogène de $A$. Cette dernière
  condition signifie que pour
  toute fonction $\a$ définie  sur la schéma de base du schéma abélien,
   à valeurs dans $\Q^\times$ et constante sur chaque composante connexe,
   $\l$ est considérée comme étant équivalente à $\a \l$. 
\item $\iota:B\drt \End(A)_\Q$ est un morphisme d'algèbres 
tel que $*$ corresponde {\`a}
l'involution de Rosati associée {\`a} $\l$. Cela signifie que si $\dag$
désigne l'involution de Rosati 
$$
\forall b\in B\;\;\; \iota (b^*)= \iota (b)^{\dag}
$$
 Une condition équivalente est
de dire que $\l$ induit un morphisme de variétés abéliennes munies
d'une action de $B$ par quasi-isogénies entre 
 $(A,\iota)$ et  $(A^\vee,\iota\circ *)$ 
\item $\eta: V\otimes \A_f \ldrt H_1(A,\A_f) [K]$ est un isomorphisme 
de $B\otimes \A_f$-modules symplectiques  
définissant une structure de niveau $K$ sur le module de Tate de $A$.
Cet isomorphisme est défini {\`a} une constante de $\A_f\*$
près puisqu'il faut fixer un isomorphisme $\A_f(1)\*\simeq \A_f\*$. 
\item La condition suivante est vérifiée : 
$$\forall b\in B\;
\det(b;Lie(A))=\det(b;V_0)
$$
Cette condition est l'analogue algébrique de la condition topologique : $\forall
h'\in X, h'$ est conjugué {\`a} $h$. 
\end{itemize}

\begin{exem}
Pla{\c c}ons nous dans le cas (A) de la section \ref{Le_cas_A}. L'action de $F$ sur $\text{Lie} (A)$ permet de décomposer 
$$
\Lie (A)=\bigoplus_{\tau\in \Hom (F,\Qb)} \text{Lie} (A)_\tau
$$
où $F$ agit sur $ \text{Lie} (A)_\tau$ via $\tau$. La condition de Kottwitz se récrit alors dans ce cadre sous la forme :
\begin{eqnarray*}
\forall \tau\in\Phi\; &\Lie (A)_\tau& \text{ est localement libre de rang } p_\tau [B:F]^{1/2} \\
\text{et } &\Lie (A)_{c\tau} & \text{ est localement libre de rang } q_\tau [B:F]^{1/2}
\end{eqnarray*}

\end{exem}

La tour de variétés 
$(\Sh_K)_K$ est munie d'une action de $G(\A_f)$ par action sur la
structure de niveau $\overline{\eta}$.

La variété de Shimura $\Sh_K$ n'est pas en général associée {\`a} la
donnée
de Shimura 
 $(G,X)$. On a en
fait la décomposition suivante définie sur $E$ (\citeg{Ko1} section 8) : 
$$
\Sh_K=\coprod_{ker^1(\Q,G)} \Sh_K(G',X)
$$
où $G'$ parcourt des formes intérieures de $G$ isomorphes {\`a} $G$ partout
localement, 
et où $\Sh_K(G',X)$ désigne le modèle canonique sur $E$ 
de la variété de Shimura
associée {\`a} la donnée de
Shimura $(G',X)$. 
 Plus précisément, si $x\in\Sh_K(\C)$ est associé {\`a}
$(A,\l,\iota,\overline{\eta})$, $x$ appartient {\`a} la composante indexée par
la classe du $G$-torseur localement trivial 
$$
Isom_{\text{B-mod.symp.}}(V,H_1(A,\Q))
$$
et $G'$ est le groupe des similitudes unitaires associé au $B$-module
symplectique $H_1 (A,\Q)$.

On remarque (\citeg{Ko1} chapitre 7)
 que dans le cas (C) 
et dans le cas (A) avec $n$ pair,  
 $G$
vérifie le principe de Hasse, c'est {\`a} dire que $\ker^1 (\Q,G)$ est réduit {\`a} un seul
élément. 
Dans ces deux cas les variétés $\Sh_K$ sont donc des variétés de
Shimura associées {\`a} la donnée de Shimura $(G,X)$. 

Il est également démontré dans la section 7 de \citeg{Ko1} que dans le cas 
(A) avec $n$ impair
l'application 
$$ 
\ker^1 (\Q,Z_G)\ldrt \ker^1 (\Q,G)
$$
est surjective. Dans ce cas, tous les groupes $G'$ sont donc en fait isomorphes {\`a}
$G$, et les variétés $\Sh_K$ sont donc une union finie de copies d'un modèle canonique de la 
variété de Shimura associée {\`a} la donnée $(G,X)$ (\cite{Del_Var_Sh},\cite{Milne_Sh}). 

\section{Systèmes locaux}

Soit $\ell$ un nombre premier différent de $p$.
Fixons définitivement un isomorphisme
$$
\iota : \Qlb \iso \C
$$

Soit $\rho$  une représentation algébrique complexe de $G$. 
  $\rho$
définit alors via $\iota$ un $\overline{\Ql}$ faisceau lisse $\LL_\rho$ sur les
différents $\Sh_K$ lorsque $K$ varie. Il est
muni d'une action de $G(\A_f)$ au sens où $\forall K_1\subset K_2$, si 
$\Pi_{K_1,K_2}:\Sh_{K_2}\drt \Sh_{K_1}$ désigne le morphisme de
  changement
de niveau, il y a des isomorphismes 
$ \Pi_{K_1,K_2}^*\LL_\rho
\drt\LL_\rho$   
et pour tout 
$g$ élément de $G(\A_f)$, si $g:\Sh_K\iso\Sh_{g^{-1}K g}$, 
 il y a des
morphismes $g^*\LL_\rho \drt \LL_\rho$, le tout vérifiant des conditions 
de compatibilité évidente. 

La variété 
$\Sh_K$ est munie d'un pro-revêtement étale $(\Sh_{K'}\ldrt
\Sh_K)_{K'\lhd K}$ de
groupe $K$. Composant l'application de projection $K\ldrt G(\Ql)$ avec
la représentation $\rho: G(\Ql)\hookrightarrow G(\Qlb) \xrig{\;^\iota \rho\;} \Gl
(V_{\;^\iota\rho})$ on obtient une représentation continue $K\ldrt  \Gl
(V_{\;^\iota\rho})$. Le système local $\LL_\rho$ est celui associé à cette
représentation. 

\section{Modèles entiers}
\subsection{Hypothèse de non ramification}

Fixons un plongement $\nu:\overline{\Q}\hookrightarrow \overline{\Qp}$ et notons $E_\nu$ le complété de $E$ correspondant. 

Soit $\O_B$ un $\Z_{(p)}$ ordre dans $B$ stable par $*$ et tel que
$\O_B\otimes\Z_p$ soit un ordre maximal dans $B_{\Qp}$.

Nous supposerons (confère \citeg{Ko1}) : 
\begin{itemize}
\item $B_{\Qp}$ est un produit d'algèbres de matrices sur des extensions
  non ramifiées de $\Qp$, c'est {\`a} dire :  
 $B$ est déployée en toutes les places de
  $F$ divisant $p$ et $F$ est non ramifié en toutes les places de $F$
  divisant $p$ 
\item Afin de s'assurer que $G_{\Qp}$ est quasidéployé nous
 supposerons qu'il existe un
réseau autodual 
  $\Lambda_0$ dans $V_{\Qp}$ 
\end{itemize}
Sous ces hypothèses le groupe $G_{\Qp}$ est non ramifié.

\subsection{Modèles}

Sous les hypothèses précédentes Kottwitz construit dans \cite{Ko1} 
 une tour $(S_{\pK})_{\pK} $
de variétés quasiprojectives lisses sur $\O_{E_\nu}$ indexée par des sous-groupes
compacts ouverts 
suffisamment petits  $\pK\subset
G(\A_f^p)$.
 Cette tour est munie d'une action de $G(\A_f^p)$. Le schéma
$S_{K^p}$  est défini comme l'espace de modules  des quadruplets
$(A,\l,\iota,\overline{\eta^p})$ où : 
\begin{itemize}
\item $A $ est un schéma abélien défini {\`a} une isogénie première {\`a} $p$ près 
\item $\l$ est une polarisation première {\`a} $p$ et $\Z_{(p)}\*$
  homogène. Cela siginifie que si $\a$ est une fonction  définie sur le schéma
  de base sur lequel est défini $A$, à valeurs dans $\Z_{(p)}\*$ et
  localement constante, alors $\l$ est considérée comme étant
  équivalente à $\a\l$. 
\item $\iota:\O_B\drt \End(A)\otimes \Z_{(p)}$ est un morphisme
d'algèbres 
 transformant l'involution 
$*$ sur $\O_B$ 
en l'involution de Rosati associée {\`a} $\l$ sur $\End
 (A)\otimes \Z_{(p)}$
\item $\eta^p:V\otimes \A_f^p \ldrt H_1(A,\A_f^p) [\pK]$ est un
  isomorphisme de $B\otimes \A_f^p$-modules symplectiques 
  définissant une structure de
  niveau $K^p$ au sens de la définition 2.5.5 de \cite{Ko1}. Cet
  isomorphisme est pris {\`a} un élément de 
  $\A_f^{p\times}$ près. 
\end{itemize}

Le groupe 
$G(\A_f^p)$ opère sur la tour $(S_{K^p})_{K^p}$ par action sur
$\overline{\eta^p}$. 

Soit un réseau autodual $\Lambda_0\subset V_{\Qp}$. Soit 
$C_0=Stab_{G(\Qp)}(\Lambda_0)$
le sous-groupe compact hyperspécial associé. Il y a alors des 
isomorphismes compatibles pour $K^p$ variant 
$$
S_{\pK}\otimes_{\O_{E_\nu}} E_\nu \iso \Sh_{C_0 \pK}\otimes_{E} E_\nu
$$
Ces isomorphismes dépendent du choix de $\Lambda_0$. 

\section{Structure de $G(\Qp)$\label{Struct}}
\subsection{Cas (A) }
 On se place dans le cadre de la section \ref{Le_cas_A}. 
Tout élément $\tau\in \Phi$ fournit un plongement 
$$
\nu \circ \tau: F \hookrightarrow \Qpb
$$
Soient $(w_i)_{i\in I}, (w_j)_{j\in J}$ les places de $F$ divisant $p$ 
associées {\`a} tous ces plongements lorsque $\tau$ parcourt $\Phi$, et où
on suppose que  
 $\forall i\in I \;
w_i\neq w_i^c$ et $\forall j\in J \; w_j= w_j^c$.
On a alors,
$$
G_{\Qp}=\prod_{i\in I} \Gl_n (F_{w_i})\times G \left (
\prod_{j\in J} GU(F_{w_j} ; n )
\right )
$$
où si $J=\emptyset, G(\emptyset)=\Qp\*$, 
 $\Gl_n(F_{w_i})$ désigne la restriction des scalaires de $F_{w_i}$
{\`a} $\Qp$ de $\Gl_n$, $GU(F_{w_j}; n)$ est la forme
intérieure quasidéployée du groupe des similitudes unitaires en n
variables sur $F_{w_j}$ vu comme groupe sur $\Qp$ (plus précisément, il
 s'agit du sous-groupe de la
restriction des scalaires de $F_{w_j}$ {\`a} $\Qp$ 
 du groupe des similitudes unitaires sur
 $F_{w_j}$ 
 ayant un facteur de similitude  dans $\Qp\*$)  
 et où le $G$ devant le produit signifie que
l'on prend le sous-groupe du produit formé des uplets ayant m{\^e}me
facteur de similitude.

L'algèbre semi-simple 
$B_{\Qp}$ est de la forme 
$$
\prod_i \left ( M_d(F_{w_i})\times M_d(F_{w_i})^{opp} \right )
\times \prod_j M_d(F_{w_j})
$$

L'équivalence de Morita permet de supposer qu'en chaque place on est
dans l'un des cas suivant :  
\subsubsection{ Cas (AL) }

 Soit $i\in I$, 
$B_{\Qp}=F_{w_i}\times F_{w_i},\; \O_{B_{\Qp}}=\O_{F_{w_i}}\times
  \O_{F_{w_i}},\;
(x,y)^*=(y,x),\;  V_{\Qp}=V_i\oplus V_i^\vee,\;
<x\oplus\varphi,x',\oplus\varphi'>
=\varphi'(x)-\varphi(x')$

Un réseau autodual est de la forme
$\Lambda_0=\Lambda\oplus\Lambda^\vee$ où $\Lambda$ est un réseau de $V_i$.

\subsubsection{Cas (AU) }

 Soit $j\in  J$, 
$B_{\Qp}=F_{w_j}$, l'involution
$*$ est égale {\`a} $\s^{[F_{w_j}:\Qp]/2}$ où $\s$ est le
Frobenius de l'extension non ramifiée $F_{w_j}|\Qp$.
$\O_{B_{\Qp}}=\O_{F_{w_j}}$,
Dans une base convenable de $V_j$ identifié alors {\`a} $F_{w_j}^n$ :
$$
\forall X,Y\in F_{w_j}^n\; \; <X,Y>=\text{Tr}_{F_{w_j}/\Qp}(\a {}^t X^*
J Y) 
$$ 
où $\a\in F_{w_j}\*$ est tel que $\a^*=-\a$,
$$
J=\left (\begin{array}{cc}
0 & I_{n/2} \\
I_{n/2} & 0 \\
\end{array}\right )
$$
si $n$ est pair  et
$$
J=\left ( \begin{array}{ccc}
0 & 0  & I_{(n-1)/2} \\
0 & 1 & 0 \\
I_{(n-1)/2} & 0 & 0 \\
\end{array} \right )
$$
si $n$ est impair.

Un réseau autodual dans le cas où $n$ est pair est alors donné par 
$\Lambda_0=\Lambda^{(1)}\oplus\Lambda^{(2)}$ où $\Lambda^{(1)},
\Lambda^{(2)}$ sont deux sous-$\O_{F_{w_j}}$-modules
totalement isotropes tels que l'accouplement  $\Lambda^{(1)}\times 
\Lambda^{(2)}\xrig{<.,.>} \Zp$ soit parfait.

Dans le cas où $n$ est impair, on peut prendre comme réseau autodual 
 $\Lambda_0=\Lambda^{(1)}\oplus \O_{F_{w_j}}.\e\oplus
\Lambda^{(2)}$ où $\Lambda^{(1)}$ et $\Lambda^{(2)}$ sont comme
précédemment et $<\e,\e>=1$. 

\subsection{Cas (C) }
$$
B_{\Qp}=\prod_{w|p} M_d (F_w)
$$
$$
G_{\Qp}=G(\prod_{w|p} GSp_n(F_w))
$$
où $Gsp$ est le sous-groupe de la restriction des scalaires de $F_w$ {\`a}
$\Qp$ du groupe des similitudes symplectiques qui est formé des éléments ayant
un facteur de similitude dans $\Qp\*$, et
où le $G$ signifie que les facteurs de similitude des différents
facteurs sont égaux .

On peut supposer par équivalence de Morita que l'on est dans le cas
suivant : $B_{\Qp}=F_w, \O_{B_{\Qp}}=\O_{F_w}, V_{\Qp}=F_w^n, <X,Y>={}^t XJ
Y$
où 
$$
J=\left ( \begin{array}{cc}
0 & -I_n \\
I_n & 0 \\
\end{array} \right )
$$
Un réseau autodual est de la forme $\Lambda_0=\Lambda^{(1)}\oplus
\Lambda^{(2)}$ où $\Lambda^{(1)},\Lambda^{(2)}$ sont des  
sous-$\O_{F_w}$-modules 
totalement
isotropes et l'accouplement
 $\Lambda^{(1)}\times \Lambda^{(2)} \xrig{<.,.>} \Zp$ est
parfait.

\chapter{Espaces de Rapoport-Zink}
\section{Isocristaux munis de structures additionnelles : le cas de
  $GL_n$}

Nous donnons ici une description explicite pour le groupe linéaire
de la classification des
isocristaux munis de structures additionnelles donnée par Kottwitz dans
\citeg{Ko2} et \citeg{Ko3}. Étant donné que nous avons essentiellement
en vu comme
application dans cette thèse  le cas de $GL_n$ nous avons mis les
autres cas en annexe. 

\subsection{$\GL_n$}\label{Isoc_GL}

Oublions momentanément les notations globales du chapitre précédent.

Soit $F|\Qp$ une extension non ramifiée de degré d. Soit  
$G=\Res_{F/\Qp} (\Gl_n)$. Le groupe algébrique $G$ se déploie sur $F$ en 
$$
G_{/F}=\prod_{i\in \Z/d\Z} \Gl_{n/F}
$$
où le groupe $\Gamma=\Gal(F|\Qp)=\Z/d\Z$ opère par permutation sur les facteurs
en translatant par un élément de $\Z/d\Z$. 

Un tore maximal de $G$ défini sur $\Qp$ est $T=(F\*)^n$. Un tore
déployé maximal est $A=(\Qp\*)^n\hookrightarrow (F\*)^n=T$. Le groupe
de Weyl absolu est $W=(\S_n)^d$. 

\begin{defi}
Le corps  $L$ est  le complété de l'extension maximale non ramifiée de $\Qp$.

L'automorphisme  $\s$ est le Frobenius géométrique de l'extension $L|\Qp$.  
\end{defi}

\begin{defi}
  On appelle isocristal un $L$-espace vectoriel $N$ muni d'une
  application bijective $\s$-linéaire $V:N\ldrt N$.

On appelle $F$-isocristal un $L$-espace vectoriel $N$ muni d'une
application bijective $\s^d$ linéaire $V:N\ldrt N$ ($\s^d$ est le
Frobenius géométrique de l'extension $L|F$).  
\end{defi}

Rappelons également que deux éléments $b_1,b_2$ de $G(L)$ sont dits $\s$-conjugués s'il existe
un $g$ dans $G(L)$ tel que 
$$
b_1= g b_2 g^{-\s}
$$ 
Rappelons alors la définition 
\begin{defi}[\cite{Ko2}]
L'ensemble $B(G)$ désigne l'ensemble des classes de $\s$-conjugaison
dans $G(L)$.
\end{defi}

L'ensemble 
$B(G)$  classifie les
isocristaux munis d'une action de $F$. A la classe de $\s$-conjugaison
de $b$ est associée la classe d'isomorphisme de l'isocristal muni d'une
action de $F$, ($F^n\otimes_{\Qp} L, b \s)$. 
 Il est bien connu que la catégorie des isocristaux munis d'une action
 de $F$ 
 est équivalente {\`a} la catégorie des $F$-isocristaux, ce qui
se traduit dans le langage de Kottwitz par le lemme de Shapiro
$B(G)\simeq B(\Gl_{n/F})$ (\citeg{Ko2}).

Un élément  $b\in B(G)$ est donné par les pentes du $F$-isocristal
associé $$\l_1=\frac{d_1}{h_1}>\dots>\l_r=\frac{d_r}{h_r}$$ où $d_i\wedge
h_i=1$, et des multiplicités $m_1,\dots,m_r$ vérifiant $\sum_i m_i h_i=n$. Le
point de Newton associé (\cite{Ko2}) est
$$
\nu_b=(\underbrace{\frac{\l_1}{d},\dots,\frac{\l_1}{d}}_{m_1
  h_1},\dots,\underbrace{\frac{\l_r}{d},\dots,\frac{\l_r}{d}}_{m_r
  h_r})\in X_*(A)_\Q=\Q^n
$$
(le dénominateur $d$ provient du fait que les pentes d'un $F$-isocristal sont $d$
fois les pentes de l'isocristal sous-jacent). Le centralisateur du
morphisme des pentes $\nu_b$ est le sous-groupe de Levi 
$$
M_b=\Res_{F/\Qp}(\Gl_{m_1 h_1})\times\dots \times \Res_{F/\Qp}(\Gl_{m_r h_r})
$$
et $b$ provient d'une classe basique de $M_b$, où l'on rappelle qu'une
classe basique est une classe possédant une unique pente.   
Rappelons que l'on note $J_b$ le groupe des automorphismes de
l'isocristal muni de ses structures additionnelles. Dans le cas
présent, la structure additionnelle consiste en l'action de $F$. $J_b$
est donc le groupe des automorphismes du $F$-isocristal associé : 
$$
J_b=\Res_{F/\Qp}\left (\Gl_{m_1}(D_{\l_1})\right )\times\dots\times \Res_{F/\Qp}\left 
( \Gl_{m_r}(D_{\l_r})\right )
$$
où $D_\l/F$ est l'algèbre centrale simple sur $F$ 
d'invariant $\l$. 

Le groupe $G^{der}$ étant simplement connexe,
$X^*(Z(\widehat{G})^\Gamma)=X_*(D)_\Gamma=B(D)$ où $D=G/G^{der}$.
Or, le déterminant induit un isomorphisme de groupes  
 $\det:G/G^{der}\iso \Res_{F/\Qp}\Gm$, et donc $B(D)=\Z$. On identifiera
donc $X^*(Z(\widehat{G})^\Gamma)$ {\`a} $\Z$. 
L'application de Kottwitz $\kappa$ définie dans la section 6 de
\cite{Ko5} (confère également le théorème 1.15 de \cite{RR}) est alors 
\begin{eqnarray*}
\kappa : B(G) & \ldrt & \Z=X^*(Z(\widehat{G})^\Gamma) \\
 b & \longmapsto & v_p(det (b))
\end{eqnarray*}

et est telle que $$\kappa(b)=\dpt{\sum_{i=1}^r m_i d_i}$$
qui est la dimension de l'isocristal ou encore le point terminal du
polygone de Newton.
\\

{\bf  Polygone de Hodge : }
Soit maintenant un cocaractère $\mu:{\Gm}_{/\overline{\Qp}}\drt
G_{/\overline{\Qp}}$. Ce cocaractère 
 est défini sur $F$ et définit un élément de 
$$
X_*(T)/W\simeq \{(x_{ij})_{i\in\Z/d\Z,1\leq j \leq n} |\; x_{ij}\in
\Z, \forall i\;\;
\forall j\leq j' \;\;\; x_{ij}\geq x_{ij'} \}
$$
sur lequel le groupe de Galois $\Gamma$ opère via $\s.(x_{ij})=(x_{i-1
  j})_{i,j}$. 

Nous supposerons $\mu$ minuscule donné, avec les notations ci dessus, par $\mu=(\mu_{ij})_{i,j}$ où 
$\forall i \; \mu_{ij}=1$ pour $1\leq j\leq a_i$ et $\mu_{ij}=0$ si
$j>a_i$, les $a_i$ étant des entiers compris entre $1$ et $n$. 

Le groupe de Galois $\Gamma=\Z/d\Z$ agit par permutation sur
$(a_i)_{i\in\Z/d\Z}$,  et le
corps de définition de la classe de conjugaison de $\mu$ a pour groupe
de  Galois le stabilisateur de $(a_i)_{i\in\Z/d\Z}\in \N^{\Z/d\Z}$
  dans $\GG=\Z/d\Z$. 

Le L-groupe connexe de $G$ est donné par 
$$
{}^L G^0=\widehat{G}=\prod_{\Z/d\Z} \Gl_n(\C)
$$ 
et le L-groupe par 
${}^LG=\widehat{G}\rtimes \Z/d\Z$ où $\Z/d\Z$ agit par permutation
cyclique des
composantes. Au cocaractère $\mu$ de $G$, Kottwitz associe un
caractère du centre du L-groupe connexe $Z(\widehat{G})$ défini par 
$\mu_1=\widehat{\mu}_{|Z(\widehat{G})^\Gamma}\in \Z$ (\cite{Ko3}
section 6.1), 
 et donné avec nos notations par 
$$
\mu_1=\sum_{i\in \Z/d\Z} a_i
$$
le point terminal du polygone de Hodge.

De même (\cite{Ko3} section 6.1), Kottwitz définit  $\mu_2$ par 
\begin{eqnarray*}
  X_*(T)_\Gamma & \xrig{\; Newt_T\;} & X_*(A)_\Q=\Q^n \\
\widehat{\mu}_{|\widehat{T}^\Gamma} & \longmapsto & \mu_2 
 \end{eqnarray*}
Avec les notations précédentes 
$$
\mu_2=\frac{1}{d} \sum_{i\in\Z/d\Z} (\underbrace{1,\dots,1}_{a_i}
,0,\dots , 0) 
$$

Kottwitz définit dans la section 6 de \cite{Ko3} le sous ensemble $B
(G,\mu)$ de $B(G)$ formé des classes de $\s$-conjugaison vérifiant le
théorème de Mazur sur la position relative des polygones de Newton et
de Hodge tel qu'il généralisé dans \cite{RR}. 
D'après ce qui précède, dans notre cas particulier cet ensemble possède la
description suivante   
\begin{eqnarray*}
B(G,\mu) &= &\{b\in B(G)\;|\; Newt(b)\leq \mu_2\in X_*(A)_\Q \text{ et }
\kappa(b)=\mu_1 \} \\
&=& \left \{ 
(\underbrace{\frac{\l_1}{d},\dots,\frac{\l_1}{d}}_{m_1
  h_1},\dots,\underbrace{\frac{\l_r}{d},\dots,\frac{\l_r}{d}}_{m_r
  h_r}) \leq 
\frac{1}{d} \sum_{i\in\Z/d\Z} (\underbrace{1,\dots,1}_{a_i}
,0,\dots , 0) \right \}
 \\ 
 &\hspace{20mm} & \text{ et  }  \sum_{i=1}^r m_i d_i=\sum_{i\in \Z/d\Z} a_i \;\;
 \\
\end{eqnarray*}
où $\leq$ désigne l'ordre usuel sur la chambre de Weyl positive
(section 2.1 de \cite{RR}). L'inégalité s'interpète comme ``le
polygone de Newton est au dessus du polygone de Hodge'', et l'égalité
comme ``les points terminaux des polygones de Newton et de Hodge sont
égaux''. 

\begin{exem}
  Lorsque $F=\Qp$, et $b\in B(G)$, $b\in B(G,\mu)$ ssi l'isocristal
  filtré $(L^n,b\s,\mu)$ est faiblement admissible au sens de
  Fontaine. 
\end{exem}

Pour un élément de $B(G)$, la condition d'appartenance à $B(G,\mu)$
est donc une condition généralisant la condition d'admissibilité faible
de l'isocristal filtré associé qui tient compte des structures additionnelles.
\\
%\begin{rema}
%Lorsque $\mu$ est défini sur $\Qp$, la condition d'appartenir {\`a}
%$B(G,\mu)$ est équivalente {\`a} celle d'admissibilité faible classique du
%$F$-isocristal associé.
%\end{rema}

{\bf Classes basique dans $B(G,\mu)$ : }

Il existe une unique classe basique dans $B(G,\mu)$. Notons la
$b_0$. Avec les notations précédentes celle-ci est
donnée par l'unique pente 
$$
\l=\frac{1}{n} \sum_{i\in \Z/d\Z} \; a_i 
$$
Si la pente $\l$ est égale à $\dpt{\frac{r}{s}}$ où $r \wedge s=1$, alors, la multiplicité de
la pente est $\dpt{\frac{n}{s}}$, et l'on a 
$$
J_{b_0}=\Res_{F/\Qp} \Gl_{\frac{n}{s}} ( D_\l) 
$$

\section{Généralités sur les espaces de Rapoport-Zink}\label{Generalites_RZ}

\subsection{Données locales de Rapoport-Zink}

\subsubsection{Donnée de type E.L. non ramifiée simple}\label{EL_RZ}

Une donnée de type E.L. non ramifiée simple $(F,V,b,\mu)$ 
consiste en la donnée
\begin{itemize}
\item
 d'une
extension finie non ramifiée $F$ de $\Qp$ 
\item
d'un $F$-espace vectoriel de
dimension finie $V$
\item d'une classe de $\s$ conjugaison $b\in B(G)$ où
$G=\Res_{F/\Qp} \Gl (V)$ 
\item 
d'une classe de conjugaison de
cocaractère minuscule $\mu :{\Gm}_{\Qpb} \ldrt G_{\Qpb}$ du type de
celles de la section \ref{Isoc_GL}.
\end{itemize}
 On suppose que ces données sont reliées entre elles par la relation 
$$
b\in B(G,\mu)
$$

Notons $I_F=\Hom_{\Qp}( F, \Qpb)$ et $n=\dim_F V$. $\mu$ est alors donné par des
couples d'entiers
$$
(p_\tau,q_\tau)_{\tau\in I_F}
$$
vérifiant $p_\tau+q_\tau=n$ (dans les notations de la section
\ref{Isoc_GL}, ces couples d'entiers sont les $(a_i,n-a_i)_{i\in \Z/
d\Z}$). 

Le corps réflex local associé $E$ est le corps de définition de la classe
de conjugaison de $\mu$. $\Gal (\Qpb |\Qp)$ opère par permutation sur
$I_F$. $\Gal ( \Qpb | E)$ est alors le stabilisateur dans
$\Gal(\Qpb|\Qp)$ de
$(p_\tau,q_\tau)$. 

\subsubsection{Donnée de type P.E.L. unitaire non ramifiée simple}
\label{PEL_RZ}

Une donnée de type P.E.L. unitaire non ramifiée  simple $(F,*,V,<.,.>,b,\mu)$ consiste
en la donnée 
\begin{itemize}
\item d'une extension finie $F$ de $\Qp$ non ramifiée munie d'une
involution non triviale $\; *$
\item d'un $F$-espace vectoriel de dimension finie
\item d'un produit hermitien symplectique $<.,.>: V\times V \ldrt \Qp$
qui est tel qu'il existe un réseau autodual dans $V$ pour ce produit
symplectique 
\item si $G$ désigne le groupe des similitudes unitaires associé,
d'une classe de $\s$ conjugaison $b\in B(G)$
\item d'une classe de conjugaison de cocaractère minuscule 
$\mu:{\Gm}_{\Qpb} \ldrt G_{\Qpb}$
\end{itemize}
On suppose que ces donnée sont reliées entre elles par la relation 
$$
c\circ \mu (z) =z 
$$
où $c$ désigne le facteur de similitude, et 
$$
b\in B(G,\mu)
$$
Le donnée d'un tel $\mu$ est équivalente {\`a} la donnée de couples
d'entiers
$$
(p_\tau,q_\tau)_{\tau\in I_F}
$$
vérifiant $p_\tau+q_\tau=n$ et  $(p_{\tau *},q_{\tau
 *})=(q_\tau,p_\tau)$ (dans la section \ref{Isoc_AU} ces entiers sont
notés $(a_i)_i$). 

 Le corps réflex associé se décrit comme précédemment.

\subsubsection{Donnée de type P.E.L. symplectique non ramifiée simple}

Elle consiste en la donnée $(F,V,<.,.>,b,\mu)$ vérifiant des
conditions analogues {\`a} celles du cas unitaire. $G$ est alors le
groupe des similitudes symplectiques. 

\subsection{Espaces de Rapoport-Zink associés}

Nous récrivons la définition 3.21 de \cite{RZ} dans nos cas non ramifiés.
 Notons
$\Eb=\widehat{\Qp^{nr}}$ le complété de l'extension
maximale non ramifiée de $\Qp$.
. Nous noterons $\bar{k}\simeq
\overline{\mathbb{F}}_p$ le corps résiduel de $\Eb$.

\subsubsection{Cas E.L. non ramifié}\label{esp_RZ_ass}

\begin{defi}
Soit $(V,F,b,\mu)$ une donnée locale de type E.L. non ramifiée simple. 
Soit $\XX,$ un groupe p-divisible sur $\bar{k}$ muni d'une action de
$\O_F$ et d'isocristal muni de son action de $F$ isomorphe {\`a}
$(V_L,b\s)$ (l'existence d'un tel $\XX$ résulte de résultats récents de
Rapoport et Kottwitz). 
Nous 
noterons $\M (b,\mu)/\spf (\O_{\Eb})$ l'espace de Rapoport-Zink
associé. C'est l'espace de  
 modules de groupes
p-divisible munis d'une action de $\O_F$, $X$, munis d'une
quasi-isogénie compatible {\`a} l'action de $\O_F$ 
$$
\rho : \XX \ldrt X
$$
et vérifiant la condition : si 
$$
\Lie (X)=\bigoplus_{\tau \in I_F} \Lie (X)_\tau
$$
est la décomposition de $\Lie (X)$ suivant l'action de $\O_F$ alors,
$$
\forall \tau\in I_F \;\; \Lie (X)_\tau \text{ est localement libre de
rang } p_\tau
$$
\end{defi}

$\M (b,\mu)$ est localement formellement de type fini sur $\spf (\O_{\Eb})$  ce qui signifie
 que localement
$\M(b,\mu)$ est de la forme 
$$
\spf(\O_{\Eb}<T_1,\dots,T_n>[[ X_1,\dots,X_m]]/I)
$$
pour un idéal $I$.

Il résulte de la théorie de Grothendieck Messing que 
$\M (b,\mu)$ est formellement lisse sur $\spf (\O_{\Eb})$.
Soit $x$ un point fermé de sa fibre spéciale. Le complété formel 
$\M(b,\mu)_{/\{x \}}^{\widehat{}}$ s'identifie {\`a} l'espace de déformation 
par isomorphismes du groupe $p$-divisible $X^{univ}_x$muni de son action de
$\O_F$ ($X^{univ}$
désignant le groupe p-divisible universel sur $\M (b,\mu)$).
Il résulte de la théorie de Grothendieck Messing que cet espace de
déformation est isomorphe au complété formel d'une certaine
Grassmanienne (les fameux modèles locaux de \cite{RZ} (définition
3.27)) en un point : 
$$
\M (b,\mu)^{\widehat{}}_{/\{ x\}} \simeq \spf \left ( \O_{\Eb} [[
X_1,\dots,X_{\sum_\tau p_\tau q_\tau} ]] \right )
$$
\vspace{0mm}

\subsubsection{Cas P.E.L. unitaire non ramifié} \label{esp_RZ_ass_PEL}

\begin{defi}
Soit $(F,*,V,<.,.>,b,\mu)$ une donnée locale de type
P.E.L. unitaire non ramifiée. Nous noterons $\M (b,\mu)/\spf
(\O_{\Eb})$ l'espace de Rapoport-Zink 
associé.  Soit $(\XX,\l)$ un un groupe $p$-divisible sur $\bar{k}$
muni d'une action de $\O_F$, d'une polarisation principale $\l$ tels
que l'isocristal polarisé de $(\XX,\l)$ muni de son action de $F$
s'identifie {\`a} $(V,b\s,<.,.>)$ (une fois de plus l'existence d'un tel
$(\XX,\l)$ devrait résulter de travaux récents de Rapoport et
Kottwitz). Alors, le schéma formel $\M (b,\mu)$ classifie les groupes p-divisibles
principalement polarisés $(X,\l')$ munis d'une action de $\O_F$, tels
que via cette action l'involution de Rosati sur $\O_F$ soit donnée par
l'involution $*$, et munis d'une quasi-isogénie compatible {\`a} l'action
de $\O_F$ : 
$$
\rho : \XX \ldrt X
$$
tels via $\rho$ $\;\l$ et $\l'$ différent, localement sur la base, d'une constante dans
$\Qp^\times$. On suppose de plus que si 
$$
\Lie (X)=\bigoplus_{\tau\in I_F} \Lie (X)_\tau
$$ 
désigne la décomposition de $\Lie (X)$ sous l'action de $\O_F$ alors 
$$
\forall \tau \;\; \Lie (X)_\tau \text{ est localement libre de rang }
p_\tau 
$$
\end{defi}

Comme précédemment dans le cas E.L., $\M (b,\mu)$ est formellement
lisse sur $\spf (\O_{\Eb})$. Son complété formel en un point fermé
de sa fibre spécial
 est le spectre formel d'une algèbre de séries
formelles en $\frac{1}{2} \sum_\tau p_\tau q_\tau$ variables. 

\subsection{Cas P.E.L. symplectique}

Étant donnée une donnée locale symplectique non ramifiée simple on peut 
définir de la m{\^e}me fa{\c c}on que dans le cas unitaire l'espace de
Rapoport-Zink associé. Nous le noterons $\M (b,\mu)$.

\subsection{Décomposition d'une donnée globale non ramifiée en produits de donnée
locales simples}\label{Dec_glo_loc}
\subsubsection{Cas (A) }

Revenons aux notations globales du chapitre précédent. La donnée
globale  $\mathcal{D}$
couplée au plongement $\nu$ fournit une donnée 
$(B\otimes \Qp, V\otimes \Qp, <.,.>, \mu_{\Qpb})$ où $\mu_{\Qpb}$ est
obtenu {\`a} partir de $\mu$ gr{\^a}ce {\`a} $\nu$. Celle ci se scinde en un
produit suivant les indices $i$ dans $I$ et $j$ dans $J$.

 Utilisant
l'équivalence de Morita on en déduit pour tout $i$ dans $I$ une donnée
de type $E.L.$ non ramifiée simple (sans sa classe de $\s$ conjugaison):
$(F_{v_i},V_i,\mu_i)$. Nous noterons $E_i$ son corps réflex. 
Avec les notations de \ref{EL_RZ},  le cocaractère 
$\mu_i$
est associé aux couples suivant d'entiers : soient $\Phi_i \subset
\Phi$ l'ensemble des $\tau$ dans $\Phi$ vérifiant $\nu\circ \tau$
induit $v_i$, alors
$$
I_{F_{v_i}} = \{ \nu\circ \tau\; |\; \tau\in \Phi_i \}
$$ 
et $\mu_i$ est associé aux 
$$
(p_{i,\tau},q_{i,\tau})_{\tau\in I_{F_{v_i}}} = (p_{\nu \circ \tau},
q_{\nu\circ \tau})_{\tau\in \Phi_i}
$$

Si maintenant $j\in J$, on obtient une donnée de type
P.E.L. non ramifiée unitaire simple 
sans sa classe de $\s$ conjugaison 
$(F_{v_j},V_j,*,<.,.>,\mu_j)$ où $*$ est la conjugaison complexe et
 $\mu_j$ se calcule de la m{\^e}me fa{\c c}on que
précédemment. Nous noterons $E_j$ son corps réflex associé. 
\\

Il y a donc une inclusion 
$$
E_\nu \subset \prod_{i\in I} E_i \prod_{j_\in J} E_j
$$

Remarquons maintenant que 
$$
B(G_{\Qp},\mu_{\Qpb})=\prod_{i\in I} B(\Gl_n(F_{w_i}),\mu_i) \times \prod_{j\in J}
B(GU(F_{w_j},n),\mu_j)
$$
En effet, l'inclusion $G_{\Qp}\hookrightarrow \prod_{i\in I} \Gl_n (
F_{w_i}) \times \prod_{j\in J} GU( F_{w_j},n)$ induit une inclusion
du membre de gauche dans  celui de droite. Si maintenant  $b=((b_i)_{i\in
  I},(b_j)_{j\in J})$ désigne un élément du membre de droite, pour
tout $j$ dans $J$ la condition d'appartenance à $B(GU(
F_{w_j},\mu_j))$ implique que $v_p (c(b_j))=1$. Quitte à
$\s$-conjuguer les $b_j$ on peut donc supposer que pour tout $j$ dans
$J$, $c(b_j)=p$. Cela montre que $b$ provient d'une classe de
$\s$-conjugaison de $G_{\Qp}$. 

Pour tout élément $b$ de $B(G_{\Qp},\mu_{\Qpb})$ notons $b=\prod_{i_\in I} b_i
\prod_{j\in J} b_j$ sa décomposition.  Si $b\in B(G,\mu)$ est donné,
la donnée globale $\mathcal{D}$ fournit donc des données locales pour
tout $i$ dans $I$ et $j$ dans $J$. 

\begin{defi}
Pour un $b\in B(G_{\Qp},\mu_{\Qpb})$ nous noterons
$(\mathcal{D}_{\Qp},b)$ la donnée locale 
de Rapoport-Zink $(B\otimes \Qp,V\otimes\Qp,<.,.>,b,\mu_{\Qpb})$ 
qui se décompose en un  produit de données simples comme ci dessus. 
\end{defi} 

\subsubsection{Cas (C)} 

La donnée se décompose en un produit de données symplectiques non
ramifiées simples indexées
par les places de $F$ divisant $p$. Pour $b\in
B(G_{\Qp},\mu_{\Qpb})\;$ on note également $(\mathcal{D}_{\Qp},b)$ la donnée
locale associée.

\subsection{Espace de Rapoport-Zink associé {\`a} $(\mathcal{D}_{\Qp},b)$ }

Rappelons rapidement la définition de l'espace de Rapoport-Zink
associé {\`a} $\mathcal{D}$ après tensorisation par $\Qp$. 

Soit donc $b\in B(G_{\Qp},\mu_{\Qpb})$, soit
$\Eb=\widehat{E^{nr}_\nu}=\widehat{\Qp^{nr}}$ le complété de l'extension
maximale non ramifiée de $E$ 
. Nous noterons $k$ le corps résiduel de $E_\nu$ et $\bar{k}\simeq
\overline{\mathbb{F}}_p$ le corps résiduel de $\Eb$. Choisissons un
groupe p-divisible $\mathbb{X}/\bar{k}$ ayant pour isocristal muni de
structures additionnelles $(V\otimes_{\Qp} L,b\s)$. Nous supposerons que
la polarisation de $\XX$ est principale c'est {\`a} dire que le réseau associé au
cristal de $\XX$ est autodual dans $V_L$.

Rappelons la définition de $\M(\mathcal{D}_{\Qp},b)$ : 

\begin{defi} Le schéma formel 
  $\M(\mathcal{D}_{\Qp},b)/\spf(\O_{\Eb})$ est le schéma formel localement formellement de
  type fini vérifiant : $\forall S/\O_{\Eb}$ un schéma sur lequel $p$
  est localement nilpotent, 
$$
\M(\mathcal{D}_{\Qp},b)(S)=\{(X,\l,\iota,\rho) \}/\sim
$$
où $X$ est un groupe p divisible sur $S$ muni d'une action $\iota$ de
$\O_B\otimes\Zp$, d'une polarisation principale $\l:X\ldrt \check{X}$ 
qui induit l'involution $*$ sur $\O_B\otimes\Zp$, et d'une
rigidification 
$$
\rho:\XX\times_k \overline{S} \ldrt X\times_S \overline{S}  
$$
(où $\overline{S}$ désigne la réduction de $S$ modulo $p$) 
qui est une quasi-isogénie compatible {\`a} l'action de
$\O_B\otimes\Zp$. On suppose de plus que
via $\rho$ les polarisations associées à $<.,.>$ sur $\XX\times_k
\overline{S}$ et $\l$ sur 
$X\times_S\overline{S}$ diffèrent localement sur $\overline{S}$ d'une
constante de $\Qp\*$. Cela signifie qu'il existe une fonction
localement constante $\a:S\ldrt \Qp^\times$ telle que le diagramme
suivant commute 
\begin{diagram}
\XX\times_{\bar{k}} \overline{S} & \rTo ^\rho & X\times_S \overline{S}
\\
\dTo & & \dTo_{\a.(\l\times Id)} \\
\check{\XX}\times_{\bar{k}}\overline{S} & \lTo^{\check{\rho}} & \check{X}\times_S\overline{S}  
\end{diagram}
 On suppose de plus que $X$ vérifie la condition de Kottwitz
 relativement {\`a} $\mu_{\Qpb}$.
Deux
quadruplets $(X,\l,\iota,\rho)$ et $(X',\l',\iota',\rho')$ sont
équivalents s'il existe un isomorphisme $f:X\iso X'$ compatible aux
actions, aux polarisations et tel que $f\circ \rho=\rho'$
\begin{diagram}
& & X\times_S\overline{S} \\
 & \ruTo^{\rho} \\
\XX\times_{\bar{k}} \overline{S} & & \dTo_{f}^{\simeq} \\
& \rdTo_{\rho'} \\
& & X'\times_S\overline{S}
\end{diagram}
\end{defi}

\begin{rema}
 Il résulte des travaux de Colmez, Fontaine, Breuil
que $\XX$ se relève en un groupe p-divisible sur $\O_{\Eb}$ dont
l'isocristal filtré est $(V_L,b\s,\mu_{\Qpb})$.
\end{rema}

\begin{rema}
  L'espace $\M(\mathcal{D}_{\Qp},b)$ est formellement lisse sur $\spf
  (\O_{\Eb})$. Puisque, comme on la verra, les espaces de
  Rapoport-Zink uniformisent dans certains cas des strates des
  variétés Shimura,
cela peut paraître en contradiction avec l'existence de
  singularités dans ces strates. Néanmoins, la fibre spéciale de  
$\M(\mathcal{D}_{\Qp},b)$ (le sous schéma réduit) est en général
  singulière et il n'y a pas de lien entre la lissité formelle d'un
  schéma formel et les singularités de sa fibre spéciale. Ce qui
  pourrait induire en erreur est l'exemple de l'espace Drinfeld. Dans
  ce cas là, les singularités de la fibre spéciale du schéma formel de
  Drinfeld (des diviseurs à croisement normaux) sont identiques aux
  singularités de la fibre spéciale des courbes de Shimura qu'ils
  uniformisent. Cela est dû au fait que ce schéma formel est adique
  sur $\spf (\O_{\Eb})$ et que donc, sa réduction modulo $p$ est un
  schéma. En général $\M(\mathcal{D}_{\Qp},b)$ n'est pas un schéma
  formel $p$-adique.
\end{rema}

\subsection{Décomposition de l'espace associé {\`a} $(\mathcal{D}_{\Qp},b)$ en
produit d'espaces simples de type E.L. et P.E.L.} \label{Dec_glob}
\subsubsection{Cas (A) }

Commen{\c c}ons par remarquer que si $\M (b,\mu)$ désigne un espace de
Rapoport-Zink de type P.E.L. unitaire non ramifié simple, celui ci est muni
d'un morphisme 
$$
\b : \M (b,\mu) \ldrt \Qp^\times/\Zp^\times
$$
mesurant le ``défaut'' de la polarisation du groupe p-divisible
universel {\`a} co{\"\i}ncider avec celle définie par $<.,.>$ sur $V_{\Qp}$. 
\\

Soit $X^{univ}$ le groupe p-divisible universel sur
$\M(\mathcal{D}_{\Qp},b)$. L'action de $\O_{B}\otimes \Zp$ sur 
$X^{univ}$, la polarisation $\l$ et $*$ 
 permettent de décomposer $X^{univ}$ en 
$$
X^{univ} = \bigoplus_{i\in I} ( X^{univ}_i \oplus \check{X}^{univ}_i)
\oplus 
\bigoplus_{j\in J} X^{univ}_j
$$

Utilisant le fait que $\O_B\otimes \Zp$ est un ordre maximal, 
on peut appliquer une ``équivalence de Morita'' {\`a} chaque
terme pour obtenir que
l'espace de Rapoport-Zink associé {\`a} $(\mathcal{D}_{\Qp},b)$ 
se décompose  de la fa{\c c}on suivante si
$J\neq \emptyset$ :
$$
\M(\mathcal{D}_{\Qp},b) = \prod_{i\in I} \M(b_i,\mu_i) \times \left (\prod_{j\in J} \M
(b_j,\mu_j)\right )^1
$$
où  
$$
\left ( \prod_{j\in J}
\M (b_j,\mu_j) \right )^1\hookrightarrow  \prod_{j\in J}
\M (b_j,\mu_j) 
$$
est l'ouvert fermé défini par l'égalité : $\forall j,j'\in J,\; \b_j=\b_{j'}$. 

Dans le cas où $J=\emptyset$, 
$$
\M(\mathcal{D}_{\Qp},b)=\prod_{i\in I} \M (b_i,\mu_i) \times \Qp\*/\Zp\*
$$

\begin{rema}
Cette décomposition en produit d'espaces plus simples
par équivalence de Morita  est un des
points clef de tout ce qui suit. Par exemple, elle permet 
dans \cite{Har4} de ramener
l'étude d'espaces de modules associés {\`a} des groupes p-divisibles de
dimension supérieur {\`a} $1$ {\`a} celle de groupes p-divisibles de dimension $1$. 
\end{rema}

\begin{exem}\label{defQp}
Supposons que pour un $i\in
I$, $\mu_i$ soit tel qu'il existe un  $\tau_0\in I_{F_{v_i}}$ tel que
$\forall \tau\neq \tau_0 \; p_{i,\tau}=0$.
Alors,  le schéma formel 
 $\M(b_i,\mu_i)$ est l'espace de modules des $\O_{F_{w_i}}^{\tau_0}$-modules
$p$-divisibles $X$ (au sens de \cite{Har4}) sur un schéma $S$ sur lequel $p$ est localement nilpotent
 munis d'une rigidification avec un  $\O_{F_{w_i}}^{\tau_0}$-module
$p$-divisible sur 
$\overline{\mathbb{F}_p}$ de $F$ isocristal associé {\`a} $b_i$ et vérifiant  
$$
\text{Lie}(X)\text{ est un }  \O_S-\text{module
libre de rang } p_{i,\tau_0} 
$$
\end{exem}

\begin{exem}\label{exemAetale}
Avec les notations précédentes, si pour un $i$ appartenant {\`a} $I$ $\;\mu_i=1$
alors, $B(\Gl_n (F_{w_i},\mu_i))$ est réduit {\`a} un seul élément $b_i$ qui est
associé {\`a} un isocristral étale et  
$$
\M (b_i,\mu_i)=\Gl_n (F_{w_i})/ C_0
$$
comme espace discret (c'est {\`a} dire union disjointe de copies de 
$\spf(\O_{\Eb})$).  

Il en est de m{\^e}me pour un $j$ appartenant {\`a} $J$. 
\end{exem}

\begin{exem}  
Avec les notations de l'exemple \ref{defQp} si $p_{i,\tau_0}=1$ et
 $b_i$ est basique, 
 $\M(b_i,\mu_i)$ paramétrise des déformations par
quasi-isogénies de $\O_{F_{v_i}}^{\tau_0}$-modules p-divisibles de dimension $1$ et de
hauteur $n$.
$$
\M(b_i,\mu_i)\simeq \spf(\O_{\Eb}[[T_1,\dots,T_{n-1}]])\times \Z 
$$
l'espace de Lubin Tate sur lequel opère $J_{b_i}=D_{\frac{1}{n}}\*$. 
\end{exem}
 
\subsubsection{Cas (C) } 

Comme ci dessus, il y a une décomposition 
$$
\M(\mathcal{D}_{\Qp},b) =\left ( \prod_{w|p} \M (b_w,\mu_w)\right )^1
$$

\section{Continuité de l'action de $J_b$}

Nous notons dans cette section $\M$ un espace de Rapoport-Zink associé
{\`a} une donnée locale simple ou bien {\`a} une donnée globale.

Le schéma formel $\M$ est muni d'une action {\`a} gauche de $J_b$ via 
$$
\forall g\in J_b\;\; \forall (X,\rho)\in \M\;\;\;\;
(X,\rho).g=(X,\rho\circ g^{-1}) 
$$

Lorsque $\M=\M (\mathcal{D}_{\Qp},\mu_{\Qpb})$, 
cette action se fait composante par composante au sens suivant : dans
le cas (A), 
$\prod_{i\in I} J_{b_i} \times \prod_{j\in J} J_{b_j}$ opère sur 
$$
\prod_{i\in I} \M (b_i,\mu_i) \times \prod_{j\in J} \M (b_j,\mu_j)
$$
et $J_b\subset \prod_{i\in I} J_{b_i} \times \prod_{j\in J}
J_{b_j}$ stabilise l'ouvert fermé $$\M(\mathcal{D}_{\Qp},b) \subset \prod_{i\in I} \M(b_i,\mu_i) \times
\prod_{j\in J} \M (b_j,\mu_j)$$

\begin{defi}
Soit $\X$ un schéma formel muni d'une action d'un groupe
localement compact 
totalement discontinu $H$. Nous dirons que $H$ agit contin{\^u}ment sur
$\mathfrak{X}$ si tout ouvert  quasicompact $U$ de $\mathfrak{X}$
 est stabilisé par un voisinage de l'identité dans $H$,
 et si $\mathcal{I}_U\subset
\O_U$ est un idéal de définition de $\O_U$, $\forall N\in\N \;$
il existe un sous-groupe compact ouvert $K_N\subset H$ suffisamment
petit tel que $\forall g\in K_N,\; g^*:\O_U/\mathcal{I}_U^N \ldrt
\O_U/\mathcal{I}_U^N$  est l'identité.
\end{defi}

\begin{prop}\label{ActC} Le groupe 
$J_b$ opère contin{\^u}ment sur $\M$.
\end{prop}
\begin{proof} On vérifie aisément que l'assertion ci dessus est
équivalente {\`a} la suivante : si $S$ est un schéma quasicompact sur
lequel $p$ est nilpotent et $(X,\rho)\in \M(S)$, il existe un
voisinage de l'identité dans $J_b$ tel que pour $g $ dans un tel
voisinage, $(X,\rho\circ g)\simeq (X,\rho)$. 

Notons  $\End(\XX)$ l'anneau des endomorphismes du groupe
$p$-divisible $\XX$ muni de ses structures additionnelles ($\l$ et
$\iota$) 
sur $\overline{\mathbb{F}_p}$. L'anneau $\End(\XX)$ est un $\Z_p$-module libre de
type fini. De plus, $\forall N>0, Id+p^N \End(\XX)$ est un
sous-groupe compact ouvert 
de $J_b$ et lorsque $N$ varie cette famille de sous-groupes forme une
base de voisinages de l'identité dans $J_b$.

Considérons la quasi-isogénie sur $\overline{S}$ 
$$
\rho\circ (Id+p^N (u\times 1)) \circ \rho^{-1}= Id+ p^N 
\rho (u\times 1) \rho^{-1} \in \End(X\times_S\overline{S})_\Q
$$  
Le schéma 
$S$ étant  quasicompact et $p$ nilpotent, par rigidité des
quasi-isogénies
$$\forall u\in \End(\XX)\;\exists N \; p^N \rho (u\times 1)\rho^{-1}
\text{ se relève en un endomorphisme de } X/S
$$
 Mais l'algèbre $\End(\XX)$ est un $\Zp$-
module de type fini. On peut donc trouver un tel $N$ uniformément
pour tous les éléments $u$ de $\End(\XX)$.
 C'est {\`a} dire, 
 $$\exists N \;\forall u\in\End(\XX) \;\; p^N \rho(u\times 1)\rho^{-1}
\text{ se relève en un endomorphisme de } X$$ Nous noterons  
$ (p^N \rho(u\times
  1)\rho^{-1})\widetilde{}$ ce relèvement.
 Comme $\End(\XX)$ est un anneau p-adique,
$$
Id+p ( p^N \rho(u\times
  1)\rho^{-1})\widetilde{} \in\Aut (X)
$$
On en déduit que $\forall u\in\End(\XX)$, $Id+p^{N+1} \rho (u\times
1)\rho^{-1}$ se relève en un automorphisme de $X$ ce qui montre que 
$$
(X,\rho)\simeq (Id+p^{N+1}u).(X,\rho)
$$
\end{proof}

\begin{rema}
Bien s{\^u}r, cette démonstration s'applique {\`a} un espace de Rapoport-Zink
général tel qu'il est défini dans \cite{RZ}. 
\end{rema}

\section{Un théorème de finitude pour l'action de $J_b$}\label{sect_theo_fin}

Dans cette section nous oublions les notations globales précédentes et
considérons des espaces associés {\`a} des donnés locales simples. 

\begin{defi}
On posera comme précédemment $L=W(\bar{k})_{\Q}$. L'automorphisme $\s$ désignera toujours un Frobenius
géométrique de l'extension non ramifiée $L|\Qp$. Le cristal associé {\`a} un groupe p-divisible  
est son cristal de Dieudonné covariant $(M,V)$ où $V$ est une application 
$\s$ linéaire, le Verschiebung.
\end{defi}

Bien que dans notre situation $L=\Eb$ nous préférons distinguer les
notations puisque le rôle de ces deux corps est totalement différent.

\subsection{Description de $\M(b,\mu)(\bar{k})$ en termes de
  cristaux}\label{Description_en_termes_cristaux}  

Soit donc une donnée locale de type E.L. ou P.E.L. non ramifiée simple.

Considérons l'isocristal $(V\otimes_{\Qp} L,b\s)$ muni de son action de $F$
et, dans le cas P.E.L.,  de
sa polarisation $$<.,.>:V_L\times V_L\ldrt \Qp(1)$$
 Un élément de $\M(b,\mu)(\bar{k})$
s'identifie alors {\`a} un cristal $M$ dans $(V_L,b\s)$ stable par l'action de
$\O_F$ tel que $M$ comme $\O_L$ réseau soit autodual {\`a} une constante
de $\Qp\*$ près (c'est {\`a} dire $\exists k\in\Z\; M^\vee=p^k M$) et tel qu'on ait
une égalité de fonctions polynomiales :
$$
\forall b\in\O_F\;\;\; \det(b;M/VM)=\det(b;V_0)
$$
Dit autrement, la dernière condition s'exprime en disant que le
polygone de Hodge arithmétique ``avec structures additionnelles'' est
fixé égal {\`a} $\mu$. Reprenons les notations de \ref{Generalites_RZ} et
détaillons les différents cas :
\\

{\bf Cas (AL) : } Reprenons les notations de \ref{Isoc_GL}. 
Posons $N=V_L$, $V=b\s$. Décomposons $N$ sous l'action de $F$ 
$$
N=\bigoplus_{i\in\Z/d\Z} N(i)
$$
où
$$
N(i)=\{ n\in N\; |\; \forall x\in F \; x.n=\s^i (x) n\}
$$
Supposons 
 $\mu$ donné par des $(a_i)_{i\in \Z/d\Z}$, $1\leq a_i \leq n$
(qui sont les $p_\tau$ dans les notations de \ref{EL_RZ}). Un élément de
 $\M(\bar{k})$ est donné par un $\O_F\otimes 
\O_L$-réseau $M\subset N$ tel que $pM\subset VM\subset M$, et si 
$M=\bigoplus_{i\in\Z/d\Z}M(i)$ où 
$$M(i)=\{m\in
M\;|\; \forall x\in\O_F\; x.m=\s^i(x) m\}$$ alors,  l'inclusion
$VM(i)\subset M(i+1)$ induit une décomposition  
$$
M/VM=\bigoplus_{i\in\Z/d\Z} M(i)/VM(i-1)
$$
où $M(i)/VM(i-1)$ est un $\bar{k}$-espace vectoriel sur lequel
le corps résiduel de $F$ agit via le Frobenius géométrique à la
puissance $i$. La condition de Kottwitz s'écrit alors avec ces
notations  
$$
\dim_{\bar{k}} M(i)/VM(i-1)=a_i
$$ 
\\

{\bf Cas (AU) : }
Prenons les notations de l'annexe \ref{Isoc_AU}. 
Notons encore $(N,V)$ l'isocristal $(V_L,b\s)$. Le
 produit symplectique sur $V$ induit une polarisation
   $<.,.>:N\times N\ldrt
\Qp(1)$. Si $$N=\bigoplus_{i\in\Z/2d\Z} N(i)$$ est la décomposition de
$N$ sous l'action de $F$ 
comme dans le cas (AL), la
condition d'{\^e}tre  un $F$-module hermitien :
  $$\forall b\in F\; <bv,v'>=<v,b^* v'>$$ s'exprime en
disant que 
$$\forall i\neq j+d\;
<N(i),N(j)>=0
$$
$$\text{et le produit }
<.,.> \text{ est parfait sur } N(i)\times N(i+d)
$$

Avec les notations de \ref{Isoc_AU}
supposons $\mu$ donné par des entiers $(a_i)_{0\leq i\leq d-1}$.
  Un élément de $\M(b,\mu)(\bar{k})$ est alors donné par un
$\O_F\otimes\O_L$-réseau $M=\bigoplus_{i\in\Z/2d\Z} M(i)$ tel que $pM\subset
VM\subset M$, $M$ est autodual {\`a} un élément de $\Qp\*$ près et 
$$\forall\;
0\leq i\leq  d\;\;\;
\dim_{\bar{k}} (M(i)/VM(i-1))=a_i
$$
\\

{\bf Cas (C) : }
L'isocristal $(N,V)$ est polarisé via  $<.,.>:N\times N\ldrt \Qp(1)$. Comme
précédemment il y a une décomposition 
$N=\bigoplus_{i\in\Z/d\Z} N(i)$.
De plus, 
$$
\forall i\neq j \;\; <N(i),N(j)>=0
$$
$$
\text{et } <.,.> \text{ est parfait sur } N(i)\times N(i)
$$ 

 Un élément de
$\M(b,\mu)(\bar{k})$ est alors donné par un $\O_F\otimes\O_L$-réseau
$M=\bigoplus_{i\in\Z/d\Z} M(i)$ tel que $pM\subset VM\subset M$, qui
est 
autodual {\`a} un élément de $\Qp\*$ près et tel que 
$$
\forall 0\leq i< d \;\; \dim_k (M(i)/VM(i-1))=n
$$
\vspace{1mm}

{\bf Action de $J_b$ }
Avec les notations précédentes, le groupe $J_b$ admet la description
suivante  
$$
J_b=\{ g\in \Aut(N,V) \; | \;\forall n,n'\; <gn,gn'>=c(g) <n,n'> \; \}
$$
Pour un élément $g$ de $J_b$ vu comme un élément de $G(L)$, l'action de $g$ sur $\M(\bar{k})$
est, en termes de cristaux $M$, l'application $M\longmapsto g.M$.
\\

\subsection{Lien avec l'immeuble de $G$} \label{Lien_imm_G}

Donnons une autre description de $\M(\bar{k})$ dans l'esprit de la section 4
de 
(\citeg{Mi1}) qui le fait appara{\^\i}tre comme un sous-ensemble de
l'immeuble de $G$ sur $L$ :
$$
\M(b,\mu)(\bar{k})=\{ g\in G(L)/K\; |\; Kg^{-1} bg^\s K=K\mu(p) K\}
$$
où $K=G(W(k))$ est un sous-groupe compact hyperspécial. 
 On vérifie cette égalité en utilisant 
le lemme 7.4 de 
\citeg{Ko1}.  Le lien avec la définition en termes de cristaux est le
suivant :
soit
$M\subset V_L$ un cristal, $M\in\M(\bar{k})$, de stabilisateur $K$
dans $G(L)$.
 A $g\in G(L)/K$ est associé le réseaux $g.M$. Lorsque $g$ vérifie
la condition ci dessus, $g.M$ est bien un cristal de polygone de Hodge
associé {\`a} $\mu$. Dans cette description, 
 $g_0\in J_b$
  opère de la fa{\c c}on suivante :
$$ gK\longmapsto g_0 g K$$

\begin{rema}
Si l'on suppose $b$ ``decent'' au sens de \citeg{RZ}, c'est {\`a} dire 
$$
\exists s\;\; (b\s)^s= (s\nu_b)(p) \s^s
$$
où $\nu_b$ désigne le morphisme des pentes alors, dans toutes les
considérations 
précédentes on peut remplacer $L$ par $\Q_{p^s}$ l'extension non
ramifiée de degré $s$ de $\Qp$.
\end{rema}

\subsection{Orbites de $J_b$ dans $\M (\bar{k})$}

  On note $(X,\rho)$ les éléments de $\M(\bar{k})$ où $X$ désigne un
  groupe p-divisible sur $\bar{k}$ muni de son action de $\O_F$ et, dans le cas
  unitaire ou symplectique, de sa polarisation principale. La
  quasi-isogénie $\rho:\XX\ldrt X$ désigne la rigidification avec le groupe
  p-divisible $\XX$ défini sur $\bar{k}$. 

\begin{lemm}\label{OrbJb}
  Deux éléments $(X_1,\rho_1),(X_2,\rho_2)\in \M(\bar{k})$ sont dans
  la même $J_b$-orbite ssi $$X_1\simeq X_2$$ comme groupes p-divisbles
  munis de structures additionnelles. 
\end{lemm}
\begin{proof}
  Soit $g\in J_b$ tel que $(X_1,\rho_1)\simeq(X_2,\rho_2\circ
  g^{-1})$. Alors, $X_1\simeq X_2$. Réciproquement, supposons qu'il y
  ait un isomorphisme de groupes p-divisibles munis de structures additionnelles
$$
f: X_1\iso X_2
$$
Considérons la quasi-isogénie $g=\rho_1^{-1}\circ f^{-1} \circ \rho_2 :
\XX\ldrt \XX$. Par définition $g\in J_b$. Le diagramme suivant commute
\begin{diagram}
  \XX & \rTo^{\rho_1} & X_1 \\
\dTo_{g^{-1}} & & \dTo^f \\
\XX & \rTo^{\rho_2} & X_2
\end{diagram}
ce qui montre que $(X_1,\rho_1)\simeq (X_2,\rho_2\circ g^{-1})$.
\end{proof}

\subsection{Action sur les composantes irréductibles}
\begin{defi} 
Si $\M$ désigne un espace de Rapoport-Zink, 
nous noterons $\Mb$ la fibre spéciale de $\M$ c'est {\`a} dire le 
sous-schéma réduit défini sur $\bar{k}$.
\end{defi}

Rappelons (\citeg{RZ}) que les composantes irréductibles de $\Mb$
sont des variétés projectives sur $\bar{k}$ et que $\Mb$ étant localement de
type fini une composante irréductible n'intersecte qu'un nombre fini
d'autres composantes. Le groupe $J_b$ permute ces composantes.
\\

\begin{theo}\label{TheoFini} Soit $\M$ un espace de Rapoport-Zink de
type E.L ou P.E.L. non ramifié simple ou encore associé {\`a} une donnée globale 
non ramifiée.   
Il y a un nombre fini d'orbites de 
 de composantes irréductibles de $\Mb$ sous l'action de $J_b$.
\end{theo}

\begin{proof}

Du point de vue de l'immeuble de $G$ sur $L$ ce théorème est un
théorème de finitude pour l'action de $J_b$ sur l'immeuble qui
généralise la proposition 2.18 de \citeg{RZ}.

D'après les remarques du début, un ouvert quasicompact de $\Mb$
intersecte un nombre fini de composantes irréductibles et réciproquement une union
finie de composantes irréductibles est quasicompacte. Le théorème est donc équivalent
{\`a} l'énoncé : il existe un ouvert quasicompact $U$ de $\Mb$ tel que
$\Mb=J_b.U$.

On vérifie en utilisant la section \ref{Dec_glob} que l'on peut se
ramener au cas des espaces associés {\`a} une donnée locale non
ramifiée simple. Nous supposerons donc que $\M$ est de ce type l{\`a}. 

%Reprenons les notations précédentes qui décrivent les éléments de $\Mb
%(\bar{k})$ en termes de réseaux $M$ dans $V\otimes L$. 

Commen{\c c}ons par rappeler le 

\begin{enonce}{Critère de quasicompacité}[corollaire 2.31 de \citeg{RZ}] 
Un ouvert $U$ de $\Mb$ est quasicompact ssi
$\exists M\in \Mb (\bar{k}) \;\exists N\in \N \;$
$$
U(\bar{k})\subset \{M'\in \M(\bar{k}), p^{N}M\subset M'\subset p^{-N} M\} 
$$
de fa{\c c}on équivalente : 
 $\exists N$, la quasi-isogénie universelle $\rho^{univ}$ soit
telle que sur $U$ $\; p^N\rho^{univ}$ et $p^N(\rho^{univ})^{-1}$ soient
des isogénies.
\end{enonce}

La démonstration consiste maintenant {\`a} exhiber un nombre fini
d'orbites dans $\Mb(\bar{k})$  (ou, dans le cas unitaire et
symplectique, dans un ensebmle plus gros : confère les définitions
suivant 
\ref{def_quasi}) 
sous l'action de $J_b$ vérifiant : il existe une constante $c$ telle
que tout point de $\Mb(\bar{k})$ soit {\`a} distance inférieure {\`a} $c$ d'un
point d'une de ces orbites.

Étant donné qu'il n'y a qu'un nombre fini de polygones de Hodge
possibles pour un groupe p-divisible muni de structures additionnelles
isogène à $\XX$, l'énoncé précédent est équivalent au même énoncé sans
la condition de Kottwitz sur le polygone de Hodge arithmétique dans la
définition de $\M(\bar{k})$. Nous ignorerons donc cette
condition. Cela signifie  que nous allons démontrer le théorème
annoncé non pas sur 
l'espace $\M$ mais sur l'espace plus gros défini de façon analogue,
sans la condition de Kottwitz sur le polygone de Hodge. 
 Nous noterons
encore $\M$ cet espace. Le lemme \ref{OrbJb} est bien sûr encore
valable pour cet espace. 
\\

Commençons par rappeler la définition suivante : 

\begin{defi}[\cite{Zink1} section 3]
  Soit $S$ un schéma de caractéristique $p$ 
et $X$ un groupe p-divisible sur $S$. Soit $\l$
  un nombre rationnel positif. On dit que $X$ est divisible par
  rapport à la pente $\l$ s'il existe des entiers naturels $r,s$ tels
  que $\dpt{\l=\frac{r}{s}}$ et tels que la quasi-isogénie 
$$
F^s p^{-r} : X \ldrt X^{(p^s)}
$$
soit une isogénie. 
\end{defi}

Nous allons utiliser le résultat suivant de Zink  : 
\begin{prop}[\citeg{Zink1} proposition 12]\label{leZi}
  Soit $K$ un corps parfait de caractéristique $p$ et $X/K$ un groupe p-divisible
de pentes $\l_1>\dots>\l_r$. 
 Il existe
  une constante $c$ ne dépendant que de la hauteur de $X$ et un groupe
  p-divisible $Y/K$ muni d'une isogénie avec $X$ de degré plus petit
que $c$ tel que $Y$ vérifie la propriété suivante : $Y$ est 
muni d'une filtration par des sous groupes p-divisibles 
$$
(0)=Y_0\subset Y_1 \subset \dots\subset Y_r=Y
$$
telle que pour tout $i$,
$Y_i$ soit divisible par rapport à la pente $\l_i$ et
$ Y_i/Y_{i-1}$ soit de pente $\l_i$.
De plus, une telle filtration se scinde. 
\end{prop}

\begin{proof}[Rappels sur la démonstration]
Rappelons que dans la démonstration, si $(M,V)$ est le cristal de $X$ et
si $U=p^{-r_1} V^{s_1}$, on pose
$$
M'=M+UM+\dots+U^{h-1} M
$$
Le point crucial de la démonstration consiste {\`a} démontrer que 
$M'$ est stable par $U$.
Une fois cela démontré 
 on peut donc décomposer  $M'$ en partie
bijective et topologiquement  nilpotente pour l'action de $U$ :
 $$M'=M'_{bij}\oplus M'_{nil}$$
Rappliquant le m{\^e}me procédé {\`a} $M'_{nil}$ pour
$U=p^{-{r_2}}V^{s_2}$ on obtient ainsi par récurrence la filtration
souhaitée. 

De plus, {\`a} la première étape $M\subset M'\subset p^{-r_1 (h-1)}
\subset p^{-h(h-1)} M$ car $r_1\leq s_1 \leq h$. Répétant ce type
d'inégalité {\`a} chaque étape, le nombre de pentes étant inférieur {\`a} $h$,
l'existence du $c$ s'en déduit.
\end{proof}

\begin{enonce}{Fait}\label{Fait_divisible}
Si le corps $K$ est algébriquement clos
et si les pentes $\l_i$ sont fixées, 
 tous les groupes $Y$ intervenant dans
l'énoncé du lemme précédent 
sont isomorphes. Si $\l_1,\dots,\l_r$ désignent les pentes du cristal 
$(M,V)$ associé et $(m_i)_{1\leq i \leq r}$ les multiplicités associées  
$$M \simeq \oplus_{i=1}^r \left (\oplus_{1}^{m_i} M^{\l_i}\right )$$
où si $\l=\frac{r}{s}, r\wedge s=1$ 
$$
M^\l=W(K)[F,V]/(V^s-p^r,FV-p)
$$
\end{enonce}

Revenons à la démonstration du théorème. 
Expliquons comment démontrer le théorème dans le cas $(AL)$ lorsque
$F=\Qp$, c'est à dire lorsqu'il n'y a pas de structure additionnelle.

 Appliquons le lemme de Zink à $\XX$. Soit donc $Y_0/\bar{k}$
un groupe p-divisible tel qu'il est donné par le lemme de Zink avec
une isogénie : $f_0:\XX\ldrt Y_0$. Soit $\widetilde{\mathcal{O}}$
la $J_b$-orbite de $(Y_0,f_0)$ dans $\M (\bar{k})$. 

 Soit maintenant $(X,\rho)\in \M(\bar{k})$. Appliquons le lemme de Zink
à $X$: soit $Y/\bar{k}$ et $f:X\ldrt Y$ l'isogénie de degré inférieur
à $c$ (où $c$ ne dépend pas de $X$ mais seulement de $\XX$ puisqu'il
ne dépend que de la hauteur du groupe p-divisible). Considérons le point $(Y,f\circ \rho)\in \M (\bar{k})$. Il est clair que
ce point est à distance inférieure à $c$ de $(X,\rho)$. De plus, 
d'après le lemme \ref{OrbJb} et \ref{Fait_divisible}, 
$(Y,f\circ \rho)\in \widetilde{\O}$. 

 On en déduit donc le théorème lorsqu'il
n'y a pas de structures additionnelles.
%Une variante immédiate dans la
%démonstration du lemme  précédent consistant {\`a} remplacer cristaux
%par $\O_F$-cristaux permet également de démontrer le théorème dans le cas
%(AL) lorsque $\mu$ est défini sur $\Qp$. 
\\

Démontrons maintenant le cas (AL) en général, c'est à dire lorsque
l'on ajoute l'action de $\O_F$. Pour cela démontrons des variantes de
\ref{leZi} et \ref{Fait_divisible}. 

\begin{lemm}\label{leZi2}
   Soit $K$ un corps parfait de caractéristique $p$ et $X/K$
   un groupe p-divisible   muni d'une action de $\O_F$
de pentes $\l_1>\dots>\l_r$. 
 Il existe
  une constante $c$ ne dépendant que de la hauteur de $X$,  un groupe
  p-divisible $Y/K$ muni d'une action de $\O_F$, et d'une isogénie
  $\O_F$ équivariante avec $X$ de degré plus petit
que $c$, tel que $Y$ vérifie la propriété suivante : $Y$ est 
muni d'une $\O_F$-filtration par des sous groupes p-divisibles munis d'une
   action de $\O_F$ 
$$
(0)=Y_0\subset Y_1 \subset \dots\subset Y_r=Y
$$
telle que pour tout $i$,
$Y_i$ soit divisible par rapport à la pente $\l_i$ et
$ Y_i/Y_{i-1}$ soit de pente $\l_i$. 
De plus, une telle filtration se scinde.  
\end{lemm}
\begin{proof}
  Soit $(M,V)$ le cristal de $X$ muni de son action de $\O_F$, c'est
  à dire d'un morphisme $\O_F\ldrt\End ( M,V)$. Reprenons le rappel de
  la démonstration du lemme \ref{leZi}. L'opérateur $U$ qui y est
  défini commute à l'action de $\O_F$, puisque $V$ commute à l'action
  de $\O_F$. Le cristal $M'$ est donc stable par l'action de
  $\O_F$. La décomposition de $M'$ en parties bijectives et
  topologiquement nilpotente sous l'action de $U$ est donc stable par l'action de
  $\O_F$.
Le reste de la démonstration est identique. 
\end{proof}

\begin{lemm}
  Soient $Y_1, Y_2$ deux groupes p-divisibles sur un corps
  algébriquement clos $K$ munis d'une action de
  $\O_F$, isogènes comme groupes munis d'une action de $\O_F$, et du
  type de ceux fournis par le lemme précédent. Alors, $Y_1\simeq Y_2$
  comme groupes p-divisibles munis d'une action de $\O_F$. 
\end{lemm}
\begin{proof}
  D'après \ref{Fait_divisible}, la catégorie $\Zp$-linéaire $\mathcal{C}$ des groupes p-divisibles
  sur  $K$ de pente $\l$ divisibles par rapport à $\l$ est
  semi-simple,  d'objet simple ayant pour cristal le cristal
  $M^\l$. L'algèbre $\End ( M^\l)$ est isomorphe à $\O_{D_\l}$, l'ordre
  maximal dans une algèbre à division d'invariant $\l$ sur $\Qp$. La
  catégorie  des objets de $\mathcal{C}$ munis d'une action de $\O_F$
  est donc équivalente à celle des $\O_{D_\l}\otimes \O_F$-modules de
  type fini sans torsion (on pourra consulter le chapitre II de
  \citeg{Weil} pour 
  la théorie des $\O_{D_\l}$ modules de type fini qui est analogue à
  celle des modules de type fini sur un anneau de valuation
  discrète). L'algèbre $\O_{D_\l}\otimes \O_{F}$ est isomorphe à $M_d
  (\O_{D'_{\l'}})$ pour un entier $d$ et un nombre rationnel $\l'$, où
  $D'_{\l'}$ désigne une algèbre à division d'invariant $\l'$ sur $F$, et
  $\O_{D'_{\l'}}$ son ordre maximal. Par équivalence de Morita la
  catégorie des  $\O_{D_\l}\otimes \O_{F}$-modules est équivalente à
  celle des $\O_{D'_{\l'}}$-modules. Or, deux $\O_{D'_{\l'}}$-modules
  de type fini sans torsion $M_1,M_2$ vérifiant $M_1\otimes \Qp\simeq
  M_2 \otimes \Qp$ comme $D'_{\l'}$-modules sont isomorphes puisqu'ils
  ont même rang. D'où le résultat lorsque $Y_1$ et$Y_2$ ont une seule
  pente. Le cas général s'en déduit puisque la décomposition de
  \ref{Fait_divisible} en somme de groupes p-divisibles isoclins  est stable par l'action de
  $\O_F$.   
\end{proof}

En utilisant ces deux lemmes,
la démonstration du cas (AL) est analogue au cas déjà démontré lorsque
$F=\Qp$ : il suffit de remplacer groupe p-divisible par groupe
p-divisible muni d'une action de $\O_F$. 
\\

Passons maintenant au cas (AU). 
Commençons avant par une définition 
\begin{defi}\label{def_quasi}
  Soit $S$ un schéma et $X$ un groupe p-divisible sur $S$. Soit $N$
  un entier. On appelle polarisation $p^N$-quasi-principale
une polarisation  $\l:X\ldrt \check{X}$ telle que $p^{N} \l^{-1}$ soit
  une isogénie.  
\end{defi}
Une polarisation $p^0$-quasi-principale est donc une polarisation
principale. 

\begin{defi}
  Soit $N$ un entier. On définit le schéma formel $\M_{(p^N)}$ sur
  $\spf (\O_{\Eb})$ comme
  étant l'espace de modules des groupes p-divisibles $X$ sur un schéma
  sur lequel $p$ est localement nilpotent, munis d'une action
  de $\O_F$, d'une polarisation $p^N$-quasi-principale $\l$, et d'une
  rigidification $\rho:\XX\ldrt X$ qui est une quasi-isogénie
  compatible à l'action de $\O_F$ et qui transforme $\l$ en un multiple
  rationnel de la polarisation sur $\XX$. On ne suppose pas la
  condition de Kottwitz vérifiée. 
\end{defi}
Il est implicite dans cette définition que $\M_{(p^N)}$ est
représentable, ce qui se déduit du théorème 2.16 de \cite{RZ} comme
dans le théorème 3.25 de \cite{RZ}. Cependant, nous n'utiliserons que
l'ensemble $\M_{(p^N)} (\bar{k})$. Le schéma formel $\M_{(p^N)}$ est muni d'une action de
$J_b$ définie de façon analogue à celle sur $\M$. Il y a une immersion
$J_b$-équivariante 
$\M=\M_{(p^0)}\hookrightarrow \M_{(p^N)}$.  

Nous allons démontrer l'existence d'un entier $N$ et d'un  nombre fini
de $J_b$-orbites
$\widetilde{\O}_1,\dots,\widetilde{\O}_r$ dans $\M_{(p^N)} (\bar{k})$
vérifiant: il existe une constante $c$ telle que tout élément de $\M(
\bar{k})$ est à distance inférieure à $c$ d'une des orbite
$\widetilde{\O}_1,\dots,\widetilde{\O}_r$. Cela démontrera le théorème
dans le cas $(AU)$. 

Commençons par une nouvelle variante du lemme \ref{leZi} : 
\begin{lemm}
  Dans l'énoncé du lemme \ref{leZi2}, supposons de plus $X$ muni d'une
  polarisation principale $\l$ compatible à l'action de $\O_F$ au sens
  où l'involution de Rosati induit  l'involution $*$ sur $\O_F$. Il
  existe alors un entier $N$ ne dépendant que de la hauteur de $X$ tel
  que le groupe p-divisible $Y$ soit muni d'une polarisation
  $p^N$-quasi-principale $\l'$ compatible à l'action de $\O_F$ et à
  l'isogénie entre $X$ et $Y$. 
\end{lemm}
\begin{proof}
Soient $M_1$ le cristal de $X$ et $M_2$ le cristal de $Y$. Soit une
isogénie 
$f:Y\ldrt X$ de degré inférieur à $c$. Elle 
 induit une inclusion $M_2\ldrt M_1$. Le cristal $M_1$ est muni d'une
 polarisation 
$<.,.>:M_1\times M_1 \ldrt \Zp (1)$ où $\Zp(1)$ est le cristal
$(W(K),V)$ avec $V$ le Verschiebung sur les vecteurs de Witt. Le
produit symplectique $<.,.>$ est parfait au niveau des $W(K)$-modules
sous-jacents aux cristaux. Si $\La$ est un $W(K)$-réseau dans le
$W(K)_\Q$ espace vectoriel associé à $M_1$, nous noterons $\La^\vee$
le réseau dual par rapport à $<.,.>$. La restriction du produit $<.,.>$ à
$M_2$ donne une polarisation $\l'$ de $Y$ telle que le diagramme
suivant commute
\begin{diagram}[size=1cm]
  X & \rTo^f  &Y \\
\dTo_\l & & \dTo^{\l'} \\
\check{X} & \lTo^{\check{f}} & \check{Y}
\end{diagram}
Étant donné que le degré de $f$ est borné par $c$, c'est à dire $[M_1
:M_2] \leq c$, il existe un entier $N$ ne dépendant que de $c$ et de
la hauteur de $X$ (et donc que de la hauteur de $X$) tel que
$p^N M_1\subset M_2\subset M_1$. Cela implique que 
$$
M_2^\vee\subset p^{-N} M_1^\vee =p^{-N} M_1\subset p^{-2N} M_2 
$$
Et donc, $p^{2N} {\l'}^{-1}$ est une isogénie. 
\end{proof}

Passons maintenant à un analogue faible de \ref{Fait_divisible} :
\begin{lemm}
  Supposons $K$ algébriquement clos. Soit $N$ un entier.
  Il n'y a qu'un nombre fini de classes d'isomophismes de groupes
  p-divisibles $Y$ sur $K$ munis d'une action de $\O_F$, d'une polarisation
  $p^N$-quasi-principale compatible à l'action de $\O_F$,  possédant une
  $\O_F$-filtration comme dans le lemme précédent et  isogènes avec
  leurs structures additionnelles à un même $X$ fixé. 
\end{lemm}
\begin{proof}
  Soient $\l_1,\dots,\l_r$ les pentes de $X$ (qui sont symétriques par
  rapport à $1/2$). Soit l'anneau $R=\prod_{i}
  \O_{D_{\l_i}}$. Cet anneau est Morita-équivalent à $\End (X)$. Il y
  a un isomorphisme $\End (\check{X})\simeq \End (X)^{opp}$ et la
  polarisation principale $\l$ induit donc une involution sur $\End
  (X)$. Celle-ci induit une involution sur $R$ permutant les facteurs
  $\O_{D_{\l}}$ et $\O_{D_{1-\l}}$. Soit $A$ l'anneau $R\otimes_{\Zp}
  \O_F$ muni de l'involution $\#$ déduite de celle sur $R$ par
  tensorisation avec celle sur $\O_F$. Considérons la catégorie
 $\Zp$-linéaire $\mathcal{C}$ des groupes p-divisibles $Y$ sur $K$, de pentes
  incluses dans l'ensemble $\{\l_1,\dots,\l_r\}$, munis d'une
  action de $\O_F$ et possédant une $\O_F$ filtration comme dans
  le lemme \ref{leZi}. Cette catégorie est semi-simple, équivalente à
  la catégorie des $A$-modules à gauche de type fini sans torsion. 
 Si $Y$ est comme dans l'énoncé, le $A$ module associé est muni d'une
  application symplectique $\#$-hermitienne associée à sa polarisation
$$
<.,.>: M\times M \ldrt \Zp
$$
La condition d'être $p^N$-quasi-principale signifie que 
$$
M\subset M^\vee\subset p^{-N} M^\vee \subset M\otimes \Qp
$$
Le résultat est donc une conséquence du lemme qui suit. 
\end{proof}

\begin{lemm}\label{finit}
  Soit $A$ une $\Zp$-algèbre qui est un $\Zp$-module libre de rang
  fini, $*$ une involution de $A$, et $B$ l'anneau $A_{\Qp}$ muni de $*\otimes 1$. 
 Soit $\mathcal{R}$ une classe d'isomorphismes de $A$-modules qui sont
  des $\Zp$ modules libres de type fini. Si $M\in \mathcal{R}$, et
  $<.,.>:M\times M \ldrt \Zp$ est un produit symplectique
  hermitien relativement à $\#$, parfait sur le $\Qp$-espace vectoriel $M_{\Qp}$,  
 nous noterons $M^\vee \subset
  M_{\Qp}$ le réseau dual.  Soit $N\in \N$. Il existe un nombre fini de
  classes d'isomorphismes de tels $A$-modules  symplectiques
  hermitiens 
  $(M,<.,.>)$ tels que la classe d'isomorphisme de $M$ soit 
  $\mathcal{R}$ et $M\subset M^\vee\subset p^{-N} M$.
\end{lemm}

\begin{proof}
Soit $n$ un entier vérifiant $n\geq 4N+2$. 
Soit $M$ un $A$-module dont la classe d'isomorphisme est 
$\mathcal{R}$. Il existe un
nombre fini de structures de $A/p^{n} A$-module symplectique hermitien 
sur
$M/p^{n}M$. Montrons que si $<.,.>_1$ et $<.,.>_2$ sont deux
structures  symplectiques hermitiennes  sur $M$ vérifiant 
$$
\forall x,y\in M\;\;\; <x,y>_1\equiv <x,y>_2\; [p^n]
$$
alors il existe un isomorphisme $(M,<.,.>_1)\simeq (M,<.,.>_2)$, ce qui
conclura. 

Soient donc deux tels produits $<.,.>_1$ et $<.,.>_2$. 
Nous devons montrer l'existence d'un $g\in \Gl_A (M)$ tel que 
$$
\forall x,y\in M\;\;\; <gx,gy>_2 = <gx,gy>_1  
$$
sachant que, relativement à $<.,.>_1$ et $<.,.>_2$, $M^\vee \subset
p^{-N} M$. 

Il existe $u\in \End_{A_{\Qp}} ( M_{\Qp})$ vérifiant
$$
\forall x,y\in M_{\Qp} \;\; <x,y>_2=<u(x),y>_1
$$
Notons $*$ l'adjonction par rapport à $<.,.>_1$ sur  $\End_{A_{\Qp}} (
M_{\Qp})$. Désormais tous les réseaux duaux seront pris par rapport à
$<.,.>_1$. 

On a l'égalité $u^*=u$ et 
\begin{eqnarray*}
  \forall x,y\in M & \;\; & <u(x),y>\equiv <x,y> \; [p^n] \\
 & \limpl &  (u-Id)(M) \subset p^n M^\vee \subset p^{n-N} M 
\end{eqnarray*}
Cherchons un $\a\in\End_A(M)$ tel que 
$$
\forall x,y\in M\;\; <(Id+p^{
\left [\frac{n}{2}\right ]+1}\a)(x),(Id+p^{\left [ \frac{n}{2} \right
]+1}\a)(y)>_2 \equiv
<x,y>_1 \;[p^{n+1}] 
$$
Un petit calcul montre que 
si $\a$ est solution de
 $$\a^*+\a=\frac{u-Id}{p^{\left [ \frac{n}{2} \right ]+1}}$$
et vérifie $\a^*(M)\subset M$, 
 $\a$ conviendra. Or, 
\begin{eqnarray*}
  (u-Id)(\La_i)\subset p^{n-N} M &\impl & 
\left (\frac{u-Id}{p^{
\left [\frac{n}{2}\right ]+1}}
\right )(M) \subset p^{n-{
\left [\frac{n}{2}\right ]-1-N}} M \subset p^{\frac{n}{2}-1-N} M \\
& \impl & \left ( \frac{u-Id}{p^{
\left [\frac{n}{2}\right ]+1}
} \right )^* (M) \subset p^{\frac{n}{2}-1-2N} M 
\end{eqnarray*}
la dernière inclusion résultant
de
$$
\forall v\in\End_A(M)\;\; v^*(M)\subset p^{-N} M
$$ 
car $\forall x\in M\; v^*(x)\in M^\vee \subset p^{-N} M$ (l'appliquer à $\dpt{v=p^{-\left [
     \frac{n}{2} \right ]+1+N}  \left (\frac{u-Id}{p^{
\left [\frac{n}{2}\right ]+1}}
\right )} $). 

Donc, comme $n\geq 4N+3$, $\dpt{\left [\frac{n}{2}\right ]-1-2N \geq 1}$ ce qui
  implique (le $1$ étant là pour le cas $p=2$) que  
$$
\a=\frac{1}{2} \left [ \left (\frac{u-Id}{p^{
\left [\frac{n}{2}\right ]+1}
} \right )+ 
\left (\frac{u-Id}{p^{
\left [\frac{n}{2}\right ]+1}
} \right )^*\right ]\in \End_A(M)
$$
convient. De plus, $Id+p^{
\left [\frac{n}{2}\right ]+1} \a \in \Gl_A(M)$.

Quitte {\`a} modifier $<.,.>_2$ par cet automorphisme on peut maintenant
supposer que $<x,y>_1\equiv <x,y>_2 \; [p^{n+1}]$ et rappliquer la
méthode de résolution par récurrence pour obtenir une suite d'éléments
de $\Gl_B(V)$ convergeant vers un $g$ tel que $<g\bullet,g\bullet>_2=<\bullet,\bullet>_1$
(étant donné que le $g$ construit précédemment appartient à $Id+ p^{
\left [\frac{n}{2}\right ]+1} \End_A(M)$, la convergence est
assurée). 
\end{proof}

Appliquons maintenant le lemme \ref{OrbJb} (encore valable pour
$\M_{(p^N)}$). On obtient donc l'existence d'un nombre fini de $J_b$-orbites
dans $\M_{(p^N)}(\bar{k})$ tel que tout point soit à distance inférieure à $c$
d'une telle orbite. 
\\

Le cas (C) se traite de la fa{\c c}on analogue au cas précédent 
en utilisant le lemme \ref{finit}.
\end{proof}

\section{Espaces de Rapoport-Zink rigides}
\label{Esp_Rap_rig}
Nous noterons dans cette section $\M$ un espace de Rapoport-Zink
associé {\`a} une donnée globale non ramifiée $\mathcal{D}$ ou bien un
espace de Rapoport-Zink associé {\`a} une donnée locale non ramifiée simple.
Si $\M$ est associé {\`a} une donnée globale $\mathcal{D}$ nous noterons  $V$
pour $V\otimes \Qp$, $B$ pour $B\otimes \Qp$ et $\O_B$ pour
$\O_B\otimes \Zp$.
Si $\M$ est associé à une donnée locale nous noterons $\O_B$ pour
$\O_F$. 
 On fixe un réseau autodual $\La_0$ dans $V$. 
$E\subset \Qpb$ désignera soit $E_\nu$ dans le cas d'un espace associé
{\`a} une donnée globale, soit le corps réflex de la donné locale dans le
cas d'un espace associé {\`a} une donnée locale.

\begin{defi}
 $\M^{an}$ désigne la fibre générique de $\M$ sur $\Eb$
au sens des espaces 
analytiques de Berkovich.

On notera $\M^{rig}$ la fibre générique de $\M$ au sens des espaces
adiques (\cite{Hu1} et l'annexe \ref{an_ad}). 
\end{defi}

\subsection{Torseur des périodes}

Si $K|\Eb$ est une extension 
 de degré fini, {\`a} $x\in \M^{an}(K)=\M(\O_K)$ est associé un
groupe p-divisible $X=X^{univ}_x$ sur $\O_K$ d'isocristal filtré 
$$
(V_L,b\s,\breve{\pi}_1(x))
$$
où $\breve{\pi}_1:\M^{an}\drt \breve{\F}^{ad}$ est la première
composante du morphisme des  périodes (\citeg{RZ})

Le module de Tate rationnel de $X$, $V_p(X)$ est une représentation
cristalline de $G_K=\Gal(\overline{K}|K)$ munie d 'une action de $B$
et d'un produit symplectique $B$-hermitien 
$$
V_p(X)\times V_p(X)\xrig{\;\; (.,.)\;\; } \Qp(1)
$$
dans lequel $T_p(X)$ est un réseau autodual.

On a alors :
$$
V_p(X)\simeq (Fil^0\Hom(V_L,B_{cris}))^{\ph}
$$
où $V_L$ est filtré par $\breve{\pi}_1(x)$ et $B_{cris}$ de
fa{\c c}on usuelle, et
$$
V_L=\Hom_{G_K}(V_p(X),B_{cris})
$$
sur lequel $b\s$ provient du Frobenius cristallin 
$\ph$ sur $B_{cris}$ 
 et $\breve{\pi}_1(x)$ de la filtration de $B_{cris}$ (confère
\citeg{Bcris}). 

On s'intéresse alors {\`a} la variation du torseur des périodes le long de
$\M^{an}$ :
$$
\M^{an}(\overline{\Eb})\ni x \longmapsto
\text{Isom}_{\text{B-mod.symp.}}
(V_p(X^{univ}_x),V)
$$
qui est un $G$-torseur non vide gr{\^a}ce {\`a} la rigidification $\rho$.

D'après \citeg{RZ} dans le cas qui nous intéresse ($G^{der}$ est
simplement connexe ) et plus généralement gr{\^a}ce {\`a} \citeg{Wint1} et aux
travaux de Colmez et Fontaine, 
\begin{theo}
  Le torseur des périodes est constant, trivial sur $\M(\overline{\Eb})$.
\end{theo}

\begin{rema}
 Dans le cas non ramifié auquel nous nous intéressons
on peut démontrer cela de fa{\c c}on encore plus directe de la fa{\c c}on
suivante : $V_p(X)$ et $V$ sont isomorphes comme $B$-modules et
contiennent tous deux un réseau autodual $T_p(X)$ et $\La_0$. Donc
d'après le lemme 7.2 de 
\citeg{Ko3} ils sont isomorphes comme modules symplectiques.
\end{rema}

\begin{rema}
Plus précisément on a le résultat suivant : si
$\mu:{\Gm}_{/\overline{\Qp}}\ldrt G_{/\overline{\Qp}}$ et $b\in B(G)$
sont tels que $(\mu,b)$ soit faiblement admissible, la classe de
cohomologie de ce $G$-torseur est 
$$
\kappa(b)-\mu_1\in H^1(\Qp,G)\simeq
\pi_0(Z(\widehat{G})^\Gamma)^D\subset X^*(Z(\widehat{G})^\Gamma)
$$
où $\kappa:B(G)\drt  X^*(Z(\widehat{G})^\Gamma)$ et
$\mu_1=\widehat{\mu}_{|Z(\widehat{G})^\Gamma}$.

En particulier, si $b\in B(G,\mu)$ (c'est {\`a} dire  $(\mu,b)$ est faiblement
 admissible en
un sens plus fort qui tient compte des structures additionnelles), 
$\kappa(b)=\mu_1$ et la classe du $G$-torseur est triviale.
\end{rema}

\subsection{Structures de niveau}

Au groupe p-divisible universel $X^{univ}/\M$ est associé un groupe
p-divisible rigide 
$$
(X^{univ})^{an}=\limi X^{univ}[n]^{an}
$$
où $X^{univ}[n]^{an}/\M^{an}$ est étale fini de degré $p^{n ht(\XX)}$. 

$T_p(X^{univ}):=(X[n]^{an})_{n\in \N}$ est donc un $\Zp$ faisceau
analytique étale localement constant sur des rev{\^e}tements étales finis,
muni d'une action de $\O_B$ et d'une polarisation 
$$
T_p(X^{univ})\times T_p(X^{univ})\ldrt \Zp(1)
$$
où $\Zp(1)$
 est le $\Zp$ faisceau analytique $(\mu_{p^n})_{n\in \N}$. Cet
 accouplement de faisceaux est parfait au sens où il induit un
 isomorphisme 
$$
T_p(X^{univ})^\vee \iso T_p(X^{univ})(-1)
$$
D'après la trivialité du torseur des périodes $\forall
x\in\M^{an}(\overline{\Eb})$ il existe un isomorphisme de $\O_B$-modules
symplectiques 
$$
\La_0\iso T_p(X^{univ}_x)
$$
Si $K\subset C_0$ est un sous-groupe compact ouvert, on peut donc
parler de  structure de niveau $K$ sur $T_p(X^{univ})$ au sens suivant
:
\begin{defi}
  Fixons un $x\in\M(\overline{\Eb})$ dans chaque composante connexe de
  $\M^{an}$. Une structure de niveau $K$ sur $T_p(X^{univ})$ est un
  isomorphisme $\forall x$ de $\O_B$-modules symplectiques 
$$
\eta:\La_0\iso T_p(X^{univ}_x)
$$
tel que via $\eta$ l'action de $\pi_1^{an}(\M^{an},x)$ se fasse {\`a}
travers $K$ agissant sur $\La_0$. Un $\eta$ étant considéré modulo
composition par un élément de $K$. On notera $\bar{\eta}$ sa classe.
\end{defi}

Bien s{\^u}r cette définition ne dépend pas du choix des
$x\in\M(\overline{\Eb})$.

Une définition équivalente est la suivante :

La classe d'isomorphisme de $T_p(X^{univ})$ comme $C_0$ torseur est
donné par un élément de 
$$\check{H}^1(\M^{an},C_0) := \limp \check{H}^1
(\M^{an},C_0\otimes \Z/p^n\Z)$$ 
Dire qu'un système local est trivial
équivaut {\`a} dire que cette classe est un cobord. De la m{\^e}me fa{\c c}on, une
trivialisation modulo $K$ d'un élément de
$$\check{Z}^1(\M^{an},C_0)=\limp \check{Z}^1(\M^{an},C_0\otimes
\Z/p^n\Z)$$ 
est un élément  $g\in \check{C}^0 (\M^{an},C_0)$ modulo
$\check{C}^0(\M^{an},K)$ tel que $$c=\text{cobord}(
g) \text{ modulo } K$$

Dit d'une autre fa{\c c}on, $\eta :\La_0\iso T_p(X^{univ}) [K]$ est donné
par une famille compatible 
$$
\eta_n:\La_0/p^n\La_0 \iso X[n]^{an} \; [K]
$$
où $\forall n \;\; \exists \{U_i\}_i$ un recouvrement étale de
$\M^{an}$ et des isomorphismes 
$$
(\eta_n)_i:\La_0/p^n\La_0 \iso  X[n]^{an}_{|U_i}
$$
tels que sur $U_i\cap U_j$, 
$$
(\eta_n)_{i|U_i\cap U_j} \circ (\eta_n)_{|U_i\cap U_j}^{-1}\in
K(U_i\cap U_j)=K^{\pi_0 (U_i\cap U_j)}
$$

$\La_0$ étant fixé, pour $K\subset C_0$ compact ouvert on définit une
tour $(\M_K)_K$ d'espaces analytiques sur $\Eb$ 
munis de morphismes de transition étales finis pour $K\subset K'$
$$
\Pi_{K,K'}:\M_{K'}\ldrt \M_K
$$
d'oubli de la structure de niveau. $\Pi_{K,K'}$ est
galoisien de groupe $K/K'$ si $K'\lhd K$.

\begin{defi}[\citeg{RZ}]
  $\M_K$ classifie les structures de niveau $K$ sur $T_p(X^{univ})$.
\end{defi}

Nous noterons $\M_K^{rig}$ les espaces adiques correspondants.

\subsection{Action se $G(\Qp)$ sur les structures de niveau}

La tour $(\M_{K})_{K\subset G(\Qp)}$ est munie d'une action de $G(\Qp)$ 
décrite rapidement dans \citeg{RZ}. Donnons en une description
détaillée.  

Si $K\subset C_0$ et  $g\in G(\Qp)$ sont tels que $g^{-1}Kg\subset
C_0$ il y a un
isomorphisme 
$$
g:\M_K\iso\M_{g^{-1}Kg}
$$
défini ainsi : si $Y$ est un espace rigide quasicompact sur $\Eb$, un
élément de $\M_K^{rig}(Y)$ est donné par un $(X,\rho)\in \M(S)$ et un
$\bar{\eta}$ tel que $(X,\rho,\bar{\eta})\in \M^{rig}_K(S^{rig})$ où
$S/\spf(\O_{\Eb})$ est un schéma formel admissible de fibre générique
$Y$ et 
$$
\eta:\La_0\iso T_p(X^{rig})\;[K]
$$
Étant donné que $g^{-1}Kg\subset C_0$, le réseau $g.\La_0$ est stable
par l'action de $K$ et donc $\eta (g.\La_0)$ est stable  par l'action
du $\pi_1$ de $Y$ sur $T_p(X^{rig})$. Il 
définit donc un groupe p-divisible rigide $X'$ muni d'une
quasi-isogénie $X\xrig{f} X'$ telle que via la composée 
$$ \ph:
\La_0\otimes \Qp \xrig{\;\simeq\;\eta\otimes 1\; } V_p(X) \xrig{\simeq
  \; f_*} V_p(X')
$$
on ait $\ph^{-1}(T_p(X'))=g.\La_0$. 

$X'$ est de la forme $X^{rig}/U$ où $U$ est un groupe rigide plat fini
sur $S$ 
et $f$ est de la forme $p^nq$ où $q:X\twoheadrightarrow X/U$ et
$n\in\Z$. D'après la version relative de la théorie de Raynaud, après
un éclatement formel admissible $\widetilde{S}\drt S$, $U$ se prolonge
en un groupe plat fini $\widetilde{U}/\widetilde{S}$. Posons 
$$
\widetilde{X}'=(X\times_S \widetilde{S})/\widetilde{U}
$$
et $\widetilde{f}:X\times_S \widetilde{S} \ldrt \widetilde{X}'$ la
  quasi-isogénie $p^n \widetilde{q}$, où
  $\widetilde{q}:X\times_S \widetilde{S}
\twoheadrightarrow X\times\widetilde{S}/\widetilde{U}$, qui prolonge $f$
sur $\widetilde{S}$. Munissons $\widetilde{X}'$ de la rigidification 
$$
\widetilde{\rho}:\XX\times_k \overline{\widetilde{S}} \xrig{\;
  \rho\times 1} X_{\widetilde{S}}\times \overline{\widetilde{S}} 
\xrig{f\times 1} \widetilde{X}'\times \overline{\widetilde{S}}
$$
ce qui nous donne un $(\widetilde{X},\widetilde{\rho})\in
  \M(\widetilde{S})$. Quant {\`a} la structure de niveau
  $\widetilde{\eta}$ sur $T_p(\widetilde{X})$, 
$$
\widetilde{\eta}=\eta\circ g
$$
au sens où le diagramme suivant commute :
\begin{diagram}
g.\La_0 && \rTo^{\eta\otimes 1}& & V_p(X) \\
\uTo^{g} &&& & \dTo \\
\La_0 & \rTo^{ \widetilde{\eta}} & T_p( \widetilde{X}') & \rInto & 
V_p(\widetilde{X}') \\
\end{diagram}
La fibre générique de $\widetilde{S}$ est $Y$ et alors 
$(\widetilde{W},\widetilde{\rho},\bar{\widetilde{\eta}})\in \M_{g^{-1}K g} (
Y)$ que l'on pose comme étant égal {\`a} $g.(X,\rho,\bar{\eta})$.
\\

\begin{rema}
Rappelons  que si $x\in \M^{an}(M)$ pour $M|\Eb$ une extension
 de degré
fini, il n'est pas nécessaire d'effectuer un éclatement formel
admissible pour prolonger notre groupe p-divisible sur $\O_M$ ce qui
simplifie la définition de l'action de $G(\Qp)$ sur les points 
$\M_K(\overline{\Eb})$. 
\end{rema}

Rappelons également que $\forall K\; \M_K$ est muni de l'action de $J_b$
prolongeant celle de $\M^{an}$,
commutant {\`a} celle de $G(\Qp)$ et compatible aux morphismes de
transition.

\begin{exem}
Avec les notations du premier premier chapitre, si $J\neq \emptyset$, si
$$
K=\prod_{i\in I} K_i \times \prod_{j\in J} K_j
$$
$$
\M_K (b,\mu)=\prod_{i\in I} \M_{K_i} (b_i,\mu_i) \times \left (
\prod_{j\in J} \M_{K_j}(b_j,\mu_j) \right )^1
$$
et l'action de $G(\Qp)$ se fait composante par composante.

Si $J=\emptyset$ et si 
$$
K=\prod_{i\in I} K_i \times K_t
$$
où $K_t$ est un sous-groupe compact ouvert de $\Zp\*$, 
$$
\M_K (b,\mu)=\prod_{in\in I} \M_{K_i} (b_i,\mu_i) \times \Qp\*/ K_t
$$
où $\Qp\*/K_t$ est l'espace de modules des structures de niveau $K_t$
sur le module de Tate $\Z_p (1)$. 
\end{exem}

\begin{exem}\label{exemtta}
Avec les notations de l'exemple précédent et en se pla{\c c}ant dans le
cadre de l'exemple 
\ref{exemAetale} (cas étale) : 
$$
\M_{K_i}(b_i,\mu_i)= G(\Qp)/K_i 
$$
De plus, $J_{b_i}=G(\Qp)$ qui agit {\`a} gauche sur cet
ensemble  et $G(\Qp)$ agit {\`a} droite via les correspondances
de Hecke ensemblistes classiques. 
\end{exem}

\subsection{Quelques propriétés des espaces de Rapoport-Zink rigides}

\begin{lemm}\label{Sep} Le schéma formel 
$\M$ est  séparé sur $\spf (\O_{\Eb})$. 
\end{lemm}
\begin{proof}
Il s'agit de montrer que le sous-schéma formel localement fermé 
\begin{diagram}
\M & \rInto^{\Delta\;\;} & \M\times \M
\end{diagram}
est une immersion fermée. Soit $(X^{univ},\rho^{univ})$ le
groupe $p$-divisible universel muni de sa rigidification
$\rho^{univ}$. Considérons les deux projections 
\begin{diagram}
\M\times \M & \pile{\rTo^{\; \; p_1} \\ \rTo _{\;\; p_2}} & \M
\end{diagram}
L'immersion 
$\Delta$ est alors défini comme étant le lieu défini par 
$$
p_1^* (X^{univ},\rho^{univ})\simeq p_2^* (X^{univ},\rho^{univ})
$$
ce qui est équivalent {\`a} dire que la quasi-isogénie 
$$ \a=
p_1^* \rho^{univ} \circ (p_2^* \rho^{univ})^{-1} : 
 p_2^* X^{univ} \ldrt  p_1^* X^{univ} 
$$
est un isomorphisme ou encore que $\a$ ainsi que $\a^{-1}$ sont des
isogénies. D'après la proposition 2.9 de \citeg{RZ} il s'agit d'une
condition fermée. 
\end{proof}

Pour comprendre l'énoncé qui suit il est conseillé de lire l'appendice
\ref{an_ad}. 

\begin{lemm}\label{PartProp}
Pour tout $K$, l'espace analytique $\M_K$ est lisse de bord vide. Le
morphisme des pentes  
$\breve{\pi}_1:\M^{an}\ldrt \breve{\F}^{ad}$ est étale.  

Traduit dans le langage des espaces adiques : l'espace adique
$\M_K^{rig}$ est lisse 
et partiellement propre et le morphisme  $\breve{\pi}_1$ est partiellement propre.
\end{lemm}
\begin{proof}
Montrons les assertions sur les espaces adiques qui sont équivalentes
{\`a} celles sur les espaces analytiques.
 
D'après le lemme précédent, $\M^{rig}$ est séparé. Pour voir
qu'il est partiellement propre il suffit donc de constater que les
composantes irréductibles de sa fibre spéciale $\Mb$ sont propres. 

Les morphismes $\Pi_{K,C_0}: \M_K^{rig} \drt \M^{rig}$ étant finis il en est donc
de m{\^e}me pour les $\M_K^{rig}$.

On sait déj{\`a} que le morphisme des pentes est étale au sens des espaces
adiques (proposition 5.17 de \citeg{RZ}). Les espaces adiques $\breve{\F}^{ad}$ et
$\M^{rig}$ 
étant partiellement propres, le mrophisme $\breve{\pi}_1$ est partiellement propre. 
\end{proof}

\begin{rema}
Ces propriétés peuvent également se déduire de l'existence de
l'uniformisation de variétés de Shimura par les $\M_K$. 
\end{rema}

\section{Cohomologie de $\M_K$}\label{cohomo_de_MK}
Reprenons les m{\^e}mes conventions sur $\M$ que celles énoncées au début de la section \ref{Esp_Rap_rig}.

\subsection{Lissité de l'action de $J_b$}

Rappelons (section 6 de \citeg{Berk2}) que si $X$ est un $k$-espace analytique le
groupe des automorphismes analytiques de $X$, $\Aut(X)$, est muni d'une
structure de groupe topologique et que par définition un groupe topologique $G$
agissant par automorphismes analytiques sur $X$ agit contin{\^u}ment si le 
morphisme induit $G\ldrt \Aut(X)$ est continu. 

On déduit alors de la proposition \ref{ActC} et du lemme 8.4 de
\citeg{Berk2} :

\begin{coro}\label{ActJ} Le groupe 
  $J_b$ agit contin{\^u}ment sur $\M_K$.
\end{coro}
\vspace{2mm}

\begin{defi}
Nous noterons $H^\bullet_c ({\M_K}\otimes_{\Eb}\C_p
,\Ql)$ la cohomologie {\`a} support compact
$\ell$-adique de l'espace analytique $\M_K\otimes\C_p$
 (confère l'annexe \ref{coho_l}).
\end{defi}

\begin{rema}
Le lecteur de désirant pas lire l'annexe \ref{coho_l} pourra prendre
comme définition de $H^\bullet_c ({\M_K}\otimes_{\Eb}\C_p
,\Ql)$ 
$$
\underset{U}{\limi} \underset{n}{\limp} H^\bullet_c
(U\otimes_{\Eb}\C_p,\Z/\ell^n \Z) \otimes \Ql
$$
où $U$ parcourt les ouverts relativement compacts de $\M_K$.
\end{rema}

 Cette espace de cohomologie est muni d'une action de
$J_b\times W_{E}$ de la fa{\c c}on suivante : l'action de $J_b$ est celle
induite par l'action sur $\M_K$, l'action du groupe d'inertie de
$\Gal(\overline{\Eb}|\Eb)$ est 
celle induite par action sur les coefficients $\C_p$, quant {\`a} l'action
d'un Frobenius $\s$ 
de $W_E$ elle est induite par la donnée de descente de
Rapoport-Zink (\citeg{RZ}) $$\a:\M_K\drt \M_K^{(\s)}$$
 La cohomologie
de $\M_K^{(\s)}\otimes_{\Eb} \C_p$ s'identifie {\`a} celle de 
$\M_K\otimes_{\Eb}\C_p$ via 
$1\times \s: \M_K\otimes \C_p\drt \M_K^{(\s)}\otimes \C_p$.  

Étant donné que la donnée de descente commute {\`a} l'action de $J_b$, 
$H^\bullet_c ({\M_K}\otimes_{\Eb}\C_p
,\Ql)$ est muni d'une action de $J_b\times W_E$. 
De plus, lorsque le niveau $K$ varie, le système des $ 
H^\bullet_c(\M_K\otimes \C_p,\Ql) 
$
est muni d'une action de $G(\Qp)\times J_b\times W_E$.  

\begin{rema} La donnée de descente $\a$ n'étant pas effective, il
s'agit vraiment d'une action de $W_E$ et non d'une action de
$\Gal(\overline{E}|E)$.
\end{rema}
\vspace{0mm}

\begin{defi}
Afin d'alléger les notations nous noterons désormais 
$$
H^\bullet_c(\M_K,\Ql):=H^\bullet_c(\M_K\otimes \C_p,\Ql)
$$
et 
$$ 
H^\bullet_c(\M_K,\Qlb):=H^\bullet_c(\M_K\otimes \C_p,\Ql)\otimes\Qlb
$$
\end{defi}
\vspace{0mm}

\begin{lemm} Il y a un isomorphisme 
  $H^\bullet_c(\M_K,\Ql)\simeq H^\bullet_c (\M_K^{rig},\Ql)$ 
\end{lemm}
\begin{proof}
C'est une conséquence du lemme \ref{PartProp} et du théorème 1.5 de
\citeg{Hu2}. 
\end{proof}

Combinant le théorème \ref{Cont} avec les corollaires \ref{ActJ} et
\ref{Quasi_alg} on obtient : 
 
\begin{coro} Pour tout $K$, 
  $H^\bullet_c(\M_K,\Ql)$ est un $J_b\times W_E$-module lisse pour
  l'action de $J_b$ et continu pour l'action de $W_E$.
\end{coro}

\begin{exem}\label{exemcoeta}
Pla{\c c}ons nous dans le cas d'une donnée locale étale (
exemple \ref{exemtta}). Alors,
$$
H^q_c(\M_K(b,\mu),\Qlb) = \left \{ {0\text{ si } q\neq 0 \atop 
 \mathcal{C}^\infty_c (G(\Qp)/K) \text{ si } q=0} \right.
$$
et donc, 
$$
\underset{K}{
\limi} H^0_c(\M(b,\mu),\Qlb) =  \mathcal{C}^\infty_c
(G(\Qp)) 
$$
où l'action de $J_b=G(\Qp)$ se fait par la représentation régulière gauche et
celle de $G(\Qp)$ par la représentation régulière droite. L'action de $W_E$
est l'action triviale. 
\end{exem}

\subsection{Un lemme d'induction }

Soit $\Delta=\Hom_\Z(X^*(G)_{\Qp},\Z)$. Le groupe $J_b$
étant une forme 
intérieure d'un sous-groupe de Levi $M$ de $G$, tout $\chi\in
X^*(G)_{\Qp}$ se restreint {\`a} $M$ et se transfert {\`a} $J_b$ en un
$\widetilde{\chi}\in X^*(J_b)_{\Qp}$. Notons 
$$
{J_b^1}=\bigcap_{\chi\in X^*(G)_{\Qp}} \ker |\widetilde{\chi} |
$$
où 
\begin{eqnarray*}
  |\widetilde{\chi} |:J_b &\drt \Z \\
x & \mapsto & v_p(\widetilde{\chi}(x)) 
\end{eqnarray*}
Il y a une application (\citeg{RZ} 3.52) 
\begin{eqnarray*}
\omega_{J_b}:J_b(\Qp) & \ldrt & \Delta \\
x &\mapsto & [\chi \mapsto \widetilde{\chi}(x)]
\end{eqnarray*}

Il y a également une application 
$$
\breve{\pi}_2:\M\ldrt \Delta
$$
$J_b$ équivariante. $\breve{\pi}_2$ est essentiellement la hauteur de la
rigidification $\rho$ (\citeg{RZ} 3.52).
Notons $\Delta'\subset \Delta$ l'image de $\breve{\pi}_2$ sur laquelle
$J_b$ agit avec un nombre fini d'orbites.
\\

\begin{defi} Pour tout $i$ élément de $\Delta'$ notons
$
  \M_K^{(i)}=\breve{\pi}_2^{-1}(i)
$
\end{defi}
On a donc 
$$
\M_K=\coprod_{i\in \Delta'} \M_K^{(i)}
$$
décomposition de laquelle on déduit :
\\

\begin{lemm} Il y a un isomorphisme
  $$
H^\bullet_c( \M_K,\Ql)\simeq \bigoplus_{\bar{i}\in\Delta'/J_b}
c-Ind^{J_b}_{{J_b^1}'} H^\bullet_c (\M_K^{(i)},\Ql)
$$ où
$$
{J_b^1}'=\bigcap_{\chi \in X^* (G)_{\Qp}} \ker ( |\widetilde{\chi}| )
$$
\end{lemm}

L'intér{\^e}t de ce lemme est que lorsque $b$
est une classe basique ${J_b^1}'=J_b^1$ a un centre compact. 

\begin{rema}
  L'isomorphisme ci dessus est un isomorphisme de $J_b$-modules. En
  fait, le groupe $G(\Qp)$ agit sur $\Delta$. De même, il y a une
  action de $W_E$ sur $\Delta$ triviale sur l'inertie de $W_E$ et asociée
  au fait que la donnée de descente de Rapoport-Zink transformant
  la rigidification $\rho$ en $\rho$ composée avec le Frobenius, elle
  change sa hauteur. Il y a alors un isomorphisme de $J_b\times
  G(\Qp)\times W_E$-modules
$$
\underset{K}{\limi} H^\bullet_c (\M_K,\Qp) \simeq  
 \bigoplus_{\bar{i}\in\Delta/J_b\times G(\Qp)\times W_E} 
c-Ind^{J_b}_{(J_b\times G(\Qp)\times W_E)^1} \underset{K}{\limi} H^\bullet_c (\M_K^{(i)},\Ql)
$$
où $(J_b\times G(\Qp)\times W_E)^1$ est le sous-groupe de $J_b\times
  G(\Qp)\times W_E$ agissant trivialement sur $\Delta$. 
\end{rema}

\begin{lemm} \label{lemme_annulation_Ext}
Soit $b$ une classe basique. 
  Si $r$ désigne le rang semi-simple de $J_b$,  pour
  une représentation lisse $\pi$ de $J_b$, pour $i>r$ 
$$
\Ext^{\; i}_{J_b-lisse} (H^\bullet_c(\M_K,\Ql),\pi)=0
$$
\end{lemm}
\begin{proof}  
$$
\Ext^i_{J_b-\text{lisse}}(c-Ind^J_{J_b^1} H^j_c(\M^{(0)},\Ql),\pi) = 
\Ext^i_{J_b^1-\text{lisse}} (H^j_c(\M^{(0)},\Ql),\pi)
$$
Notons $J_b^0=\bigcap_{\chi \in X^*(M)} \ker (\widetilde{\chi})$
qui est un groupe algébrique sur $\Qp$ {\`a} centre fini.
 $J_b^0(\Qp)\lhd J_b^1$ et
$J_b^1/J_b^0(\Qp)$ est un groupe compact 
 ce qui implique que si $\rho_1,\rho_2$
sont deux représentations lisses de $J_b^1$, 
$$
\Ext^i_{J_b^1-lisse}(\rho_1,\rho_2)\simeq \Ext^i_{J_b^0-lisse}
(\rho_1,\rho_2)^{
J_b^1/J_b^0}
$$
(c'est une conséquence de l'exactitude du foncteur $H^0(
J_b^1/J_b^0(\Qp),-)$ pour les coéfficients discrets).  
Dans \citeg{Sch}, il est construit pour toute représentation lisse d'un
groupe algébrique $H$ sur $\Qp$ ayant un caractère central $\chi$ sur
$Z_H^0$, la composante connexe neutre du centre,  une résolution
projective dans la catégorie des représentations lisses ayant $\chi$
comme caractère central sur $Z_H^0$ de longueur inférieure {\`a} $r$.

D'où le résultat.
\end{proof}

\subsection{Finitude}

\begin{prop} \label{Jfin} Supposons que $b$ soit la classe
basique. Alors, 
 $\forall \d\in \Delta,$
$$
H^q_c(\M_K^{(\d)},\Ql)
$$
 est un $J_b^1$-module de type fini. De m{\^e}me 
$$
H^q_c(\M_K,\Ql)
$$
est un $J_b$-module de type fini.
\end{prop}

\begin{proof}

On vérifie en utilisant le théorème  \ref{TheoFini} qu'il y a un
nombre fini d'orbites de 
composantes irréductibles de $\M^{(\d)}$ sous l'action de $J_b^1$. Si
$Z_1,\dots,Z_t$ sont des représentants de ces orbites, soit 
$$
U=Z_1^{an}\cup\dots\cup Z_t^{an}
$$
le tube au dessus des $Z_i$, un ouvert analytique dans $\M^{an}$
. Notons encore $U$ pour $\Pi_{K,C_0}^{-1}(U)\subset \M_K$. 

 Le corollaire \ref{aa} implique que $U$ s'identifie {\`a} un tube au
dessus d'un fermé dans 
l'espace analytique 
associé {\`a} une variété algébrique propre sur $\O_{\Eb}$. 
 Il résulte alors du théorème 3.3 (ii) de 
\citeg{Hu2} que 
$$
\dim_{\Qlb}(H^q_c(U,\Qlb))<+\infty
$$ 
Soit $K\subset {J_b^1}$ le stabilisateur de $U$. C'est un sous-groupe
compact ouvert. Considérons le recouvrement 
$(g.U)_{\bar{g}\in  J_b^1/K}$ de $\M_K^{(\d)}$. 

Il lui est associé une suite spectrale de cohomologie de $\check{C}$ech {\`a}
support compact 
(II.\ref{check}) concentrée en $p\leq 0,\;\dim (\Sh) \geq 
 q\geq 0$ :
$$
E^{pq}_1=\bigoplus_{
\bar{\a}\subset J_b^1/K \atop
|\a|=-p+1 } H^q_c(U(\a),\Ql) \limpl H^{p+q}_c(\M_K^{(0)},\Ql)
$$
où  $\dpt{U(\a)=\bigcap_{\bar{g}\in \a} g.U}$. 

Cette suite spectrale est $J_b^1$ équivariante où 
$$
\forall g\in J^1_b\;\; g_! :  H^q_c(U(\a),\Ql) \iso H^q_c(g.U(\a),\Ql)
$$
Si $\a\subset J_b^1/K$ notons 
$$
K_\a=\bigcap_{\bar{g}\in \a} g K g^{-1} 
$$
Les espaces $H^q_c(U(\a),\Ql)$ sont des $K_\a$-modules lisses puisque
$K_\a$ agit contin{\^u}ment sur $U(\a)$ d'après le corollaire
$\ref{ActJ}$. 

Récrivons $E^{pq}_1$ sous la forme suivante :
$$
E^{pq}_1=\bigoplus_{[\bar{\a}]\in J_b^1\bc\left (  J_b^1/K \right )^{-p+1}}
c-Ind^{J_b^1}_{K_\a} H^q_c (U(\a),\Ql)
$$
\vspace{0mm}

Montrons maintenant que 
$$
\#\{[\bar{\a}]\in J_b^1\bc \left ( J_b^1/K \right )^{-p+1} \; |\; U(\a)\neq \emptyset \}<+\infty
$$

Pour cela, remarquons que si $A$ est une union finie de composantes
irréductibles de $\Mb^{(\d)}$, $\{g\in J_b^1\; |\; g.A\cap A\neq
\emptyset\}$ est compact et que donc 
$$
\Omega=\{g\in J_b^1\; |\; g.U\cap U\neq \emptyset \}
$$
est compact et contient $K$.

Si $[\a]=[(\bar{g}_0,\dots,\bar{g}_{-p})]\in J_b^1\bc\left (  J_b^1/K\right
)^{-p+1}$ est tel que $U(\a)\neq \emptyset$ alors, $\forall i\neq
j\;\; g_i^{-1}g_j\in K\bc \Omega $ qui est fini.

Modulo l'action de $J_b^1$ {\`a} gauche sur les $(-p+1)$-uplets on peut
supposer que $g_{-p}\in K$ et donc $\forall i \; \bar{g_i}\in K\bc
\Omega$ qui est fini. D'où la finitude de l'ensemble. 
\\

On conclu donc que $E^{pq}_1$ est une somme finie d'induites compactes
de représentations de dimension finie et est donc une représentation de
type fini de $J_b^1$. D'après \citeg{BDKV} la catégorie des
$J_b^1$-modules  lisses est localement noethérienne (i.e. tout objet de type
fini est noethérien). On en déduit que $H^{p+q}_c(\M^{(\d)},\Ql)$
possède une filtration finie $F^pH^{p+q}_c(\M^{(\d)},\Ql)$ {\`a} quotients
de type fini et est donc lui m{\^e}me de type fini. 
\end{proof}

\begin{coro}\label{Extf}
 Soit $b$ la classe basique.  
  Pour toute  représentation admissible $\pi$ de $J_b$,
  pour tous $K,p,q$
$$
\dim_{\Qlb} \text{Ext}^{\; p}_{J_b\text{-lisse}} \left ( H^q_c
  (\M_K,\Qlb),\pi \right ) \; <+\infty
$$
\end{coro}
\begin{proof}
 La réciprocité de Frobenius donne un isomorphisme :
$$
\text{Ext}^{\; p}_{J_b\text{-lisse}} \left ( H^q_c
  (\M_K,\Qlb),\pi \right ) \simeq \bigoplus_{\bar{\d}\in \Delta/J_b } 
\text{Ext}^{\; p}_{J_b^1\text{-lisse}} \left ( H^q_c
  (\M_K^{(\d)},\Qlb),\pi_{|J_b^1} \right )
$$
La représentation 
 $\pi_{|J_b^1}$ est encore 
admissible. On conclut
alors gr{\^a}ce {\`a} la proposition précédente et le lemme qui suit en posant
$\dpt{H=\bigcap_{\chi\in X^* (G) } \ker ( \chi )}$ car 
$$
\text{Ext}^{\; p}_{J_b^1\text{-lisse}} \left ( H^q_c
  (\M_K^{(\d)},\Qlb),\pi_{|J_b^1} \right )
=\text{Ext}^{\; p}_{H(\Qp)\text{-lisse}} \left ( H^q_c
  (\M_K^{(\d)},\Qlb),\pi_{|J_b^1} \right )^{J_b^1/H(\Qp)}
$$
\end{proof}

\begin{lemm}
Soit $H$ un groupe $p$-adique semi-simple, $\pi_1$ une représentation
lisse
de type fini de $H$ et $\pi_2$ une représentation admissible. Alors, 
$$
\forall i \;\; \dim (\text{Ext}^i_{H\text{-lisse}}(\pi_1,\pi_2))<\infty
$$  
\end{lemm}
\begin{proof}
$\pi_1$ étant de type fini il existe une surjection 
$$ \a: 
\bigoplus_{i\in I} c-ind^H_{K_i}\rho_i \twoheadrightarrow \pi_1
$$
où $I$ est fini, les $K_i$ sont des sous-groupes compacts ouverts et
les $\rho_i$ de dimension finie. La représentation de droite étant
elle m{\^e}me de type fini, $\ker (\a)$ est de type fini (locale noethérianité 
de la catégorie des représentations lisses de $H$ (\cite{BDKV})) et on
peut rappliquer le processus pour construire ainsi par récurrence une
résolution projective $P^\bullet \drt \pi_1$ de $\pi_1$ telle que
$\forall q \;P^q$ soit une représentation de type finie de $H$. On
conclut aussitôt puisque si $\rho$ est de type fini et $\pi_2$
admissible alors $\Hom_H(\rho,\pi_2)$ est de dimension finie. En
effet, un
morphisme $H$ équivariant de $\rho$ dans $\pi_2$ est déterminé par
l'image d'un nombre fini de vecteurs de $\rho$ engendrant $\rho$. Or
une telle famille finie de vecteurs est contenue dans $\rho^V$ pour un
sous-groupe compact ouvert $V$ et donc d'image dans $\pi_2^V$ qui est de dimension finie. 
\end{proof}

\begin{rema}
 Dans \citeg{Sch} ce lemme est démontré d'une autre
fa{\c c}on mais en supposant de plus que $\pi_1$ est admissible. Cette
restriction est donc inutile.
\end{rema}

\begin{coro}
Soit $b$ la classe basique.
Supposons $J_b$ anisotrope modulo son centre
. Alors, 
$$\forall \d\in \Delta \; \forall q\;\;\;\; 
\dim_{\Qlb} H^q_c(\M_K^{(\delta)},\Qlb) <+\infty
$$
\end{coro}

\begin{rema}
Supposons la classe $b$ basique.
On ne peut espérer mieux que la proposition \ref{Jfin} dans le cas où
$J_b$ n'est pas anisotrope modulo son centre. En effet, on pourrait
penser que $H^\bullet_c(\M_K^{(\delta)},\Qlb)$ est un $J_b^1$-module de
longueur finie. Mais cela est faux comme le montre l'exemple étale
(exemple \ref{exemcoeta}).   
\end{rema}

\chapter{Uniformisation des variétés de Shimura de type P.E.L.}
\label{uniformi_PEL}

L'uniformisation des variétés de Shimura se décompose en plusieurs
étapes :
\begin{itemize}
\item stratification de la fibre spéciale selon la classe d'isogénie du
groupe $p$-divisible muni de structures additionnelles
\item raffinement d'une strate en classes d'isogénies $\phi$ de
variétés abéliennes munies de structures additionnelles 
\item classes d'isomorphismes dans une classe d'isogénie :
\begin{itemize}
\item en $p$ gr{\^a}ce {\`a} l'espace $\M_{K_p}$
\item hors $p$ ``tout est étale'' et il s'agit d'une description de
réseaux munis de structures de niveau  c'est {\`a} dire d'éléments de $G(\A_f^p)/K_p$ 
\end{itemize} 
\end{itemize}

Certaines de ces étapes apparaissent dans le comptage des points des
variétés de Shimura sur les corps finis. 

\section{Stratification par la classe d'isogénie du groupe p-divisible}

Reprenons les notations globales du premier chapitre.

Soit $\Sb_{K^p}=S_{K^p}\times_{\O_{E_\nu}} k$ où $k$ désigne ici le corps
résiduel de $E_\nu$. Soit $\mathcal{A}/\Sb$ le schéma abélien universel.

Gr{\^a}ce au théorème de spécialisation des cristaux de Grothendieck
généralisé dans \citeg{RR}, le schéma $\Sb$ est stratifiée par le polygone de
Newton du cristal muni de structures additionnelles 
$R^1 f_{cris *} \O_{\mathcal{A}}$. 

Plus précisément, si $b\in B(G_{\Qp},\mu_{\Qpb})$ 
$$
\Sb (b) =\{ x\in \Sb(\kb)\; |\;
(H_{1,cris}(\mathcal{A}_x,W(\kb)_\Q),F)
\simeq (V_{\Qp}\otimes W(\kb)_\Q, b\otimes \s) \}
$$
(l'isomorphisme étant pris au sens des isocristaux munis de structures
additionnelles). L'ensemble $\Sb(b)$ est l'ensemble des points géométriques 
sous-jacents {\`a} un sous-schéma localement fermé réduit de $\Sb$ encore noté
$\Sb(b)$. On a la stratification 
$$
\Sb=\coprod_{b\in B(G_{\Qp},\mu_{\Qpb})} \Sb(b)
$$
et si $P$ est un polygone fixé 
$$
\coprod_{Newt(b)\geq P} \Sb(b) 
$$
est fermé dans $\Sb$.
\\

Deux strates se distinguent particulièrement :
\begin{itemize}
\item la strate basique qui est fermée et est associée {\`a} la
  classe basique de $B(G_{\Qp},\mu_{\Qpb})$. Son polygone de Newton est maximal
\item la strate $\mu$-ordinaire associée {\`a} l'image via $B(T)\ldrt
  B(G_{\Qp})$ de $\widehat{\mu}_{|\widehat{T}^\Gamma}\in B(T)$ de polygone
  de Newton minimal dans $B(G_{\Qp},\mu_{\Qpb})$. Elle est ouverte et dense dans
  $\Sb$ (\citeg{Wed})
\end{itemize}

Citons :
\begin{conj} On a l'équivalence suivante 
  $\Sb(b)\neq \emptyset \ssi b\in B(G,\mu)$ 
\end{conj}
(La conjecture est pour l'implication de la droite vers la gauche,
l'implication de la gauche vers la droite étant la généralisation du
théorème de Mazur). Il semble que des travaux récents de Rapoport et
Kottwitz permettent de démontrer cette conjecture. Nous montrerons de
tout fa{\c c}on plus tard que la strate basique est non vide, ce qui
est suffisant pour les applications que nous avons en vue dans cette
thèse. 

\section{Uniformisation de Rapoport-Zink}\label{Unif}

Fixons un point base $x\in \Sb(\kb)$. Le théorème de Serre Tate permet
d'uniformiser $S$ ``verticalement'' au dessus de $x$ c'est {\`a} dire  sur un
voisinage formel de $x$ :
\begin{center}
$
S^\wedge_{\{x\}} \simeq$: \{ déformations par isomorphismes du groupe 
p-divisible muni de structures additionnelles $\mathcal{A}_x[p^\infty]$ \}.
\end{center}
Nous voulons uniformiser une plus grande partie de la strate qu'un
point, c'est pourquoi on déforme le groupe p-divisible non plus par
isomorphismes mais par quasi-isogénies gr{\^a}ce {\`a} l'espace $\M(\mathcal{D}_{\Qp},b)$.

Soit $\phi$ la classe d'isogénie du triplet $(\mathcal{A}_x,\l,\iota)$
et 
$$
I^\phi=\Aut(\mathcal{A}_x,\l,\iota)
$$
un groupe réductif sur $\Q$. Notons $b\in B(G_{\Qp},\mu_{\Qpb})$ la
classe associée {\`a} $\phi$ et $\M=\M(\mathcal{D}_{\Qp},b)$ l'espace de Rapoport-Zink associé.

A $\phi$ est alors associé un ensemble $\St$ 
 de sous-schémas fermés
projectifs de $\Sb(b)_{\bar{k}}$  (confère le théorème 6.23 de
 \cite{RZ} pour la définition de $\St$ où cet ensemble est noté $\mathcal{T}$) tel que 
$$
\St(\kb)=\{ y\in \Sb(b)(\kb)\; |\; \text{ la classe d'isogénie de }
(\mathcal{A}_y,\l,\iota)\in\phi\; \}
$$

Le théorème 6.23 de \cite{RZ} affirme alors qu'il y a 
 un isomorphisme de schémas formels sur $\spf(\O_{\Eb})$ : 
$$\Theta: 
I^\phi(\Q) \backslash \left ( \M \times G(\A_f^p)/K^p \right )
\iso  \left ( S_{K^p}\otimes_{\O_E} \O_{\Eb} \right )^\wedge_{/\St}
$$
où l'application d'uniformisation est définie en tordant le triplet 
$
(\mathcal{A}_x,\l,\iota)
$
par une déformation par quasi-isogénie de
$(\mathcal{A}_x[p^\infty],\l,\iota)$ et en tordant la structure de
niveau $\overline{\eta^p}$ par un élément de $G(\A_f^p)/K^p$.

L'action de $I^\phi (\Q)$ sur $\M$ se fait {\`a} travers $J_b$ par action
de $I^\phi(\Q)$ par quasi-isogénies sur
$(\mathcal{A}_x[p^\infty],\l,\iota)$ ce qui donne une injection 
$$
I^\phi (\Q)\hookrightarrow J_b
$$
 et sur $G(\A_f^p)$ via l'action
d'un élément de $I^\phi (\Q)$ sur
le module de Tate $H_1(\mathcal{A}_x,\A_f^p)$ ce qui donne une
injection 
$$
I^\phi(\Q) \hookrightarrow G(\A_f^p)
$$

Lorsque $K^p$ varie les différents isomorphismes d'uniformisation sont
compatibles et commutent {\`a} l'action de $G(\A_f)$ (et donc aux
correspondances de Hecke hors $p$) en un sens évident. 

\subsection{Quelques propriétés de l'isomorphisme d'uniformisation}\label{recrit}

Soit $(y_i)_{i\in I}$ un ensemble de représentants des orbites de
$I^\phi (\Q)$ dans $G(\A_f^p)/K^p$, c'est {\`a} dire 
$$
G(\A_f^p)/K^p = \coprod_{i\in I} I^\phi (\Q). y_i
$$
Pour tout $i\in I$ notons $\GG_i=Stab_{I^\phi (\Q)} y_i$ que l'on
verra comme un sous-groupe de $J_b$. 

L'isomorphisme d'uniformisation se récrit alors : 
$$
\coprod_{i\in I} \GG_i\bc \M \iso \left ( S_{K^p}\otimes
\O_{\Eb}\right )^{\widehat{\;}}_{/\St}
$$

\begin{defi}
  Appelons $\mathcal{C}$ l'ensemble des sous-groupes $\Gamma\subset J_b$ de la forme
  $I^\phi(\Q)\cap (J_b\times K^p)$ pour $K^p\subset G(\A^p_f)$ tel
  qu'il existe un $g\in G(\A_f^p), gK^p g^{-1}$ soit suffisamment petit.
\end{defi}
Les groupes $\Gamma\in \mathcal{C}$ sont exactement les groupes $\GG_i$ 
 qui interviennent dans
les morphismes d'uniformisation ci dessus lorsque $K^p$ et
 $(y_i)_i$ varient.

\begin{lemm}\label{leG} Les éléments de $\mathcal{C}$ vérifient :
\begin{itemize}
\item
  Tout élément $\GG$ dans $\mathcal{C}$ est discret sans
  torsion dans $J_b$. 
\item
$\forall \Gamma_1,\Gamma_2\in\mathcal{C}\;\;$ $\Gamma_2$ et $\Gamma_2$
sont commensurables. 
\item
$\forall A\subset J_b$ fini il existe $\Gamma\in\mathcal{C}$ tel que
$\Gamma\cap A=\emptyset$. 
\end{itemize}
\end{lemm}
\begin{proof}
 Le fait que $\Gamma$ soit discret et sans torsion  est
démontré lors de la démonstration du théorème 6.23 de  \citeg{RZ}.

La commensurabilité résulte de la commensurabilité des $K^p\subset
G(\A_f^p)$ compacts ouverts.

Quant {\`a} la dernière propriété, c'est une conséquence du fait que si
$\ell\neq p, I^\phi(\Q)\hookrightarrow G(\Ql)$ et $G(\Ql)$ est séparé
au sens où pour tout sous-ensemble fini
 $A$ de $\subset G(\Ql)$ il existe un sous-groupe
ouvert de $G(\Ql)$ ne rencontrant pas $A$.
\end{proof}

Rappelons également :

\begin{prop}[\citeg{RZ}] Pour tout ouvert quasicompact $U$ dans $\M$,
  pour tout groupe $\GG$ dans $\mathcal{C}$,  l'ensemble 
 $\{\g\in \Gamma\; |\; U.\g\cap U\neq\emptyset\;\}$ est  
  fini, et $\Gamma\setminus \{Id\}$ n'a pas de point fixe dans $\M$.
\end{prop}

Le groupe
$\Gamma$ agit donc en quelque sorte de fa{\c c}on proprement discontinue et
sans points fixes sur $\M$, comme c'est le cas pour les groupes
arithmétiques uniformisant les points complexes des variétés de
Shimura et agissant sur les domaines symétriques hermitiens
$X=G(\R)/K_\infty$.

\begin{coro}\label{aa} Tout ouvert quasicompact $U$ dans $\M$ 
  est isomorphe au 
  complété formel d'une variété quasi-projective sur $\O_{\Eb}$ le long d'un
  sous-schéma fermé de sa fibre spéciale.
\end{coro}
\begin{proof}
 D'après la proposition précédente et le lemme \ref{leG}
il existe $\Gamma \in\mathcal{C}$ tel que $\forall g\in\Gamma\;\;
U.\gamma\cap U=\emptyset$. 

On peut insérer $\Gamma$ dans des $(\Gamma_i)_{i\in I}$ uniformisant
$S$. Considérons alors le composé 
$$
U\hookrightarrow \M \drt \GG\bc \M \hookrightarrow  \coprod_{i\in I}
\Gamma_i\bc \M \iso
S^\wedge_{/\St}
$$
qui induit une immersion ouverte de $U$ dans $S^\wedge_{/\St}$. L'ouvert $U$
étant quasicompact, la fibre spéciale de l'ouvert de $S^\wedge_{/\St}$
image par ce composé ne rencontre qu'un nombre fini d'éléments de
$\St$. Cet ouvert  
 s'identifie donc au complété formel de $S$ le long d'un ouvert d'une
union finie d'éléments de $\St$, i.e. un ouvert d'un fermé de
$\Sb$. Le schéma $S$ étant quasiprojectif on en déduit le résultat.
\end{proof}

\begin{coro} Pour tout  $\GG$ élément de $\mathcal{C}$ le morphisme 
   $\M\ldrt \GG\bc\M$ est étale
  (i.e. adique et formellement étale) et le morphisme d'uniformisation
  est étale.
\end{coro}
\begin{proof}
 Il faut montrer que le morphisme  $\M\drt \GG\bc\M$ est adique
puisqu'il est clairement formellement étale. Il suffit de montrer que
ce morphisme est adique en restriction {\`a} tout ouvert
quasicompact de $\M$. Soit $U$ un ouvert quasicompact de $\M$. Choisissons
un sous groupe $\Gamma'\subset \Gamma$ tel que $\Gamma'\in\mathcal{C}$ et $\forall
\gamma\in\Gamma'\setminus\{ Id \}\;\; \g.U\cap
U=\emptyset$. Décomposons l'application $U\ldrt \M/\Gamma$ en 
\begin{diagram}
  U & \rInto^j & \M & \rTo^p & \M/\Gamma' \\
& && & \dTo^\pi \\
& && & \M/\Gamma
\end{diagram}
Si l'on choisit $\Gamma'\subset \Gamma$ tel que $\Gamma$ soit de
la forme $I^\phi(\Q)\cap (J_b\times K^p)$ et $\Gamma'$ de la forme
$I^\phi(\Q)\cap {K^p}'$ pour ${K^p}'\subset K^p$, on a un diagramme
commutatif 
\begin{diagram}
  \M/\Gamma' & \rInto & (S_{{K^p}'})_{/\St}^\wedge \\
\dTo^\pi & & \dTo \\
\M/\Gamma & \rInto &  (S_{K^p})_{/\St}^\wedge
\end{diagram}
L'application de droite étant étale on en déduit que $\pi$
l'est. Le morphisme $p\circ j$ étant une immersion ouverte cela
conclu la démonstration. 
\end{proof}

\begin{lemm}\label{plgt}
  L'inclusion $I^\phi(\A_f^p)\subset G(\A_f^p)$ est un homéomorphisme
  sur son image.
\end{lemm}
\begin{proof}
 L'existence d'un isomorphisme de $B\otimes \A_f^p$-modules symplectiques  
$$
\eta:
V\otimes \A_f^p\iso H_1(\AA_x,\A_f^p)
$$
implique que la $\widehat{\Z}^p$ structure sur $G(\A_f^p)$ est la m{\^e}me
que celle induite sur $I^\phi(\A_f^p)$.
\end{proof}

\begin{rema}
 Se donner un isomorphisme $\eta$ comme dans la
démonstration précédente est plus fort que se donner $\forall \ell\neq
p$ un isomorphisme $V\otimes\Ql\iso H_1(\AA_x,\Ql)$. C'est ce qui fait
marcher le lemme précédent.
\end{rema}

\begin{lemm}\label{coco} Pour tout élément $\GG$ dans $\mathcal{C}$,
 le groupe  $\GG$ est cocompact dans $I^\phi (\Q_p)$.
\end{lemm}

\begin{proof} Le groupe réductif 
 $I^\phi$ est anisotrope modulo son centre et $I^\phi (\R)$
est compact modulo son centre. On en déduit que si $K'\subset I^\phi
(\A_f)$ est un sous-groupe compact ouvert l'ensemble
$$
I^\phi(\Q) \backslash I^\phi (\A_f) / K'
$$
est fini.

L'inclusion 
$I^\phi (\A_f)\hookrightarrow J_b\times G (\A_f^p)$ étant continue, si $K$
est un sous-groupe compact ouvert de $J_b\times G(\A_f^p)$, le groupe
$K\cap I^\phi 
(\A_f)$ est compact ouvert dans $I^\phi (\A_f)$. 

Si $\GG$ est un élément de $\mathcal{C}$, il existe un sous groupe
compact ouvert $K^p$ dans $G(\A_f^p)$
tel que $\GG=I^\phi (\Q)\cap K^p$. Si  $K_p'$ est un sous-groupe
compact ouvert de $I^\phi (\Qp)$ et si l'on pose  $K'=K'_p.(K^p\cap
I^\phi(\A_f^p))$, alors l'inclusion
$$
\GG\backslash I^\phi (\Qp) /K'_p \hookrightarrow I^\phi(\Q)\backslash
I^\phi (\A_f) /K'
$$
montre que l'ensemble de gauche est fini. 
\end{proof}

\subsection{Le cas de la strate basique }\label{str_bas}

Supposons dans cette section que $b=b_0$ est la classe basique de
$B(G_{\Qp},\mu_{\Qpb})$. 

 D'après \citeg{RZ}(6.34) on a alors :
\begin{itemize}
\item L'ensemble $\{\phi\;|\; b(\phi)=b_0\}$ est fini
\item $\forall \phi\;\; b(\phi)=b_0$, $I^\phi$ est une forme
  intérieure de $G$
\item $I^\phi(\Qp)=J_b$, $\forall \ell\neq p\;\; I^\phi(\Ql)=G(\Ql)$
  et donc $I^\phi (\A_f)=J_b\times G (\A_f^p)$ (comme groupes topologiques
  d'après \ref{plgt}) 
\item $I^\phi (\R)$ est la forme intérieure compacte modulo le centre
  de $G(\R)$
\item Le lemme \ref{coco} implique alors que $\forall \GG\in \mathcal{C}\;\; \GG$ est cocompact dans $J_b$ 
\end{itemize}

Dans le cas (A) et (C) de \citeg{Ko1} que nous considérons on a en
fait mieux :
\begin{prop}[\citeg{Ko1},\citeg{RZ}]
Le nombre de classes d'isogénie $\phi$ intervenant dans la strate
basique est égal a $|\ker^1 (\Q,G)|$.

De plus, pour toutes ces classes d'iosgénie $\phi$ les groupes $I^\phi$ sont isomorphes sur
$\Q$ et il existe des isomorphismes compatibles avec les isomorphismes 
$I^\phi (\A_f)\simeq J_b\times G(\A_f^p)$.
\end{prop}
\begin{proof}
Il résulte du lemme 6.28 de \citeg{RZ} que les classes d'isogénie de la strate
basique sur un corps fini $\mathbb{F}_{p^r}$ sont associées {\`a} un m{\^e}me
triplet $(\g_0,\g,\delta)$ de \citeg{Ko1}. Le reste de la proposition se
déduit de \citeg{Ko1}.
\end{proof}

\begin{coro}
On a alors une uniformisation de la strate basique :
$$
\coprod_{\ker^1(\Q,G)} I^\phi (\Q) \bc \left ( \M \times G(\A_f^p)/K^p\right  ) \iso 
\left (S_{K^p}\right )^\wedge_{/\Sb(b_0)}
$$
où $I^\phi (\Q)\bc G(\A_f^p)/K^p$ est fini.
\end{coro}

\section{Uniformisation rigide}

Revenons au cas d'un $b$ quelconque. 

\begin{lemm}\label{quotrigfo} Pour tout $\GG$ dans 
 $\mathcal{C}$ 
il y a un isomorphisme 
$$\GG\bc \M^{an}\iso \left ( \GG\bc
  \M
  \right )^{an} $$
\end{lemm}

\begin{proof}
 Si $U$est un ouvert quasicompact de $\M$, soit
$\GG'\lhd \GG$, $\GG'\in\mathcal{C}$ choisi tel que $\forall
\g'\in\GG'\setminus\{Id\}\;
\g'.U\cap U=\emptyset$ et tel que le morphisme $\GG'\bc \M\drt \GG\bc \M$ soit étale fini
galoisien de groupe $\GG/\GG'$.

 Notons $p':\M\drt \GG'\bc \M$, $p:\M\drt \GG\bc\M$ les projections et 
 $V=p^{-1}(p(U))=\dpt{\cup_{\g\in \GG} U.\g}$.

Étant donné que $\GG'$ est d'indice fini dans $\GG$, 
$p'(V)$ est quasicompact. Le morphisme $p'(V)\drt p(U)$ est galoisien de groupe
$\GG/\GG'$ ce qui implique que 
$$
p(U)^{an}=(\GG/\GG')\bc p'(V)^{an}
$$
De plus, $p'(V)=\GG'\bc V$ au sens où $p'(V)$ est le schéma formel
recollé des $(\g.U)_{\g\in\GG}$ le long des $\g_1.U\cap \g'\g_2.U$
pour $\g_1,\g_2\in\GG$ et $\g'\in\GG'$. Donc $p'(V)^{an}$ est l'espace
analytique recollé des domaines analytiques fermés $\g.U^{an}$ le long
des fermés analytiques $\g_1.U^{an}\cap \g'\g_2.U^{an}$. Cela
implique que $\GG'\bc p'(V)^{an}=V^{an}$.

Au final, 
$$
p(U)^{an}=(\GG/\GG')\bc p'(V)^{an}=(\GG/\GG')\bc (\GG'\bc
V^{an})=\GG\bc V^{an}
$$

Écrivant $\M$ comme union croissante de tels $U$ on obtient le
résultat.
\end{proof}

Soit maintenant $\phi$ une classe d'isogénie comme dans la section
\ref{Unif} et $\St$ 
la famille de sous-schémas fermés de $\Sb_{\bar{k}}$ associée. 
Notons $S^\wedge_{K^p}$ le complété p-adique du schéma
$S_{K^p}\times_{\O_E}\O_{\Eb}$ et $\Sh_{C_0 K^p}^{an}$ l'espace
analytique sur $\Eb$ 
associée {\`a} la variété algébrique $\Sh_{C_0 K^p}\otimes\Eb$. Il y a
une immersion 
$$
(S^\wedge_{K^p})^{an}\hookrightarrow \Sh_{K^pC_0}^{an}
$$
qui 
est un domaine analytique fermé égal {\`a} tout $
\Sh_{K^pC_0}^{an}$ lorsque $S$ est propre ce qui est le cas si $\End_B(V)$
est une algèbre {\`a} division sur $F$ (\citeg{Ko3}).

Le morphisme de spécialisation 
$$
sp: (S^\wedge_{K^p})^{an}\ldrt \Sb_{K^p}
$$
est défini sur ce domaine analytique. Si le schéma $S$ n'est pas propre les
points de $\Sh_{K^pC_0}^{an}\setminus (S^\wedge_{K^p})^{an}$ se
spécialisent sur le bord de la fibre spéciale d'une compactification
de $S$. Dit en d'autres termes, les points géométriques de 
$(S_{K^p}^\wedge)^{an}$ sont les points géométriques $x$ de $
\Sh_{K^p C_0}^{an}$ où la variété abélienne $\mathcal{A}_x$ 
muni de ses structures additionnelles a bonne réduction.

\begin{defi}
Nous noterons 
$$
\Sh_{C_0 K^p}^{an}(\phi)=\left ((S_{K^p})^\wedge_{/\St}\right )^{an}
$$
la fibre spéciale du schéma formel complété de $S$ le long de
$\St$.
\end{defi}

C'est le tube au dessus de $\St$ au sens suivant :

$\forall Z\in\St\;\; (S_{/Z}^\wedge)^{an}=sp^{-1}(Z)$ est un
domaine analytique ouvert dans $(S^\wedge)^{an}$. L'espace analytique 
$\Sh_{K^pC_0}^{an}(\phi)$ est alors l'espace analytique recollé des
ouverts $sp^{-1}(Z),Z\in\St$ le long des $sp^{-1}(Z_1)\cap
sp^{-1}(Z_2),\; Z_1,Z_2\in\St$. On a donc que
$\Sh_{C_0K^p}^{an}(\phi)$ est l'ouvert analytique de $(S^\wedge)^{an}$
union des $sp^{-1}(Z),Z\in\St$. 

\begin{rema}
En fait, on a mieux car les $Z\in\St$ sont des variétés projectives et on
en déduit que $\Sh_{C_0K^p}^{an}(\phi)$ n'est pas seulement ouvert
dans $(S^\wedge )^{an}$ mais aussi dans $\Sh_{C_0 K^p}^{an}$.
\end{rema}

\begin{rema}\label{rem_u_pas_tout}
A part pour la strate basique, 
$$
\coprod_{\phi,b(\phi)=b} \Sh_{C_0 K^p}^{an}(\phi)\subsetneq sp^{-1}(\Sb(b))
$$
Par exemple, un point de $sp^{-1}(\Sb(b))$ se spécialisant sur le
point générique d'une composante irréductible de $\Sb(b)$ n'appartient
pas au membre de gauche pour $b$ différent de la classe basique.
\end{rema}

\begin{defi}
  Si $K_p\subset C_0$, $K=K_pK^p$ nous noterons 
$$
\Sh_K^{an}(\phi)=\Pi_{C_0K^p,K}^{-1}(\Sh_{C_0K^p}^{an}(\phi))
$$
un ouvert analytique de $(\Sh_K)_\eta^{an}$.
\end{defi}
On a alors le théorème suivant combiné de l'uniformisation formelle et
du lemme \ref{quotrigfo} :  

\begin{theo}\citeg{RZ}
Pour $K=K_pK^p$ variant
il y a des isomorphisme  compatibles d'espaces analytiques sur $\Eb$ 
$$
I^\phi(\Q)\backslash \left ( \M_{K_p}\times G(\A_f^p)/K^p \right ) \iso
\Sh^{an}_K (\phi) 
$$
\end{theo}

\begin{coro}\label{Quasi_alg} Pour tout 
  $K_p$, $\M_{K_p}$ est un espace analytique quasi-algébrique
  (annexe \ref{def_quasi_alg}). 
\end{coro}

\begin{coro} Pour tout $\GG$
  dans $\mathcal{C}$ le morphisme $\M_K  \ldrt \GG\bc \M_K$ est un
  isomorphisme local d'espace analytiques. 
\end{coro}

\begin{proof}
 Soit $x\in \M^{an}$. D'après le corollaire précédent, $x$ possède une base de voisinages formée de domaines affino{\"\i}des. De plus, si $sp(x)\in Z$ une composante irréductible de $\Mb$, $Z^{an}$ est un voisinage de $x$ vérifiant 
$$
\{\g\in\GG\; | \;\g. Z^{an}\cap Z^{an} \} \text{ est fini }
$$
$x$ possède donc une base de voisinages
affino{\"\i}des $\mathcal{B}$ tels que $\forall V\in\mathcal{B}\;$
$\{\g\in \GG\; |\; \g.V\cap V\neq \emptyset \}$ soit fini. Rappelons
également que 
$\GG\setminus\{Id\}$ n'a pas de points fixes dans  $\M^{an}$.

Soit donc $V\in \mathcal{B}$ fixé et 
$$\{\g_1,\dots,\g_d\}=\{\g\in\GG\setminus \{Id \}\;|\; \g.V\cap V\neq \emptyset\}$$
Alors,
$$
\forall i\in\{1,\dots,d\} \;\; \bigcap_{W\in \mathcal{B}\; W\subset
  V } W\cap \g_i.W =\emptyset
$$
car $x$ n'est pas un point fixe de $\g_i$. Par compacité des
$W\in\mathcal{B}$, 
$$\forall i,\exists W_i\in \mathcal{B}, \;\; W_i\cap \g_i. W_i=\emptyset 
$$
Alors, $W'=\bigcap_{i=1}^d W_i$ est un domaine affino{\"\i}de voisinage de
$x$ vérifiant $\forall \g\in \GG\setminus \{Id\} \;\; W'\cap \g. W'
=\emptyset$ et $W'\hookrightarrow \GG\bc \M^{an}$.
\end{proof}

\begin{coro}\label{iso_loc_unif}
Le morphisme d'uniformisation 
$$
 \M_{K_p}\times G(\A_f^p)/K^p \ldrt \Sh_K^{an}(\phi)
$$
est un isomorphisme local.
\end{coro}

%%% Local Variables: 
%%% mode: latex
%%% TeX-master: t
%%% End: 

\part{Une suite spectrale d'Hochschild-Serre pour l'uniformisation de
Rapoport-Zink} 
\chapter[Suite spectrale de Hochschild-Serre]
{Suite spectrale de Hochschild-Serre}

\section{Suite spectrale de Cohomologie de \v{C}ech équivariante}
\label{check}

Nous adoptons les notations de l'appendice \ref{coho_l} sur la
cohomologie $\ell$-adique des espaces analytiques. 

Soit $k$ un corps valué complet non archimédien. 

\begin{defi}
  Nous dirons qu'un $k$-espace analytique est de dimension finie sur $k$
  si $\dim (X) + cd_\ell(k)<+\infty$. 
\end{defi}

\begin{prop}\label{check_equiv}
  Soit $X$ un espace analytique de dimension finie sur $k$ muni d'une
  action d'un groupe $G$ à gauche, soit $(U_i)_{i\in I}$ un recouvrement ouvert
  localement fini de $X$ stable par l'action de $G$ et soit $\F\in
  \Lpf\Xet$ muni d'une action de $G$ compatible à celle sur $X$. Il
  existe  alors un complexe borné $(\I^\bullet)$ 
de faisceaux de $\La$-modules sur
  $X_{\text{ét}}$, muni d'une action de $G$ (compatible à celle sur
  $X$)
comme complexe de
  faisceaux, 
et un complexe double 
$$
C^{p,q}=\bigoplus_{\a\subset I \atop |\a|=-p+1} \GG_c(U(\a),\I^q)
$$
où $\dpt{U(\a)=\bigcap_{i\in \a} U_i}$, 
où $C^{p,q}\ldrt C^{p,q+1}$ est induit par $\I^q\ldrt \I^{q+1}$, qui est
muni d'une action de $G$ via les morphismes 
$$
\forall g\in G\;\;\; g_!: 
\GG_c(U(\a),\I^q)\ldrt \GG_c (U(\a).g,\I^q)
$$
tel que la suite spectrale $G$ équivariante associée à la filtration
par les ligne
soit 
$$
E^{pq}_0=C^{pq}\;\;\;\; E^{pq}_1=\bigoplus_{\a\subset I\atop
  |\a|=-p+1} H^q_c(U(\a),\F) \limpl H^{p+q}_c (X,\F)
$$ 
où l'action sur $H^{p+q}_c(X,\F)$ est celle induite par l'action de
$G$ sur $\F$.
\end{prop}

\begin{proof}
 Par définition une action de $G$ sur $\F$ est la donnée
de morphismes 
$$
\forall g\in G\;\;\; g^*\F \xrig{\; \a_g\; } \F
$$
se composant de façon naturelle :
$$\forall g_1,g_2\in G \;\;
\a_{g_2}\circ g_2^*\a_{g_1}=\a_{g_1g_2}$$
 et vérifiant $\a_e=Id$. De tels morphismes sont nécessairement des
 isomorphismes et  $\a_{g^{-1}} : (g^{-1})^*\F \drt \F$ est
l'adjoint  
l'inverse de $\a_g$ via l'égalité $(g^{-1})^*\F= g_*\F$.

Si $U$ est un ouvert de $X$ et $g$ un élément de $G$, le morphisme 
$$
g_! : \R \GG_c ( U, \F) \ldrt \R \GG_c (g.U,\F)
$$
est défini de la façon suivante : $\a_{g^{-1}}$ définit par
adjonction un morphisme $u:\F\drt (g^{-1})_* \F$   et $g_!$ est la
composé 
\begin{diagram}
\R \GG_c (U,\F) & \rTo^{ \R\GG_c (u) } & \R\GG_c (U, (g^{-1})_* \F) & 
\rTo^\sim & \R\GG_c (g.U,\F)
\end{diagram}
Cela définit bien une action à gauche au sens où $(g_1 g_2)_!= g_{1!}
g_{2!}$.
\\

Si $\I^\bullet$ est un complexe d'injectifs dans $\Lf\Xet$
représentant $\mathbb{R}\pi_* (F)$, l'action de $G$ sur $\F$ induit  une
action de $G$ sur $\mathbb{R}\pi_* (\F)$ dans $\D^+(X_{\text{ét}},\La)$ et donc une
action sur le complexe $\I^\bullet$ à homotopie près au sens où
$\forall g\in G\;\;$ on a un morphisme de complexes 
$\b_g:g^* \I^\bullet \ldrt\I^\bullet$ tel que $\b_{g_1 g_2}$ soit homotope
à $\b_{g_2}\circ g_2^*\b_{g_1}$. 

Pour avoir une vraie action nous allons nous placer dans une catégorie
de faisceaux plus petite.

\begin{defi}
  Si $Y$ est un espace analytique  muni d'une action de $G$, 
  nous noterons $\LpGf_{/Y_{\text{ét}}}$, resp. $\LGf_{/Y_{\text{ét}}}
$ la catégorie des faisceaux $\La$-adiques, resp. de $\La$-modules sur
  $Y_{\text{ét}}$ munis d'une action de $G$ compatible à celle sur
  $Y$. Nous noterons $\D(Y_{\et},\La_\bullet-G)$, $\D(Y_{\et},\La-G)$ les
  catégories dérivées associées.
\end{defi}

Les catégories 
$\LpGf_{/Y_{\text{ét}}}$ et $\LGf_{/Y_{\text{ét}}}$ 
 sont des catégories abéliennes $\La$ linéaires. Si l'action de $G$
 sur $Y$ est triviale, elles coïncident avec
 $\La_\bullet[G]-\mathcal{F}\text{sc}_{/Y_{\text{ét}}}$, resp. 
$\La[G]-\mathcal{F}\text{sc}_{/Y_{\text{ét}}}$, les catégories de faisceaux de
 $\La[\GG]$-modules.  
\\ 

 \begin{defi} Notons
\begin{eqnarray*}
   \iota_\bullet:\LpGf_{/Y_{\text{ét}}} &\hookrightarrow & 
 \Lpf_{/{Y_{\et}}} \\
 \iota:\LGf_{/Y_{\text{ét}}} &\hookrightarrow  &\Lf_{/Y_{\et}} 
 \end{eqnarray*} 
les plongements canoniques, et 
\begin{eqnarray*}
  \widetilde{\pi}_* : \LpGf_{/Y_{\text{ét}}} &\ldrt & \LGf_{/Y_{\text{ét}}}
  \\
   (\F_n)_n & \longmapsto & \limp \F_n
\end{eqnarray*}
 \end{defi}

 \begin{lemm} Les catégories 
$\LpGf_{/Y_{\text{ét}}}$ et $\LGf_{/Y_{\text{ét}}}$
    possèdent suffisamment d'injectifs, les foncteurs $\iota_{\bullet}$ et $\iota$
    sont exactes et envoient suffisamment d'injectifs sur des
    injectifs.
 \end{lemm}
\begin{proof}
Les foncteurs 
  $\iota_{\bullet}$ et $\iota$ possèdent un adjoint à droite
:
$$
\mathcal{G}\longmapsto Ind^G_1 \mathcal{G}=\prod_{g\in G} g^* \mathcal{G}
$$
Soit alors $\mathcal{H}$ un objet de $\LpGf_{/Y_{\text{ét}}}$ (resp. 
de  $\LGf_{/Y_{\text{ét}}}$). 
Si $\iota(\mathcal{H})\hookrightarrow \I$ où $\I$ est un injectif
(resp. $\iota_{\bullet} \mathcal{H} \hookrightarrow \I$), on a une injection
$\mathcal{H}\hookrightarrow Ind^G_1(\I)$ et le membre de droite est un
injectif puisque $Ind^G_1$ possède un adjoint à gauche. On en déduit
que nos deux catégories possèdent suffisamment d'injectifs. 
De plus $\iota(Ind^G_1 (\I))$ (resp. $\iota_{\bullet} (Ind^G_1 (\I))$)
est un injectif et donc tout objet de l'une de nos deux catégories 
 se plonge dans un
injectif d'image un objet injectif par $\iota$
(resp. $\iota_{\bullet}$), d'où la seconde assertion.
\end{proof}

\begin{lemm}
On a une égalité 
  $$\iota\circ \mathbb{R}\widetilde{\pi}_*=\mathbb{R}\pi_* \circ \iota_{\bullet}$$
\end{lemm}
\begin{proof}
 Cela résulte de l'égalité $\iota\circ \widetilde{\pi}_*=\pi_*\circ
\iota_{\bullet}$, du lemme précédent et de
l'exactitude de $\iota$.
\end{proof}

Appliquons maintenant cela à $\F\in \LpGf\Xet$ : 
$$
\iota(\mathbb{R}\widetilde{\pi}_* \F) =\mathbb{R} \pi_* \iota_{\bullet} (\F)
$$
Il existe donc un complexe d'injectifs $\I^\bullet$ dans $\Lf\Xet$ où
$\I^q=0$ pour $q<<0$ et $\I^\bullet$ est muni d'une action de $G$ tel
que $\I^\bullet \simeq \mathbb{R}\pi_* (\F)$ comme objets de la catégorie
dérivée munis d'une action de $G$. 

Rappelons que $\mathbb{R}\GG_c=\mathbb{R}\GG_!\circ \mathbb{R}\pi_*$ (lemme \ref{Rpp})
 et que donc la
cohomologie du complexe de $\La[G]$-modules $\GG_c(X,\I^\bullet)$
calcule $H^\bullet (X_{\et},\F)$ muni de son action de $G$. 

$X$ étant de dimension finie sur $k$, $H^n_c(X,\F)=0$ pour $n>>0$ et
donc, quitte à remplacer $\I^\bullet$ par $\tau_{\leq d} \I^\bullet$
pour $d>>0$ on peut supposer que le complexe $\I^\bullet$ est borné. 
\\

\begin{lemm}
\label{resU}
  Si $f:U\ldrt X$ est étale, la cohomologie du complexe
$\GG_c(U,\I^\bullet_{|U})$ calcule 
$H^\bullet_c(U,\F_{|U})$. C'est en particulier le cas si $U$ est un
ouvert de $X$.
\end{lemm}
\begin{proof}
 Nous noterons 
\begin{eqnarray*}
  f^*: \Lpf\Xet & \ldrt & \Lpf_{/Y_{\et}} \\
(\mathcal{G}_n)_n & \longmapsto & (f^*\mathcal{G}_n)_n
\end{eqnarray*} qui est exact et
possède un adjoint à gauche $(\mathcal(\G_n))_n\longmapsto 
(f_!\mathcal{G}_n)_n$. 

$\pi_*\circ f^*=f^*\circ \pi_*$, donc $\mathbb{R}\pi_*\circ f^*=f^*\circ
\mathbb{R}\pi_*$.
\end{proof}

Nous pouvons maintenant appliquer SGAIV,XVII 6.2.10  pour obtenir des
résolutions simpliciales: 
$$
\forall\bullet \leq 0 \;\;\;\; J^{\bullet,q} =\bigoplus_{\a\subset I\atop
  |\a|=-\bullet+1 } {e_{\a}}_! e_{\a}^*\I^q \ldrt \I^q
$$
où $e_\a:U(\a)\hookrightarrow X$.  Ces résolutions simpliciales sont
celles associées au couple de foncteurs adjoints :
$$
(j_!,j^*): (\coprod_{i\in I} U_i)_{\et}^{\widetilde{}} \ldrt X_{\et}^{\widetilde{}} 
$$
où $j:\coprod_{i\in I} U_i \drt X$
(confère \cite{Illusie1} par exemple).

$J^{\bullet\bullet}$ est un complexe double muni d'une action de $G$ : 
$$
J^{p,q}=\bigoplus_{\a\subset I\atop
  |\a|=-p+1 } J^{p,q}_\a
$$
et l'action de $G$ se décompose sur les composantes en des morphismes
$$\forall g\in G\;\;\;\; g^*J^{p,q}_{g.\a} \ldrt J^{p,q}_\a$$

Le complexe double
 $$
C^{p,q}=\GG_! (X,J^{p,q})=\bigoplus_{\a\subset I\atop |\a|=-p+1} 
\GG_!(U(\a),\I^q_{|U(\a)})
$$
convient alors car d'après le lemme \ref{resU} on a 
$$
E^{pq}_1=\bigoplus_{\a\subset I\atop |\a|=-p+1} H^q_c(U(\a),\F_{|U(\a)})
$$  
\end{proof} 

On peut montrer (nous n'aurons pas besoin de cette proposition) 
  en utilisant des résolutions par des faisceaux
discrets comme dans  \citeg{Har2} :
\begin{prop}
  Sous les mêmes hypothèses que précédemment, supposons de plus $X$
  quasi-algébrique, $\F$ localement constant, $\forall i \; U_i$ est
  un ouvert distingué et l'action de $G$ sur $X$ est continue. On peut
  alors trouver un complexe double comme précédemment : 
 $$C^{p,q}= \bigoplus_{\a\subset
  I\atop |\a|=-p+1}
 C^{p,q}_\a$$ tel que $\forall \a\;\; C^{p,q}_\a$ soit un
  $\text{Stab}_G(U(\a))$-module lisse.
\end{prop}
\vspace{1mm}

\section
{Suite spectrale pour l'action d'un
groupe discret}\label{suit_disc}

Soit $X$ un espace analytique de dimension finie sur $k$, $\GG$ un
groupe discret agissant sur $X$ de telle manière que $\forall x\in
X\;\;\exists U$ un voisinage de $x$ tel que $\forall \g\in\GG\setminus
\{Id\} \;\;\g. U\cap U=\emptyset$. Le faisceau  $\GG\bc X$ sur le grand
site étale de $k$ est alors représentable par un espace analytique
obtenu par recollement de tels ouverts (utiliser 
\citeg{Berk1} 1.3.2,1.3.3).

 Soit une extension de degré fini $L|\Ql$, $V$ un
$L$-espace vectoriel de dimension finie et $\rho:\GG\ldrt \Gl(V_\rho)$ une
représentation continue pour $\GG$ muni de la topologie profinie et
$\Gl(V_\rho)$ de la topologie $\ell$-adique. 

Le morphisme 
$p:X\ldrt \GG\bc X$ étant un isomorphisme local il est associé à $\rho$
un
faisceau étale $L$-adique localement constant $\F_\rho$ 
 sur $\GG\bc X$ : si
$M\subset V_\rho$ est un $\O_L$ réseau stable par $\GG$, $\forall n\;\;
\rho_n:\GG \ldrt \Gl(M/\varpi_L^n M)$ se factorise par un sous-groupe
d'indice fini $\GG_n\lhd \GG$ et $\F_\rho=(\F_n)_n\otimes_{\O_L} L$
où $\F_n$ est trivial sur $\GG_n\bc X$ (qui est étale fini au dessus de
$\GG\bc X$) associé à la représentation 
$$\pi^{an}_1(\GG\bc X)\ldrt \GG_n\bc \GG \xrig{\;\;\bar{\rho}_n\;\;} \Gl(M/\varpi_L^n M)
$$
\vspace{0mm}
\begin{theo}\label{Hosch} 

\begin{enumerate} \item
Soit $\F$ un faisceau $L$-adique étale sur $\GG\bc X$. Il y a une
suite spectrale convergente 
$$
E^{pq}_2= \Ext^{\; p}_\GG \left ( H^{-q}_c (X, p^* \F), 1\right )
\limpl H^{-(p+q)}_c ( \GG\bc X, \F)^*
$$
\item
  Il y a une suite spectrale convergente 
$$
E^{pq}_2=\Ext^{\; p}_\GG \left ( H^{-q}_c(X,L),\check{\rho}\right ) 
\limpl H^{-(p+q)}_c(\GG\bc X,\F_\rho)^*
$$
\end{enumerate}
\end{theo}
\vspace{3mm}
\begin{proof}

Commençons par remarque que 2) se déduit de 1). En effet, il y a un
isomorphisme 
de faisceaux munis d'une action de $\GG$
$$
p^* \F_{\rho} \simeq \underline{V}_\rho
$$
où $\underline{V}_\rho$ désigne le faisceau constant muni de l'action
de $\GG$ via $\rho$. Donc 
$$
H^{-q}_c ( X,p^*\F)\simeq H^{-q}_c ( X,L) \otimes \rho
$$
comme $\GG$-modules. Si $M$ et $N$ sont deux $L[\GG]$-modules il y a une formule
d'adjonction (puisque $L$ est un corps) 
$$
\Ext^{p}_\GG ( M\otimes \rho, N)\simeq \Ext^{p}_\GG ( M,
\check{\rho}\otimes N)
$$
Qui montre donc que 1) implique 2).
\\

Passons donc maintenant à la démonstration de 1). 

La démonstration consiste à écrire un complexe de cohomologie
de $\check{C}ech$ à support compact dont la cohomologie calcule la cohomologie à
support compact de $X$ à coefficients dans $p^* \F$ 
dans la catégorie dérivée bornée des complexes de $L[\GG]$
modules,  et à identifier
l'image par le foncteur $\R\Hom_\GG (\bullet,1)$  de cet objet avec le
dual d'un
complexe de cohomologie de $\check{C}$ech à support compact dont la
cohomologie est 
la cohomologie à support compact de $\GG\bc X$ à coefficients dans $\F$. 
\\

 Soit $(U_i)_{i\in I}$ un recouvrement ouvert de $X$ par
des ouverts vérifiant 
$$
\forall i\neq j\; U_i\neq U_j\text{ et }\forall
\g\in\GG\setminus\{Id\}\;\;\; U_i.\g\cap U_i=\emptyset
$$

Nous travaillerons parfois directement avec des faisceaux $L$-adiques 
et des faisceaux de $L$-modules 
pour ne
pas alourdir les notations. Le lecteur effrayé pourra travailler avec
des faisceaux $\O_L$ adiques comme précédemment 
 et tensoriser par $L$ dès qu'il verra un
$H^\bullet_c$.

Soit $\mathcal{G}\in {\O_L}_\bullet-\mathcal{F}\text{sc}\Xet$  tel que
$\F=\mathcal{G}\otimes L$. Considérons 
$$
\R \pi_* \mathcal{G} \in \mathbb{D}^+ ( (\GG\bc X)_{\et},\O_L)
$$
et soit $\mathcal{I}^\bullet$ un complexe de faisceaux de $L$-modules
injectifs représentant $(\R\pi_* \mathcal{G})\otimes L$. Tronquons
$\mathcal{I}^\bullet$ de telle manière que si l'on note encore
$\mathcal{I}^\bullet$ le tronqué, la cohomologie du complexe $\GG_!
(\GG \bc X, \mathcal{I}^\bullet)$ calcule $H^\bullet_c ( \GG\bc X,
\F)$. 
$p$ étant étale, $p^* \mathcal{I}^\bullet$ est un complexe de
$L$-modules  injectifs qui est muni d'une action de $\GG$. 
De plus, d'après la démonstration du lemme \ref{resU} 
$$
p^* (\R\pi_* \mathcal G) = \R \tilde{\pi}_* ( p^* \F) 
$$
où rappelons que $\tilde{\pi}$ est l'équivalent de $\pi$ mais pour les
faisceaux munis d'une action de $\GG$ compatible à celle sur
$X$. Donc, la cohomologie du complexe
$\GG_! (X,p^*\mathcal{I}^\bullet)$  muni de l'action de $\GG\;$  
 déduite de celle sur $p^*\mathcal{I}^\bullet$ 
calcule
 $H^\bullet_c ( X, p^* \F)$ comme $\GG$-module.

Considérons le complexe double de cohomologie de $\check{C}$ech à
support compact 
(confère la démonstration de \ref{check_equiv}) associé à
$p^*\I^\bullet$ et au 
recouvrement $(U_i.\g)_{i\in I,\g\in \GG}$ de $X$ :
\begin{eqnarray*}
  C^{p,q} &=& \bigoplus_{\a\subset I\times \GG \atop |\a|=-p+1}
  \GG_!(U(\a),p^*\I^q) \\
\end{eqnarray*}
où $\dpt{U(\a)=\bigcap_{(i,\g)\in \a} U_i.\g}$. 

Le groupe $\GG$ opère sur $C^{\bullet\bullet}$ via les morphismes 
$$
\forall \g\in\GG\;\;\; \GG_!(U(\a),p^*\I^q)
\xrig{\;\; \g_!\;\;} \GG_!(U(\g.\a),\I^q)
$$
 où pour
un $\g$ l'application $\a\mapsto
\g.\a$ est définie grâce à l'hypothèse $\forall i\neq j \; U_i\neq U_j$.

On a alors l'isomorphisme de $\GG$-module 
$$
H^n(\text{Tot}^\oplus C^{\bullet\bullet})\simeq H^n (X,p^*\F)
$$

Posons pour un entier $p$ 
$$\mathcal{P}_p=\{\a\subset I\times\GG\;|\; |\a|=-p+1\;\}
$$ sur lequel $\GG$ agit à gauche via la seconde composante. Remarquons
maintenant : 

 {\bf Identité fondamentale  : } 
 $$
C^{p,q}= \bigoplus_{\overline{\a}\in \GG\bc \mathcal{P}_p}
  \text{c-Ind}^\GG_1 \; \GG_!(U(\a),p^*\I^q) 
$$
%Considérons le complexe double
%$$
%\Hom_\GG(C^{\bullet\bullet},1)\simeq \Hom_\GG \left ( \bigoplus_{\overline{\a}\in
% \GG\bc \mathcal{P}_\bullet} \text{c-Ind}^\GG_1 \; \GG_!(U(\a),\J^\bullet)
%,\check{\rho} \right )
%$$

Il résulte de cette identité que $\text{Tot}^\oplus
C^{\bullet\bullet}$ est un complexe de projectifs dans la catégorie
des $L[\GG]$-modules et que donc, si $F$ désigne le foncteur
contravariant exact à gauche 
$\Hom_\GG(-,1)$, 
$$
\R^n F(\text{Tot}^\oplus C^{\bullet\bullet}) \simeq
H^{-n}(F(\text{Tot}^\oplus C^{\bullet\bullet})) 
$$
où $\R^nF$ désigne l'hypercohomologie de $F$. Quant à la seconde suite
spectrale d'hypercohomologie elle donne : 
$$
E^{pq}_2=R^p F (H^{-q}(\text{Tot}^\oplus C^{\bullet\bullet}))\limpl \R^{p+q} F
(\text{Tot}^\oplus C^{\bullet\bullet}) 
$$
or, on sait déjà que 
$$
R^pF(H^q(\text{Tot}^\oplus C^{\bullet\bullet}))
=\Ext^{\;p}_\GG \left
  ( H^{-q}_c(X,L), 1 \right )
$$

Il reste donc à montrer que $H^n ( F (\text{Tot}^\oplus C^{\bullet\bullet}))
 \simeq
H^n_c(\GG\bc X,\F)$. 
\\

Le complexe $F(\text{Tot}^\oplus C^{\bullet\bullet})$ est le complexe
simple associé au complexe double $\Hom_\GG( C^{\bullet\bullet},1)$.
Montrons que $\Hom_\GG(C^{\bullet\bullet},1)$ est isomorphe au dual du
complexe double de cohomologie de $\check{C}$ech à support compact associé à
$\I^\bullet$ et au 
recouvrement $(p(U_i))_{i\in I}$ de $\GG\bc X$ qui calcule 
$H^\bullet_c (\GG\bc X,\F)$, ce qui conclura. 

Notons $\mathcal{C}_p$ un ensemble de représentants de
$\mathcal{P}_p$ modulo $\GG$. D'après l'identité fondamentale et la
réciprocité de Frobenius pour les induites compactes 
$$
\Hom_\GG (C^{p,q},1)=\left (\bigoplus_{\a\in\mathcal{C}_p} 
\GG_! (U(\a),p^*\I^q)  \right )^*
$$
C'est le complexe double associé au complexe de complexes
cosimpliciaux dual du complexe de complexes sipmliciaux 
associé aux applications 
$$
\forall \a \in \mathcal{C}_{p+1} \forall \b\in \mathcal{C}_p \;\; f_{\a,\b}:
\GG_! (U(\a),p^* \I^q) \ldrt \GG_! ( U (\b),p^* \I^q)
$$
où $f_{\a,\b}=\g_!$ s'il existe $\g\in \GG \; \g .\a\subset \b$ et
$f_{\a,\b}=0$ sinon.

Le complexe de cohomologie de $\check{C}$ech à support compact associé
à $\I^\bullet$ et 
$(p(U_i))_{i\in I}$ est :
$$
B^{p,q}=\bigoplus_{\b\subset I\atop |\b|=-p+1} \; \GG_! \left ( \bigcap_{i\in
  \b} p(U_i),\I^q\right )
$$
On vérifie en utilisant que  $\forall i\;\forall \g\in \GG\setminus\{Id\}\;\;
U_i.\g\cap U_i=\emptyset$ que si 
\begin{eqnarray*}
\xi:\mathcal{C}_p &\ldrt & \{ \b\subset
I\; |\; |\b|=-p+1\; \} \\
\{(i_0,\g_0),\dots,(i_{-p},\g_{-p})\} &\longmapsto &
\{i_0,\dots,i_{-p} \}
\end{eqnarray*}
alors, 
$$
\forall \b \subset I \;\;\;
\bigcap_{i\in \b} p(U_i)=\coprod_{\a\in \mathcal{C}_p\atop \xi(\a)=\b} p(U(\a))
$$
et que donc 
$$
B^{p,q}=\bigoplus_{\a\in \mathcal{C}_p} \GG_!(p(U(\a)),\I^q)
$$
or $p_{|U(\a)}:U(\a)\iso p(U(\a))$ et donc 
$$
\GG_!(p(U(\a)),\I^\bullet)\iso \GG_!(U(\a),p^*\I^\bullet)
$$
d'où l'isomorphisme 
$$
\Hom_\GG ( C^{\bullet,\bullet},1)\simeq B^{\bullet,\bullet}
$$
(on vérifie facilement, par exemple au niveau des complexes simpliciaux, que
les applications de bord sont les mêmes)
\end{proof}

\begin{rema}{(\it Indépendance du choix du recouvrement)}\label{ind_recouv}

 La suite
spectrale construite dans le théorème précédent semble à priori
dépendre du choix d'un recouvrement $(U_i)_{i\in I}$ de $X$ vérifiant
$\forall i\in I \;\forall \g\in \GG\setminus \{ Id \}\;\; U_i\cap \g . U_i
=\emptyset$. Si $\mathcal{U}$ désigne un tel recouvrement notons
$E^{pq}_r (\mathcal{U})$ la suite spectrale associée.  
Si $\mathcal{U}$ et $ \mathcal{V}$ sont deux
tels recouvrements on peut toujours trouver un troisième recouvrement
$\mathcal{W}$ plus fin que $\mathcal{U}$ et $\mathcal{W}$ vérifiant
ces hypothèses. 
Si de plus $\mathcal{V}$ est plus fin que $\mathcal{U}$
il a un  isomorphisme de suites spectrales
(pour $r\geq 2$) 
\begin{diagram}[size=1cm]
E^{pq}_r (\mathcal{U}) & \rTo^\sim & E^{pq}_r (\mathcal{V})
\end{diagram}
associé (il s'agit d'un morphisme naturel défini au niveau des
complexes de cohomologie de $\check{C}$ech à support compact qui
induit un isomorphisme puisque c'est le cas au niveau des
$E^{pq}_2$). Ces isomorphismes se composent de façon naturel.
On a donc un système inductif de suites spectrales $(E^{pq}_r
(\mathcal{U}))_{\mathcal{U}}$ où $\mathcal{U}$ parcourt des
recouvrements vérifiant les hypothèses ci dessus et où 
 $\mathcal{U} \leq \mathcal{V}$ si
$\mathcal{V}$ est plus fin que $\mathcal{U}$. Les morphismes de
transition sont des isomorphismes et si l'on veut une définition
canonique ne dépendant pas du choix d'un recouvrement on peut poser
$$
E^{pq}_r= \underset{\mathcal{U}}{
\limi} E^{pq}_r (\mathcal{U})
$$
\end{rema}

\begin{rema}\label{remeq}
  La suite spectrale est naturelle en $X,\GG,\rho$ au sens suivant :
  si $X',\GG',\rho'$ vérifient les mêmes hypothèses que $X,\GG,\rho$ et si 
$$
f:X\ldrt X',\;\; g:\GG\ldrt \GG'
$$
où $f$ est étale fini et $g$ un morphisme de groupe 
sont tels que $$\forall \g\in \GG \;\forall x\in X\;\;
f(x.\g)=f(x).g(\g)$$ et
$$
\rho=\rho'\circ g
$$ 
alors, il y a un morphisme induit :$\bar{f}:X/\GG\ldrt X'/\GG'$ qui
induit un morphisme $\bar{f}^*\F_{\rho'}\ldrt \F_\rho$ qui induit un
morphisme 
${}^t\bar{f}_!:H^{-(p+q)}_c(X'/\GG',\F_{\rho'})^*\ldrt
H^{-(p+q)}_c(X/\GG,\F_\rho)^*$.

Il y a alors un morphisme de suite spectrales $E^{pq}_r(X',\GG',\rho')
\ldrt E^{p,q}_r(X,\GG,\rho)$ qui induit sur la limite ${}^t \bar{f}_!$
et tel que
le morphisme sur les $E^{pq}_2$ soit le morphisme 
$$\Ext^{\; p}_{\GG'} (H^{-q}_c(X',L),\check{\rho}') 
\ldrt \Ext^{\; p}_\GG (H^{-q}_c(X,L),\check{\rho})$$ qui 
se déduit des deux applications $f_!:H^p_c(X,L)\ldrt H^p_c(X',L)$ et
$\check{\rho}'\ldrt \check{\rho}$ par fonctorialité des Ext.

Ces  morphismes de suites spectrales se composent
naturellement au sens où l'application $(X,\GG,\rho)\mapsto E^{pq}_r$
est un foncteur de la catégorie des triplets vérifiant les hypothèses
ci dessus et munis de morphismes du type de ceux ci dessus dans la
catégorie des suites spectrales.
\end{rema}
\vspace{1mm}

Nous utiliserons en fait la variante suivante du théorème précédent :

\begin{theo}\label{Hoschp} Sous les hypothèses précédentes, 
il y a une suite spectrale $\Gal (\bar{k}| k)$ équivariante 
$$
E^{pq}_2 = \text{Ext}^p_\GG ( H^{-q}_c ( X\hat{\otimes}_k \hat{\bar{k}}),
\check{\rho}) \limpl H^{-(p+q)}_c (\GG\bc X\hat{\otimes}_k \hat{\bar{k}}, \F_\rho)^*
$$
\end{theo}
\begin{proof}
Reprenons les notations de la démonstration du théorème
précédent. Appliquons la démonstration de ce théorème au recouvrement 
$(U_i\hat{\otimes} \hat{\bar{k}})_{i\in I}$ de
$X\hat{\otimes}\hat{\bar{k}}$ et au complexes de faisceaux obtenus à
partir de $\mathcal{I}^\bullet$ par pull-back
vers $X\hat{\otimes}_k \hat{\bar{k}}$. L'équivariance de la suite
spectrale obtenue ne pose alors pas de problème. 
\end{proof}

Les remarques \ref{ind_recouv} et \ref{remeq} s'appliquent bien sûr
également au théorème \ref{Hoschp}. 

\begin{enonce}{Variantes}
\begin{itemize}
\item Au lieu de nous placer dans le cadre des coefficients
$\ell$-adiques considérons des $\overline{\mathbb{F}}_\ell$ coefficients. La
démonstration précédente s'adapte aussitôt (et est même plus simple
puisqu'on n'a pas besoin d'utiliser la machinerie du foncteur $\R\pi_*$)
pour montrer l'existence 
d'une suite spectrale du même type. Le point clef étant de remarquer
que si $M$ est un $\overline{\mathbb{F}}_\ell$-espace vectoriel alors $c-Ind^\GG_1
M$  est un $\GG$-module projectif, comme c'était le cas pour les $L$ coefficients. 
\item Par contre, si l'on prend comme coefficients un anneau artinien
local $R$
 de corps résiduel un corps de caractéristique $\ell$, la démonstration
 ne marche pas en général puisque les $H^{p}_c( U_i, R)$ ne sont pas
 forcément des $R$-modules projectifs et les induites compactes de ces
 $R$-modules ne sont donc pas forcément des $\GG$-modules projectifs. La
 démonstration marche encore si tout point de $X$ possède une base de
 voisinages formée d'ouverts $U$ vérifiant $H^\bullet_c(U,R)$ est un
 $R$-module  libre.
\end{itemize}
\end{enonce}

\section[Une suite spectrale d'Hochschild-Serre pour l'uniformisation 
de R.Z.]{Une suite spectrale de Hochschild-Serre pour l'uniformisation 
de Rapoport-Zink}

La suite spectrale établie dans cette section n'est pas seulement
valable dans le cadre des variétés de Shimura de type P.E.L.
considérées dans la première partie mais pour
n'importe quelle variété de Shimura possédant une uniformisation par
des espaces rigides ``raisonnables''. 

\subsection{Spéculations}

L'idée de l'existence d'une suite spectrale telle qu'elle est énoncée
 dans le théorème qui suit est due à Michael Harris
(\cite{Har1} 
). Dans \cite{Har1} M.Harris démontre l'existence de cette suite
 spectrale dans le cas de l'espace de Drinfeld uniformisant des
 variétés de Shimura différentes des nôtres. Néanmoins, la
 démonstration que nous allons donner s'applique à n'importe quel type
 d'uniformisation. 

Reprenons les notations globales du premier
chapitre. Soit donc $b\in B(G_{\Qp},\mu_{\Qpb})$ et $\phi$ une classe d'isogénie
associée à $b$.

\begin{defi}
$$
G^\phi (\A_f)= J_b \times G(\A_f^p)
$$
\end{defi}
Il y a donc une inclusion 
$$
I^\phi (\A_f) \hookrightarrow G^\phi (\A_f)
$$
via l'action d'un automorphisme de $\phi$ sur l'homologie
cristalline en $p$ et $\ell$-adique pour $\ell\neq p$.

Il y a alors des isomorphismes lorsque $K=K_p K^p$ varie  
\begin{eqnarray*}
I^\phi (\Q) \bc \left ( \M_{K_p}\times G(\A_f^p)/K^p\right ) & \iso& 
\left (\M_{K_p}\times I^\phi (\Q) \bc G^\phi (\A_f)/K^p \right ) /J_b
\\
\; [x, y K^p] & \longmapsto & [x, I^\phi (\Q) y K^p ] 
\end{eqnarray*}
où, contrairement au membre de gauche, $J_b$ agit à droite sur
$\M_{K_p}$ dans le membre de droite. 
 Récrivons donc 
l'isomorphisme d'uniformisation rigide  sous la forme 
$$
\left (\M_{K_p}\times I^\phi (\Q) \bc G^\phi (\A_f)/K^p \right ) /J_b
\iso \Sh_K^{an} (\phi)
$$
De ce point de vue là, l'uniformisation de Rapoport-Zink est une
uniformisation p-adique. 

Pour simplifier supposons que $\rho=1$. 
Supposons l'existence d'une suite spectrale du type Hochschild-Serre
associée au  quotient p-adique précédent :
$$ E^{pq}_2=
\text{Ext}^{\; p}_{J_b-\text{lisse}} \left ( H^{-q}_c (\M_{K_p}\times
I^\phi(\Q) \bc G^\phi (\A_f)/K^p, \Qlb), 1 \right ) \limpl H^{-(p+q)}_c
(\Sh_K^{an} (\phi), \Qlb)^*
$$
Alors, 
$$
 H^q_c (\M_{K_p}\times
I^\phi(\Q) \bc G^\phi (\A_f)/K^p, \Qlb) = H^q_c (\M_{K_p},\Qlb)
\otimes \Qlb [I^\phi(\Q) \bc G^\phi (\A_f)/K^p]
$$
comme $J_b$-module. Et donc, 
$$
E^{pq}_2 =\text{Ext}^{\; p}_{J_b-\text{lisse}} \left ( H^q_c
(\M_{K_p},\Qlb) , \underbrace{\Homf ( I^\phi (\Q) \bc G^\phi
(\A_f),\Qlb)^{K^p}_{J_b
\text{lisse}}}_{(\mathcal{A}^\phi_1)^{K^p}}\right )
$$
où $\mathcal{A}^\phi_1$ est un espace de ``formes automorphes sur
$G^\phi$ `` de type $1$ à l'infini.

Le théorème qui suit donne une démonstration de ces spéculations. 
La démonstration que nous donnons est valable dans un cadre beaucoup
plus général que celui des variétés de Shimura que nous
considérons (par exemple dans le cadre de la conjecture 4.2 de \cite{ECM}).
 Elle est différente de celle de \cite{Har1} qui ne se
généralise pas au cadre général.

\subsection{La suite spectrale}

Rappelons que nous notons 
$$
H^\bullet(\M_K,\Qlb)=H^\bullet(\M_K\otimes \C_p,\Ql)\otimes\Qlb
$$

Fixons $b\in B(G_{\Qp},\mu_{\Qpb})$ et $\phi$ une classe d'isogénie dans la strate
indexée par $b$ (confère le chapitre \ref{uniformi_PEL}). 

Nous allons définir un espace de formes automorphes sur $G^\phi(\A_f)$.

\begin{defi}
  $\mathcal{A}^\phi_\rho=\{f:G^\phi(\A_f)\ldrt V_\rho \; |\; f \text{
  est } G^\phi(\A_f) \text{ lisse à droite et } \\ \forall \g\in I^\phi(\Q)\;
  \;
f(\g\bullet)=\rho(\g).f(\bullet) \}$
\end{defi}
\vspace{2mm}

Pour $\g\in I^\phi(\Q), \rho(\g)$ signifie que $\g$ est vu comme un
élément de $G(\Qlb)$ via $I^\phi(\Q)\hookrightarrow
G(\Ql)\hookrightarrow G(\Qlb)$. 

Via l'isomorphisme fixé $\Qlb\iso \C$
on a 
$$
\mathcal{A}^\phi_\rho=\{f:G^\phi(\A_f)\ldrt V_\rho \;|\; f \text{ est
  lisse et } \forall\g\in I^\phi(\Q)\;\; f(\g\bullet)=\rho(\g).f(\bullet)\}
$$
où ici $\rho(\g)$ est défini via $I^\phi(\Q)\hookrightarrow
I^\phi(\C)\hookrightarrow G(\C)$. 

$A^\phi_\rho$ est l'espace des ``formes automorphes sur $G^\phi$'' de
type $\rho$ à l'infini. Il est muni d'une action lisse de
$G^\phi(\A_f)$ par translations à droite.

Si $K^p\subset G(\A_f^p)\subset G^\phi(\A_f)$ nous noterons 
$(\AA^\phi_\rho)^{K^p}$ l'espace des formes invariantes par $K^p$
i.e. de niveau inférieur à $K^p$. 

\begin{rema}
Soit $\mathcal{A} (I^\phi)_\rho$ l'espace des formes automorphes sur
$I^\phi$ se transformant via $\check{\rho}_{|I^\phi (\R)}$ à l'infini.
Il y a alors un isomorphisme 
$$
\mathcal{A}^\phi_\rho \simeq Ind^{G^\phi (\A_f)}_{I^\phi (\A_f)}
\mathcal{A} (I^\phi)_\rho
$$
où $Ind$ désigne l'induite lisse.
\end{rema}

Soit $N=\dim(\M^{an})=\dim(\Sh)$. Soit $K=K_pK^p$ un niveau.
\\

\begin{defi} Posons 
$$
H^\bullet(\Sh^{an}_K(\phi)
,\LL_\rho)=H^{2N-\bullet}_c(\Sh_K^{an}(\phi)
,\LL_{\check{\rho}})^*(-N)  
$$
Le dual de Poincaré algébrique. 
\end{defi}

\begin{rema}
Cette définition est complètement formelle puisqu'à part pour la strate
basique, ces espaces de cohomologie sont de dimension infinie. 

Néanmoins, dans le cas de la strate basique, ces groupes peuvent
s'interpréter comme l'aboutissement d'une suite spectrale des cycles
évanescents $\ell$-adique (II.\ref{cycles_l_adiques})
\end{rema}

\vspace{4mm}

\begin{theo}\label{sui_spe_hhhh}
  Il y a un système compatible de suites spectrales $G(\A_f)\times
  W_{E_\nu}$ équivariant  
$$
E^{pq}_2(K_pK^p)=\Ext^{\; p}_{J_b\text{-lisse}}\left 
(H^{2N-q}_c(\M_{K_p},\Qlb)(N),(\AA^\phi_\rho)^{K^p}\right )
\limpl H^{p+q}(\Sh^{an}_K(\phi),\LL_\rho^{an})
$$
\end{theo}

Quelques remarques s'imposent avant la démonstration :
$G(\A_f)=G(\Qp)\times G(\A_f^p)$ et $\forall g\in G(\Qp)$ il y a un
morphisme 
$$
H^{2N-q}_c(\M_{g^{-1}K_p g},\Qlb)\xrig{\; (g^{-1})_!\;} H^{2N-q}_c(\M_{K_p},\Qlb)
$$
qui induit par fonctorialité des $\Ext$ un morphisme $E^{pq}_2(K_p K^p) 
\drt E^{pq}_2(g^{-1}K_p g K^p)$.
\\
Dire que la suite spectrale est $G(\Qp)$ équivariante signifie que
ce morphisme est le morphisme associé au niveau des groupes $E^{pq}_2$ à un
 morphisme de suites spectrales
$$
E^{pq}_r (K_p K^p) 
\drt E^{pq}_r(g^{-1}K_p g K^p)
$$
 qui
induit dans la limite le morphisme usuel 
$$
H^{p+q}(\Sh_K(\phi),\LL_\rho^{an})\ldrt H^{p+q}(\Sh_{g^{-1}Kg}(\phi),\LL_\rho^{an})
$$
Cela signifie également que ces morphismes de suites spectrales se
composent de façon naturelle. 
De même, $G(\A_f^p)$ opère sur $\AA^\phi_\rho$ et donc induit un
morphisme pour $g\in G(\A_f^p)$ 
$$
(\AA^\phi_\rho)^{K^p}\ldrt (\AA^\phi_\rho)^{g^{-1}Kg}
$$
qui induit le morphisme  $E^{pq}_2(K)\ldrt E^{pq}_2(K_p g^{-1}K^p g)$.
La $G(\A_f^p)$ équivariance est alors prise au même sens que pour
$G(\Qp)$. Quant à l'action de $W_{E_\nu}$ elle se fait sur chaque
$E^{pq}_2 (K)$ et ne pose pas de problème d'interprétation. 

Reste à expliquer ce que l'on appelle système compatible. Soient donc
$K'_p \subset K_p$ et ${K^p}'\subset K^p$. La compatibilité signifie
qu'il y a un morphisme de suites spectrales 
$$
E^{pq}_r (K') \ldrt E^{pq}_r (K)
$$
qui est induit au niveau des termes $E^{pq}_{2}$ par les morphismes
$J_b$
équivariants 
$$
(\Pi_{K_p,K'_p })_! :
H^q_c (\M_{K_p'},\Qlb) \ldrt H^q_c (\M_{K_p},\Qlb)
$$
et l'inclusion
$$
(\mathcal{A}^\phi_\rho )^{K^p} \hookrightarrow (\mathcal{A}^\phi_\rho )^{{K^p}'} 
$$
couplés à la fonctorialité de Ext. 
Bien s\^ur, ces morphismes de transition sont compatibles à l'action
de $G(\A_f)\times W_{E_\nu}$. 

\begin{proof} 
Prenons les notations suivantes pour l'uniformisation rigide : 
$$\Theta : 
I^\phi (\Q) \bc \left ( \M_{K_p}\times G(\A_f^p)/K^p \right )\iso
\Sh_K^{an} (\phi) 
$$
Soient $(y_i)_{i\in I}$ des éléments de $G(\A_f^p)/K^p$ vérifiant 
$$
G(\A_f^p)/K^p = \coprod_{i\in I} I^\phi (\Q). y_i 
$$
et notons pour tout $i$ dans $I$ $\;\; \GG_i=Stab_{I^\phi (\Q)}
(y_i)$. Nous verrons les groupes $\GG_i$ comme des sous-groupes de
$J_b$ via le plongement $I^\phi (\Q) \hookrightarrow J_b$. A un tel
choix est associé un isomorphisme 
$$
\Xi_{(y_i)_i} :
\coprod_{i\in I} \GG_i\bc \M_{K_p} \iso I^\phi (\Q) \bc \left (
\M_{K_p} \times G(\A_f^p)/K^p\right ) 
$$
Notons $\rho_i:\GG_i \ldrt \Gl (V_\rho)$ le morphisme continu (pour
$\GG_i$ muni de la topologie profinie et $V_\rho$ la topologie
$\ell$-adique) défini par le composé 
\begin{diagram}
\GG_i & \rInto & I^\phi ( \Q) & \rInto & G(\Ql) & \rTo^\rho & \Gl (V_\rho)
\end{diagram}

\begin{lemm}\label{fscth}
$$
  \Xi_{(y_i)_i}^*(\Theta^*\LL_\rho) =\coprod_{i\in I} \F_{\rho_i}
$$
où $ \F_{\rho_i}$ est le système local associé à $\rho_i$ (cf. le début de
la section \ref{suit_disc}). 
\end{lemm}

\begin{proof}
 Il suffit de regarder ce que donne une variation du
niveau $K^p$ à droite sur la décomposition de gauche.
\end{proof}
 
 Nous allons avoir besoin des deux lemmes suivants : 
\begin{lemm}\label{leind} Soit $i\in I$.
  Soient $\pi_1$ une représentation lisse de $J_b$ et $\pi_2$ une
  représentation lisse de $\GG_i$.
$$
\Ext^\bullet_{J_b-\text{lisse}}(\pi_1,Ind^{J_b}_{\GG_i} \pi_2)\simeq \Ext^\bullet_\GG(\pi_1,\pi_2)
$$
où les induites $Ind$ sont des induites lisses.
\end{lemm}

\begin{proof}
 Cela résulte de la réciprocité de Frobenius usuelle pour les groupes
$p$-adiques et de 
ce que $\GG_i$ étant discret,
$\Ext_{\GG_i-lisse}=\Ext_{\GG_i}$.
\end{proof}

\begin{lemm}\label{leprod}
  Soient $\pi$ et $(\rho_i)_{i\in I}$ des représentations lisses de
  $J_b$.
Il y a un isomorphisme canonique :
$$
\prod_{i\in I} \Ext^\bullet_{J_b-\text{lisse}}(\pi,\rho_i)\simeq
\Ext_{J_b-\text{lisse}}(\pi,\left ( \prod_{i\in I} \rho _i \right )^{\text{lisse}}) 
$$
où $(\prod_i \rho_i)^{lisse}$ est la partie lisse de la représentation produit.
\end{lemm}

\begin{proof}
 Soit $P_\bullet \twoheadrightarrow \pi$ une résolution
projective de $\pi$ dans la catégorie des $J_b$-modules lisses.
\begin{eqnarray*}
  \prod_{i\in I} \Ext^p_{J_b-\text{lisse}} (\pi,\rho_i) & =&
  \prod_{i\in I} H^p
  (\Hom_{J_b} (P_\bullet,\rho_i)) \\
 &=& H^p \left (\prod_{i\in I} \Hom_{J_b} (P_\bullet,\rho_i)\right ) \\
 &=& H^p \left ( \Hom_{J_b} (P^\bullet,\prod_{i\in I} \rho_i ) \right ) \\
\end{eqnarray*}
or $ \Hom_{J_b} (P^\bullet,\prod_{i\in I} \rho_i )= 
 \Hom_{J_b} (P^\bullet,(\prod_{i\in I} \rho_i)^{lisse} )$
car $P^\bullet$ étant lisse l'image d'un morphisme de $J_b$-module dans
 $\prod_i \rho_i$ atterrit dans $(\prod_i \rho_i)^{lisse}$. 

D'où le résultat.
\end{proof}

Pour $i\in I$ nous considérons la suite spectrale du théorème
\ref{Hoschp} (rappelons qu'a fin de ne pas alourdir les notations on
ne note pas les extensions des scalaires à $\Cp$) : 
$$
E^{pq}_2(i)=\text{Ext}^{\; p}_{\GG_{i}} \left (
  H^{-q}_c(\M_{K_p},\Qlb),\rho_{i}\right ) \limpl H^{-(p+q)}_c
  (\GG_i\bc \M,\F_{\check{\rho}_i})^* 
$$
Le produit de ces suites spectrales est une suite spectrale convergente 
$$
\prod_{i\in I} E^{pq}_2 (i) \limpl \prod_{i\in I}
 H^{-(p+q)}_c (\M/\GG_{i,K^p},\F_{\check{\rho}_{i,K^p}})^*
$$
L'aboutissement de cette suite spectrale se récrit  
\begin{eqnarray*} &&
\prod_{i\in I}  H ^{-(p+q)}_c  (\M/\GG_{i,K^p},\F_{\check{\rho}_{i,K^p}})^*
\iso  H^{-(p+q)}_c (\coprod_{i\in I} \M/\GG_{i,K^p},\coprod_{i\in I}
\F_{\check{\rho}_{i,K^p}})^* \\  &&\underset{\sim}{\xrig{\; (\Xi_{(y_i)_i })_*
\;}} H^{-(p+q)}_c\left ( I^\phi (\Q)\bc \left ( 
\M_{K_p}\times G(\A_f^p)/K^p \right ),\LL_{\check{\rho}} \right)^*
\underset{\sim}{\xrig{\;\;\Theta_*\;\;}} H^{-(p+q)}_c
(\Sh_K^{an}(\phi),\LL_{\check{\rho}})^* \\
 &\;& =  H^{p+q+2N} (\Sh_K^{an} (\phi),\LL_{\rho}) (N)
\end{eqnarray*}
Quant aux premiers termes, grâce au lemme  \ref{leind}
$$
\forall i\in I\;\forall p,q\;\; \text{Ext}^{\; p}_{\GG_i}
(H^{-q}_c (\M_{K_p},\Qlb),\rho_i) \simeq \text{Extp}^{\;
p}_{J_b\text{-lisse}} ( H^{-q}_c (\M_{K_p},\Qlb),
\text{Ind}^{J_b}_{\GG_i} \rho_i)  
$$
et donc grâce au lemme \ref{leprod}
$$
\prod_{i\in I} E^{pq}_2 (i) \simeq 
\Ext^p_{J_b\text{-lisse}}\left (H^{-q}_c(\M_{K_p},\Qlb),\left (\prod_i
  \rho_i 
\right )^{\text{lisse}}\right) 
$$

L'ingrédient suivant est le lemme 
\begin{lemm} Il y a un isomorphisme de $J_b$-modules
$$\zeta_{(y_i)_i} : 
\left ( \prod_{i\in I} Ind_{\GG_i}^J \rho_i
 \right )^{lisse} \iso 
 (\AA^\phi_\rho )^{K^p}
$$
\end{lemm}

\begin{proof} 
Définissons $\zeta_{(y_i)_i}$. Soit $(f_i)_{i\in I} \in \prod_{i}
Ind_{\GG_i}^J \rho_i$. On a donc, $\forall i\in I \; f_i : J_b\ldrt
V_\rho$ et est tel que $\forall \g \in \GG_i \; f_i
(\g\bullet)=\rho_i (\g) . f_i (\bullet)$. Définissons
$f=\zeta_{(y_i)_i} ( (f_i)_i)$ : $\left
(\mathcal{A}^\phi_\rho\rho\right )^{K^p}$ s'identifie aux fonctions 
de $J_b \times G(\A_f^p)/K^p $ à valeurs dans $V_\rho$ se transformant
via $\rho$ par l'action à gauche de $I^\phi (\Q)$. 
Soit  $(a,b)\in J_b\times G(\A_f^p)/K^p$. $\exists ! i\in I\; \exists
\g \in I^\phi (\Q) \;\; b=\g.y_i$. Posons alors 
$$
f(a,b)=\rho (\g). f_i (\g^{-1} a)
$$
$f(a,b)$ est bien défini car si $b=\g'.y_i$ où $\g'\in I^\phi (\Q)$
alors ${\g'}^{-1} \g \in \GG_i$ et donc 
\begin{eqnarray*}
\rho ( \g'). f_i ({\g'}^{-1} a)) &=& \rho (\g'). f_i ( ( {\g'}^{-1} \g
)  
\g^{-1} a ) \\ 
&=&  \rho (\g') \rho_i ({\g'}^{-1} \g) f_i (\g^{-1} a) = \rho (\g) f_i
(\g^{-1} a) 
\end{eqnarray*}
Par définition de $f$, $\forall \g \in I^\phi (\Q) \; f (\g
\bullet)=\rho (\g). f(\bullet)$. La lissité de $f$ est claire quant à
l'action de $G(\A_f^p)$ puisqu'elle est invariante par $K^p$, quant à
la lissité vis à vis de l'action de $J_b$ elle résulte de l'existence
d'un sous-groupe compact ouvert $C$ dans $J_b$ vérifiant $\forall c\in
C\; \forall i\in I \;\; f_i (\bullet c)= f_i (\bullet)$. 

Définissons l'inverse de $\zeta_{(y_i)_i}$. Si $f\in \left (
\mathcal{A}^\phi_\rho \right )^{K^p}$, $f: J_b \times
G(\A_f^p)/K^p\ldrt V_\rho$, pour un $i$ dans $I$ et $a$ dans $J_b$
posons 
$$
f_i (a)= f (a,y_i)
$$
$f\mapsto (f_i)_i$ est bien l'inverse de $\zeta_{(y_i)_i}$. 
\end{proof}

On obtient donc au final une suite spectrale $W_{E_\nu}$ équivariante 
$$
E^{pq}_2((y_i)_{i\in I}) = \text{Ext}^{\; p}_{J_b-\text{lisse}} (
H^{-q}_c(\M_{K_p},\Qlb),(\mathcal{A}^\phi_\rho)^{K^p}) \limpl 
H^{2N+p+q} ( \Sh_K^{an} (\phi) ,\LL_{\check{\rho}}) (N) 
$$
où $E^{pq}_2 ((y_i)_i)$ et l'aboutissement ne dépendent pas du choix
$(y_i)_i$ mais où les autres termes en dépendent à priori.

La suite spectrale annoncée s'obtient en translatant verticalement
cette suite spectrale de $2N$ et en la tordant par $\Qlb (-N)$. 
\\

Montrons qu'en fait cette suite spectrale est indépendante du choix
des $(y_i)_i$ au sens suivant :
\begin{lemm}\label{lemme_indep}
  Supposons choisis des $(y'_j)_{j\in J},\; y'_j\in G(\A_f)/K^p$ tels que 
$$
G(\A_f)/K^p = \coprod_{j\in J} I^\phi (\Q). y'_j
$$
Il y alors un isomorphisme naturel de suites spectrales 
$E^{pq}_r ((y_i)_i)\drt E^{pq}_r((y'_j)_j)$ induisant l'identité sur
$E^{pq}_2$ et l'aboutissement :
\begin{diagram}
  E^{pq}_2((y_i)_{i\in I}) &=&  \text{Ext}^{\; p}_{J_b-\text{lisse}} (
H^{-q}_c(\M_{K_p},\Qlb),(\mathcal{A}^\phi_\rho)^{K^p}) \limpl  &
H^{-(p+q)}_c (\Sh_K^{an} (\phi), \LL_{\check{\rho}}) ^*
\\
\dTo^{Id} &&& \dTo_{Id} \\
  E^{pq}_2((y'_j)_{j\in J}) &=& \text{Ext}^{\; p}_{J_b-\text{lisse}} (
H^{-q}_c(\M_{K_p},\Qlb),(\mathcal{A}^\phi_\rho)^{K^p}) \limpl &
H^{-(p+q)}_c ( \Sh_K^{an} (\phi),\LL_{\check{\rho}})^*
\end{diagram}
La naturalité étant prise au sens où si l'on fait un troisième choix
$(y''_k)_k$ comme ci dessus 
le morphisme de suites spectrales $E^{pq}_r((y_i)_i)\drt
E^{pq}_r((y''_k)_k)$ est le composé de $E^{pq}_r((y_i)_i)\drt 
E^{pq}_r((y'_j)_j)$ et $E^{pq}_r((y'_j)_j)\drt  E^{pq}_r((x''_k)_k)$
\end{lemm}
\begin{proof}

Quitte à réindexer on peut supposer que $J=I$ et que $\forall i\in I
\; y'_i\in I^\phi (\Q). y_i$. Fixons pour tout $i$ dans $I$ un $\g_i$
dans $I^\phi (\Q)$ vérifiant $y'_i=\g_i. y_i$. On a donc 
$$
\forall i\in I \; \GG'_i=Stab_{I^\phi (\Q)} y'_i=\g_i \GG_i \g_i^{-1}
$$
On définit comme précédemment $\rho'_i$. 
Notons pour tout $i$ dans $I$ $\; E^{pq}_r (i)'$ la suite spectrale
associée à $(\M_{K_p},\GG'_i,\rho_i')$.

Considérons l'isomorphisme de suites spectrales
 associé au morphisme de triplets 
$(X,\GG,\rho)\ldrt (X',\GG',\rho')$ par la remarque \ref{remeq} où
 $X=X'=\M_{K_p}$, le morphisme 
 $X\drt X'$ est l'action de $\g_i\in J_b$,
 $\GG=\GG_{i},\GG'=\GG'_{i}$ et le morphisme de groupes $\GG\drt \GG'$
 est la conjugaison par $\g_i^{-1}$,
  et le morphisme  $\rho\drt\rho'\circ \text{int}_{\g_i^{-1}}$ est la
 conjugaison par $\rho (\g_i)$. A un tel choix de $(\g_i)_i$ est donc
 associé un isomorphisme de suites spectrales  
 \begin{diagram}
   E^{pq}_r(i)  & \limpl &
   H^{p+q}_c(\GG_i\bc\M_{K_p},\F_{\check{\rho}_i} ) \\
\dTo_{\simeq}^{\a_{(\g_i)_i}} & & \dTo_{\simeq } \\
 E^{pq}_r(i)' & \limpl &
   H^{p+q}_c(\GG'_i\bc\M_{K_p},\F_{\check{\rho'}_i} ) 
 \end{diagram} 
Le morphisme $\a_{(\g_i)_i}$ induit au niveau des $E^{pq}_2$ un morphisme
$$
\text{Ext}^{\; p}_{\GG_i}
   (H^q_c(\M_{K_p},\Qlb),\rho_i) \ldrt  \text{Ext}^{\; p}_{\GG'_i}
   (H^q_c(\M_{K_p},\Qlb),\rho'_i)
$$
qui se décrit en utilisant la fonctorialité de Ext.

Le produit de ces isomorphismes de suites spectrales induit un
isomorphisme  de suites spectrales convergentes 
$$\prod_{i\in I} E^{pq}_r(i)\ldrt \prod_{i\in I}
E^{pq}_r(i)'$$
Rappelons qu'au choix des $y_i$ on a associé deux isomorphismes 
$$
\prod_{i\in I} E^{pq}_2 (i)\iso \text{Ext}^{\; p}_{J_b\text{-lisse}} (
H^q_c (\M_{K_p},\Qlb), \left (\mathcal{A}^\phi_\rho\right )^{K^p}) 
$$
et 
$$
\prod_{i\in I} H^{-(p+q)}_c (\GG_i \bc \M_{K_p}, \F_{\check{\rho}_i} )^*
\iso H^{-(p+q)}_c (\Sh_K^{an} (\phi),\LL_{\check{\rho}})^*
$$
Il en est de même pour $y'_i$. On en déduit l'existence de
l'isomorphisme de suites spectrales annoncé.

Afin de montrer que ce morphisme induit
bien l'identité sur les termes associés à $r=2$ il faut montrer le
lemme suivant 
\begin{lemm}
  Les diagrammes suivants commutent :
  \begin{diagram}
 & & \prod_i E^{pq}_2(i) \\
& \ruTo^\simeq \\
    \text{Ext}^{\; p }_{J_b\text{-lisse}}(
    H^{-q}_c(\M_{K_p},\Qlb),(\mathcal{A}^{\phi}_\rho)^{K^p}) & &\dTo \\
& \rdTo^\simeq \\
 & & \prod_i E^{pq}_2(i)'
  \end{diagram}
  \begin{diagram}
    \prod_i H^{p+q}_c(\GG_i\bc \M,\F_{\rho_i}) \\
 &  \rdTo^\simeq \\
\dTo_\simeq^{\prod_i {\g_i}_!} & & H^{p+q}_c (I^\phi (\Q)\bc (
  \M_{K_p}\times G(\A_f^p)/K^p)
,\Theta^*\LL_\rho) \\
 & \ruTo^\simeq \\
 \prod_i H^{-(p+q)}_c(\GG'_i\bc \M,\F_{\rho'_i})
  \end{diagram}
\end{lemm}

\begin{proof}
 Pour le deuxième diagramme c'est clair puisque le diagramme
suivant commute :
\begin{diagram}
  \coprod_i \GG_i\bc \M \\
& \rdTo^\simeq \\
\dTo^{\prod_i \g_i} & & I^\phi(\Q)\bc (\M_{K_p}\times G(\A_f^p)/K^p) \\
& \ruTo^\simeq \\
\coprod_i \GG'_i\bc \M 
\end{diagram}
Quant au premier cela se ramène aisément (il suffit de suivre les
différents isomorphismes) à montrer que le diagramme suivant commute :
\begin{diagram}
   & & \prod_i \text{Ind}^J_{\GG_i} \rho_i \\
& \ruTo^{\zeta_{(y_i)_i}}   \\
(\mathcal{A}_\rho^\phi)^{K^p} && \dTo_{\prod_i \Delta_i} \\
& \rdTo^{\zeta_{(y'_i)_i}} \\
 & & \prod_i \text{Ind}^J_{\GG'_i} \rho'_i
\end{diagram}
où $\Delta_i:Ind^J_{\GG_i}\rho_i \iso Ind^J_{\GG'_i} \rho'_i$ est
induit par $(\text{int}_{g_i},\rho(\g_i)):(\GG_i,\rho_i)\drt
(\GG'_i,\rho'_i)$. La commutativité est laissée au lecteur.
\end{proof}

Remarquons que l'isomorphisme de suites spectrales ainsi construit
dépend du choix des $\g_i$ qui ne sont pas uniques. Néanmoins, si l'on
fait un autre choix de $\g_i$ l'isomorphisme de suites spectrales
induit encore l'identité sur le premier terme de la suite
spectrale. Or deux morphismes de suites spectrales coïncidant sur les
termes $E^{pq}_2$ coïncident. L'isomorphisme ne dépend donc pas du
choix des $\g_i$.

Finalement il reste à vérifier la naturalité des morphismes de suites
spectrales construits mais cela est clair par le même argument que
précédemment 
puisque sur les termes $E^{pq}_2$
tous donnent l'identité.
Cela finit la démonstration de l'indépendance des du choix des $y_i$. 
\end{proof}

Afin d'avoir une définition intrinsèque ne dépendant pas du choix des 
$(y_i)_i$ nous poserons 
$$
E^{pq}_r = \underset{(y_i)_i}{\limp} E^{pq}_r ((y_i)_i)
$$
où les isomorphismes de transition sont les isomorphismes construits
dans le lemme \ref{lemme_indep}. Il résulte également de ce lemme qu'il y a
un isomorphisme  
$$
E^{pq}_2 \iso\text{Ext}_{J_b-\text{lisse}} ( H^{-q}_c (\M_{K_p},\Qlb), (\mathcal{A}^\phi_\rho)^{K^p})
$$ 
\\

Reste à vérifier l'équivariance de la suite spectrale.

Nous allons faire varier le groupe $K=K_pK^p$, c'est pourquoi nous
préciserons les notations précédentes en notant 
$$
E^{pq}_r (K, y_i )
$$
la suite spectrale notée auparavant $E^{pq}_r (i)$, 
$$
E^{pq}_{r} (K, (y_i)_i)
$$
la suite spectrale construite précédemment associée à des $(y_i)_i$
vérifiant 
$$
G(\A_f^p)/K^p=\coprod_{i\in I} I^\phi (\Q). y_i
$$
et 
$$
E^{pq}_r (K)=\underset{(y_i)_i}{\limp} E^{pq}_r (K, (y_i)_i)
$$

Commençons par la $G(\Qp)$ équivariance. 
Soit $g\in G(\Qp)$ et $(y_i)_{i\in I}$. 
Puisque l'action de $G(\Qp)$ commute à celle de
$J_b$, pour
tout $i$ dans $I$  
 $g^{-1}$ induit un isomorphisme de triplet (\ref{remeq})  
$$
(\M_{g^{-1}K_p g},\GG_i,\rho_i)\iso (\M_{K_p},\GG_i,\rho_i)
$$ 
d'où d'après la remarque \ref{remeq} un isomorphisme de suites spectrales 
$$
E^{pq}_r (K_p,y_i) \iso E^{pq}_r (g^{-1} K_p g,y_i)
$$
qui induit le morphisme annoncé sur les termes $E^{pq}_2$, c'est à
dire le transposé du morphisme 
$$
H^{-q}_c(\M_{g^{-1}K_p g},\Qlb)\xrig{\; (g^{-1})_!\;} H^{-q}_c(\M_{K_p},\Qlb)
$$
 Prenant le
produit de ces morphismes lorsque $i$ varie on obtient un isomorphisme 
$$
E^{pq}_r (K,(y_i)_i) \iso E^{pq}_r (g^{-1}K_p g K^p,(y_i)_i)
$$
Si de plus on fait un autre choix de $(y'_i)_i$, le diagramme suivant
est commutatif
\begin{diagram}[size=1cm]
 E^{pq}_r (K, (y_i)_i) & \rTo^\sim &  E^{pq}_r (g^{-1}K_p g K^p,
 (y_i)_i) \\
\dTo^\simeq & & \dTo^\simeq \\
 E^{pq}_r (K, (y'_i)_i) & \rTo^\sim &  E^{pq}_r (g^{-1}K_p g K^p,
 (y'_i)_i)
\end{diagram} 
où les applications verticales sont celles construites dans le lemme
\ref{lemme_indep}. En effet, ce diagramme de suites spectrales commute
au niveau des termes $E^{pq}_2$, il commute donc. 
On en déduit en prenant la limite un isomorphisme 
$$
E^{pq}_r ( K) \iso E^{pq}_r (g^{-1}K_p g K^p)
$$

Passons à la $G(\A_f^p)$ équivariance. Soit donc $g\in G(\A_f^p)$
et $(y_i)_i$ comme précédemment. 
 $g$
induit un isomorphisme
\begin{eqnarray*}
\a: G(\A_f^p)/K^p & \iso & G(\A_f^p)/ (g^{-1} K^p g) \\
 x K^p & \mapsto & xg (g^{-1}K^p g) 
\end{eqnarray*}
$\a$ commute à l'action de $I^\phi (\Q)$ et donc, comme $(y_i)_i$ est tel que
$$
G(\A_f^p)/K^p=\coprod_{i\in I} I^\phi (\Q). y_i
$$
alors 
$$
G(\A_f^p)/(g^{-1} K^p g)= \coprod_{i\in I} I^\phi (\Q). \a(y_i)
$$
Remarquons que pour  $i\in I\;\;$ $\GG_i=Stab_{I^\phi (\Q)} y_i=
Stab_{I_\phi (\Q)} \a (y_i)$. Il y a donc des égalités
$$
E^{pq}_r (K,y_i)= E^{pq}_r (K_p g^{-1} K^p g, \a (y_i))
$$
Prenant le produit sur $i$ on obtient l'égalité 
$$
E^{pq}_r (K, (y_i)_i) =E^{pq}_r ( K_p g^{-1} K^p g, (\a(y_i))_i)
$$
Pour montrer que ce morphisme induit bien l'application attendue sur
les termes $E^{pq}_2$ il faut montrer le lemme suivant dont la
démonstration ne pose aucun problème 
\begin{lemm}
Le diagramme suivant commute 
\begin{diagram}
( \mathcal{A}^\phi_\rho)^{K^p} & \rTo^{\zeta_{(y_i)_i}} & \prod_{i\in
I} \text{Ind}^{J_b}_{\GG_i} \rho_i \\
\dTo^\simeq & & \dEqualto \\
 ( \mathcal{A}^\phi_\rho)^{g^{-1}K^p g} & \rTo^{\zeta_{(\a(y_i))_i}} & \prod_{i\in
I} \text{Ind}^{J_b}_{\GG_i} \rho_i 
\end{diagram}
où l'application verticale de gauche est induite par la translation d'une fonction
de $\mathcal{A}^\phi_\rho$ par $g$. 
\end{lemm}
Si maintenant $(y'_i)_i$ est un autre choix, le diagramme suivant
commute 
\begin{diagram}[size=1cm]
 E^{pq}_r (K, (y_i)_i) & \rTo^\sim &  E^{pq}_r (K_p g^{-1} K^p g,
 (\a(y_i))_i) \\
\dTo^\simeq & & \dTo^\simeq \\
 E^{pq}_r (K, (y'_i)_i) & \rTo^\sim &  E^{pq}_r (K_p g^{-1} K^p g,
 (\a(y'_i))_i)
\end{diagram}
où les applications verticales sont celles construites dans le lemme
\ref{lemme_indep}. En effet, il résulte de la commutativité du
diagramme précédent que ce diagramme commute au niveau des termes
$E^{pq}_2$. Il commute donc. Passant à la limite on obtient donc 
un isomorphisme 
$$
E^{pq}_r (K) \iso E^{pq}_r ( K_p g^{-1} K^p g)
$$

Passons maintenant à la compatibilité lorsque $K$ varie. Soit donc
$K_p'\subset K_p$ un sous-groupe compact ouvert et $K'=K'_p K^p$. Soit $(y_i)_i$ comme
précédemment. Il y a alors pour tout $i$ dans $I$ un morphisme de
triplet 
$$
(\M_{K'_p},\GG_i,\rho_i)\ldrt (\M_{K_p},\GG_i,\rho_i)
$$
qui induit un morphisme de suite spectrale $$E^{pq}_r (K,y_i)\ldrt
E^{pq}_r (K',y_i)$$ Prenant le produit des ces morphismes on obtient
un morphisme de suite spectrale 
$$
E^{pq}_r (K, (y_i)_i) \ldrt E^{pq}_r (K', (y_i)_i) 
$$
dont il est clair qu'il induit bien le morphisme attendu au niveau des
termes $E^{pq}_2$ et qu'il est compatible aux isomorphismes du lemme
\ref{lemme_indep} lorsque l'on fait varier $(y_i)_i$. 

Soit maintenant ${K^p}'\subset K^p$ un sous-groupe compact ouvert et
$K'=K_p {K^p}'$. La projection
$$
pr:G(\A_f^p)/K^p \ldrt G(\A_f^p)/ {K^p}'
$$
est compatible à l'action à gauche de $I^\phi (\Q)$.
 Soit $(y_i)_i$
tel que 
$$
G(\A_f^p)/K^p=\coprod_{i\in I} I^\phi (\Q).y_i
$$
et $(x_j)_j$ des éléments de $G(\A_f^p)/{K^p}'$ vérifiant 
$$
G(\A_f^p)/{K^p}'=\coprod_{j\in J} I^\phi (\Q).x_j
$$
et $\forall j\in J\; \exists i\in I \; pr (x_j)=y_i$. Pour $j\in J$ notons $\b
(j)\in I$ l'élément vérifiant $pr (x_j)=y_{\b (j)}$. Pour tout $j\in
J$, $\GG_j \subset GG_{\b (j)}$ d'où un morphisme de triplets
$$
(\M_{K_p},\GG_j,\rho_j)\ldrt (\M_{K_p},\GG_{\b(j)},\rho_{\b (j)})
$$
qui induit un morphisme de suites spectrales 
$$
E^{pq}_r (K, y_{\b (j)}) \ldrt E^{pq}_r (K', y_j)
$$
On obtient donc en assemblant ces morphismes un morphisme de suites spectrales 
$$
E^{pq}_r (K, (y_i)_i) \ldrt E^{pq}_r (K', (x_j)_j)
$$
Il résulte du lemme suivant dont la démonstration ne pose pas de problème  que ce morphisme est bien le morphisme
attendu au niveau des termes $E^{pq}_2$.

\begin{lemm}
Le diagramme suivant commute 
\begin{diagram}
(\mathcal{A}^\phi_\rho)^{K^p} & \rTo^{\zeta_{ (y_i)_i}} &
\prod_{i\in I} \text{Ind}^{J_b}_{\GG_i} \rho_i \\
\dInto & & \dTo \\
(\mathcal{A}^\phi_\rho)^{{K^p}'} & \rTo^{\zeta_{ (x_j)_j}} &
\prod_{j\in J} \text{Ind}^{J_b}_{\GG_j} \rho_j 
\end{diagram}
où les morphismes de gauches sont induits par les égalités 
$\text{Res}^{\GG_j}_{\GG_{\b(j)}} \rho_{\b(j)}= \rho_j$.
\end{lemm}
Si maintenant $(y'_j)_j$ et $(x'_j)_j$ sont comme ci dessus, il y a un
diagramme commutatif 
\begin{diagram}
E^{pq}_r (K, (y_i)_i) & \rTo & E^{pq}_r (K', (x_j)_j) \\
\dTo^\simeq & & \dTo^\simeq \\
E^{pq}_r (K, (y'_i)_i) & \rTo & E^{pq}_r (K', (x'_j)_j)
\end{diagram}
où les application verticales sont celles construites dans le lemme
\ref{lemme_indep}. En effet, c'est déjà le cas au niveau des termes
$E^{pq}_2$. 

On en déduit un morphisme 
$$
E^{pq}_r (K) \ldrt E^{pq}_r(K')
$$
\\

Les actions de $G(\Qp)$ et de $G(\A_f^p)$  construites commutent et se
composent naturellement 
puisque c'est la cas au niveau des termes $E^{pq}_2$. De même, ces
actions commutent aux morphismes de changement de niveau construits ci
dessus. 
\end{proof}
\vspace{4mm}

\begin{coro}
  Il y a une suite spectrale $G(\A_f)\times W_{E_\nu}$ équivariante 
$$
E^{pq}_2=\underset{K}{\limi}\;
\Ext^{\; p}_{J_b-\text{lisse}}\left 
(H^{2N-q}_c(\M_{K_p},\Qlb)(N),(\AA^\phi_\rho)^{K^p}\right )
\limpl \underset{K}{\limi}\;
H^{p+q}(\Sh^{an}_K(\phi),\LL_\rho^{an})
$$
\end{coro}

\chapter[Matsushima $p$-adique]{Application à la strate basique : ``une formule de Matsushima
$p$-adique''}

Soit $b=b_0$ la classe basique de $B(G_{\Qp},\mu_{\Qpb})$. Supposons $\Sb (b_0)\neq
\emptyset$ (ce que nous montrerons plus tard).  

Si  la classe d'isoggénie $\phi$ induit la classe $b_0$, soit $I^\phi$
la forme intérieure de 
$G$ associée. Le groupe réductif $I^\phi$ étant anisotrope modulo son
centre, le $\C$-espace vectoriel 
 $\AA^\phi_\rho$ est l'espace des
formes automorphes sur $I^\phi$ de type à l'infini $\check{\rho}'$ où  
\begin{diagram}
  \rho':I^\phi (\R) & \rInto & I^\phi (\C) =G(\C) & \rTo^\rho & \Gl (V)
\end{diagram}
au sens où
$$
\mathcal{A}^{\phi}_{\rho}=\Hom_{I^\phi (\R)} (\check{\rho}',\mathcal{A} (I^\phi ))
$$
Nous noterons encore $\rho=\rho'$.

Supposons maintenant $\rho$ irréductible et soit $\mathcal{A}(I^\phi)$
l'espace des représentations automorphes de $I^\phi$ (ce n'est pas l'ensemble
des classes d'isomorphismes : on tient compte des multiplicités).

Le groupe de Lie 
$I^\phi(\R)$ étant anisotrope modulo son centre, pour un sous groupe
compact ouvert $K^p$ de $G(\A_f^p)$ 
$$
(\mathcal{A}^\phi_\rho)^{K^p}=\bigoplus_{\Pi\in \mathcal{A}(I^\phi)
  \atop \Pi_\infty=\check{\rho}} \Pi_p\otimes (\Pi^p)^{K^p}
$$

 Il y un isomorphisme pour tous $p$ et $q$ 
$$
\text{Ext}^{\;
  p}_{J_b\text{-lisse}}\left ( H^q_c(\M_{K_p}),(\mathcal{A}^\phi_\rho)^{K^p}\right )
\simeq \bigoplus_{\Pi\in \mathcal{A}(I^\phi)\atop \Pi_\infty =\check{\rho}} 
\text{Ext}^{\; p}_{J_b\text{-lisse}} \left (
  H^q_c(\M_{K_p}),\Pi_p\right )\otimes (\Pi^p)^{K^p}
$$
où $H^q_c(\M_{K_p})=H^q_c ( \M_{K_p},\Qlb)$. 
En effet, $H^q_c(\M_{K_p})$ étant un $J_b$-module de type fini (\ref{Jfin}), les
$\text{Ext}$ se calculent par des résolutions de $H^q_c(\M_{K_p})$ dont chaque
terme est une somme finie de termes du type 
$c-Ind^{J_b}_{K_0}\pi$ où $K_0$ 
compact ouvert et $\pi$ de dimension finie. Or,
$$
Hom_{J_b}(c-Ind^{J_b}_{K_0} \pi, (\mathcal{A}^\phi_\rho)^{K^p})
\simeq \Hom_{K_0}(\pi, (\mathcal{A}^\phi_\rho)^{K^p})
=\Hom_{K_0}(\pi,  (\mathcal{A}^\phi_\rho)^{K_1K^p})
$$
où $K_1\subset K_0$ est compact ouvert tel que $\pi(K_1)=\{ Id\}$. 
$(\mathcal{A}^\phi_\rho)^{K_1K^p}$ est de dimension finie
(admissibilité de l'espace des formes automorphes) et donc 
$$
\#\{ \Pi\in\mathcal{A}(I^\phi)\; |\; \Pi_\infty=\rho \;\; \Pi^{K_1
  K^p}\neq (0)\}<+\infty
$$
On en déduit que le calcul des Ext ne fait intervenir qu'un nombre
fini de $\Pi$ d'où le résultat.

On en déduit alors en utilisant le lemme I.\ref{lemme_annulation_Ext}
ainsi que la proposition I.\ref{Jfin}
(confère \citeg{Har4} pour la définition du groupe de
Grothendieck généralisé  $\text{Groth}(G(\A_f)\times W_{E_\nu})$ ) :
\vspace{4mm}

\begin{coro}
  Il y a une égalité dans $\text{Groth}(G(\A_f)\times W_{E_\nu})$ :
\begin{eqnarray*}
&\;& 
\sum_{i,j \atop { \Pi\in\mathcal{A}(I^\phi) \atop \Pi_\infty =\check{\rho}}}
(-1)^{i+j} \left [ \underset{K_p}{\limi}
  \text{Ext}^{\;i}_{J_b\text{-lisse}}
 \left ( H^j_c(\M_{K_p}\otimes \C_p,\Qlb)(N),\Pi_p\right ) \right ] 
\otimes 
\left [ \Pi^p \right ] \\
&=& \sum_i (-1)^i \left [ \underset{K}{\limi} H^i(
  \Sh_K^{an}(\phi),\LL_\rho^{an}) \right ] 
\end{eqnarray*}
\end{coro}

La décomposition 
$$
\Sh_K^{an}=\coprod_{\ker^1(\Q,G)} \Sh_K(G',X)^{an}
$$
induit une décomposition 
$$
sp^{-1} (\Sb(b_0))=\coprod_{\ker^1(\Q,G)} \Sh_K^{an}(G',X)\cap sp^{-1}(\Sb(b_0))
$$
où $ \Sh_K^{an}(G',X)\cap sp^{-1}(\Sb(b_0))$ est un ouvert de
$\Sh_K(G,X)^{an}$. 

\begin{defi}
Dans la décomposition précédente nous noterons 
$$
\Sh_K(G',X)^{an}_{\text{basique}}= sp^{-1}(b_0)\cap \Sh_K (G',X)^{an}
$$
\end{defi}

Rappelons (section \ref{str_bas}) que dans le cas que nous considérons
les groupes $G'$ sont tous isomorphes à $G$, d'où un isomorphisme 
des modèles canoniques sur le corps réflex $E(G,X)$ :
$$
\Sh_K(G,X)\iso \Sh_K(G',X)
$$
Les points géométriques dans $\overline{E}$ 
de $\Sh_K(G,X)\cap sp^{-1}(b_0)$ se décrivent comme étant associés à
des triplets $(A,\l,\iota)$ dont l'isocristal associé est
basique. L'isomorphisme ci dessus induit donc un isomorphisme
d'ouverts analytiques  
$$
\Sh_K(G,X)_{\text{basique}} \iso \Sh_K(G',X)_{\text{basique}}
$$
Donc, dans $\Groth (G(\A_f)\times W_E)$ on a :
$$
[H^\bullet_c (sp^{-1} (b_0),\LL_\rho^{an})]= |\ker^1(\Q,G) |. 
[H^\bullet_c (\Sh_K(G,X)_{\text{basique}},\LL_{\rho}^{an})]
$$

Étant donné (\ref{str_bas}) que tous les groupes $I^\phi$ sont
isomorphes et qu'il y  a $|\ker^1(\Q,G)|$ classes $\phi$ dans la
strate basique, en sommant l'égalité du corollaire précédent sur tous
les $\phi$ intervenant dans la strate basique et en divisant par 
$|\ker^1(\Q,G)|$ on obtient le corollaire suivant :

\begin{coro}\label{coro_Matsushima}
Choisissons la classe d'isogénie $\phi$ intervenant dans la strate basique. 
Il y a une égalité dans $\text{Groth}(G(\A_f)\times W_{E_\nu})$ :
\begin{eqnarray*}
&\;& 
\sum_{i,j \atop { \Pi\in\mathcal{A}(I^\phi) \atop \Pi_\infty =\check{\rho}}}
(-1)^{i+j} \left [ \underset{K_p}{\limi}
  \text{Ext}^{\;i}_{J_b\text{-lisse}}
 \left ( H^j_c(\M_{K_p}\otimes \C_p,\Qlb)(N),\Pi_p\right ) \right ] 
\otimes 
\left [ \Pi^p \right ] \\
&=& \sum_i (-1)^i \left [ \underset{K}{\limi} H^i(
  \Sh_K(G,X)^{an}_{\text{basique}},\LL_\rho^{an}) \right ] 
\end{eqnarray*}
\end{coro}

%%% Local Variables: 
%%% mode: latex
%%% TeX-master: "thesef"
%%% End: 

\part{Contribution de la cohomologie de la strate basique}
\chapter{Formule de Lefschetz sur la fibre spéciale}
 
Le but de ce chapitre est d'établir la formule des traces
\ref{Lefschetz_speciale} pour la trace de certaines fonctions dans
l'algèbre de Hecke globale fois un élément du groupe de Weil local $W_{E_\nu}$
sur la cohomologie des strates des variétés de Shimura de type
P.E.L. non ramifiées en $p$.   

L'un des points clef qui nous permettra d'utiliser cette formule est
que les termes intervenant dedans s'expriment comme un
produit de la trace d'une fonction de l'algèbre de Hecke en $p$ sur un
espace de cohomologie fois une intégrale orbitale hors $p$.
Cette factorisation des termes locaux dans la formule des traces de
Lefschetz est une conséquence des trois résultats suivants :
\begin{itemize}
\item Le théorème de Fujiwara \ref{theoFuji} qui sous certaines
conditions permet d'écrire ces
termes locaux comme la trace de certaines fonctions sur la fibre des
cycles évanescents.  
\item Le théorème de comparaison de Berkovich qui montre que cette
fibre ne dépend que du complété formel de la variété de Shimura en un
point fixe.
\item Le théorème de Serre-Tate (ou plus généralement l'uniformisation
de Rapoport-Zink) qui montre que ce complété formel ne dépend que de
la classe d'isomorphisme du groupe p-divisible en un point fixe et qui
permet donc de mettre en facteur la trace de notre fonction dans
l'algèbre de Hecke en $p$ sur un espace de déformation de groupes p-divisibles.
\end{itemize}

Pour appliquer ces trois résultats nous avons besoin de spécialiser
les correspondances de Hecke sur la fibre spéciale de modèles entiers
de nos variétés de Shimura. Malheureusement, à part en niveau
parahorique en $p$, on ne dispose pas de bon modèles entiers de ces
variétés. L'idée consiste alors à considérer le spécialisé au sens de
\cite{Fuj} de l'image
directe sur la fibre générique vers un niveau parahorique d'une
correspondance de Hecke. Cependant, nous avons ensuite besoin de
montrer que les termes locaux obtenus ne dépendent finalement que de
la situation originelle sur la fibre générique.
Afin de faire la lien avec la seconde partie nous avons également
besoin d'exprimer nos résultats sous forme analytique rigide. 
 C'est pourquoi nous développons
un formalisme de spécialisation des correspondances cohomologiques d'un
point de vue analytique rigide et montrons que l'isomorphisme de
Berkovich (\cite{Berk3}) est compatible à ce formalisme.

\section{Rappels sur les correspondances cohomologiques }\label{Rappels_coho}

Ce qui suit est principalement tiré des articles \citeg{Lef},
\citeg{Pink} et \citeg{Fuj}.

\subsection{Définition}

Soit $S$ le spectre d'un corps. Considérons un diagramme de schémas
sur $S$ 
\begin{diagram}[size=1cm]
 & &  C \\
 &  \ldTo^{c_1}(2,4) & \dTo^{c} & \rdTo(2,4)^{c_2} \\
 & & X_1 \times X_2 \\
 & \ldTo_{p_1} & & \rdTo_{p_2} \\
X_1 & & & & X_2
\end{diagram}
où $p_1$ et $p_2$ désignent les deux projections de $X_1\times X_2$ sur
$X_1$ et $X_2$.

 Nous supposerons toujours dans la suite que $c$ est
propre.

\begin{defi}
Si $(\F_1,\F_2)\in \D^b_c (X_1,\Qlb)\times \D^b_c(X_2,\Qlb)$ on note 
$$
\Cohf_C (\F_1,\F_2)=c_* \Homf(c_1^*\F_1,c_2^! \F_2)\in \D^b_c
(X_1\times X_2, \Qlb)
$$
le faisceau des correspondances cohomologiques à support dans $C$.

Nous noterons 
$$
\Coh_C(\F_1,\F_2)=H^0 (X,\Cohf_C(\F_1,\F_2))=\Hom (c_1^*\F_1,c_2^!\F_2)
$$
Lorsque $C=X_1\times X_2$ on notera $\Cohf$, resp. $\Coh$ pour $\Cohf_C$,
resp. $\Coh_C$
\end{defi}

Rappelons que la formule d'induction implique qu'il y a un
isomorphisme canonique :
$$
c_*c^! \Cohf(\F_1,\F_2)\iso \Cohf_C (\F_1,\F_2)
$$
L'adjonction $c_* c^!\drt Id$ donne une flèche d'oubli du support 
$$
\Cohf_C(\F_1,\F_2) \ldrt \Cohf(\F_1,\F_2) 
$$
et par application de $H^0(X_1\times X_2,-)$ : 
$$
\Coh_C(\F_1,\F_2) \ldrt \Coh(\F_1,\F_2) 
$$

\subsection{Seconde définition }

La formule de K{\"u}nneth pour $\Homf$ (\citeg{Del-finitude}) implique qu'il
existe un isomorphisme canonique en $\F_1$ et $\F_2$ :
$$
D \F_1 \boxtimes \F_2 \iso \Cohf (\F_1,\F_2)
$$

où $D\F=\Homf (\F, K_X)$ où $K_X$ désigne le faisceau dualisant sur $X$.

\subsection{Exemples}

\begin{itemize}
\item Si $\Delta_X\hookrightarrow X\times_S X$ désigne la diagonale de
$X\times X$, on a 
$$
\Coh_{\Delta_X} (\F_1,\F_2)\iso \Hom (\F_1,\F_2)
$$
En particulier, pour tout $\F$ on notera $\Delta$ la correspondance
cohomologique sur $\F$ associée à $Id\in\Hom (\F,\F)$.
\item Plus généralement, si $f:X_2\drt X_1$ et $C=\GG_f\subset
X_1\times X_2$ désigne le graphe de $f$ alors
$$
\Coh_C(\F_1,\F_2)\simeq \Hom (f^*\F_2,\F_1)
$$
\item A un cycle $C\subset X\times X$
 est associé une correspondance cohomologique de $\Ql$ dans
$c_2^!\Ql(-d)[2d]$ pour un $d\in \Z$ (cf. \citeg{Lef})
\item Soit $X_0/\Fq$ et $\F\in \mathbb{D}^+(X_0,\Qlb)$. Notons $X=X_0\times_{\Fq}
\Fqb$. Soit $Fr:X_0\ldrt X_0$ l'endomorphisme de Frobenius
relatif à $\Fq$. Il y a un isomorphisme canonique 
$$
Fr^*\F\iso \F
$$
d'où une correspondance sur $\F$ à support dans le graphe de
$Fr$. Par extension des scalaires cette correspondance donne une
correspondance appelée correspondance de Frobenius sur
$\overline{\F}=pr_1^*\F$ où $pr_1:X_0\times_{\Fq}
\Fqb \drt X_0$. Elle est donnée par 
$$
(Fr\times Id)^* \overline{\F} \ldrt \overline{\F}
$$
\end{itemize}

\subsection{Image directe}

Considérons un diagramme commutatif 

\begin{diagram}
 & & C \\
 & \ldTo^{c_1} & \dTo^{f'} & \rdTo^{c_2} \\
X_1 & & & & X_2 \\
\dTo^{f_1} &  & D & & \dTo^{f_2} \\
& \ldTo^{d_1} & & \rdTo^{d_2} \\
Y_1 & & & & Y_2
\end{diagram}

On peut alors définir (\citeg{Lef}) un morphisme 
$$
f'_* \Cohf_C (\F_1,\F_2) \ldrt \Cohf_D (f_{1!}\F_1,f_{2*}\F_2 )
$$
qui induit un morphisme image directe
$$
\Coh_C(\F_1,\F_2)\ldrt \Coh_D (f_{1!}\F_1,f_{2*}\F_2)
$$
noté $u\mapsto f_* u$.

\begin{rema}
Supposons $Y_1=X_1, Y_2=X_2$ et $D=X_1\times X_2$. On retrouve alors
le morphisme d'oubli du support.
\end{rema}

Lorsque $Y_1=Y_2=D=S$ on obtient l'action d'une correspondance
cohomologique sur la cohomologie. Plus précisément, il s'agit d'un
morphisme 
$$
\Coh_C(\F_1,\F_2)\ldrt \Hom_{\D^b_c(\Qlb)} (R\GG_c(X_1,\F_1),R\GG (X_2,\F_2))
$$
et l'on notera $u\mapsto u_*$ le morphisme associé dans ce cas là.

\subsection{Un cas particulier}

Supposons $c_1$ propre. Dans ce cas là, le morphisme du cas précédent
associé à l'action d'une correspondance cohomologique se factorise par 
$R \GG_c (X_2,\F_2)\drt R \GG (X_2, \F_2)$ et possède une description
simple.

Soit en effet $u:c_1^*\F_1 \ldrt c_2^! \F_2$. Par adjonction, la donnée
d'un tel $u$ équivaut à celle de $\tilde{u}:c_{2!} c_1^* \F_1 \ldrt
\F_2$. $u_*$ est alors  définie par le composée 
\begin{diagram}
R\GG_c (X_1,\F_1) & \rTo^{c_1^*} & R \GG_c (C,c_1^*\F_1) & \rTo^\sim &
R\GG_c(X_2,c_{2!} c_1^*\F_1) & \rTo^{R\GG_c(X_2,\tilde{u})} & R \GG_c (X_2,\F_2)
\end{diagram}
où la première application est définie car $c_1$ est propre.

Le fait que $c_1$ propre implique que la correspondance cohomologique
agit sur la cohomologie à support compact se comprend mieux du point
de vue de la section \ref{Ext_zero}.

\subsection{Restriction}

Soient $V_1\subset X_1$ et $V_2\subset X_2$ deux sous-schémas
localement fermés. Supposons que $c_2^{-1}(V_2)\subset
c_1^{-1}(V_1)$, d'où une ``sous-correspondance'' de $c$ 
\begin{diagram}
c_2^{-1}(V_2) & \rTo^{(c'_1,c'_2)} & V_1\times V_2
\end{diagram}
Soit donc $u\in \Coh_C (\F_1,\F_2)$ et $\tilde{u}:c_{2!}c_1^*\F_1 \drt
\F_2$ le morphisme associé par adjonction.

Considérons la ligne suivante 
\begin{diagram}
c'_{2!}{c'_1}^* (\F_{1|V_1})& \rEqualto & c'_{2!} \left (
(c_1^*\F_1)_{|c_2^{-1}(V_2)}\right ) & \rTo^\sim & (c_{2!}c_1^* \F_1)_{|V_2} &
\rTo^{\tilde{u}_{|V_2}} & \F_{2|V_2}
\end{diagram}
où l'isomorphisme du milieu est un isomorphisme de changement de base. 
Elle définit par adjonction une correspondance restreinte. On en
déduit un morphisme de restriction
$$
\Coh_C (\F_1,\F_2) \ldrt \Coh_{c_2^{-1}(V_2)}(\F_{1|V_1},\F_{2|V_2})
$$
Dans le cas où $X_1=X_2$, $V_1=V_2$ nous dirons que $V_1$ est stable
par $c$. 

\subsection{Extension par zéro}\label{Ext_zero}

Supposons $c_1$ propre. Supposons nous donné des
compactifications $j_1:X_1\hookrightarrow\overline{X_1}$,
$j_2:X_2\hookrightarrow\overline{X_2}$ ainsi qu'une compactification
$j:C\hookrightarrow \overline{C}$ de $C$
compatible aux deux précédentes, au sens où il y a un diagramme
commutatif 
\begin{diagram}
 C & \rInto^j & \overline{C} \\
\dTo^c & & \dTo^{(\overline{c}_1,\overline{c}_2)=\overline{c}} \\
X_1\times X_2 & \rInto^{j_1\times j_2} & \overline{X_1} \times
\overline{X_2} 
\end{diagram}
Il existe alors un isomorphisme canonique (\citeg{Fuj} lemme 1.3.1) 
$$
j_* \Cohf_C(\F_1,\F_2) \iso \Cohf_{\overline{C}} (j_{1!} \F_1,j_{2!} \F_2)
$$
qui induit un isomorphisme 
$$
\Coh_C (\F_1,\F_2) \iso \Coh_{\overline{C}}(j_{1!} \F_1,j_{2!} \F_2)
$$
noté $u\mapsto j_! u$ et qui se décrit ainsi : $\forall u\in \Coh_C
(\F_1,\F_2)$, $u:c_1^* \F_1 \drt c_2^!\F_2$. $c_1$ étant propre il
existe un morphisme 
$$
\overline{c_1}^*j_{1!} \F_1 \ldrt j_! c_1^* \F_1
$$
défini par adjonction.

Alors, 
\begin{diagram}
j_!u : & \overline{c_1}^*j_{1!} \F_1 & \rTo & j_! c_1^* \F_1 & \rTo^{j_!
(u )} & j_! c_2^! \F_2 & \rTo & \overline{c_2}^! j_{2!} \F_2
\end{diagram}
où le dernier morphisme est défini par adjonction.
\\

L'extension par zéro est compatible aux images directes au sens où
``$j_! f_* u =\overline{f}_* j_! u$''.

En particulier, l'action de $u$ sur la cohomologie à support compact
s'interprète comme l'action de $j_! u$ : le diagramme suivant commute 
\begin{diagram}
 R\GG_c (X_1,\F_1) & \rTo^{u_*} & R\GG_c (X_2,\F_2) \\
\dTo^{\wr} & & \dTo_{\wr} \\
R\GG (\overline{X_1},j_{1!}\F_1) & \rTo^{(j_!u)_*} & R\GG
(\overline{X_2},j_{2!} \F_2)
\end{diagram}

\subsection{Accouplement des correspondances cohomologiques }

Nous supposerons désormais que $X_1=X_2=X$. Considérons un diagramme
cartésien 
\begin{diagram}
& & E=C\times_X D \\
 & \ldTo & \dTo^e & \rdTo \\
C & & &&  D \\
& \rdTo_{c} & &  \ldTo^d \\
&& X\times_S X
\end{diagram}
i.e. $E$ est le support intersection de $C$ et $D$. 

Il y a alors un morphisme (\citeg{Lef} paragraphe 4) d'accouplement des
correspondances cohomologiques 
$$
\Cohf_C(\F_1,\F_2)\otimes \Cohf_D (\F_2,\F_1) \ldrt e_* K_E
$$
qui induit une forme bilinéaire 
\begin{eqnarray*}
\Coh_C(\F_1,\F_2) \otimes \Coh_D (\F_2,\F_1) &\ldrt & H^0 (E,K_E) \\
u\otimes v & \longmapsto & <u,v>
\end{eqnarray*}

\subsection{Termes locaux de Lefschetz} 

Il y a un isomorphisme 
$$
H^0 (E,K_E) \iso \bigoplus_{\b\in \pi_0 (E) } H^0 (\b,K_\b)
$$
et si $\b\in \pi_0 (E)$, $\dim \b=0$, alors $K_\b=\Qlb$ et donc 
$H^0(\b,K_\b)\iso\Qlb$.
On note alors $<u,v>_\b\in \Qlb$ l'image de $<u,v>$ par les deux
applications ci dessus. 

Ces termes locaux sont invariants par localisation étale. Cependant,
en général, ils ne dépendent pas que de la fibre des faisceaux en des
points fixes.

\subsection{Formule de Lefschetz Verdier Grothendieck}

Elle dit simplement que l'accouplement des correspondances
cohomologiques est compatible aux images directes propres. 

Plus précisément, étant donné un diagramme commutatif dont les faces
supérieures et inférieures sont cartésiennes 
\begin{diagram}[size=8mm]
 & & C& & \lTo && E \\
& \ldTo  & \vLine &&& \ldTo  & \dTo_b \\
X\times X &&  \lTo &\HonV & D   \\
\dTo~{f=f_1\times f_2} &&  \dTo & & \dTo \\
  & & C' &  \lTo & \VonH &  \hLine  & E' \\
 & \ldTo &&  &  & \ldTo \\
Y \times Y && \lTo & & D'
\end{diagram}
où $f_1$ et $f_2$ sont propres. Alors, le diagramme suivant commute 
\begin{diagram}
\Coh_C(\F_1,\F_2)\otimes \Coh_D(\F_2,\F_1) & \rTo & H^0 (E,K_E) \\
\dTo^{f_{1*}\otimes f_{2*} } & & \dTo_{b_*} \\
\Coh_{C'}(f_{1*}\F_1,f_{2*}\F_2)\otimes
\Coh_{D'}(f_{2*}\F_2,f_{1*}\F_1) & \rTo & H^0(E',K_{E'})
\end{diagram}
soit : $$<f_*u,f_*v>=b_*<u,v>$$
\\

En particulier, si $X$ est propre sur $S$, si $(u,v)\in
\Coh_C(\F_1,\F_2)\times \Coh_D(\F_2,\F_1)$ et si $\dim E=0$ alors, 
$u_*\in \Hom( R\GG (X_1,\F_1),R\GG (X_2,\F_2))$ et $v_*\in 
\Hom( R\GG (X_2,\F_2),R\GG (X_1,\F_1))$ et 
$$
\text{Tr} (v_*\circ u_*)=\sum_{\b \in \pi_0 (E)} <u,v>_\b
$$
\\

Supposons de plus que $\F_1=\F_2=\F$ et que $v=\Delta$ la
correspondance diagonale. 
\begin{defi}
Nous noterons
$$\text{Fix}(c)=C\times_{X\times X} \Delta_X$$
et si $\beta\in \pi_0 (\text{Fix} (c))$, $\dim \b=0$ ($\b$ est un point fixe
isolé) 
$$
\text{Loc}_\b (u)=<u,\Delta>_\beta \in \Qlb 
$$
\end{defi}

Ainsi, lorsque tous les points fixes sont isolés 
$$
\text{tr}(u, R\GG (X,\F))=\sum_{\b\in \pi_0(\text{Fix}(c))}
\text{Loc}_\b (u)
$$

\subsection{Termes locaux na{\"\i}fs} 

Soit $u\in \Coh_C (\F,\F)$ et soit $\b\in \pi_0 (\text{Fix} (c)), \dim
\b=0$. Soit $x=c_1 (\b)=c_2 (\b)\in X$ et supposons $c_2$ quasi-fini
au voisinage de $x$. 

Le morphisme 
$u:c_1^*\F \drt c_2^! \F$ induit par adjonction $\tilde{u}:c_{2!} c_1^* \F \drt \F$, et
$c_2$ étant quasi-fini au voisinage de $x$ 
$$
(c_{2!}c_1^*\F)_x \simeq \bigoplus_{y\in c_2^{-1} (x)} (c_1^*\F)_y 
$$
D'où un morphisme $\F_x = (c_1^*\F)_\b \hookrightarrow
(c_{2!}c_1^*\F)_x$ qui composé avec $\tilde{u}_x$ donne un
endomorphisme 
$$
u_\b\in \End (\F_x)
$$
\begin{defi}
On pose 
$$
\text{Loc.na{\"\i}f}_\b (u)=\text{Tr} (u_\b)\in \Qlb$$
\end{defi}

\subsection{Interprétation des termes locaux na{\"\i}fs comme de vrais
termes locaux}

Bien que $\text{Loc}_\b (u)$ soit invariant par localisation étale, il
ne dépend pas que de la limite $\F_{c_1(\b)}$ i.e. en général 
$\text{loc}_\b \neq \text{Loc.na{\"\i}f}_\b$. 
\\

Dans \citeg{Pink} il est démontré que si $\b$ vérifie les hypothèses
précédentes alors le calcul de $\text{Loc}_\b (u)$ se ramène à la
situation suivante :
$
j:X\setminus \{ x \}\hookrightarrow X
$
est stable par $c$. Et qu'alors, 
$$
\text{Loc}_\b (u)= \text{Loc}_\b (j_! j^* u)+\text{Loc.na{\"\i}f}_\b (u)
$$
Le premier terme dans le membre de droite s'interprétant donc comme
un terme d'erreur par rapport au terme local na{\"\i}f. 

Le calcul de $\text{Loc}_\b (u)$ se ramène donc au cas où la fibre du
faisceau est nulle.

\subsection{Exemples}

\begin{itemize}
\item Si $X$ est lisse, $\F$ localement constant et l'intersection en
$\b$ de $C$ et $\Delta$ est transverse 
$$
\text{Loc}_\b (u)= \text{Loc.na{\"\i}f}_\b (u) 
$$
\item Si $X$ est lisse, si $u$ est la correspondance sur le faisceau
constant associé à un cycle, si $\b$ est isolé, l'application classe
de cycle transformant le produit d'intersection des cycles en
cup-produit,
$$
\text{Loc}_\b (u)=(C.\Delta)_\b
$$

\end{itemize}

\subsection{Un lemme sur les termes locaux na{\"\i}fs}

Ce lemme nous servira plus tard dans la démonstration du théorème
\ref{Lefschetz_speciale}.

\begin{lemm}\label{ramene}
Soit 
\begin{diagram}[size=8mm,labelstyle=\scriptstyle]
C & \rTo^{(c_1,c_2)} & X\times X \\
\dTo^{f'} & & \dTo^{f\times f} \\
D & \rTo^{(d_1,d_2)} & Y\times Y
\end{diagram}
un morphisme de correspondances où les $c_i,d_i (i=1,2)$ sont finis et
où $f'$ et $f$ sont étales finis. Supposons que $\dim
\text{Fix}(D)=0$, et que 
$$
\forall \G \;\forall u\in \Coh_D (\G)\;\; \forall \b\in \text{Fix}(D) \;\;\;
\text{Loc}_\b (u)=\text{Loc.na{\"\i}f}_\b (u)
$$
Alors, $\dim \text{Fix}(C)=0$, et $$\forall \F\;\forall u\in
\Coh_C(\F) \;\forall \b'\in  \text{Fix}(C)\;\;\; \text{Loc}_\b
(u)=\text{Loc.na{\"\i}f}_\b (u) 
$$
\end{lemm}
\begin{proof}
L'assertion $\dim \text{Fix}(C)=0$ 
 résulte du fait que $f'$ envoie les points fixes de $C$
sur ceux de $D$ 
et de la finitude des fibres de $f'$. 

Quant à la seconde assertion elle résulte facilement de l'invariance
par localisation étale des termes locaux de Lefschetz appliquée à
$f_*\F$ muni de la correspondance cohomologique $f_* u$ et de la
définition des termes locaux na{\"\i}fs. 
\end{proof}

\section{Le théorème de Fujiwara}

Le théorème qui suit, démontré par Fujiwara dans \citeg{Fuj}, avait
été conjecturé par Deligne. Des résultats partiels avaient été obtenus
par Zink (\citeg{Zink_Lef}) et Pink (\citeg{Pink}).

\begin{theo}\label{theoFuji}
Soit \begin{diagram}
C_0 & \rTo^{(c_1,c_2)} & X_0
\end{diagram} une correspondance définie sur $\Fq$. Soit $k|\Fq$ et notons 
$X=X_0\times_{\Fq} k$, $C=C_0\times_{\Fq} k$ et $Fr$ la correspondance
de Frobenius sur $X_0$ étendue à $X$.

Supposons $c_2$ quasi-fini et notons 
$$
\forall \b\in C\;\;\;\; d(\b)=[k(\b):k(c_2(\b))]_{insep.}
. \text{long} \left ( \O_{c_2^{-1}(c_2(\b)),\b}\right )
$$
Alors, $\forall N\in \N \;\; \forall \b\in \text{Fix}(C.Fr^N)$,
$q^N>d(\b) \limpl \b$ est un points fixe isolé et 
$$
\forall \F \; \forall u\in \Coh_{C.Fr^N}(\F,\F) \;\;\; \text{Loc}_\b
(u)=\text{Loc.na{\"\i}f}_\b (u)
$$
\end{theo}

\section{Spécialisation des correspondances cohomologiques
algébriques}

Soit $S=\spec (\O)$ un trait de point générique $\eta$ et de point
fermé $s$. Nous noterons $\La$ un anneau de coefficients de torsion
première à $p$ ou bien $\Qlb$. 

Nous noterons $\Rpsi$ (ou souvent, afin de na pas alourdir
les notations, $\rpsi$) le foncteur des cycles évanescents. De m{\^e}me,
pour un morphisme $f$ au dessus de $S$ nous noterons $f_{\bar{s}}$ sa
fibre spéciale étendue à $\bar{s}$ et $f_\eta$ sa fibre générique
(mais souvent, lorsqu'il n'y a pas d'ambigu{\"\i}té, nous noterons
$f$ pour $f_{\bar{s}}$ et $f_\eta$). 

\subsection{Règles de calculs sur les cycles évanescents}\label{RegPsi}

\subsubsection{Image directe}

Soit $f:X\ldrt Y$ au dessus de $S$. Il y a alors un morphisme naturel
(\citeg{Cycles} 2.1.7.)  
$$
\rpsi f_* \ldrt f_* \rpsi
$$
qui est un isomorphisme lorsque $f$ est propre. Ce morphisme est bien
s{\^u}r compatible à la composition :
\begin{diagram}[size=9mm]
\rpsi ( (fg)_*) & \rEqualto & \rpsi (f_* g_*) & \rTo & f_* \rpsi (g_*)
\\
\dTo & & & & \dTo \\
(fg)_* \rpsi & \rEqualto & & & f_* g_* \rpsi 
 \end{diagram}

est un diagramme commutatif.

\subsubsection{Image inverse}
Pour $f:X\ldrt Y$ il y également un morphisme naturel (\citeg{Cycles} 2.1.7.)
$$
f^* \rpsi \ldrt \rpsi f^* 
$$
qui est un isomorphisme lorsque $f$ est lisse. Comme ci dessus ce
morphisme est compatible à la composition vis à vis de l'égalité  $(f
g)^*=g^* f^*$.  

\subsubsection{Compatibilité au changement de base}

\begin{lemm}\label{RpsiBC}
Soit 
\begin{diagram}[size=9mm]
X' & \rTo^{f'} & Y' \\
\dTo^{g'} & & \dTo^g \\
X & \rTo^f & Y 
\end{diagram}
un diagramme sur $S$. Pour $\F\in \D^+ (Y',\Lambda)$ considérons le
morphisme de changement de base 
$$
f^* g_* \F \xrightarrow{\; \a (\F)\; } g'_* {f'}^* \F
$$
Le diagramme suivant commute :
\begin{diagram}[size=9mm]
f^* g_* \rpsi (\F) & \rTo^{\a (\rpsi (\F))} & g'_* {f'}^* \rpsi (\F) \\
\uTo & & \dTo \\
f^* \rpsi (g_* \F) & & g'_*\rpsi ({f'}^*\F) \\
\dTo & & \uTo \\
\rpsi (f^*g_* \F) & \rTo^{\rpsi (\a (\F))} & \rpsi (g'_* {f'}^* \F)
\end{diagram}
où les applications verticales sont déduites des fonctorialités
précédentes par image directe et image inverse.
\end{lemm}
\begin{proof}
Elle ne pose aucun problème.
\end{proof}

\subsection{Spécialisation}
Considérons un diagramme au dessus de $S$ : 
\begin{diagram}[size=8mm]
 & & C \\
 & \ldTo^{c_1} & & \rdTo^{c_2} \\
X_1 & & & & X_2
\end{diagram}
où comme d'habitude $c:C\ldrt X_1\times X_2$ est propre. 

Fujiwara a défini dans \citeg{Fuj} un morphisme de spécialisation des
correspondances cohomologiques 
$$
\Coh_{C_\eta} (\F_1,\F_2) \ldrt \Coh_{C_{\bar{s}}} (\Rpsi
(\F_1),\Rpsi (\F_2))
$$
La définition de ce morphisme utilise l'existence de l'isomorphisme
$\rpsi D \simeq D \rpsi$ qui se transpose mal au cadre de la géométrie
rigide. Afin d'avoir une description simple qui se transporte dans le
cadre de la géométrie rigide nous allons décrire le morphisme de
spécialisation dans un cadre plus restreint.
Nous ferons pour cela l'hypothèse suivante :
\\

{\it Hypothèse : } $c_1$ et $c_2$ sont propres.
\\

Soit donc $u\in \Coh_{C_\eta}(\F_1,\F_2)$, $u:c_{1\eta}^* \F_1 \drt
c_{2\eta}^! \F_2$. $c_2$ étant propre il lui est associé 
$\tilde{u}: c_{2\eta *} c_{1\eta }^* \F_1 \drt \F_2$.
D'où 
$$
\Rpsi (\tilde{u}) : \Rpsi (c_{2\eta *} c_{1\eta }^*\F_1) \ldrt \Rpsi (\F_2)
$$
Appliquant les diverses fonctorialités de $\rpsi$ et la propreté de
$c_2$ on obtient alors la ligne suivante : 
\begin{diagram}
c_{2\bar{s}_*} c_{1\bar{s}}^*\Rpsi (\F_1) & \rTo &
c_{2\bar{s}*}\Rpsi (c_{1\bar{s}}^*\F_1) & \lTo^\sim & 
\Rpsi (c_{2\bar{s}*}c_{1\bar{s}}^*\F_1) & \rTo^{ \Rpsi (\tilde{u})} &
\Rpsi (\F_2) 
\end{diagram}
qui fournit par
 adjonction un $u_{\bar{s}}\in \Coh_{C_{\bar{s}}}(\Rpsi
(\F_1),\Rpsi (\F_2))$ dont on peut vérifier qu'il co{\"\i}ncide avec celui
défini par Fujiwara (nous n'utiliserons pas ce fait par la suite).  
\\

La proposition qui suit est démontrée par Fujiwara (proposition 1.6.1
de \citeg{Fuj}). Néanmoins nous aurons besoin de démontrer son
analogue pour les correspondances rigides. Sa démonstration ne
s'adaptant pas dans le cadre rigide, nous  donnons une démonstration
(sous nos hypothèses plus restrictives) 
qui n'utilise pas la dualité $D$ et le fait que $\rpsi$ commute à
cette dualité. 

\begin{rema}
  Si $\mathfrak{X}$ désigne un schéma formel sur $S$ de la forme
  $\widehat{X}_{/V}$ où $X$ est un schéma  séparé de type fini sur
  $S$, et $V$ un sous schéma localement fermé de sa fibre spéciale
  alors, il résulte de l'isomorphisme de comparaison de Berkovich
et de la commutation des cycles évanescents algébriques à la dualité  que
  si $\F$ désigne un faisceau étale constructible de torsion première
  à $p$ sur $X_{\et}$,  il y a un isomorphisme 
$$
D\Rth ( \F^{an}_{|\X^{an}}) \simeq \Rth ( D( \F^{an}_{|\X^{an}}) )  
$$
(on utilise également la commutation de $D$ au foncteur $\F\mapsto
\F^{an}$, confère \cite{Berk5}). Néanmoins, cet isomorphisme dépend à
priori du choix de la réalisation de $\X$ comme le complété formel
d'un $X$ le long d'un $V$ et n'est donc pas formulé de façon
intrinsèque. 
\end{rema}

\begin{prop}\label{speima}
Soit un diagramme sur $S$
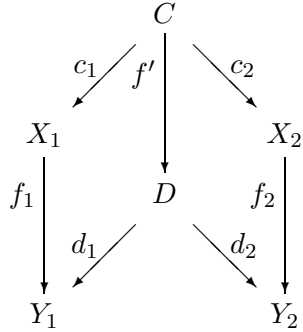
 \begin{diagram}[size=8mm]
 & & C \\
 & \ldTo^{c_1} & \dTo^{f'} & \rdTo^{c_2} \\
X_1 & & & & X_2 \\
\dTo^{f_1} &  & D & & \dTo^{f_2} \\
& \ldTo^{d_1} & & \rdTo^{d_2} \\
Y_1 & & & & Y_2
\end{diagram}
où $f_1,f_2,c_1,c_2,d_1,d_2$ sont propres. Alors, le diagramme suivant commute
:
\begin{diagram}[size=9mm]
\Coh_{C_{\eta}} (\F_1,\F_2) & \rTo & \Coh_{C_{\bar{s}}} (\Rpsi
(\F_1),\Rpsi (\F_2)) \\
\dTo^{f_{\eta *}} & & & \rdTo{f_{\bar{s} *}} \\
 & & & & \Coh_{D_{\bar{s}}}(f_{1\bar{s} *}\Rpsi (\F_1), f_{2\bar{s} *} \Rpsi (\F_2))
\\
 & & & \ruTo^\simeq \\
\Coh_{D_{\eta}} ( f_{1 \eta *}\F_1, f_{2\eta * } \F_2) & \rTo &
\Coh_{D_{\bar{s}}} (\Rpsi (f_{1\eta * }\F_1), \Rpsi (f_{2 \eta *}\F_2)
\end{diagram}

\end{prop}
\begin{proof}
Soit $u\in \Coh_{C_\eta} (\F_1,\F_2)$. 
Appliquons le lemme \ref{RpsiBC} au carré de gauche et au morphisme de
changement de base $d_1^*f_{1*} \ldrt f'_* c_1^* $. On obtient un
diagramme commutatif 
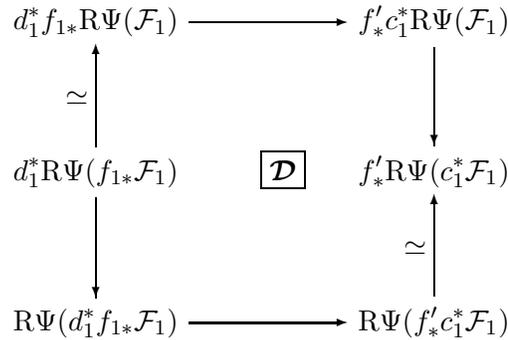
\begin{diagram}[size=1cm]
d_1^* f_{1*} \rpsi (\F_1) & \rTo && & f'_* c_1^* \rpsi (\F_1) \\
\uTo^\simeq & && & \dTo \\
d_1^* \rpsi (f_{1*} \F_1) & & \boxed{\boldsymbol{\mathcal{D}}} & & f'_*\rpsi
(c_1^*\F_1) \\
\dTo & & && \uTo^\simeq \\ 
\rpsi (d_1^* f_{1*} \F_1) & \rTo && & \rpsi (f'_* c_1^*\F_1) 
\end{diagram}

Considérons maintenant le diagramme 

\begin{diagram}[size=1cm]
d_{2*} d_1^* f_{1*} \rpsi (\F_1) & \rTo && & d_{2*} f'_* c_1^* \rpsi
(\F_1) & \rEqualto & f_{2*} c_{2*}c_1^* \rpsi (\F_1) &
\rTo^\b & f_{2*} \rpsi (\F_2) \\
\uTo^\simeq & && & \dTo & & \dTo & \ruTo(2,4)~{\scriptstyle{\a}} & \uTo^\simeq  \\ 
d_{2*}d_1^* \rpsi (f_{1*} \F_1) & & 
\boxed{d_{2*}\boldsymbol{\mathcal{D}}} & & d_{2*}f'_*\rpsi 
(c_1^*\F_1) & \rEqualto & f_{2*} c_{2*} \rpsi (c_1^*\F_1)   \\
\dTo & & & &  \uTo^\simeq & & \uTo^\simeq  \\
d_{2*} \rpsi (d_1^* f_{1*} \F_1) & \rTo && & d_{2*}\rpsi (f'_*
c_1^*\F_1)  & & f_{2*} \rpsi (c_{2*} c_1^*\F_1)  \\
\uTo^\simeq & &&& \uTo^\simeq & & \uTo^\simeq \\
\rpsi (d_{2*} d_1^* f_{1*} \F_1) & \rTo & & & \rpsi (d_{2*} f'_*
c_1^*\F_1) & \rEqualto & \rpsi (f_{2*} c_{2*} c_1^* \F_1) & \rTo &
\rpsi (f_{2*} \F_2) \\
\end{diagram}
\\

Le carré en haut à gauche est obtenu en appliquant $d_{2*}$ à
$\mathcal{D}$. Il commute donc. 
L'application $\a$ est égale à $f_{2*} \rpsi (\tilde{u})$ (où
$\tilde{u} : c_{2*} c_1^* \F_1 \drt \F_2$). L'application $\b$ est
définie à partir de $\a$ de telle manière que le triangle supérieur commute. 

L'autre triangle (ou plut{\^o}t trapèze) et tous les autres rectangles
 dans le diagramme commutent 
 car ils sont tous associés à des cas de
fonctorialité par image directe ou inverse des cycles évanescents
 (cf. \ref{RegPsi}).  

On en déduit donc que le grand rectangle extérieur  commute. 

Le résultat s'en déduit car en parcourant la ligne du haut on obtient
l'image directe de la correspondance spécialisée tandis qu'en joignant
les deux sommets du haut en passant par les trois autres cotés on
obtient la spécialisation de l'image directe. 

\end{proof}

\begin{coro}
Si $X_1/S$ et $X_2/S$ sont propres alors la suite spectrale des cycles
évanescents est équivariante pour l'action de $u_{\bar{s}}$ sur le
terme initial et de $u$ sur l'aboutissement : le diagramme suivant
commute
\begin{diagram}[size=8mm]
R\GG (X_{1\bar{s}},\Rpsi (\F_1)) & \rTo^\sim & R\GG
(X_{1\bar{\eta}},\bar{\F}_1) \\
\dTo^{(u_{\bar{s}})_*} & & \dTo^{\bar{u}_*} \\
  R\GG (X_{2\bar{s}},\Rpsi (\F_2)) & \rTo^\sim & R\GG
(X_{2\bar{\eta}},\bar{\F}_2)
\end{diagram}
\end{coro}

\section{Cycles évanescents analytiques rigides}
\subsection{Rappels} \label{Rappels_Berk}

Soit $\X$ un schéma formel localement formellement de type fini sur $\spf(\O)$ de
fibre générique $\X^{an}$ et de fibre spéciale $\X_s$. 

Il y a un isomorphisme de sites :
$$
(\X_s)_{\et} \iso 
\X_{\et} 
$$
Si de plus $U/\X$ est étale, $U^{an}/\X^{an}$ est quasi-étale. On a
donc un diagramme de sites 
\begin{diagram}[size=1cm]
(\X_s)_{\et} & \rEqualto & \X_{\et} & \lTo^\nu & \X^{an}_{q-\et} &
\rTo^\mu & \X^{an}_{\et}
\end{diagram}
\\

Rappelons (m{\^e}me si nous ne l'utiliserons pas dans la suite) que
$\X^{an}_{q-\et}\simeq \X^{rig}_{\et}$ et que $\mu^* :
\widetilde{\X^{an}_{\et}} \ldrt \widetilde{\X^{an}_{q-\et}}$ est
pleinement fidèle d'image essentielle les faisceaux surconvergents sur
$\X^{rig}_{\et}$. 
\\

Berkovich pose alors (\citeg{Berk2},\citeg{Berk3}) 
$$
\Theta =\nu_* \mu^* : \widetilde{(\X^{an})_{\et}} \ldrt \widetilde{(\X^{s})_{\et}}
$$

Si $\hat{\bar{\eta}}$ désigne le complété de $\bar{\eta}$ il y a
également un
foncteur 
$$
\Theta_{\bar{\eta}} :  \widetilde{(\X^{an})_{\et}} \ldrt
\widetilde{(\X^{an}\hat{\otimes}{\hat{\bar{\eta}}})_{\et}} \ldrt
\widetilde{(\X_{\bar{s}})_{\et}} 
$$
obtenu grâce à l'extension des scalaires de $\eta$ à
$\hat{\bar{\eta}}$.

Soit $\Lambda$ un anneau local artinien de torsion première à
$p$. Nous considérerons 
$$
\Rth : \D^+ ( \X^{an},\Lambda) \ldrt \D^+ ( \X_{\bar{s}}, \Lambda)
$$
$\Rth (\F)$ est muni d'une action de $\Gal (\bar{\eta}| \eta)$ compatible à
celle sur $\X_{\bar{s}}$. 

Rappelons le théorème de Berkovich 

\begin{theo}[\citeg{Berk3}] Soit $X/S$ un schéma de type fini et
$Y\subset X_s$ un sous-schéma localement fermé de sa fibre
spéciale. Pour tout faisceau étale abélien constructible $\F$ de torsion première
à $p$, $\forall q\in \N$ il existe un isomorphisme canonique :
$$
\text{R}^q \Psi_{\bar{\eta}} (\F)_{|Y}\iso \text{R}^q \Theta_{\bar{\eta}} (\widehat{\F}_{/Y})
$$
où $\widehat{\F}_{/Y}$ désigne la restriction de $\F^{an}$ sur
$X_\eta^{an}$ au domaine analytique $(\widehat{X}_{/Y})^{an} \hookrightarrow
X_\eta^{an}$. 
\end{theo}

Rappelons comment est construit ce morphisme. Commen{\c c}ons par construire un
morphisme 
$$
\Psi_{\bar{\eta}} (\F)_{|Y} \ldrt \Theta_{\bar{\eta}} (\widehat{\F}_{/Y})
$$
$\Psi_{\bar{\eta}} (\F)_{|Y}$ est le faisceau associé au préfaisceau qui
à un $U/Y_{\bar{s}}$ étale associe $\underset{\mathcal{D}}{\limi}
\bar{\F} (V_\eta)$ où la limite inductive est prise selon les
diagrammes $\mathcal{D}$ suivants 
\begin{diagram}[size=8mm]
U & & \rTo & & V \\
\dTo & & & & \dTo^{\text{étale}} \\
Y_{\bar{s}} & \rInto & X_{\bar{s}} & \rTo & \bar{X} 
\end{diagram}

$ \Theta_{\bar{\eta}} (\widehat{\F}_{/Y})$ est le faisceau qui à $U/Y_{\bar{s}}$
étale associe $\GG(\tilde{U}^{an},\bar{\F}^{an})$ où
$\tilde{U}/\widehat{X}_{/Y}$ est l'unique relèvement étale de $U$ et
$\tilde{U}^{an}/X_{\hat{\bar{\eta}}}^{an}$ est quasi-étale (rappelons
que si $A\xrig{p} X^{an}$ est quasi-étale et $\G$ est un faisceau  sur $X^{an}_{\et}$,
$\GG(A_{\et},p^* \G)= \GG(A_{q-\et},\mu^* \G)$). 

Soit donc maintenant $\mathcal{D}$ un diagramme comme précédemment et
$s\in \GG(V_\eta,\bar{\F})$. Le diagramme $\mathcal{D}$ et la
propriété de relèvement unique des morphismes étales via à vis des
immersions nilpotentes implique qu'il existe un unique morphisme
\begin{diagram}[size=9mm]
\tilde{U} & \rTo^a & \widehat{V}_{/Y} \\
\dTo & \ldTo \\
\widehat{X}_{/Y} 
\end{diagram}
d'où un morphisme 
\begin{diagram}[labelstyle=\scriptstyle]
\GG(V_\eta,\bar{\F}^{an}) & \rTo^{\text{restriction}} &
\GG(\widehat{V}_{/Y}^{an}, \bar{\F}^{an}) & \rTo^{a^*} & \GG
(\tilde{U}^{an},\bar{\F}^{an}) 
\end{diagram}
L'image de $s$ par ce composé est l'image par le morphisme que l'on
voulait décrire au niveau des préfaisceaux. Il induit le morphisme
cherché au niveau des faisceaux. 
\\

Si $\a:Y_{\bar{s}}\hookrightarrow X_{\bar{s}} $ et $ \b :
(\widehat{X}_{/Y})^{an} \subset (X_\eta)^{an}$ il y a donc un morphisme 
$$
\a^* \Psi_{\bar{\eta}} \ldrt \Theta_{\bar{\eta}} \b^*
$$
D'où un morphisme  
$$
\a^* \Rpsi \ldrt \Rth \b^*
$$
Et le théorème de Berkovich implique qu'il s'agit d'un isomorphisme
dans la catégorie dérivée $\D^+_c(X_{\bar{s}},\Lambda)$.

\subsection{Comparaison avec les cycles évanescents adiques}

  A un schéma formel $\X/S$ séparé localement formellement de type fini est
  associé un espace adique $\X^{rig}$ noté $\tilde{d} (\X)$ par
  Huber (confère l'appendice \ref{an_ad}). Celui ci est muni d'un
  morphisme de spécialisation $\tilde{\l}: \X\ldrt \X_s$.
  Ce morphisme induit un morphisme de topos
$$
(X^{rig})_{\et} \ldrt \X_{\et}
$$
Avec les notations précédentes, les faisceaux étales abéliens sur $ X^{rig}$
s'identifient aux faisceaux sur le site quasi-étale de $\X^{an}$ et ce
morphisme de sites n'est rien d'autre que le morphisme $\nu$ défini précédemment. 
Huber
  définit alors les cycles évanescents adiques  (confère par exemple 
\cite{Hu3} 3.12) pour  un faisceau $\F$ sur $(\X^{rig})_{\et}$ comme
étant $\R \tilde{\l}_* \F$. Du point de vue analytique, les cycles
évanescents adiques ne sont  donc rien d'autre que les $\R \nu_* \F$. 

Les cycles évanescents adiques et analytiques coïncident au sens où
l'on a avec les identifications précédentes entre le site quasi-étale
analytique et le site étale adique l'égalité 
$$
\forall \F\in \D^+ (\X^{an},\La) \;\; \rth (\F) = \R\tilde{\l}_* (
\mu^* \F) 
$$
Il suffit pour cela de montrer que $\R (\nu_* \mu^*)= \R\nu_* \circ
\mu^*$. Or, il résulte par exemple du théorème 3.3 de \cite{Berk2}
que $\mu^*$ envoi les faisceaux mous sur des faisceaux
$\nu_*$-acycliques, d'où le résultat. 

Nous n'utiliserons que le point de vue des cycles évanescents de
Berkovich, mais faisons 
remarquer au lecteur que Huber possède des théorèmes de comparaison
et de finitude plus généraux que ceux déduits du théorème de
comparaison de Berkovich (proposition 3.11
et proposition 3.15 de \cite{Hu3}).

\subsection{Fonctorialité de $\Rth$ et de l'isomorphisme de
comparaison}

\subsubsection{Image directe}

Soit $f:\X \ldrt \Y$ un morphisme de schémas formels formellement de
type fini sur $S$. Il y a alors un isomorphisme  naturel
(\citeg{Berk3}, corollaire 2.3 (ii))
$$
\Rth f^{an}_* \iso f_{\bar{s}*} \Rth
$$
D'où en particulier une suite spectrale des cycles évanescents pour
tout $\X$ : 
$$
R\GG (\X_{\bar{s}},\Rth (\F)) \simeq R\GG( \X^{an}\hat{\otimes}
\hat{\bar{\eta}} , \bar{\F})
$$

\begin{lemm}{\label{Corig1}}
Soit $f:X\ldrt Y$ un morphisme de schémas de type fini sur $S$. Soient
$W_0\subset Y_s$ et $W_1\subset X_s$ deux sous-schémas fermés tels que 
$W_{1 red} \subset f^{-1}(W_0)_{red}$. Le morphisme  $f$ induit donc un morphisme
$\hat{f} : \widehat{X}_{/W_1} \ldrt \widehat{Y}_{/W_0}$. Soit $\F\in
\D^+ ( X_\eta ,\Lambda)$. Il y a alors un diagramme commutatif 
\begin{diagram}[width=2cm]
\Rpsi (f_{\eta *} \F)_{|W_0} & \rTo & \Rth \left ( (\hat{f}^{an})_* 
\widehat{\F}_{/W_1} \right ) \\
\dTo & & \dTo^\simeq \\
\left ( f_{\bar{s}*} \Rpsi (\F) \right )_{|W_0} & \rTo &
\hat{f}_{\bar{s} *} \Rth \left ( \widehat{\F}_{/W_1} 
\right )  
\end{diagram}
\end{lemm}
\begin{proof}

Le diagramme s'insère dans le grand diagramme suivant dont on doit
montrer qu'il est commutatif 

\begin{diagram}[width=7mm,height=15mm,labelstyle=\scriptstyle]
\Rpsi (f_{\eta *} \F)_{|W_0} & \rTo^{ A \simeq} & \Rth (
\widehat{(f_{\eta *}\F)}_{/W_0} ) & \rTo^{\rth(\a) \simeq} & 
\Rth \left ( (f_{\eta *}^{an}\F^{an} )_{|(\widehat{Y}_{/W_0})^{an}}
\right ) & \rTo^{\rth (\b) } & \Rth \left ( \hat{f}^{an}_* (
\F^{an}_{|\widehat{X}_{/W_1}^{an}})\right )   \\ 
\dTo & && &&& \dTo^{\simeq} \\
\left ( f_{\bar{s}*} \Rpsi (\F) \right )_{|W_0} && \rTo^\gamma & 
(f_{\bar{s}|W_1})_* \left ( \Rpsi (\F)_{|W_1} \right ) & \rTo^{B \simeq} && 
\hat{f}_{\bar{s} *} \Rth \left ( \F^{an}_{|(\widehat{X}_{/W_1})^{an}}
\right ) 
\end{diagram}

 où $\a$ est induit par l'isomorphisme de changement de
base associé au diagramme de sites (corollaire 7.5.3 de \cite{Berk1}) 
\begin{diagram}[size=7mm]
X_{\eta,\et}^{an} & \rTo & X_{\eta,\et} \\
\dTo & & \dTo \\
Y_{\eta,\et}^{an} & \rTo & Y_{\eta,\et} 
\end{diagram}
$\b$ et $\gamma$ sont des applications de changement de base associées
à des restrictions, $A$ et $B$ sont induites par l'isomorphisme de
Berkovich, et les deux applications verticales sont induites par la
fonctorialité de $\rpsi$, resp. $\rth$, vis à vis de l'image directe.

Soient $a:W_0\hookrightarrow Y$ et $b:
(\widehat{X}_{/W_1})^{an}\hookrightarrow X_\eta^{an}$. On doit
vérifier que deux applications naturelles 
\begin{diagram}
a^*\Rpsi R f_{\eta * } & \pile{ \rTo \\ \rTo } & R f_{\bar{s} * } \Rth b^*
\end{diagram}
co{\"\i}ncident. Or, $b$ étant une quasi-immersion, $b^*$ envoie les
faisceaux mous sur des faisceaux mous qui sont eux m{\^e}mes
$\Theta_{\bar{\eta}}$ acycliques. De plus, $\Theta_{\bar{\eta}}$
envoie les faisceaux mous sur des faisceaux flasques (proposition 2.2
de \cite{Berk3}).
 Donc,
$$
R( f_{\bar{s}*} \Theta_{\bar{\eta}} b^* )= R f_{\bar{s} * } \Rth b^*
$$
Il est également clair que $R ( a^* \Psi_{\bar{\eta}} f_{\eta *} )=
a^* \Rpsi R f_{\eta *}$.
On vérifie alors aussit{\^o}t que les deux  applications naturelles sont
induites par dérivation des deux applications naturelles 
\begin{diagram}
 a^* \Psi_{\bar{\eta}} f_{\eta *} & \pile{ \rTo \\ \rTo } &
 f_{\bar{s}*} \Theta_{\bar{\eta}} b^* 
\end{diagram}
Dont on doit donc montrer qu'elles sont égales. Il faut donc montrer
que le diagramme commute sans foncteur dérivés et pour un faisceau
$\F$. Il suffit pour cela de démontrer que le diagramme commute au
niveau des préfaisceaux définissant les cycles évanescents. 
Soit donc $U/W_{0\bar{s}}$ étale et $s\in \GG( U, \Psi_{\bar{\eta}} (
f_{\eta *} \F)_{|W_0})$. Cette section est donnée par un germe de
diagramme 
\begin{diagram}[size=8mm]
  U & \rTo & && V \\
\dTo & & && \dTo_{\text{étale}} \\
W_{0\bar{s}} & \rInto & Y_{\bar{s}} & \rInto & \bar{Y} 
\end{diagram}
et une section 
$$
t\in \GG( V_{\bar{\eta}}, \overline{f_{\eta *} \F}) = \GG (
V_{\bar{\eta}}, f_{\bar{\eta} *} \bar{\F}) = \GG ( f_{\bar{\eta}}^{-1}
( V_{\bar{\eta}}), \bar{\F}) =\GG (f^{-1} (V)_{\bar{\eta}}, \bar{\F})
$$
La section image de $s$ par les trois application horizontales du haut
est une section $t'\in \GG ( W_{0\bar{s}}, \Theta_{\bar{\eta}} (
\hat{f}^{an}_* ( 
\F^{an}_{|\widehat{X}_{/W_1}^{an}}))$ qui se décrit ainsi :
comme dans le \ref{Rappels_Berk} il y a des diagrammes 
\begin{diagram}[size=8mm]
  \widetilde{U} & \rTo & \widehat{V}_{/W_0}& &   \widetilde{U}^{an} & \rTo & \widehat{V}_{/W_0}^{an} \\
\dTo & \ldTo  & & \text{ et } & \dTo & \ldTo  \\
Y_{/W_0}  & & & & Y_{/W_0}^{an}
\end{diagram}
Il y a donc un morphisme $\xi:(f^{an})^{-1} ( \widetilde{U}^{an}) \ldrt
(f^{an})^{-1} ( \widehat{V}^{an}_{/W_0})$. 
La section $t$ donne une section $u$  de $\GG( (f^{an})^{-1} (
\widetilde{U}^{an}), \bar{\F}^{an})$ par les applications 
\begin{diagram}[size=8mm]
\GG ( f^{-1}(V)_{\bar{\eta}}, \bar{\F}) & \rTo & \GG ( f^{-1}
(V)_{\hat{\bar{\eta}}}^{an}, \bar{\F}^{an}) & \rTo & \GG (
(f^{an})^{-1} ( \widehat{V}^{an}_{/W_0}),\bar{\F}^{an}) & \rTo^{\xi^*}
& \GG(  (f^{an})^{-1} ( \widetilde{U}^{an}),\bar{\F}^{an})
\end{diagram}
or, on a l'égalité
$$  \GG(  (f^{an})^{-1} ( \widetilde{U}^{an}),\bar{\F}^{an})= \GG (
\widetilde{U}^{an} ,  \hat{f}^{an}_* (
\F^{an}_{|\widehat{X}_{/W_1}^{an}}))$$ 
l'image de cette section $u$ par l'application de spécialisation 
$$
\GG( \widetilde{U}^{an},  \hat{f}^{an}_* (
\F^{an}_{|\widehat{X}_{/W_1}^{an}})) \ldrt \GG ( U, \Theta_{\bar{\eta}} (
\hat{f}^{an}_* ( 
\F^{an}_{|\widehat{X}_{/W_1}^{an}}))
$$
est $t'$. De plus, $\GG( U, \hat{f}_{\bar{s} *} \Theta_{\bar{\eta}} (
\F^{an}_{|(\widehat{X}_{/W_1})^{an}} ) = \GG ( f_{|W_1}^{-1}(U),\Theta_{\bar{\eta}} (
\F^{an}_{|(\widehat{X}_{/W_1})^{an}} ))$, et il y a des diagrammes 
\begin{diagram}[size=8mm]
  U_1=(f_{|W_1})^{-1} (U) & &  \rTo && f^{-1} (V)  &&&&  \widetilde{U_1} &
  \rTo & f^{-1} (V)_{/W_1} \\
\dTo & & & & \dTo_{\text{étale}} & & \text{  }& & \dTo & \ldTo  \\
W_{1\bar{s}} & \rInto & X_{\bar{s}} &\rInto &  \bar{X} && && \widehat{X}_{/W_1}
\end{diagram}

 Par définition de
l'application de
commutation naturelle des cycles évanescents analytiques à l'image
directe,   l'image de $t'$ par
l'application verticale de droite dans le grand diagramme est l'image de
$u$ par la composée 
$$
\GG( (f^{an})^{-1} (
\widetilde{U}^{an}), \bar{\F}^{an}) \ldrt \GG ( \widetilde{U_1}^{an},
\bar{\F}^{an}) \ldrt \GG ( U_1 , \Theta_{\bar{\eta}} ( \hat{f}^{an}_* ( 
\F^{an}_{|\widehat{X}_{/W_1}^{an}})))
$$
où la dernière application est l'application de spécialisation
associée aux deux diagrammes ci dessus.  Cela donne une description de
l'image de $t$ en parcourant le grand diagramme par le haut puis par
l'application verticale de droite. 

Décrivons l'image de $t$ de l'autre manière. 

L'image de $t$ par l'application verticale de gauche composée avec
l'application 
$\g$ dans le grand diagramme est donnée par son image par la composée
composée  
$$ 
\GG( f^{-1} (V)_{\bar{\eta}}, \bar{\F})  \ldrt \GG ( U_1,
\Psi_{\bar{\eta}} ( \F))= \GG ( U, (f_{\bar{s}|W_1})_*   \Psi_{\bar{\eta}} ( \F)_{|W_1})
$$
où l'application du milieu est l'application de spécialisation
associée au diagramme de gauche  ci dessus.
Par définition de l'isomorphisme de Berkovich, l'image de $t$, par le
second chemin, dans $\GG(U,(\hat{f}_{\bar{s}})_* \Theta_{\bar{\eta}} (
\F^{an}_{|\widehat{X}_{/W_1}^{an}})) = \GG (U_1, \Theta_{\bar{\eta}}
(\F^{an}_{|\widehat{X}_{/W_1}^{an}})) $ est l'image de $t$ par la
composée 
$$
\GG( f^{-1} (V), \bar{\F}) \ldrt \GG ( f^{-1}(V)^{an}, \bar{\F}^{an})
\ldrt \GG ( \widetilde{U_1}^{an}, \bar{\F}^{an}) \ldrt \GG ( U_1,
\Theta_{\bar{\eta}} ( \F^{an}_{|\widehat{X}_{/W_1}^{an}}))
$$
où la dernière application est l'application de spécialisation. 

Le fait que ces deux descriptions donnent la même image est alors clair
une fois que l'on mes bout à bout les différents morphismes décrits. 
\end{proof}

\begin{coro}
Dans le lemme précédent, supposons $Y=S$ et  $f$ propre. Il y a alors un
morphisme Galois équivariant de suite spectrale des cycles évanescents 
\begin{diagram}[size=9mm]
H^p( X_{\bar{s}}, R^q\Psi_{\bar{\eta}} \F ) & \rImplies & H^{p+q} (X_{\bar{\eta}},\F)
\\
\dTo & & \dTo \\
H^p(  W_{1\bar{s}} , R^q\Psi_{\bar{\eta}}  \F) & \rImplies & 
H^{p+q} ((\widehat{X}_{/W_1})^{an}\hat{\otimes}\hat{\bar{\eta}}, \F^{an})
\end{diagram}
\end{coro}

\subsubsection{Image inverse}

Soit $f:\X \ldrt \Y$. Il y a alors un morphisme naturel 
$$
f^* \rth \ldrt \rth f^*
$$
défini par adjonction grâce à l'isomorphisme $f_* \rth f^* \simeq \rth
f_*   f^* $ et l'application d'adjonction $Id\drt f_* f^*$. 

Le morphisme se décrit ainsi au niveau de la cohomologie : une section
de  $ f^*
R^q\Theta_{\bar{\eta}}\F$ est donnée étale localement par  un
$U/\X_{\bar{s}}$ étale, par un
diagramme 
\begin{diagram}[size=8mm]
U & \rTo& \X_{\bar{s}} \\
\dTo & & \dTo^{f} \\
V  & \rTo & \Y_{\bar{s}}
\end{diagram}
où $V/\Y_{\bar{s}}$ est étale, et une section $t\in H^q(\tilde{V}^{an},
\F^{an})$.
Le diagramme se relève en 
\begin{diagram}[size=8mm]
\tilde{U}^{an} & \rTo & \X^{an} \\
\dTo & & \dTo^{f} \\
\tilde{V}^{an}  & \rTo & \Y^{an}
\end{diagram}
et l'image réciproque de $t$ par l'application verticale de gauche
donne un élément de $H^q(\tilde{U}^{an},f^*\F)$, ce qui conclut la
description de notre morphisme.
\\

\begin{lemm}\label{Corig2}
Pla{\c c}ons nous dans le cadre du lemme \ref{Corig1}. Il y a alors un
diagramme commutatif 
\begin{diagram}[width=2cm]
\left ( f^* \Rpsi (\F) \right)_{|W_1} & \rTo & \hat{f}_{\bar{s}}^* \Rth
(\widehat{\F}_{/W_0}) \\
\dTo & &   \dTo \\ 
\Rpsi (f^*\F)_{|W_1} &
\rTo & \Rth \left (\hat{f}^{an *} (\widehat{\F}_{/W_0}) \right )
\end{diagram}
où les deux applications verticales sont celles associées à la
fonctorialité de $\rpsi$ et $\rth$ pour $f^*$.
\end{lemm}
\begin{proof}

Le diagramme s'insère dans le diagramme suivant 

\begin{diagram}[size=1cm]
\left ( f^* \Rpsi (\F) \right)_{|W_1} & \rEqualto & (f_{|W_1})^* \left (
\Rpsi (\F)_{|W_0}\right ) & \rTo^{A\sim} & (\hat{f}_{\bar{s}})^* \Rth
(\widehat{\F}_{/W_0}) \\
\dTo & & & &  \dTo \\ 
\Rpsi (f^*\F)_{|W_1} & \rTo^{B\sim} & \Rth (\widehat{f^*\F}_{/W_1}) &
\rEqualto & \Rth \left (\hat{f}^{an *} (\widehat{\F}_{/W_0}) \right )
\end{diagram}
où $A$ est induit par l'isomorphisme de Berkovich et $B$ est
l'isomorphisme de Berkovich. 

On doit montrer que deux morphismes 
\begin{diagram}
\left ( f^* \Rpsi (\F) \right)_{|W_1} & \pile{ \rTo \\ \rTo } &
 \Rth \left (\hat{f}^{an *} (\widehat{\F}_{/W_0}) \right )
\end{diagram}
co{\"\i}ncident.
Soient $a:W_1\hookrightarrow X_s$ et $b:(\widehat{Y}_{/W_0})^{an}
\hookrightarrow Y_\eta^{an}$. On a donc deux morphismes 
\begin{diagram}
a^* f^* \rpsi & \pile{\rTo \\ \rTo} & \rth \hat{f}^* b^*
\end{diagram}
Rappelons que $\rth= (R \nu_*)\circ \mu^*$.
 Se donner deux tels morphismes équivaut donc par adjonction à
se donner deux morphismes 
\begin{diagram}
\nu^* a^* f^* \rpsi &  \pile{ \rTo \\ \rTo } & \mu^* \hat{f}^* b^*
\end{diagram}
Ce qui revient à se donner deux morphismes 
\begin{diagram}
\nu^* a^* f^* \Psi &  \pile{ \rTo \\ \rTo } & \mu^* \hat{f}^* b^*
\end{diagram}
ou encore deux morphismes
\begin{diagram}
 a^* f^* \Psi &  \pile{ \rTo \\ \rTo } & \Theta \hat{f}^* b^*
\end{diagram}
On s'est donc ramené à vérifier la commutativité du diagramme dans le
cas où il n'y a plus de foncteurs dérivés, ce qui est facile en
utilisant la définition de l'isomorphisme de Berkovich.
\end{proof}

\subsubsection{Compatibilité au changement de base}

Idem. cas algébrique.
 
\section{Correspondances cohomologiques analytiques rigides}

$\Lambda$ désigne un anneau de coefficients de torsion première à
$p$. 

\subsection{Définition}

Soit 
\begin{diagram}[size=8mm]
 & & C \\
 & \ldTo^{c_1} & & \rdTo^{c_2} \\
X_1 & & & & X_2
\end{diagram}
un diagramme d'espaces analytiques de Berkovich sur $\eta$, et
$(\F_1,\F_2)\in \D^+( X_1,\Lambda)\times\D^+( X_2,\Lambda)$.

\begin{defi}
$$
\Coh_C(\F_1,\F_2)=\Hom (c_{2*}c_1^*\F_1,\F_2)
$$
\end{defi}

\subsection{Restriction à un domaine analytique}

Soient $U_1\subset X_1, U_2\subset X_2$ deux domaines analytiques
localement fermés tels que 
$c_2^{-1}(U_2) \subset c_1^{-1} (U_1)$. Considérons la correspondance 
\begin{diagram}[size=1cm]
& & c_2^{-1} (U_2) \\
& \ldTo^{c'_1} & & \rdTo^{c'_2} \\
U_1 & &&& U_2
\end{diagram}

Soit $u\in \Coh_C (\F_1,\F_2)$. Considérons le diagramme 
\begin{diagram}[size=8mm]
c_2^{-1}(U_2) & \rInto^b & C \\
\dTo^{c'_2} & & \dTo^{c_2} \\
U_2 & \rInto^a & X_2 
\end{diagram}
duquel il résulte un isomorphisme de changement de base $a^*c_{2*}
\iso c'_{2*} b^*$. 

La restriction de $u$ est alors définie par la composition suivante 
\begin{diagram}
res(u):& c'_{2*} {c_1'}^* ( \F_{1|U_1}) = c'_{2*} b^* c_1^* \F_1 & 
\lTo^\sim & a^* c_{2*} c_1^* \F_1 & \rTo^{a^* u} & a^* \F_2 = \F_{2|U_2}
\end{diagram}

$res (u)\in \Coh_{c_2^{-1}(U_2)}(\F_{1|U_1},\F_{2|U_2})$ et cela
définit un morphisme 
$$
\Coh_C (\F_1,\F_2) \ldrt \Coh_{c_2^{-1} (V_2)} (\F_{1|U_1},\F_{2|U_2})
$$

\subsection{Image directe}

Soit un diagramme d'espaces analytiques 

\begin{diagram}[size=1cm]
 & & C \\
 & \ldTo^{c_1} & \dTo^{f'} & \rdTo^{c_2} \\
X_1 & & & & X_2 \\
\dTo^{f_1} &  & D & & \dTo^{f_2} \\
& \ldTo^{d_1} & & \rdTo^{d_2} \\
Y_1 & & & & Y_2
\end{diagram}

et soit $u\in \Coh_C (\F_1,\F_2)$. On définit $f_* u$ de la fa{\c c}on
suivante 
\begin{diagram}
f_*u: & d_{2*} d_1^* f_{1*} \F_1 & \rTo & d_{2*} f'_* c_1^* \F_1 & 
\rEqualto & f_{2*} c_{2*} c_1^* \F_1 & \rTo^{f_{2*}u} & f_{2*}\F_2
\end{diagram}
où la première application est déduite du morphisme de changement de
base associé au losange de gauche. 
\\

On a donc une application :
$$
\Coh_C (\F_1,\F_2) \ldrt \Coh_D (f_{1*}\F_1,f_{2*}\F_2)
$$
D'où en particulier une action sur la cohomologie :
$$
\Coh_C (\F_1,\F_2) \ldrt \Hom ( R\GG (X_1,\F_1), R\GG (X_2,\F_2))
$$
qui est définie ainsi :
\begin{diagram}
u_*: & R\GG (X_1,\F_1) & \rTo^{c_1^*} & R\GG (C,c_1^*\F_1) & \rTo^\sim
& R\GG (X_2,c_{2*}c_1^*\F_1 ) & \rTo & R\GG (X_2,\F_2)
\end{diagram}

Ces images directes sont compatibles à la restriction en un sens
facile à définir. 

\subsection{Analytification d'une correspondance algébrique}

Soit 
$$
C\xrig{ (c_1,c_2) } X_1\times X_2
$$
une correspondance algébrique avec $c_1$ propre. 

Soit $u\in \Coh_C (\F_1,\F_2)$ une correspondance cohomologique
algébrique où $(\F_1,\F_2)\in \D^+_c (X_1,\La)\times \D^+_c
(X_2,\La)$. Il lui est associé un morphisme 
$\tilde{u}:c_{2*}c_1^*\F_1\ldrt \F_2$. $u^{an}$ est alors définie
ainsi :
\begin{diagram}[size=18mm]
u^{an} : & c_{2*}^{an} c_1^{an *} \F_1^{an} & \rEqualto & c_{2*}^{an}(
c_1^*\F_1)^{an} & \lTo^\sim & (c_{2*} c_1^* \F_1)^{an} &
\rTo^{\tilde{u}^{an}} & \F_2^{an}
\end{diagram}
où l'isomorphisme du milieu des l'isomorphisme de changement de base
du corollaire 7.5.3 de \cite{Berk1}. 
D'où un morphisme
$$
\Coh_C(\F_1,\F_2) \ldrt \Coh_{C^{an}}(\F_1^{an},\F_2^{an})
$$

On vérifie facilement que ce morphisme d'analytification est
compatible aux images directes propres. 

\section{Spécialisation des correspondances cohomologiques analytiques
rigides}
\subsection{Définition}

Soit 
\begin{diagram}[size=1cm]
 & & \mathfrak{C} \\
& \ldTo^{c_1} & & \rdTo^{c_2} \\
\X_1 & & && \X_2
\end{diagram}

un diagramme de schémas formels localement formellement de type fini
sur $S$, c'est à dire une correspondance formelle. Soit $u\in \Coh_{\CC^{an}} (\F_1,\F_2)$. La correspondance
spécialisée $u_{\bar{s}}$ est définie ainsi :
\begin{diagram}
u_{\bar{s}} : & c_{2*} c_1^* \Rth (\F_1) & \rTo & c_{2*} \Rth (c_1^* \F_1) &
\lTo^\sim & \Rth (c_{2*}c_1^*\F_1) & \rTo^{\rth(u)} & \Rth (\F_2)
\end{diagram}
où l'on a utilisé les diverses fonctorialités de $\rth$. Cela définit
 un morphisme : 
$$
\Coh_{\CC^{an}} (\F_1,\F_2) \ldrt \Coh_{\CC_{\bar{s}}} (\Rth(\F_1),\Rth(\F_2))
$$

\subsection{Compatibilité à l'image directe}

Considérons la proposition \ref{speima} dans le cadre rigide en ne
faisant aucune hypothèse de propreté sur les morphismes. 
Sa démonstration   s'adapte alors automatiquement au
cas rigide en 
changeant $\Psi$ en $\Theta$ et en remarquant que toutes les
hypothèses de propreté sont inutiles puisque $\Theta$, contrairement à
$\Psi$, commute toujours aux images directes (il faut bien s{\^u}r
utiliser les différentes fonctorialités établies pour $\Rth$ pour
adapter la démonstration). 

On en déduit donc la compatibilité entre image directe quelconque et
spécialisation dans le cas analytique rigide.

\section[Compatibilité entre spécialisations rigides et
algébriques]{Compatibilité entre les spécialisations algébriques et 
analytiques rigides} 

Soit une correspondance sur algébrique $S$
\begin{diagram}
C & \rTo^{(c_1,c_2)} & X_1\times X_2
\end{diagram}
où l'on suppose que $c_1$ et $c_2$ sont propres. 
Soient $V_1\subset X_{1s}$ et $V_2\subset X_{2s}$ deux sous-schémas fermés
de leur fibre spéciale. Supposons 
$$
c_2^{-1} (V_2)\subset c_1^{-1} (V_1)
$$
d'où une correspondance sur $s$ 
\begin{diagram}
c_2^{-1}(V_2) & \rTo^{(c'_1,c'_2)}  & V_1\times V_2
\end{diagram}
et une correspondance formelle 
\begin{diagram}
\CC & \rTo^{(\hat{c}_1,\hat{c}_2)}  & \X_1\times \X_2
\end{diagram}
où $\CC=\widehat{C}_{/c_2^{-1}(V_2)}$, $\X_1=\widehat{X_1}_{/V_1}$, 
 $\X_2=\widehat{X_2}_{/V_2}$.

\begin{prop} 
Soient $(\F_1,\F_2)\in \D^+_c (X_1,\La)\times \D^+_c (X_2,\La)$. 
Le diagramme suivant commute 
\begin{diagram}[labelstyle=\scriptstyle] 
\Coh_{C_\eta} (\F_1,\F_2) & \rTo^{\text{analytification}} &
\Coh_{C_\eta^{an}}(\F_1^{an},\F_2^{an}) & \rTo^{\text{restriction}} &
\Coh_{\CC^{an}}(\widehat{\F_1}_{/V_1},\widehat{\F_2}_{/V_2}) \\
\dTo^{sp} & & & &  \dTo^{sp} \\
\Coh_{C_{\bar{s}}} (\Rpsi (\F_1),\Rpsi (\F_2)) & \rTo^{\text{restriction}} &
\Coh_{c_2^{-1}(V_2)} (\Rpsi (\F_1)_{|V_1},\Rpsi (\F_2)_{|V_2}) &
\rTo^\sim & \Coh_{\CC_{\bar{s}}} (\Rth (\widehat{\F_1}_{/V_1}),\Rth(\widehat{\F_2}_{/V_2}))
\end{diagram}
\end{prop}
\begin{proof}
Soit $u\in \Coh_{C_\eta} (\F_1,\F_2)$ et notons $\tilde{u}:c_{2*}c_1^*
\F_1\ldrt \F_2$ le morphisme associé par adjonction. 
La démonstration consiste  à appliquer la compatibilité de
l'isomorphisme de comparaison avec les 
deux cas de fonctorialité pour les cycles évanescents.

Appliquons le lemme \ref{Corig2} à $\F_1$ et 
\begin{diagram}[size=1cm]
c_2^{-1}(V_2) & \rInto & C \\
\dTo^{c'_1} & & \dTo^{c_1} \\
V_1 & \rInto & X_1 
\end{diagram}
On obtient alors un diagramme 
\begin{diagram}[width=8mm,height=1cm]
{c'_1}^*\left ( \rpsi (\F_1)_{|V_1}\right ) && \rTo^\sim && (\hat{c}_{1\bar{s}})^* \rth
(\widehat{\F_1}_{/V_1}) \\
\dTo && \boxed{\mathcal{D}} & & \dTo \\
\rpsi  (c_1^*\F_1)_{|c_2^{-1} (V_2)} && \rTo && \rth ( (\hat{c}_1^{an})^*
\widehat{\F_1}_{/V_1}) 
\end{diagram}

Considérons maintenant le diagramme suivant 
\begin{diagram}[width=8mm,height=1cm]
c'_{2*}{c'_1}^*\left ( \rpsi (\F_1)_{|V_1}\right ) && \rTo^\sim && c'_{2*}(\hat{c}_{1\bar{s}})^* \rth
(\widehat{\F_1}_{/V_1}) \\
\dTo && \boxed{c'_{2*}\mathcal{D}} & & \dTo \\
c'_{2*}\rpsi  (c_1^*\F_1)_{|c_2^{-1} (V_2)} && \rTo &&
\hat{c}_{2\bar{s}*}\rth ( (\hat{c}_1^{an })^* 
\widehat{\F_1}_{/V_1}) \\
\uTo^\simeq & &&& \uTo^\simeq \\
\left ( c_{2*} \rpsi (c_1^*\F_1)\right )_{|V_2} &&& &
\hat{c}_{2\bar{s}*}\rth \left ( (\widehat{c_1^*\F_1})_{/V_1} \right )
\\
\uTo^\simeq &&& & \uTo^\simeq \\
\rpsi( c_{2*}c_1^* \F_1)_{|V_2}  && \rTo && \rth \left (
\hat{c}_{2*}^{an} (\widehat{c_1^*\F_1})_{/V_1} \right ) \\
 & \rdTo(4,2) &&& \dTo^\simeq \\
 &&&& \rth\left (\widehat{c_{2*}c_1^*\F_1}_{/V_2}\right ) &
 \rTo^{\rth(\tilde{u}^{an}_{|\mathfrak{X}_2^{an}})} & \rth ( \widehat{\F_2}_{/V_2})
\end{diagram}
où le deuxième rectangle commute grâce au lemme \ref{Corig1}. 

Le résultat s'en déduit, car les deux éléments de $\Coh_{\mathfrak{C}^{an}} (
\rth ( \widehat{\F_1}_{/V_1}),\rth(\widehat{\F_2}_{/V_2}))$ associés
se déduisent de ce diagramme en le parcourant du sommet en haut à
droite vers l'élément en bas à droite des deux fa{\c c}ons possibles en
suivant les bords du diagramme.  
\end{proof}

\section{Résumé des différentes fonctorialités}

Supposons nous donné un morphisme de correspondances sur $S$ 
\begin{diagram}
C & \rTo^{(c_1,c_2)} & X_1\times X_2 \\
\dTo^{f'} &  & \dTo_{f_1\times f_2} \\
D & \rTo^{(d_1,d_2)} & Y_1\times Y_2 
\end{diagram}

Soient $W_1\subset Y_1, W_2\subset Y_2$ deux fermés tels que
$d_2^{-1} (W_2) \subset d_1^{-1} (W_1)$. Soient $V_i=f_i^{-1} (W_i)$,
$i=1,2$ et supposons que $c_2^{-1} (V_2) \subset c_1^{-1} (V_1)$.  

Soient $\mathfrak{C}=\widehat{C}_{/c_2^{-1}(V_2)},
\X_i=\widehat{X}_{/V_i}$, $i=1,2$. Soient 
$\mathfrak{D}=\widehat{D}_{/d_2^{-1}(W_2)},
\mathfrak{Y}_i=\widehat{Y_i}_{/W_i}$, $i=1,2$. Il y a donc un
morphisme de correspondances formelles
\begin{diagram}
\CC & \rTo^{(\hat{c}_1,\hat{c}_2)} & \X_1\times \X_2 \\
\dTo^{\hat{f}'} &  & \dTo_{\hat{f}_1\times \hat{f}_2} \\
\mathfrak{D} & \rTo^{(\hat{d}_1,\hat{d}_2)} & \mathfrak{Y}_1\times
\mathfrak{Y}_2  
\end{diagram}
Qui induit un morphisme de correspondances analytiques. 

Soient $(\F_1,\F_2)\in \D^+_c (X_1,\Lambda)\times \D^+_c (X_2, \Lambda)$. 
Supposons que $c_1,c_2,d_1,d_2,f_1,f_2$ sont propres.  

Le diagramme de la page suivante est alors commutatif.
\newpage

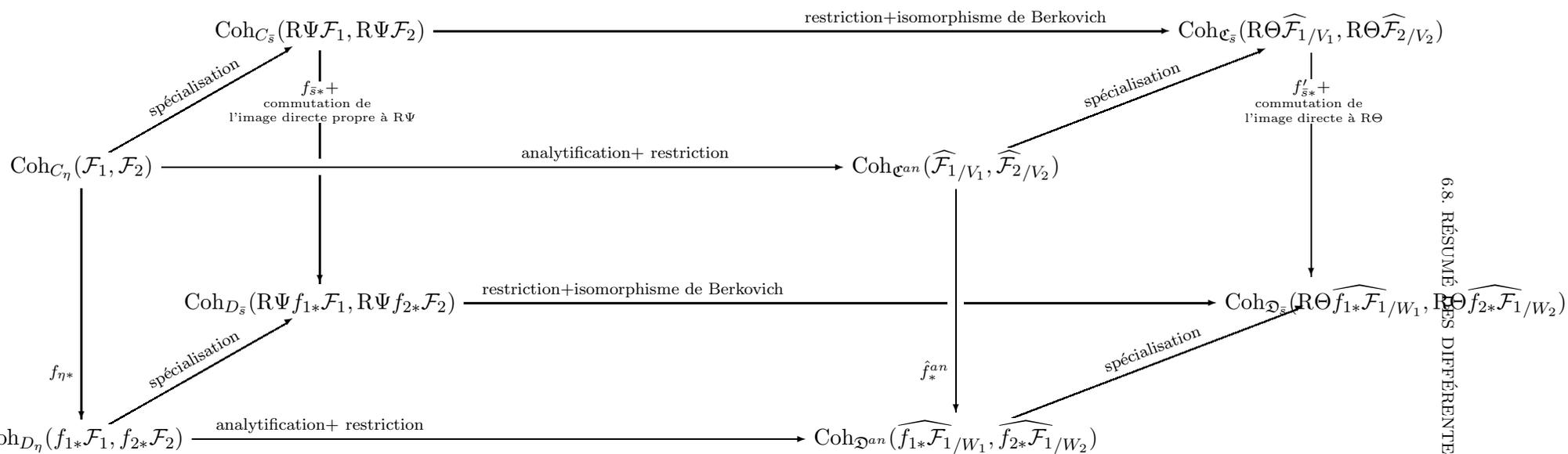
\newsavebox{\diag}
\savebox{\diag}{  
\begin{diagram}[PostScript=RadicalEye,landscape,labelstyle=\scriptstyle,width=4cm,heigth=5cm] 
  & \Coh_{C_{\bar{s}}} (\rpsi \F_1 ,\rpsi \F_2) &&
 \rTo^{\text{restriction+isomorphisme de Berkovich}}  & \Coh_{\CC_{\bar{s}}}
 (\rth \widehat{\F_1}_{/V_1},\rth \widehat{\F_2}_{/V_2}) \\
   \ruTo(1,2)^{\text{spécialisation}} & \vLine~{f_{\bar{s}*} + \atop {\text{commutation de } \atop \text{
 l'image directe propre à } \rpsi}} 
 & &  \ruTo(1,2)^{\text{spécialisation}} &
 \dTo~{f'_{\bar{s}*} + \atop {\text{commutation de } \atop \text{
 l'image directe à } \rth}}  \\
 \Coh_{C_\eta} (\F_1,\F_2) & \HonV & \rTo^{\text{analytification+
 restriction}} & \Coh_{\CC^{an}}
 (\widehat{\F_1}_{/V_1},\widehat{\F_2}_{/V_2}) \\
  &  \dTo  & &  \\
  &  \Coh_{D_{\bar{s}}} (\rpsi f_{1*}\F_1, \rpsi f_{2*}\F_2) &
 \hLine^{\;\;\text{   restriction+isomorphisme de Berkovich}} &
 \VonH & \rTo  \Coh_{\mathfrak{D}_{\bar{s}}} (\rth
 \widehat{f_{1*}\F_1}_{/W_1}, \rth
 \widehat{f_{2*}\F_1}_{/W_2}) \\
 \dTo^{f_{\eta *}}  \ruTo(1,2)^{\text{spécialisation}} & & &
 \dTo^{\hat{f}^{an}_*} \ruTo(1,2)^{\text{spécialisation}} \\
\Coh_{D_\eta} (f_{1*}\F_1,f_{2*}\F_2) & \rTo^{\text{analytification+
 restriction}} & & \Coh_{\mathfrak{D}^{an}} (
 \widehat{f_{1*}\F_1}_{/W_1},  \widehat{f_{2*}\F_1}_{/W_2}) 
\end{diagram} 
}

\begin{picture}(100,100)(0,5)
\put(650,-530){\usebox{\diag}}
\end{picture}

\newpage 

\section{Des coefficients de torsion aux coefficients $\ell$-adique}
\label{tors_ell}

Oublions dans cette section les notations précédentes et reprenons les
notations de l'appendice \ref{coho_l} : en particulier $k$ est un corps valué
complet, $\Lambda$ un anneau de valuation discrète complet dont nous
noterons $K_\La$ le corps des fractions. Soit $X$ un $k$-espace
analytique ou bien un espace adique sur $\Spa (k,k^0)$. En général
les groupes de cohomologie étale $\La$-adique 
$$
H^\bullet (X,\La)
$$
ne sont pas définis : il n'y a aucune définition intrinsèque qui
fasse que ces groupes vérifient de bonnes propriétés. Nous allons voir
cependant que dans certains cas, lorsque certaines propriétés de
finitude sont satisfaites, on peut donner une définition ad-hoc des
$H^\bullet(X,\La)$ satisfaisante.
\begin{defi}
  Soit $q\in \N$. Soit $\F=(\F_n)_n\in \Lpf\Xet$. 
Supposons que $(H^q (X,\F_n))_{n\in \N}$ soit un
  système A.R. $\La$-adique de type fini. On posera alors
$$ 
H^q (X,\F)=\underset{n\in \N}{\limp} H^q (X,\F_n)
$$ 
un $\La$-module de type fini.
\end{defi}

Cette définition est justifiée par la proposition suivante :

\begin{prop}\label{Poincare_adique}
  Supposons $X$ lisse sur $k$ (ou sur $\Spa(k,k^0)$ dans le cas
  adique) de dimension pure $d$
. Soit $\F=(\F_n)_n\in \Lpf\Xet$ localement constant. Soit
  $\check{\F}$ le système local dual. Supposons que $\forall q\in \N\;$ les
  systèmes 
$$
\left (H^q_c(X\hat{\otimes}\hat{\bar{k}},\F_n)\right )_{n\in \N}\text{ et } 
\left (H^q(X\hat{\otimes}\hat{\bar{k}},(\check{\F})_n)\right )_{n\in \N}
$$
sont A.R. $\La$-adiques de type fini et que l'application 
$$
H^q_c (X\hat{\otimes}\hat{\bar{k}},\F) \ldrt \underset{n\in \N}{\limp}
H^q_c (X\hat{\otimes}\hat{\bar{k}},\F_n) 
$$
est un isomorphisme. 

Il y a alors un isomorphisme canonique  de dualité de Poincaré 
$$
\left ( H^q_c (X,\F)\otimes_\La K_\La\right )^* \iso H^{2d-q}
(X,\check{\F} )(d)\otimes_\La K_\La
$$
\end{prop}

\begin{proof}
  Nous allons appliquer le formalisme de Ekedahl \cite{Eke} dont nous
  utiliserons librement les notations. Soit $f:X\hat{\otimes}\hat{\bar{k}}\ldrt
  \mathcal{M}(\hat{\bar{k}})$ (resp. $\Spa(
  \hat{\bar{k}},\hat{\bar{k}}^0)$ dans le cas adique). Il y a alors un
  foncteur (section 5 de \cite{Eke})
$$
\R f_! : \D (\mathcal{T}-\underline{\La}) \ldrt \D ( \mathcal{S}-\La)
$$
où $\mathcal{T}$ est le topos $\widetilde{X_{\et}}$ et 
$\mathcal{S}$  le topos ponctuel
$\widetilde{\mathcal{M}\left (\hat{\bar{k}}\right )_{\et}}$ (resp. $Spa (\hat{\bar{k}},\hat{\bar{k}}^0)_{\et}^{\widetilde{\;}}$)
vérifiant (confère par exemple la fin de l'énoncé du (iii) du théorème 6.3 de
\cite{Eke})
$$
\forall n\in \N \;\; \La/\mathfrak{m}^n \overset{\mathbb{L}}{\otimes}_{\La}
\left ( \R f_! \F\right )= \R\GG_c ( X\hat{\otimes}\hat{\bar{k}},\F_n)
\in \D^+ ( ( X\hat{\otimes}\hat{\bar{k}})_{\et},\La/\mathfrak{m}^n)
$$
(Ce foncteur est à distinguer du foncteur $\R\GG_c (X,-)$ introduit
dans l'appendice \ref{coho_l}) 
Le topos ponctuel $\mathcal{S}$ vérifie les hypothèses du paragraphe
7 de \cite{Eke}.  $\R f_!$ est de dimension cohomologique finie
($R^i f_!$ s'annule en dimension supérieure à $2d$ pour les
coefficients de torsion ( \cite{Hu1} 5.3.8, \cite{Berk1}
5.3.11). Donc, d'après les hypothèses faites 
$$
\R f_! \F\in \D^b_c ( \mathcal{S}-\La)
$$
D'après le théorème 7.1 de \cite {Eke} il y a une équivalence de
catégories 
\begin{diagram}
\D^b_c ( \mathcal{S}-\La) & \pile{ \rTo^{\R\pi_*} \\
  \lTo_{\mathbb{L}\pi^*} } & \D^b_{\text{type fini}} ( \La-\text{mod})
\end{diagram}
où $\pi_*$ est le foncteur $\limp$,  $\pi^*$ le foncteur $M\mapsto
(M\otimes  \La/\mathfrak{m}^n )_{n\in \N}$, et $\R\pi_*$ et $\mathbb{L}\pi^*$
sont quasi-inverses l'un de l'autre.
Donc, 
\begin{eqnarray*}
  \R f_! \F & \simeq & \mathbb{L} \pi^* ( \R \pi_* \R f_! \F) \\
 & = & \left ( \La/\mathfrak{m}^n \overset{\mathbb{L}}{\otimes}_\La \R
 \pi_* \R f_! \F \right )_{n\in \N} \\
 & \simeq & \left ( \R\GG_c ( X\hat{\otimes} \hat{\bar{k}},\F_n)\right
 )_{n\in \N}
\end{eqnarray*}
Soit $M^\bullet = \R\pi_* (\R f_! \F )$. La suite spectrale 
$$
E^{ij}_2 = \R ^i \pi_* \left ( (H^j_c ( X\hat{\otimes}\hat{\bar{k}},
  \F_n)_{n\in \N} \right ) \limpl H^{i+j} (M^\bullet )
$$
vérifie $E^{ij}_2=0$ si $i\neq 0$ car par hypothèse le système 
$(H^j_c ( X\hat{\otimes} \hat{\bar{k}},\F_n))_{n\in \N}$ est
A.R. $\La$-adique. Donc, 
$$
H^j (M^\bullet) \simeq \underset{n\in \N}{\limp} H^j_c (X,\F_n)
$$
Comme dans le paragraphe 6 de \cite{Eke}, utilisant les propriétés de
finitude de la dimension cohomologique de $f_!$, on peut définir 
$f^! \La\in \D (\mathcal{T}-\underline{\La})$. Il résulte du théorème
de dualité de Poincaré pour les coefficients de torsion ( théorème
7.3.1 de \cite{Berk1} dans le
cas analytique et théorème 7.5.3 de \cite{Hu1} dans le cas adique) que 
$f^! \La \simeq \underline{\La} (d)[2d]$.
En effet, le morphisme trace pour les coefficients de torsion $\R f_!
f^* (-)(d)[2d]\drt (-)$ s'étend naturellement en une application de
foncteurs de $\D^b( \mathcal{S}-\La)$ dans $\D^b(\mathcal{S}-\La)$ qui
définit par adjonction un morphisme $f^* (-) (d) [2d]\drt f^! (-)$
dans $\D^b (\mathcal{T}-\La)$. Ce morphisme est un isomorphisme
puisqu'après application du foncteur conservatif
$-\overset{\mathbb{L}}{\otimes} \La/\mathfrak{m}$ on obtient
l'isomorphisme dans le cas de torsion (grâce au 6.3. (iii) de
\cite{Eke}). 
 La dualité de Verdier $\La$-adique 
(isomorphisme 6.1 de 
\cite{Eke}) donne alors un
isomorphisme 
$$
\R\Hom ( \R f_! \F, \La_\bullet)
\simeq \R\Hom ( \F , \La_\bullet (d) [2d])
$$  
Le premier terme s'identifie via l'équivalence de catégories ci dessus à 
$$
\R\Hom_\La ( M^\bullet, \La) 
$$
quant au second à 
$$
\R f_* ( \check{\F}(d)[2d])
$$
qui, comme précédemment pour $\R f_!$, est égal à 
$$
\left ( \R \GG (X,\F_n)\right )_{n\in\N}
$$
Grâce à l'hypothèse que $\forall
j,\; (H^j(X,\check{\F}_n))_n$ est un système  A.R. $\La$-adique de type fini,
si $$N^\bullet=\R\pi_* \R f_* \check{\F} (d)[2d]$$
alors 
$$
\forall j\;\; H^j (N^\bullet)= \underset{n\in\N}{\limp} H^{2d+j} ( X,\check{F}) (d)
$$
On obtient (toujours grâce à l'équivalence du théorème 7.1 de
\cite{Eke}) donc 
$$
\R\Hom_\La (M^\bullet, \La) \simeq N^\bullet
$$
dans la catégorie dérivée bornée des complexes de $\La$-modules à
cohomologie de type fini.
 Si $W$ est un $\La$-module de type fini, 
$$
Ext^{i}_{\La} ( W,\La) = \left \{ \Hom_{\La} ( W_{\text{torsion}},
K_\La/\La) \text{ si } i=1 \atop 
0 \text{ si } i\geq 2 \right.
$$
La suite spectrale 
$$
E^{ij}_2 = Ext_{\La}^i ( H^{-j} (M^\bullet), \La) \limpl H^{i+j} ( N^\bullet)
$$
dégénère donc en des suites exactes pour tout $i$ 
\begin{diagram}[size=8mm]
  0 & \rTo &  \Hom_{\La} ( H^{2d-i+1}_c ( X,
  \La)_{\tiny\text{torsion}}, K_\La/\La )  & \rTo & H^i (
  X,\check{\F})(d) & \rTo &  \Hom_\La ( H^{2d-i}_c ( X, \F),\La) & \rTo 0 
\end{diagram}
qui donne le résultat voulu après tensorisation par $K_\La$ puisque le
terme de gauche est de torsion.
\end{proof}

\begin{rema}
  Il est sous entendu que les isomorphismes ci dessus sont Galois équivariants. 
\end{rema}

\begin{prop}\label{cycles_l_adiques}
  Les hypothèses de la proposition précédente sont vérifiées lorsque
  $X$ et $\F$ sont de la forme suivante :
  \begin{itemize}
  \item Il existe un schéma séparé de type fini $Y/k^0$, un faisceau
    étale algébrique constructible $\G$ sur $Y_\eta$ et un ouvert $U$
    de la fibre spéciale de $Y$, $Y\times_{k^0} k^0/k^{00}$, tels que
    $X$ s'identifie à la fibre générique au sens des espaces adiques
    du complété formel de $X$ le long de $U$ 
$$
(\widehat{Y}_{/U})^{rig}\subset (Y_\eta)^{rig}
$$ 
et
$\F$ s'identifie à la restriction de $\G^{rig}$ à cet ouvert (et $X$
est lisse et $\F$ est localement constant). 
\item Il existe un schéma séparé de type fini $Y/k^0$,  un faisceau
    étale algébrique constructible $\G$ sur $Y_\eta$ et, dans le cas
    analytique, un fermé propre
    $V$ de
    la fibre spéciale de $Y$ tels que $X$ s'identifie au tube au
    dessus de $V$ dans $(Y_{\eta})^{an}$, $( \widehat{Y}_{/V})^{an}$,
    et $\F$ s'identifie à la restriction de $\G^{an}$. Dans le cas adique
    mêmes conditions mais 
    $V$ est seulement supposé fermé.   
  \end{itemize}
Soit de plus $\X$ le schéma formel $\widehat{X}_{/U}$ dans le premier
cas et $\widehat{X}_{/V}$ dans le second. Soit $(\Rth (\F_n))_{n\in
  \N}$ le système de faisceaux des cycles proches analytiques rigides. Alors, 
$$
(\Rth ( \F_n))_{n} \in \D^b_c ( \X_{\bar{s}}, \La)
$$  
la catégorie dérivée bornée à cohomologie constructible $2$-colimite
des catégories $(\D^b_c ( \X_{\bar{s}}, \La/\mathfrak{m}^n))_{n\in \N}$, et il y a
une suite spectrale des cycles évanescents de $\La$-modules de type
fini convergente 
$$
E^{pq}_2= H^p ( \X_{\bar{s}}, (\Rth ( \F_n))_n
)=\underset{n\in\N}{\limp} H^{p+q} ( \X_{\bar{s}}, \Rth
( \F_n)) \limpl H^{p+q} ( 
X\hat{\otimes} \hat{\bar{k}}, \bar{\F}) 
$$
\end{prop}
\begin{proof}
  Commençons par les assertions sur les cycles évanescents. Le fait
  que $ \forall n\;
\Rth ( \F_n) \in \D^b_c ( \X_{\bar{s}}, \La/\mathfrak{m}^n)
$ résulte du théorème de comparaison de Berkovich et du théorème de
  de constructibilité de Deligne (\cite{Del-finitude}) pour les cycles
  évanescents algébriques. La formule de Künneth pour les cycles
  évanescents algébriques (appendice de \cite{Illusie2}) couplée au
  théorème de comparaison de Berkovich montre que
  $$
\forall n\;\; \Rth (\F_{n+1}) \overset{\mathbb{L}}{\otimes} \La/\mathfrak{m}^n
 \simeq  \Rth ( \F_n)
$$ 
(et cet isomorphisme est canonique). On en déduit que
$$
(\Rth (\F_n))_n \in \D_c^b ( \X_{\bar{s}},\La)
$$
 Il résulte également de \cite{Del-finitude}
  et de l'annulation des $\text{R}^i\Psi_{\bar{\eta}}$ pour $i$ grand  
  que $(\Rth (\F_n))_n$ est un système A.R. $\La$-adique (utiliser le
  lemme 12.5 de \cite{FK} ou bien \cite{Jou2} 5.3.1). Le théorème de
  finitude de la cohomologie 
  étale des faisceaux constructibles sur un corps de
  \cite{Del-finitude} couplé au théorème 5.3.1 de \cite{Jou2}  montre alors que
  $\forall p\;\forall q$ le $\La_\bullet$-module 
$$
(H^p( \X_{\bar{s}}, \Rth^q ( \F_n)))_{n\in \N}
$$
est A.R. $\La$-adique. Pour tout $n$ il y a une suite
spectrale des cycles évanescents  
$$
E^{pq}_2 (n)=
H^p ( \X_{\bar{s}}, \Rth ( \F_n))\limpl H^{p+q} ( X\hat{\otimes}
\hat{\bar{k}},\bar{\F}_n) 
$$
On obtient donc une suite spectrale convergente dans la catégorie des
$\La_\bullet$-modules 
$$ E^{pq}_2=
(E^{pq}_2 (n))_{n\in \N} \limpl \left (H^{p+q} (  X\hat{\otimes} 
\hat{\bar{k}},\bar{\F}_n )\right )_{n\in \N}
$$
Les premiers termes $E^{pq}_2$ sont dans la catégorie des
$\La_\bullet$-modules A.R. $\La$-adiques de type fini. Cette catégorie
est une sous-catégorie abélienne pleine stable par extension de celle
des $\La_\bullet$-modules (\cite{Jou2} théorème 5.2.3 et proposition
5.2.4). Le fait 
qu'elle soit abélienne implique que  $\forall p\; \forall q \;
\forall r \geq 2\;\; (E^{pq}_r (n))_{n\in \N}$ est 
A.R. $\La$-adique de type fini. $\forall p\; \forall q \;$ l'aboutissement $ (H^{p+q} (  X\hat{\otimes} 
\hat{\bar{k}},\bar{\F}_n ))_{n\in \N}$ est alors muni d'une filtration
finie dont les gradués sont A.R. $\La$-adique de type fini. La stabilité par
extensions de la catégorie des $\La_\bullet$-modules
A.R. $\La$-adiques  montre donc que cet aboutissement est 
A.R. $\La$-adique. Considérons maintenant le foncteur 
$\limp$ qui va de la catégorie des $\La_\bullet$-modules dans celle
des $\La$-modules. En restriction à la catégorie des $\La_\bullet$-modules
A.R. $\La$-adiques de type fini il est exact d'image des $\La$-modules de
type fini. Appliquant ce  foncteur exact à la suite spectrale
$(E^{pq}_2 (n))_{n\in \N}$ on obtient donc une suite spectrale
convergente de $\La$-modules de type fini 
$$
E^{pq}_2= H^p ( \X_{\bar{s}}, (\Rth ( \F_n))_n
)=\underset{n\in\N}{\limp} H^{p+q} ( \X_{\bar{s}}, \Rth
( \F_n)) \limpl \underset{n\in \N}{\limp} H^{p+q} ( 
X\hat{\otimes} \hat{\bar{k}}, \bar{\F_n}) 
$$
ce qui termine la démonstration de la partie cycles évanescents de la
proposition.  

Passons maintenant à la première partie. Dans tous les cas, on a déjà
vérifié l'hypothèse concernant $(H^q ( X\hat{\otimes}
\hat{\bar{k}},\check{F_n}))_{n\in \N}$ de la proposition
\ref{Poincare_adique}. Reste donc les assertions concernant la
cohomologie à support compact. Elles r\'esultent dans le cas des espaces
adiques du théorème 3.1 de \cite{Hu2} pour l'espace adique
associé au complété formel
le long de l'ouvert $U$ et du théorème 3.3 de \cite{Hu2} dans le cas
de l'espace adique associé au
complété formel le long du fermé $V$. Dans le cas où $V$ est propre,
l'énoncé pour l'espace analytique de Berkovich associé au complété
formel le long $V$ est une conséquence du fait que $V$ propre implique
que l'espace adique est partiellement propre, et que donc sa
cohomologie à support compact coïncide avec celle de l'espace
analytique associé (proposition 1.5 de \cite{Hu2}).   
\end{proof}

\begin{defi}
  Soit $C\xrig{ \; (c_1,c_2)\;} X_1\times X_2$ une correspondance
  d'espaces analytiques où $c_1$ et $c_2$ sont finis. 
  Soient $\F_1,\F_2\in \Lpf\Xet$. On note $\Coh_C (\F_1,\F_2)$ les
  correspondances cohomologiques à support dans $C$ que l'on pose
  étant égal à l'ensemble des systèmes compatibles de correspondances
  dans $\prod_{n\in \N} \Coh_C ( (\F_1)_n, (\F_2)_n)$. 
\end{defi}

Comme précédemment, une correspondance cohomologique agit sur la
cohomologie $\La$-adique lorsque celle ci est bien définie
par simple passage à la limite. De plus, dans le cas de la fibre
générique d'un schéma formel, une correspondance spécialisée (égale
par définition au système des correspondances spécialisées) sur les
cycles proches analytiques rigides fournit une suite spectrale
des cycles proches équivariante comme précédemment. Cela se déduit par
passage à la limite projective des suites spectrales obtenues dans le
cas de torsion grâce à la proposition précédente.

\begin{rema}
  Soit $(c_1,c_2):C\ldrt X_1\times X_2$ une correspondance analytique
  où $c_1$ et $c_2$ sont étales finis. 
  Via la dualité de Poincaré de la proposition \ref{Poincare_adique},
  l'action d'une correspondance cohomologique $\La$-adique sur la
  cohomologie à support compact doit être définie comme étant l'action
  de la correspondance duale si l'on veut que l'isomorphisme de
  Poincaré soit équivariant. 
\end{rema}

\begin{enonce}{Application}
  La représentation locale fondamentale de \cite{Har4} s'exprime comme
  le dual de Poincaré de la cohomologie $\ell$-adique à support compact
  des espaces de Lubin-Tate. 
\end{enonce}

Désormais nous utiliserons librement des cycles évanescents
$\ell$-adiques sans forcément citer la proposition
\ref{cycles_l_adiques}. Cette proposition contient tous les éléments
nécessaires (et même plus) pour justifier le passage des coefficients
de torsion aux coefficients $\ell$-adiques lorsque cela sera nécessaire.

\section{Formule des traces générale}

\begin{theo}
Soit $D \xrightarrow{\; (d_1,d_2) \; } Y_1\times Y_2$ une
correspondance sur $S$ où $D,Y_1$ et $Y_2$ sont propres sur $S$.

Soit $C \xrightarrow{\; (c_1,c_2) \; } X_1\times X_2$ une
correspondance sur $\eta$, où $c_1$ et $c_2$ sont finis. 

Soit un morphisme de correspondances sur $\eta$ :
\begin{diagram}
C & \rTo & X_1\times X_2 \\
\dTo^{f'} &  & \dTo_{f_1\times f_2} \\
D_\eta & \rTo & Y_{1\eta}\times Y_{2\eta} 
\end{diagram}
où $f_1,f_2$ sont étales finis. 

Soient $V_1\subset Y_{1s}, V_2\subset Y_{2s}$ deux fermés vérifiant 
$d_2^{-1}(V_2)\subset d_1^{-1}(V_1)$. Soient $U_1=(f_1^{an})^{-1} (sp^{-1}
(V_1)), U_2=(f_2^{an})^{-1} (sp^{-1} (V_2))$ les ouverts de $X_1^{an},
X_2^{an}$ associés. Supposons $c_2^{an -1} (U_2) \subset c_1^{an
-1}(U_1)$.

Soient $\F_1,\F_2$ deux $\Qlb$ faisceaux lisses sur $X_1$,
resp. $X_2$. Soit $u\in\Coh_C(\F_1,\F_2)$. Soit $v=(f_* u)_{\bar{s}} \in \Coh
(\Rpsi (f_{1*} \F_1),\Rpsi (f_{2*} \F_2))$ la spécialisation de l'image
directe de $u$. Il y a alors un morphisme de suites spectrales des
cycles évanescents commutant à l'action de $\Gal (\bar{\eta} |\eta)$ :
\begin{diagram}
H^p \left (V_{1\bar{s}}, R^q\Psi_{\bar{\eta}} (f_{1*}\F_1)\right ) & 
\rImplies & H^{p+q}\left ( U_1\hat{\otimes} \hat{\bar{\eta}},\bar{\F}_1
\right ) \\ 
\dTo^{v_*} & & \dTo^{u^{an}_*} \\
H^p \left (V_{2\bar{s}}, R^q\Psi_{\bar{\eta}} (f_{2*}\F_2)\right ) & 
\rImplies & H^{p+q}\left ( U_2\hat{\otimes} \hat{\bar{\eta}},\bar{\F}_2
\right )
\end{diagram}
où la cohomologie $\ell$-adique des ouverts rigides $U_1,U_2$ est
définie comme étant la Poincaré duale de la cohomologie à support
compact et où 
$$ u_*^{an} :
H^i (  U_1\hat{\otimes} \hat{\bar{\eta}},\bar{\F}_1)
 \xrig{ (c_{1}^{an})^*_{|c_2^{an-1}(U_2)}} H^i ( (c_2^{an})^{-1} (U_2) 
\hat{\otimes} \hat{\bar{\eta}},c_1^*\bar{\F}_1)  
\xrig{(c_{2}^{an})_*} H^i ( U_2 
\hat{\otimes} \hat{\bar{\eta}},c_{2*} c_1^*\bar{\F}_1)  
$$
$$\hspace{-8cm}\xrig{ H^i (\tilde{u})}
 H^i ( U_2 
\hat{\otimes} \hat{\bar{\eta}},\bar{\F}_2)  
$$
où $\tilde{u}:c_{2*}c_1^* \F_1\ldrt \F_2$ est associé à $u$ par adjonction.
\end{theo}

\begin{proof}
Pour $i=1,2$ complétons le diagramme 
\begin{diagram}[size=8mm]
&  & X_i \\
 & & \dTo \\
Y_i &  \lInto & Y_{i\eta} 
\end{diagram}
en un diagrame 
\begin{diagram}[size=8mm]
X_i^0 & \lInto & X_i \\
\dTo & & \dTo \\
Y_i &  \lInto & Y_{i\eta} 
\end{diagram}
où $X_i^0$ est un schéma de type fini sur $S$ et 
 $X_i^0/Y_i$ est fini. En particulier,
 $X_i^0$ est propre sur $S$.
 Soit le
diagramme :
\begin{diagram}[heigth=8mm,width=25mm]
 & & C \\
 & & \dTo \\
X_1^0 \times_{Y_1}  D \times_{Y_2} X_2^0 & \lInto &
 X_1\times_{Y_{1\eta}} D_\eta \times_{Y_{2\eta}} X_2 
\end{diagram}
Et, comme ci dessus, soit $C^0$ un modèle de type fini de $C$ sur $S$ 
et fini au dessus de $X_1^0 \times_{Y_1}  D
\times_{Y_2} X_2^0$. Ce dernier schéma étant propre sur $S$, $C^0$ est
propre sur $S$ 
 (cependant, et c'est la un point important qui nous a
poussé à développer notre formalisme des correspondances cohomologiques
dans la catégorie dérivée, on ne peut pas forcément trouver de modèle
$C^0$ tel que  
$C^0/X_1^0$ et $C^0/X_2^0$ soient finis). 

On a donc étendu la correspondance $C\drt X_1\times X_2$ sur $\eta$
en une correspondance $C^0\drt X_1^0\times X_2^0$ sur $S$. On a
également étendu le morphisme de correspondances $f$ en un morphisme
de correspondances 
sur $S$. Tous les
morphismes étendus sont propres. 

Remarquons maintenant que l'hypothèse $c_2^{an-1}(f_2^{an -1}
(U_2))\subset c_1^{an-1} (f_1^{an -1} (U_1))$ implique, grâce à la
surjectivité du morphisme de spécialisation, qu'il en est de m{\^e}me sur
la fibre spéciale (au niveau des schémas réduits ce qui est suffisant
pour nos besoins cohomologiques étales) pour les fermés
$V_1$ et $V_2$ : $c_{2s}^{-1} (f_{2s}^{-1} (V_2))_{\text{red}} \subset c_{1s}^{-1}
(f_{1s}^{-1} (V_1))_{\text{red}}$. 

Le théorème résulte donc de la compatibilité des différentes
fonctorialités pour les cycles évanescents démontrées dans les section
précédentes. 
\end{proof}

Nous supposerons désormais que $k(s)$ est le corps fini $\Fq$. Nous
noterons $\eta^{nr}$ le complété de l'extension maximale non ramifiée
de $\eta$, $\s\in \Gal (\eta^{nr} |\eta )\simeq \Gal (\bar{s} |s )$
le Frobenius géométrique et $W_\eta\subset \Gal ( \bar{\eta} | \eta )$
le groupe de Weil. Si $\tau\in W_\eta$ nous noterons $v(\tau)$
l'entier vérifiant $\tau_{|\eta^{nr}}=\s^{v(\tau)}$. 

Soit $U$ un espace analytique sur $\eta^{nr}$. Pour tout $\tau\in
W_\eta$ il y a un isomorphisme canonique 
$$
\left ( U\hat{\otimes}_{\eta^{nr}} \hat{\bar{\eta}}\right )^{(\tau)}
\simeq 
U^{(\s^{v(\tau)})} \hat{\otimes}_{\eta^{nr}}  \hat{\bar{\eta}}
$$

Si $U$ est un domaine analytique dans un espace analytique de la forme
$X\hat{\otimes} \eta^{nr}$ pour un espace analytique $X$ sur $\eta$, $
U^{(\s^{v(\tau)})}$ est le domaine analytique de $X\hat{\otimes}_\eta {\eta^{nr}}
$ dont les points rigides na{\"\i}fs (ceux du spectre maximal) sont l'image
par $\s^{v(\tau)}$ de ceux de $U$.

De l'isomorphisme ci dessus on déduit un morphisme (qui n'est pas au
dessus de $ \hat{\bar{\eta}}$) d'espaces analytiques pour
tout $\tau$ 
\begin{diagram}[labelstyle=\scriptstyle]
 U\hat{\otimes}_{\eta^{nr}} \hat{\bar{\eta}} & \rTo^{Id\times \tau} &
 U^{(\s^{v(\tau)})} \hat{\otimes}_{\eta^{nr}}  \hat{\bar{\eta}} 
\end{diagram}
Si maintenant $\F$ est un faisceau sur $X_{\et}$, $(Id\times
\tau)^*\bar{\F}\simeq \bar{\F}$ où $\bar{\F}$ désigne le pull-back de
$\F$ à $X\hat{\otimes} \hat{\bar{\eta}}$. On en déduit pour tout
$\tau$ dans $W_\eta$ un isomorphisme 
$$
H^q \left ( U^{(\s^{v(\tau)})}
\hat{\otimes}_{\eta^{nr}}  \hat{\bar{\eta}}, \bar{\F}\right )
 \xrig{\;\;\; \tau^*\;\; } 
H^q \left (  U\hat{\otimes}_{\eta^{nr}} \hat{\bar{\eta}},
\bar{\F}\right )
$$

Supposons de plus qu'il existe un schéma formel $\Y$ sur $S$, un
morphisme $f:X\ldrt \Y^{an}$ et  un point fermé $y\in \Y_{\bar{s}}$
tels que $U=f^{-1} (sp^{-1} ( y))$. Alors, si $v(\tau)\geq 0$  
$$
U^{(\s^{v(\tau)})}=f^{-1} ( sp^{-1} (Fr^{v(\tau)} (y)))
$$
et on en déduit donc un isomorphisme lorsque $v(\tau)\geq 0$ 
$$
 H^q \left ( f^{-1} ( sp^{-1} (Fr^{v(\tau)} (y)))
\hat{\otimes}_{\eta^{nr}}  \hat{\bar{\eta}}, \bar{\F}\right )
\xrig{\;\; \tau^*\;\; }
H^q \left (  f^{-1}(sp^{-1}(y))\hat{\otimes}_{\eta^{nr}} \hat{\bar{\eta}},
\bar{\F}\right )
$$

\begin{theo}\label{Fortree}
Pla{\c c}ons nous dans le cadre du théorème précédent lorsque $X_1=X_2=X$,
$Y_1=Y_2=Y$, $\F_1=\F_2=\F$, $V_1=V_2=V$ et $U_1=U_2=U$. Il y a donc
un
diagramme au dessus de $\eta$ :
\begin{diagram}[width=15mm,labelstyle=\scriptstyle]
C & \rTo^{(c_1,c_2)} & X\times X \\
\dTo^{f'} &  & \dTo_{f\times f} \\
D_\eta & \rTo^{(d_{1\eta},d_{2\eta})} & Y_{\eta}\times Y_{\eta} 
\end{diagram}

 Supposons de plus que
$d_{2\bar{s}}$ est fini. Il existe alors un entier $N$ ne dépendant
que du degré de $d_{2\bar{s}}$ tel que $\forall u\in \Coh_{C_\eta}
(\F), \forall \tau\in W_\eta \; v(\tau)\geq N $
$$
\text{Tr}\left ( u\times\tau; R\GG (X_{\bar{\eta}}, \bar{\F} )\right ) = 
\sum_{y\in D_{\bar{s}}(\bar{k}) \atop Fr^{v(\tau)}(d_1(y))=d_2 (y)} 
\text{Tr} \left ( u\times\tau ; R\GG \left ( f^{an-1} (sp^{-1}
(d_2(y)))\hat{\otimes} 
\hat{\bar{\eta}}, \bar{\F} \right)\right )
$$ 
et 
$$
\text{Tr}\left ( u\times\tau; R\GG (U\hat{\otimes} \hat{\bar{\eta}}, \bar{\F}^{an} )\right ) = 
\sum_{y\in D_{\bar{s}}(\bar{k}) \atop { d_2 (y)\in V \atop
Fr^{v(\tau)}(d_1(y))=d_2 (y)}} 
\text{Tr} \left ( u\times\tau ; R\GG \left ( f^{an-1} (sp^{-1}
(d_2(y)))\hat{\otimes} 
\hat{\bar{\eta}}, \bar{\F} \right)\right )
$$ 
où l'action de $u\times \tau$ dans le membre de droite est le composé
des  morphismes 
\begin{diagram}[size=10mm]
 R\GG  ( f^{an-1} (sp^{-1}
(\underbrace{Fr^{v(\tau)} d_1(y)}_{d_2(y)}))\hat{\otimes} 
\hat{\bar{\eta}}, \bar{\F} ) 
 & \rTo^{\tau^*} &
R\GG \left ( f^{an-1} (sp^{-1}
(d_1(y)))\hat{\otimes} 
\hat{\bar{\eta}}, \bar{\F} \right) &\rTo^{c_1^*} \\ 
R\GG\left  ( {f'}^{an-1}(sp^{-1}(y))\hat{\otimes} 
\hat{\bar{\eta}}, \bar{\F}\right ) & 
\rTo^{c_{2*}} &  R\GG \left ( f^{an-1} (sp^{-1}
(d_2(y)))\hat{\otimes} 
\hat{\bar{\eta}}, \bar{\F} \right) 
\end{diagram}
\end{theo}
\begin{proof}
Les cycles évanescents $\Rpsi ( f_* \F)$ sont munis d'une structure de 
faisceau de Weil et donc d'une correspondance de Frobenius associée à
$\tau$ :
$$
Fr^{v(\tau) *} \rpsi (f_*\F) \iso \tau^* \rpsi (f_*\F) \ldrt \rpsi (f_*\F)
$$
 Cette correspondance  commute aux
correspondances cohomologiques spécialisées puisque ces dernières sont
définies sur $s$. 
On peut donc appliquer le théorème de Fujiwara 
(\ref{theoFuji}) à la correspondance cohomologique spécialisée
composée avec la correspondance de Frobenius ci dessus associée à
$\tau$.  

Plus précisément, dans le premier cas on applique Fujiwara à tout
$X_s$. Dans le second cas, d'après le théorème précédent, la trace
sur la cohomologie de $U$ est égale à celle sur la cohomologie des
cycles évanescents sur V, on applique donc le théorème de Fujiwara à
$V$. 

La fibre des cycles évanescents en un point fixe s'identifie à la
cohomologie du tube au dessus de ce point fermé. 
Le calcul des termes locaux na{\"\i}fs s'effectue alors en utilisant en
utilisant les deux fonctorialités de l'isomorphisme de
comparaison de Berkovich \ref{Corig1} et \ref{Corig2} ainsi que la
compatibilité de cet isomorphisme à l'action de Galois. 
\end{proof}

\section[Modèles en niveau parahorique]{Modèles des variétés de
  Shimura et des espaces de 
Rapoport-Zink en niveau parahorique} 

Nous avons introduit dans le premier chapitre des modèles entiers de
nos espaces en niveau compact hyperspécial. Nous aurons néanmoins
besoin, pour des raisons techniques liées à la théorie des types, des
modèles définis dans \cite{RZ} en niveau parahorique en $p$. Ces
modèles ne nous servirons que comme intermédiaires de
démonstration. Nous ne rappellerons donc pas leur définition. 
Nous reprenons les notations globales de la première
partie. $\mathcal{D}$ désignera donc une donnée de type P.E.L. non ramifiée.

\subsection{Sous-groupes parahoriques et leurs normalisateurs}

Soit $\LL$ une multichaîne polarisée de réseaux dans $V_{\Qp}$
(\cite{RZ} définition 3.14). 

Notons 
$$
K_\LL=\{ g\in G (\Qp) \; |\; \forall \Lambda \in \LL\;
g.\Lambda=\Lambda \}
$$
le sous-groupe parahorique associé. Le normalisateur dans $G(\Qp)$ de
$K_\LL$ est alors défini par 
$$
N_\LL = \{ g\in G (\Qp) \; |\; \forall \Lambda \in \LL \;
g.\Lambda\in \LL \}
$$
Si $\LL=(\Lambda_i)_{i\in \Z}$ est uniforme au sens où $i\mapsto
[\Lambda_i : \Lambda_{i+1}]$ est constant alors 
$$
N_\LL= \{ g \in G(\Qp)\; |\; \exists j\in \Z \;
g.\Lambda_i=\Lambda_{i+j} \}
$$
Nous n'aurons besoin dans la suite que de ce cas là et nous
supposerons donc que toutes nos multichaînes sont uniformes.

\begin{exem}
Pla{\c c}ons nous dans le cas (A) lorsque $J=\emptyset$. Par équivalence de
Morita, la donnée d'une telle multichaîne $\LL$ est alors équivalente
à la donnée pour tout $i$ dans $I$ d'une multichaîne 
$\LL_i$ dans $F_{v_i}^n$. Pour un $i$ fixé la donnée d'une telle
multicha{\^\i}ne est alors équivalente (modulo l'action de  $\Gl_n
(F_{v_i})$) à la donnée d'un entier $d$ divisant $n$. Une cha{\^\i}ne
associée à $d$ est alors de la forme 
$$
\dots \subset \underbrace{<e_1,\dots,e_n>}_{\Lambda_0} \subset \dots
\subset \underbrace{<\varpi_{F_{v_i}}^{-1} e_1,\dots,\varpi_{F_{v_i}}^{-1}
e_{kd},e_{kd+1},\dots,e_n>}_{\Lambda_k} \subset \dots \subset
\varpi_{F_{v_i}}^{-1} \Lambda_0\subset \dots 
$$ 
et le groupe $K_\LL$ se décrit ainsi par blocs de taille $d\times d$ : 
$$
K_\LL= \left \{ \left (
\begin{matrix}
\Gl_d (\O_{F_{v_i}}) & \dots & M_d (\O_{F_{v_i}}) \\
\vdots & \ddots \\
\varpi_{F_{v_i}} M_d (\O_{F_{v_i}}) & \dots &  \Gl_d (\O_{F_{v_i}}) 
\end{matrix}\right )\right \}
$$
quant à $N_\LL$, c'est le groupe engendré par $K_\LL$ et l'élément
suivant 
$$\left (
\begin{matrix}
0 & & & \varpi_{F_{v_i}}^{-1} I_d \\
I_d & \ddots & & 0 \\
  & \ddots & 0 & \vdots \\
  & & I_d & 0
\end{matrix}\right )
$$

\end{exem}

\subsection{Variétés de Shimura}

Dans \cite{RZ}, chapitre 6, il est défini des modèles entiers des
variétés de Shimura de type P.E.L. que nous considérons sur $E_{\nu}$
en niveau parahorique, $\Sh_{K_\LL K^p}$. Nous noterons $S_{K_\LL K^p}$
ces modèles entiers. Lorsque $K^p$ varie ils forment une tour de
variétés quasiprojectives sur $\spec (\O_{E_\nu})$ munie d'une action
de $G(\A_f^p)$. 

\begin{exem}
 Soit $\La_0$ un réseau autodual dans $V_{\Qp}$ de stabilisateur $C_0$
  dans $G(\Qp)$. 
  Supposons que $\LL$ soit la plus petite multichaîne contenant $\La_0$. Alors, $K_\LL=C_0$ est compact hyperspécial
  et $S_{K_\LL 
  K^p}$ est le modèle entier défini précédemment en niveau hyperspécial.
\end{exem}

Si $\LL$ contient un réseau autodual $\La_0$ alors, $K_\LL\subset C_0$
et il y a un morphisme  de tours :
$$
(S_{K_\LL K^p} )_{K^p} \ldrt ( S_{C_0 K^p})_{K^p}
$$
étendant celui défini sur la fibre générique.

\begin{rema}
  En général ce morphisme n'est pas fini !
\end{rema}

On ne supposera pas que $\LL$ contient un réseau autodual, ce qui nous
permettra, par exemple, d'inclure parmi les $K_\LL$ des sous-groupes
compacts maximaux des groupes de similitudes unitaires non
hyperspéciaux. 

\subsection{Espaces de Rapoport-Zink}

A $b\in B (G_{\Qp}, \mu_{\Qpb})$ est associé un espace de Rapoport-Zink
$$
\M_{\LL} ( \mathcal{D}_{\Qp}, b)/ \spf (\O_{Eb_\nu})
$$
muni d'une action de $J_b$ (\cite{RZ} définition 3.21). Comme
précédemment, cet espace se scinde en produit d'espaces de
Rapoport-Zink associés à des données simples. 

Si $K_p\subset K_\LL$, il y a une tour d'espaces rigides $\M_{K_p}$ au
dessus de $\M_\LL^{an}=\M_{K_\LL}$ qui coïncide avec celle définie
auparavant lorsque $K_p\subset K_\LL\cap C_0$. 

\subsection{Extension de l'action du normalisateur d'un parahorique
aux modèles entiers}\label{ext_para_lu}

\subsubsection{Variétés de Shimura}

L'action de $G (\Qp)$ sur la tour $(\Sh_{K})_K$ permet de définir un
morphisme de groupes trivial sur $K_\LL$ 
$$
N_\LL\ldrt \Aut ( \Sh_{K_\LL K^p} )
$$
Le but de cette section est d'étendre ces automorphismes aux modèles
entiers ci dessus, c'est à dire définir un morphisme 
$$
N_\LL \ldrt \Aut ( S_{K_\LL K^p} )
$$
redonnant celui d'avant sur la fibre générique et commutant à l'action
de $G(\A_f^p)$. 
\\

Soit donc $g\in N_\LL$. Soit $T$ un schéma sur $\spec
(\O_{E_\nu})$. Un élément de $S_{K_\LL K^p} (T)$ est donné par un
triplet $(A,\l,\bar{\eta^p})$ où $A$ est une $\LL$-multichaîne de
variétés abéliennes sur $T$, $\l$ une polarisation de cette
multichaîne et $\eta^p$ une structure de niveau hors $p$, le tout
vérifiant certaines conditions (\cite{RZ} chapitre 6). 

Posons  $g.(A,\l,\bar{\eta^p})= (A',\l',\bar{\eta^p})\in S_{K_\LL K^p}
(T)$  où  avec les
notations du chapitre 6 de \cite{RZ} : 
$$
\forall \La\in\LL \;\;\; A'_\La=A_{g.\La} \text{ muni de son action de } \O_B
$$
$$
\forall \La_1,\La_2\in \LL\;\; \La_1\subset \La_2 \;\;\;
\rho'_{\La_1,\La_2}=\rho_{g.\La_1, g.\La_2}
$$
Posons pour tout $a$ dans $B^\times$ normalisant l'ordre $\O_B\otimes
\Z_{(p)}$ 
$$
\theta'_a=\theta_a 
$$

Définissons $\l'$. On doit définir pour tout $\La\in \LL$ une
quasi-isogénie 
$$
\l'_{\La}: A'_{\La} \ldrt \left ( A'_{\La^\vee}\right )^\vee
$$
telle que 
$$
\left ( \rho'_{\La^\vee,\La} \right )^\vee \circ \l'_{\La} : A'_{\La}
\ldrt \left ( A'_{\La} \right )^\vee
$$
soit une polarisation de $A'_{\La}$. 
Posons $\l'_{\La}$ comme étant égal à la quasi-isogénie composée 
\begin{diagram}
  A'_\La=A_{g.\La} & \rTo^{\l_{g.\La}} & \left ( A_{(g.\La)^{\vee}}
  \right )^\vee & \rTo^{( \rho_{ (g.\La)^\vee, g. \La^\vee})^\vee} &
  \left ( A_{g.\La^\vee} \right )^\vee = \left (A'_{\La^\vee} \right )^\vee
\end{diagram}
Alors, 
\begin{eqnarray*}
  \left ( \rho'_{\La^\vee, \La} \right )^\vee \circ \l'_{\La} &=& 
\left ( \rho_{g.\La^\vee, g.\La} \right )^\vee \circ \left ( \rho_{
    (g.\La)^\vee, g.\La^\vee} \right )^\vee \circ \l_{g.\La} \\
 &=& \left ( \underbrace{\rho_{
    (g.\La)^\vee, g.\La^\vee} \circ  \rho_{g.\La^\vee,
    g.\La}}_{\rho_{(g.\La)^\vee,g.\La}} \right )^\vee \circ \l_{g.\La}
\end{eqnarray*}
qui est donc bien une polarisation. 
Il reste à vérifier que $\l': A \ldrt A'$ est bien un morphisme de
multichaînes, mais cela ne pose pas de problème.

Le morphisme ainsi défini vérifie bien les propriétés voulues. 

\subsubsection{Espaces de Rapoport-Zink}

Posons $\M_\LL= \M_\LL (\mathcal{D}_{\Qp},b)$.
On définit comme ci dessus pour les variétés de Shimura un morphisme 
$$
N_\LL \ldrt \Aut ( \M_\LL ) 
$$
trivial sur $K_\LL$, commutant à l'action de $J_b$ et à la donnée de
descente de Rapoport-Zink, et étendant le morphisme défini sur la
fibre générique. 

Sa définition est similaire à celle ci dessus en remplaçant variétés
abéliennes par groupes p-divisibles. 

\subsection{Uniformisation}

Soit $\phi$ une classe d'isogénie comme dans la première partie.
Comme dans la première partie, il y a des isomorphismes de tours
lorsque $K^p$ varie (\cite{RZ} théorème 6.23)
$$
I^\phi(\Q) \bc \left ( \M_\LL \times G(\A_f^p)/K^p \right ) \iso
(S_{K_\LL K^p} \otimes \O_{\Eb_\nu})^{\wedge{\;}}_{/\St}
$$
compatibles à la donnée de descente de Rapoport-Zink et commutant à
l'action de $N_\LL$ définie ci dessus sur les deux membres. 

\section[Formule des traces pour la cohomologie d'une strate]
{La formule des traces pour la cohomologie des strates d'une
variété de Shimura de type P.E.L.} 

Soit $\mathcal{D}$ une donnée globale de type P.E.L. non ramifiée en
$p$. 
Fixons une multichaîne polarisée $\LL$ dans $V_{\Qp}$. Soient $K_\LL$
le groupe parahorique associé et $N_\LL$ son normalisateur.

\subsection{Action de Galois sur les domaines analytiques de
$\M_{K_p}$}

Soit $K_p\subset K_\LL$ un sous-groupe compact ouvert. Rappelons que les
espaces analytiques $\M_{K_p}$ sont définis sur
$\Eb=\widehat{\Qp^{nr}}$. Soit 
$$
\a:\M_{K_p}\ldrt \M_{K_p}^{(\s )}
$$
la donnée de descente de Rapoport-Zink. Soit $U\subset \M_{K_p}$ un
domaine analytique. $\a (U)\subset  \M_{K_p}^{(\s )}$ est alors un
domaine analytique et si $pr: \M_{K_p}^{(\s )}\drt  \M_{K_p}$ désigne
la projection, $pr (\a(U))$ est un domaine analytique dans $ \M_{K_p}$
que nous noterons $U^{(\a)}$. Il y a alors un morphisme (qui n'est pas
défini au dessus de $\Eb$)
$$
U \ldrt U^{(\a)}
$$
D'où un morphisme $\forall \tau\in W_E\; v(\tau)\geq 0$
\begin{diagram}
H^q \left (U^{(\a^{v(\tau)})}\hat{\otimes}\Cp,\Ql\right ) &
\rTo^{\tau^*} & H^q \left ( U\hat{\otimes} \Cp,\Ql\right )
\end{diagram}
(lorsque la cohomologie $\ell$-adique de $U$ est bien définie ce qui
est par exemple le cas si $U$ est le tube au dessus d'un sous-schéma
localement fermé quasicompact de $\Mb$, grâce à \ref{aa} de la première
partie et \ref{cycles_l_adiques}).

\subsection{Tubes et leur cohomologie}

$\M_\LL(\bar{k})$ est muni d'une action d'un Frobenius grâce à
$\a$. Lorsque $K_\LL=C_0$ cette
action co{\"\i}ncide avec celle décrite dans \citeg{Mi1} lorsqu'on voit
$\M_\LL(\bar{k})$ comme sous-ensemble de l'immeuble de $G$ sur $L$. Nous
noterons 
$$
\Phi:\M_\LL (\bar{k})\ldrt \M_\LL (\bar{k})
$$
ce Frobenius.

Ainsi, avec les notations précédentes, si $sp$ désigne le morphisme de
spécialisation associé à $\M_\LL$ 
$$
sp^{-1}(\Phi y)=sp^{-1} (y)^{(\a)}
$$
Les morphismes de changement de niveau commutant à la donnée de
descente on en déduit que 
$$
\Pi_{K_p,K_\LL}^{-1}\left ( sp^{-1}(\Phi.y) \right )= \Pi_{K_p,K_\LL}^{-1}
\left ( sp^{-1} (y) \right )^{(\a)}
$$ 
   
\begin{defi}
$\forall y\in \M (\bar{k}) \; \forall K_p\subset K_\LL$,
$$\M_{K_p}(y)=\Pi_{K_p,K_\LL}^{-1} \left ( sp^{-1} (y)\right )$$ un ouvert
analytique de bord vide dans $\M_{K_p}$ (qui est étale fini au dessus
du disque ouvert  $\M_{C_0}(y)$).
\end{defi}

$H^q( \M_{K_p}(y) \hat{\otimes} \Cp,\Ql)$ est alors défini comme étant
le Poincaré dual de la cohomologie à support compact. C'est un
$\Ql$-espace  vectoriel de dimension finie muni d'une action lisse de
$Stab_{J_b} (y) \times I_{E_\nu}$, $K_\LL$ agissant quant à lui sur la limite lorsque
$K_p$ varie dans $K_\LL$.
\\

Soit donc maintenant $f_p\in \mathcal{H}(G(\Qp))$ de la forme $f'_p*
\delta_z$ où $supp (f'_p)\subset K_\LL$ et $z\in N_\LL$, $\tau\in
W_{E_\nu}$ vérifiant $v(\tau)\geq 0$ et $\g\in J_b$. Supposons que $y\in \M_\LL
(\bar{k})$ vérifie 
$$
z\g\Phi^{v(\tau)}.y=y
$$
Il y a alors des morphismes bien défini pour $K_p$ variant  
\begin{diagram}
H^\bullet ( \M_{K_p}(y),\Ql) & \rTo^{\; z^*\circ \g^* \circ \tau^* \;} & 
H^\bullet (\M_{K_p}(y),\Ql)
\end{diagram}
qui composés avec l'action de $f'_p$ permettent de définir 
$$
\text{Tr}\left ( f_p\times \g\times \tau ; [H^\bullet (\M
(y),\Ql)]\right ) \in \Ql
$$

\subsection{Le théorème}

\begin{defi}
  Soit $b\in B(G_{\Qp},\mu_{\Qpb})$. Pour $K=K_p K^p$ avec $K_p\subset
  C_0$, notons $\Sh^{rig,\geq b}_{K}\subset \Sh_K^{rig}$ le tube au sens des espaces
  adiques au dessus du fermé 
$$
\bigcup_{b'\in B(G_{\Qp},\mu_{\Qpb}) \atop \text{Newton} (b') \geq
  \text{Newton} (b)} \Sb (b')
$$
lorsque $K_p=C_0$, et son image réciproque par les morphismes de
changement de niveau sinon. 
\end{defi}

\begin{defi}
 
 Nous noterons 
$$
[H^\bullet ( \Sh^{rig,\geq b}, \LL_\rho) ] = \sum_i (-1)^i
[\underset{K}{\limi} H^i ( \Sh^{rig,\geq b}_K, \LL_\rho)] \in \Groth (
G(\A_f)\times W_{E_\nu})
$$
qui est bien défini d'après la proposition \ref{cycles_l_adiques} et
coïncide avec la même définition en remplaçant adique par analytique
lorsque $S_K$ est propre ou bien lorsque $b$ est la classe basique. 
\end{defi}

Il résulte alors de la proposition \ref{cycles_l_adiques} que 

\begin{prop}
  Si $b=b_0$ la classe basique, il y a un isomorphisme canonique 
$$
H^\bullet ( \Sh^{rig,\geq b_0}_K,\LL_\rho) \iso H^\bullet ( \Sh^{an}
(b_0),\LL_\rho) 
$$
où le membre de gauche est celui définit ci dessus et celui de droite
celui défini formellement dans la seconde partie comme dual de Poincaré algébrique 
de la cohomologie à support compact du tube au dessus du fermé propre
$\Sb (b_0)$. 
\end{prop}

\begin{theo}\label{Lefschetz_speciale}
Supposons les variétés de Shimura $S_K$ propres sur $\O_{E_\nu}$. 
Soit $b\in B(G,\mu)$. 
  $f_p\in \mathcal{H}(G(\Qp))$ de la forme $f'_p*
\delta_z$ où $supp (f'_p)\subset K_\LL$ et $z\in N_\LL$. Soit
$\widetilde{K} \subset G(\A_f^p)$ un compact ouvert. Il existe alors
un $N\in\N$ tel que 
$$
\forall f^p\in \mathcal{H} (G(\A_f^p))\;\; \text{supp}(f^p)\subset
\widetilde{K}, \;\forall \tau\in W_E \;\; v(\tau)\geq N  
$$
\begin{eqnarray*} &\; & 
\text{Tr} \left ( f_p\otimes f^p\times \tau ;[H^\bullet
(\Sh^{rig,\geq b},\LL_\rho)]\right ) 
 = \sum_{\phi\atop \text{Newt}(b(\phi)) \geq \text{Newt} (b)} \hspace{1mm}\sum_{\{\g\}\in \{I^\phi
(\Q)\}} \hspace{4mm} \sum_{[y]\in
I^\phi_\g (\Q) \bc \text{Fix}(z\times \g\times \Phi^{v(\tau)};
\M_\LL (\bar{k}))} \\  &\; & \text{vol}\left ( I^\phi_\g (\Q)_y\bc I^\phi_\g
(\A_f^p)\right ) \; \text{Tr} \rho (\g) \; \; \text{Tr}\left ( f_p\times
\g\times\tau \;\; ;\; [H^\bullet (\M (y),\Ql)]\right ) \; \O_\g (f^p) 
\end{eqnarray*}
où $I^\phi_\g (\Q)_y = Stab_{I^\phi_\g (\Q)} (y)$.
\end{theo}

\begin{proof}
Toute fonction $f^p\in \mathcal{H} (G(\A_f))$ telle que $
\text{supp}(f^p)\subset\widetilde{K}$ se décompose en combinaison
linéaire de fonctions caractéristiques d'ensembles de la forme $K^p
g^p K^p$ où $K^p$ est un sous-groupe compact ouvert et où $K^p
g^p K^p\subset \widetilde{K}$. Quitte à raffiner une telle
décomposition on peut de plus supposer les $K^p$ suffisamment
petit. On peut donc supposer $f^p=1_{K^p g^p K^p}$ avec $K^p$
suffisamment petit. On peut également supposer que $f'_p= 1_{K_p g_p
  K_p}$ où $K_p\subset K_\LL$ et $g_p\in K_\LL$. 

Considérons le diagramme suivant de correspondances
de Hecke : 
\begin{diagram}
   & & & Sh_{g_p^{-1} K_pg_p\; (g^{p})^{-1} K^p g^p} \\
 & \ldTo & & \dTo && \rdTo \\
\Sh_{K_p K^p} & && (S_{K_\LL\; (g^{p})^{-1}K^p g^p })_\eta & && \Sh_{K_p K^p}  \\
\dTo & &\ldTo &  & \rdTo && \dTo \\
(S_{K_\LL K^p})_\eta & & && && (S_{K_\LL K^p})_\eta
\end{diagram}

où la correspondance du bas est définie sur $\O_{E_\nu}$. Tordons ces
correspondances par l'isomorphisme associé à $z$. Et 
appliquons donc le théorème \ref{Fortree} à ce diagramme. 
\\

Mais pour obtenir l'énoncé annoncé il faut tout d'abord remarquer la chose
suivante :
 lorsque $K^p$ diminue et $g^p$ est tel que $K^p g^p K^p\subset
\widetilde{K}$ le degré des correspondances de Hecke associées $[K^p:
K^p\cap g^p K^p g^{p-1}]$ ne reste pas borné. On ne peut donc pas
appliquer directement le théorème de Fujiwara si l'on veut obtenir la
borne uniforme $v(\tau)\geq N$ annoncée lorsque $\text{supp}(
f^p)\subset \widetilde{K}$. 
Néanmoins, fixons un sous-groupe compact ouvert suffisamment petit $K'$
tel que $K'.\widetilde{K}. K'= \widetilde{K}$ et restreignons nous aux
$K^p$ contenus dans $K'$. Pour un tel $K^p$ et $g^p$ tel que $K^p g^p
K^p \subset \widetilde{K}$, il y a un morphisme de correspondances de
Hecke sur $k$ 
\begin{diagram}
\Sb_{K^p\cap g^p K^p g^{p-1}} & \rTo & \Sb_{K^p}\times \Sb_{K^p} \\
\dTo & & \dTo \\
\Sb_{K'\cap g^p K' g^{p-1}} & \rTo & \Sb_{K'} \times \Sb_{K'}
\end{diagram}
Et on peut appliquer le lemme \ref{ramene} à ces correspondances
étendues à $\bar{k}$ puis tordues par Frobenius. Maintenant,  le degré
de la correspondance du bas est donné par $[K' : K'\cap g^p K'
g^{p-1}]$. Or, $K'\bc \widetilde{K} /K'$ est fini et $g^p\in
\widetilde{K}$. On en déduit que ce degré reste borné. Et on peut donc
appliquer le théorème de Fujiwara sous la forme du théorème \ref{Fortree}
uniformément lorsque $\text{supp}(f^p)$ reste dans $\widetilde{K}$.  
\\

Pour tout $b'\in B(G_{\Qp},\mu_{\Qpb})$, 
$\Sb (b')( \bar{k})$ est une union
disjointe selon les classes d'isogénie $\phi$ de sous-ensembles $\St
(\bar{k})$ stables par les correspondances de Hecke. Les points fixes de
$Fr^{v(\tau)}\times \text{Hecke}(K^pg^pK^p)$ sont donc une union
disjointe de points fixes associés aux différentes de $\phi$, d'où la
première somme dans la formule annoncée. 

De plus, 
$$
\St(\bar{k})\simeq I^\phi (\Q) \bc \left ( \M_\LL (\bar{k})\times
G(\A_f^p)/K^p\right ) 
$$
où $Fr$ agit via $[\Phi\times Id]$, les correspondances de Hecke de
fa{\c c}on usuelle et $z$ via $[z\times Id]$. Appliquant le lemme 5.3 de
\citeg{Mi1} on obtient 
$$
\text{Fix} \left ( Fr^{v(\tau)}\times K^p g^p K^p\times z ; \St
(\bar{k})\right ) = \coprod_{\{\g\} \in \{ I^\phi (\Q) \}} 
I^\phi_\g (\Q) \bc \text{Fix} \left( \Phi^{v(\tau)} \times z \times \g ;\M_\LL
(\bar{k})\right ) \times \text{Fix} \left ( \g\times  K^p
g^p K^p \right ) 
$$

Fibrons cet ensemble de points fixes par les éléments de 
$$
I^\phi_\g (\Q)\bc \text{Fix} \left (\Phi^{v(\tau)}\times z \times \g ;  \M_\LL
(\bar{k})\right )
$$
Le cardinal de la fibre associée à la classe de y est
$$
\text{vol} \left ( I^\phi_\g (\Q)_y\bc I^\phi_\g (\A_f^p)\right )\;\;
\O_\g (f^p)  
$$ 
De plus, en tout point $\b$ de la fibre en $[y]$, 
$$
\widehat{\Sb}_{/\b}\simeq \M_{\LL/y}^{\widehat{\;}}
$$
Et les rev{\^e}tements rigides associés sont donnés par l'uniformisation
rigide des variétés de Shimura. On en déduit que le terme local
associé dans le théorème \ref{Fortree} en un tel point $\b$ est 
$$
 \text{Tr}\left ( f_p\times
\g\times\tau \;\; ;\; [H^\bullet (\M (y),\Ql)]\right ) \;\text{Tr}
\rho (\g)
$$
où le terme $\text{Tr}\rho (\g)$ provient du fait que $\LL_\rho^{an}$
restreint aux tubes ci dessus est le faisceau constant
$\underline{V_\rho}$. 
\end{proof}

\begin{rema}\label{finitude_classes_conj}
Dans la formule précédente, lorsque $\widetilde{K}$,$v(\tau)$ et $z$ sont fixés,
indépendamment de $f^p$ le nombre de $\phi$ et de classes de
conjugaison $\{\g\}$ contribuant de fa{\c c}on non nulle est fini.
Pour le voir, posons $f^p=1_{\widetilde{K}}$. Ecrivons $\widetilde{K}$
comme une union disjointe d'un nombre fini de classes doubles 
$$
\widetilde{K} =\coprod_{a\in A} K^p g_a K^p
$$
pour un sous groupe compact ouvert $K^p$ de $G(\A_f^p)$. 
Appliquons la méthode de la démonstration à chacune des
correspondances associée à ces doubles classes pour compter le nombre
de points fixes des correspondances associées (pas la trace des
correspondances cohomologiques). On obtient une formule comme dans la
démonstration du théorème, sans les traces des correspondances
cohomologiques. 
 La finitude du nombre
de points fixes implique alors qu'il existe un nombre fini de
$(\phi,\{\g\})$ tels que $$ \text{Fix} \left( \Phi^{v(\tau)} \times z \times \g ;\M_\LL
(\bar{k})\right )\neq \emptyset$$ 
et $\O_{\gamma} ( 1_{\widetilde{K}}) \neq 0$. Ce dernier point étant
équivalent à dire que la classe de $G(\A_f^p)$ conjugaison de $\gamma$ rencontre
le compact $\widetilde{K}$. L'assertion s'en déduit. 
\end{rema}

%%% Local Variables: 
%%% mode: latex
%%% TeX-master: "thesef"
%%% End: 

\chapter[Formule de Lefschetz sur la fibre générique]
{Formule de Lefschetz sur la fibre générique}\label{premsect} 

Dans ce chapitre, nous redémontrons un cas particulier très simple de
la formule d'Arthur pour la trace des correspondances Hecke sur la
cohomologie d'intersection des variétés de Shimura, telle qu'elle est
réinterprétée topologiquement dans \cite{Gor_Mac1}. En faisant une
hypothèse simple (la correspondance de Hecke est à support dans les
éléments réguliers en une place), on s'aperçoit que tous les points
fixes sont des points C.M.. Si l'on se restreint aux variétés de
Shimura de type P.E.L. considérées dans la première section, nous
montrons  (\ref{lehteoint}) que 
si un tel point fixe se spécialise sur la strate indexée par un $b\in
B (G_{\Qp},\mu_{\Qpb})$ alors la classe de conjugaison $\g$ dans
$G(\Q)$ associée à ce point fixe est conjuguée dans $G(\Qp)$ à un
élément du centralisateur du morphisme des pentes (un sous-groupe de
Levi de $G(\Qp)$). En particulier, si $\g$ ne se spécialise pas sur la
strate basique, $\g$ n'est pas elliptique. Bien que nous n'utiliserons
pas ce fait par la suite nous avons inclus le théorème \ref{lehteoint}
car il fournit une motivation pour l'énoncé ``la cohomologie des
strates non basiques est induite en $p$'' (théorème \ref{coho_ind_en_p}). En
effet, si l'on disposait d'une formule des traces rigides pour la
trace d'un opérateur de Hecke sur la cohomologie à support compact du
tube rigide sur la 
strate basique faisant intervenir une somme sur des points fixes naïfs 
+ des termes au bord, cela montrerait que la distribution somme sur les points fixes
naïfs associée aux strates non basiques est à support dans les
éléments non elliptiques. 

\section{Généralités sur les points fixes sur la fibre générique}

Soit $(\Sh_K)_{K\subset G(\A_f)}$ une variété de Shimura associée à une donnée de
Shimura
$(G,X)$. 
La variété algébrique 
$\Sh_K$ est définie sur le corps reflex associé $E$ et les
correspondances de Hecke également. La correspondance associée à $KgK$
où $g\in G(\A_f)$ et $K$ est un sous-groupe compact ouvert de
$G(\A_f)$ est :
\begin{diagram}
  && \Sh_{gKg^{-1}\cap K} & \rTo^g & \Sh_{K\cap g^{-1}Kg} \\
& \ldTo^\Pi && & & \rdTo^\Pi \\
\Sh_K & & &&&& \Sh_K
\end{diagram}

\begin{lemm}
  Toute fonction $f\in \mathcal{H} (G(\A_f))$ est combinaison linéaire de
  fonctions caractéristiques $1_{KgK}$ telles que la correspondance de
  Hecke associée soit un cycle dans $\Sh_K\times \Sh_K$ :
$$
\Sh_{K\cap gKg^{-1}} \hookrightarrow \Sh_K\times \Sh_K
$$
\end{lemm}
\begin{proof}
  Si $K_1,K_2$ sont deux sous-groupes compacts ouverts de $G(\A_f)$
  contenus dans un même sous-groupe compact ouvert suffisamment petit
  alors
$$
\Sh_{K_1\cap K_2}  \hookrightarrow \Sh_{K_1}\times \Sh_{K_2}
$$
Il suffit donc de montrer que $f$ est combinaison linéaire de
fonctions $1_{KgK}$ où $K$ et $gKg^{-1}$ sont suffisamment petits. Or,
$f$ est combinaison linéaire de fonctions $1_{KgK}$ avec $K$ aussi
petit que l'on veut et $KgK\subset \text{supp} (f)$ ($\impl g\in
\text{supp} (f)$). $\text{supp}(f)$ étant compact, lorsque $K$
``tend'' $\{1\}$, uniformément pour $g\in\text{supp} (f)$, $gKg^{-1}$
``tend'' vers $\{ 1\}$.
\end{proof}

Étant donné que nous nous intéressons à la trace d'une telle fonction
$f$ dans la cohomologie de variétés de Shimura, ce lemme justifie que
nous ne nous restreignions qu'à de telles correspondances de Hecke.

\begin{enonce}{Hypothèse} 
  Nous supposerons dans ce chapitre que les correspondances de Hecke
  sont toutes définies par des cycles. 
\end{enonce}

De telles correspondances sont donc définies par un cycle lisse de $\Sh_K\times
\Sh_K$. 

Si de plus $\rho$ est donné, la double classe $KgK$ définit une correspondance
cohomologique sur $\LL_\rho$ à support dans ce cycle :
$$
c_2^!\LL_\rho = g^*\Pi^*\LL_\rho =g^*\LL_\rho \ldrt \LL_\rho =c_1^*\LL_\rho
$$
Ceci définit des morphismes d'algèbres compatibles pour $K$ variant :
$$
\mathcal{H}(G,K)\ldrt \Coh ( \LL_\rho,\LL_\rho)
$$
Et si 
\begin{eqnarray*}
\pi : \mathcal{H} (G,K) &\ldrt& \End ( H^\bullet ( \Sh_K,\LL_\rho)) \\
 f & \longmapsto & \int_{G(\A_f)} f(g) \rho (g) dg
\end{eqnarray*}
où $\rho$ est l'action de $G(\A_f)$ sur la cohomologie, alors le
diagramme suivant commute 
\begin{diagram}
  \mathcal{H}(G,K) && \rTo &&  \Coh ( \LL_\rho,\LL_\rho) \\
& \rdTo & & \ldTo \\
 & &  \End ( H^\bullet ( \Sh_K,\LL_\rho)) 
\end{diagram}
où l'application verticale de droite est définie par l'action d'une
correspondance cohomologique sur la cohomologie (\ref{Rappels_coho}). 
\\

Nous allons tout d'abord nous intéresser à la sous-variété des points
fixes de cette correspondance :
$$
\Fix=(\Delta\cap \Sh_{K\cap gKg^{-1}})_{\text{red}}
$$
où $\Delta$ désigne la diagonale de $\Sh_K$. Pour cela, nous allons
décrire ses points complexes :
$$
\Fix(\C)=\Delta(\C)\cap \Sh_{K\cap gKg^{-1}}(\C)=\{ x\in 
\Sh_{K\cap gKg^{-1}}(\C)\;|\; \Pi(x)=\Pi(g.x) \}
$$
Utilisons l'uniformisation complexe : 
$$
\Sh_K(\C)=G(\Q)\bc(X\times G(\A_f)/K)
$$
Dans laquelle $G(\Q)$ opère librement sur $X\times G(\A_f)/K$ (car
$K$ est suffisamment petit). On utilisera également le fait que le
morphisme 
$$
p:X\times G(\A_f)/K\ldrt G(\Q)\bc (X\times G(\A_f)/K)
$$
est un isomorphisme analytique local : tout $(x,yK)\in X\times
G(\A_f)/K$ est tel qu'il existe $\Omega$ voisinage de $x$ dans $X$ tel
que $p:\Omega\times \{yK\} \iso p(\Omega \times \{yK\})$ soit un isomorphisme
analytique. 

La correspondance de Hecke associée à $KgK$ s'écrit alors avec ces
notations : 
\begin{diagram}
  &&& [x,y(gKg^{-1}\cap K)] \\
& \ldMapsto & &&& \rdMapsto \\
[x,yK]  &&&&&& [x,ygK]
\end{diagram}
et 
$$
\Fix(\C)=\{ [x,y(gKg^{-1}\cap K)] \;|\; [x,yK]=[x,ygK]\;\}
$$
que l'on peut voir, par hypothèse, comme un sous-ensemble de $\Sh_K
(\C)$ : 
$$
\Fix(\C)=\{ [x,y K)] \;|\; [x,yK]=[x,ygK]\;\}
$$

Or, 
\begin{eqnarray*}
  [x,yK]=[x,ygK] &\lssi & \exists \g\in G(\Q)\;\; \g.x=x\text{ et }
  \g yK=ygK \\
&\lssi& \exists\g\in G(\Q)\cap ygKy^{-1} \;\; \g.x=x
\end{eqnarray*}
$x,yK$ étant fixés un tel $\g$ est unique (action libre de $G(\Q)$ sur
$X\times G(\A_f)/K$) et est semi-simple elliptique dans $G(\R)$ car
$\g\in\text{Stab}_{G(\R)}(x)$.
De plus, si $\g'\in G(\Q)$, alors  
$$
\left\{ 
  {
    \g.x=x \atop
\g.yK= ygK } \right.
\lssi
\left\{ 
  {
    (\g'\g{\g'}^{-1})\g'.x=\g'.x \atop
  (\g'\g{\g'}^{-1})\g' yK=\g'ygK
  }\right.
$$
Et donc $\g'\g{\g'}^{-1}$ est l'unique élément de $G(\Q)$ vérifiant les
égalités de droite.

On a donc démontré :
\vspace{2mm}

\begin{lemm}
  Il y a une partition naturelle 
$$
\Fix(\C)=\coprod_{\{\g\}\in \{G(\Q)\}_{ss}\atop \g \text{ elliptique
    dans }G(\R)} \Fix(\C)_{\{\g\}}
$$
où $\{G(\Q)\}_{ss}$ désigne les classes de conjugaison 
semi-simples de $G(\Q)$ et 
$$
\Fix(\C)_{\{\g\}}=\{[x,yK]\; |\; \exists \tilde{\g}\in \{\g\} \;\;
\tilde{\g}.x=x \text{ et } \tilde{\g}yK=ygK \;\}
$$
\end{lemm}

\begin{rema}
  On aurait également pu déduire cette décomposition du   lemme 5.3 de
\citeg{Mi1} comme dans la démonstration de la formule des traces sur
la fibre spéciale (théorème \ref{Lefschetz_speciale}). 
\end{rema}

\begin{lemm}
  Si la classe de conjugaison $\{\g\}$ est associée à un point fixe de
  $KgK$ alors $\{\g\}$ est 
  associé à tout point de la composante connexe de ce point dans
  $\Fix(\C)$. On a donc une partition 
$$
\pi_0(\Fix(\C))=\coprod_{\{\g\}\in\{G(\Q)\}_{ss}\atop \g \text{ elliptique
    dans }G(\R)
} \pi_0(\Fix(\C))_{\{\g\}}
$$
\vspace{3mm}
\end{lemm}
\begin{proof}
 Il suffit de démontrer que tout $[x,yK]\in
\Fix(\C)_{\{\g\}}$ possède un voisinage $U$ dans $\Sh_K(\C)$ tel que
$U\cap \Fix(\C)\subset \Fix(\C)_{\{\g\}}$ ce qui montrera que
$[x,yK]\mapsto \{\g\}$ est localement constante donc constante sur
chaque composante connexe. 

Fixons $[x,yK]$ dans $\Fix (\C)$ tel que $\g.x=x$ et $\g yK =y gK$ où
$\g\in G(\Q)$. 

Si $\Omega$ est un voisinage compact de $x$ dans $X$, 
$$
\{g\in G(\R)\; |\; g.\Omega\cap \Omega \neq\emptyset \;\} 
$$
est un compact de $G(\R)$ (écrire $X$ sous la forme $G(\R)/K_\infty$). 

Si de plus $\Omega$ est petit $p$ induit un isomorphisme 
 $\Omega \times yK\iso p(\Omega\times
\{yK\})\subset \Sh_K(\C)$ où $p(\Omega\times \{yK\})$ est un voisinage de
$[x,yK]$ dans $\Sh_K(\C)$. 

Fixons un $\Omega$ vérifiant les deux conditions ci dessus : $\Omega$
est  compact
et petit. Si
$(x',yK)\in \Omega \times \{yK\}$ est tel que $[x',yK]=p(x',yK)\in \Fix (\C)$ alors
soit $\g'\in G(\Q)$ l'unique élément tel que 
$$
\left\{ {\g'.x'=x' \atop 
\g'.yK =ygK } \right.
$$
On a donc, $\g'.\Omega\cap \Omega \neq\emptyset \text{ et } \g'\in
ygKy^{-1} \cap G(\Q)$ 

Par discrétude de $G(\Q)$ dans $G(\A)$, 
$$
A_\Omega=
\{ \g''\in G(\Q) \; |\; \g''.\Omega \cap \Omega \neq \emptyset \text{
  et } \g''\in ygKy^{-1} \; \} 
$$
est fini. Or, lorsque $\Omega$ parcourt des ensembles vérifiant les
conditions ci dessus,
$$
\bigcap_\Omega A_\Omega =\{\g \}
$$
Et donc, quitte à rétrécir $\Omega$ on peut supposer qu'il contient un
seul élément. On pose alors $U=p(\Omega \times yK)$. 
\end{proof}

\begin{prop}
  Tout élément de $\pi_0(\Fix(\C))_{\{\g\}}$ est une composante
  connexe de l'image d'une variété de Shimura associée à $G_\g^0$ où $G_\g$
  désigne le
  centralisateur de $\g$. 
\end{prop}
\begin{proof}

Soit $[x,yK]\in \Fix(\C)_{\{\g\}}$ où $\g.x=x\text{ et }
\g.yK=ygK$. Le point $x$ correspond  à un morphisme  $h_x:\mathbb{S}\ldrt G_{/\R}$, et 
$$
\g.x=x\lssi h_x:\mathbb{S}\ldrt (G_\g^0)_{\R}
$$
Le couple 
$(G_\g^0,h_x)$ définit alors une variété de Shimura de domaine
hermitien associé 
$$ Y= 
G_\g^0(\R)/C_{G_\g^0(\R)} (h_x)\hookrightarrow X
$$
Posons $K'=y K y^{-1} \cap G_\g^0 (\A_f)$. 
Le morphisme  $(G_\g^0,h_x)\ldrt (G,h_x)$ induit
alors un morphisme analytique 
\begin{eqnarray*} q: 
\Sh_{K'} (G_\g^0,h_x)(\C) &\ldrt  &
\Sh_K(G,X)(\C) \\  
\; [ x_1, y_1 K' ] &\longmapsto &  [x_1, y_1 y K ]
\end{eqnarray*}
qui est un revêtement fini au dessus de son image. 

De plus, 
$$
[x,yK]\in \im (q), \text{ et } \im (q)\subset \Fix (\C)_{\{\g\}}
$$
Considérons maintenant le lemme suivant 

\begin{lemm}
  Soient $x\in X$ et $g_0\in G(\R)$ tels que $g_0.x=x$. Alors, 
\begin{eqnarray*}
\{ x'\in X\; |\; g_0.x'=x' \; \} &=& \{ g.x \; |\; g\in
C_{G(\R)}(g_0)\;\} \\
&\simeq & C_{G(\R)}(g_0)/ (C_{G(\R)}(g_0)\cap K_\infty) 
\end{eqnarray*}
\end{lemm}
 
\begin{proof}
Soit $K_\infty =Stab_{G(\R)} (x)$. Soit 
$$\mathfrak{g}=\mathfrak{k}\oplus\mathfrak{p}
$$
 une décomposition de Cartan de $\text{Lie}(G(\R))$ où
 $\mathfrak{k}=\text{Lie}(K_\infty)$. L'application 
 \begin{eqnarray*}
 exp: \mathfrak{p} &\iso & X \\
Z &\longmapsto & exp(Z) .K_\infty \\
0 &\longmapsto & x 
 \end{eqnarray*}
est un difféomorphisme. 

Si $x'=\exp(Z)K_\infty\in X$ où $Z\in \mathfrak{p}$, 
$$
g_0.x'=x'  \lssi  \underbrace{
(g_0 \exp(Z)g_0^{-1})}_{\exp(\text{Ad}(g_0)(Z))}.K_\infty
=\exp(Z).K_\infty \\
$$
car $g_0\in K_\infty$. De plus, $g_0\in K_\infty \impl 
\text{Ad}(g_0)(\mathfrak{p})\subset \mathfrak{p}$ et donc,
$\text{Ad}(g_0)(Z)\in\mathfrak{p}$ et l'égalité 
$\exp(\text{Ad}(g_0)(Z)).K_\infty
=\exp(Z).K_\infty$ implique $\text{Ad}(g_0)(Y)=Y$,
 qui implique elle même que $g_0$ commute à $\exp(Z)$. Donc :
\begin{eqnarray*}
\{ x'\in X\; |\; g_0.x'=x' \; \} &=& \{ g.x \; |\; g\in
C_{G(\R)}(g_0)\;\} \\
&=& C_{G(\R)}(g_0)/ (C_{G(\R)}(g_0)\cap K_\infty) 
\end{eqnarray*}
\end{proof}

Appliquons ce lemme à $x$ et $g_0=\g$. On trouve que l'ensemble des
points fixes de $\g$ dans $X$ est $Y$. Le résultat s'en déduit.
\end{proof}

\begin{coro}
  $\Fix$ est lisse sur $E$ (en particulier ses composantes connexes 
sont irréductibles). 
\end{coro}
\vspace{1mm}

\begin{exem} Dans la décomposition
$$
\pi_0(\Fix(\C))=\dpt{\coprod_{\{\g\}\in\{G(\Q)\}_{ss}} \pi_0(\Fix(\C))_{\{\g\}}
}$$ 
\begin{itemize}
\item
Les classes de conjugaison elliptiques correspondent aux composantes
de points fixes compactes. 
\item Les classes de conjugaison régulières correspondent aux points
  fixes isolés (et elles sont automatiquement elliptiques dès qu'un
tel point fixe existe) 
\end{itemize}
\end{exem}

\section{Points fixes isolés}

L'idée sous-jacente à cette section est que tout vecteur tangent dans
l'espace tangent en un point fixe qui est fixe infinitésimalement
donne lieu à une géodésique dans l'espace hermitien $X$ qui est formée de
points fixes, et que donc en tout point fixe isolé l'intersection est
transverse. 
\\

Les points fixes isolés correspondent comme on l'a  vu aux points
fixes associés à des classes de conjugaison $\{\g\}$ où $\g$ est
régulier c'est à dire $G_\g^0$ est un tore (qui est automatiquement
elliptique si un tel point existe). On a donc en particulier :

\begin{lemm}
Les points fixes isolés sont des points C.M. de $\Sh_K$.
\end{lemm}

Nous aurons besoin du lemme suivant pour appliquer une formule des
traces de Lefschetz. 

\begin{lemm}
  Si $\xi\in \Sh_K(\C)$ est un point fixe isolé de $KgK$ alors, en
  $\xi$, $\Delta$ et $\Sh_{gKg^{-1}\cap K}\hookrightarrow \Sh_K\times
  \Sh_K$ s'intersectent transversalement.
\end{lemm}

\begin{proof} Soit 
$$\xi\in\Fix(\C)_{\{\g\}},\; \xi=[x,yK], \;\g.x=x\text{ et } \g yK=ygK
$$
Soit $\Omega$ voisinage de $x$ dans $X$ tel que $p:\Omega\times \{yK\}
\iso p(\Omega\times \{yK\} )$ soit un isomorphisme et $\Fix(\C)\cap
p(\Omega\times yK)=\{\xi\}$.

On a alors que $p:\Omega\times \{y(gKg^{-1}\cap K)\} \iso p(\Omega \times \{y
(gKg^{-1}\cap K)\})\subset \Sh_{gKg^{-1}\cap K}$.

L'inclusion $\Sh_{gKg^{-1}\cap K}\hookrightarrow \Sh_K\times \Sh_K$
s'identifie alors dans la carte locale $\Omega \times \{ yK\} \simeq \Omega$ à :
\begin{eqnarray*}
  \Omega &\ldrt & \Omega\times\Omega \\
x' &\longmapsto & (x',\g.x')
\end{eqnarray*}
Et $\Delta$ à 
\begin{eqnarray*}
   \Omega &\ldrt & \Omega\times\Omega \\
x' &\longmapsto & (x',x')
\end{eqnarray*}

Les intersections sont donc transverses ssi $d(\g)_x:T_x \Omega \ldrt
T_x\Omega $ ne possède pas la valeur propre $+1$. Si
$K_\infty=\text{Stab}_{G(\R)} (x),
\mathfrak{g}=\mathfrak{k}\oplus\mathfrak{p}$ est une décomposition de
Cartan de $\text{Lie}(G(\R))$ telle que
$\mathfrak{g}=\text{Lie}(K_\infty)$ alors, 
\begin{diagram}
  T_x\Omega & \rTo^{d(\g)_x} & T_x\Omega \\
\uTo_{\simeq } & & \uTo_{\simeq} \\
\mathfrak{p} & \rTo^{\text{Ad}(\g)_{|\mathfrak{p}}} & \mathfrak{p}
\end{diagram}
où $\g\in K_\infty$. L'espace propre associé à la valeur propre $+1$
s'identifie donc à $$\text{Lie}(C_{G(\R)}(\g))\cap \mathfrak{p}$$ Mais
$C_G(\g)^0$ étant un tore tel que $C_{G(\R)}(\g)^0\subset K_\infty$, 
$\text{Lie}(C_{G(\R)}(\g))\subset \mathfrak{k}$ et donc 
$$
\text{Lie}(C_{G(\R)}(\g))\cap \mathfrak{p}=0$$

\end{proof}

On déduit également de la démonstration précédente :

\begin{coro} \label{Lef_comp}
  Si $\xi\in \Fix(\C)_{\{\g\}}$ est isolé alors 
$$
\text{Lef}_\xi (KgK,\LL_\rho)=\text{Tr} (\rho(\g))
$$
\end{coro}

\section{Le cas des variétés de Shimura de type P.E.L.}

\subsection{Un lemme sur la conjugaison stable}

\begin{lemm}\label{lemstab}
 Soit $G/\Qp$ un groupe réductif tel que $G^{der}$ soit simplement connexe. 
  Soit $\g\in G(\Qp)$ régulier stablement conjugué à un élément de
  $M(\Qp)$ où $M$ est le Levi d'un parabolique propre défini
sur $\Qp$. Alors, $\g$ est conjugué dans $G(\Qp)$ à un élément de $M(\Qp)$.
\end{lemm}

\begin{proof} Supposons que $\g$ soit stablement conjugué à $\g'\in M
(\Qp)$ :
$$
\g'= g \g g^{-1} \text{  où } g\in G(\Qpb)
$$
$\g$ et $\g'$ étant réguliers, et $G^{der}$ simplement connexe, les
centralisateurs $C_G(\g)$, resp. $C_G(\g')$, sont des tores maximaux
sur $\Qp$ : 
$T$,
resp. $T'$.

L'application définie sur $\Qpb$ 
 $$
\text{int}_g : T \iso T'
$$
est en fait un isomorphisme sur $\Qp$ (les centralisateurs de deux
éléments stablement conjugués sont formes intérieures l'un de l'autre
via cette application, mais étant abéliens, le cocycle intérieur est
trivial : il s'agit de $c_\s=int_{g g^{-\s}}$ où $g g^{-\s}\in T$). Elle
induit donc un isomorphisme 
$$
\text{int}_g : T_{d} \iso (T')_d
$$
(où l'indice $d$ signifie le sous-tore déployé maximal ). 

Rappelons maintenant que si $T$ est un tore maximal de $G$ 
sur $\Qp$, le Levi
d'un  parabolique défini sur $\Qp$ minimal parmi ceux contenant $T$
est $ C_G (T_d)$ (cf. Borel, Algebraic groups)
(par exemple, $T$ est elliptique ssi
$T_d\subset Z_G$ ssi $T$ n'est contenu dans aucun parabolique
propre). 

Dans notre cas, l'isomorphisme $\text{int}_g : T_{d} \iso (T')_d$
implique que $C_G ((T')_d) =  g C_G (T_d) g^{-1}$ et que donc, $\g$ et
$\g'$ sont dans deux paraboliques sur $\Qp$ conjugués dans
$G(\Qpb)$. Mais (Borel, Algebraic groups) deux paraboliques sur $\Qp$ conjugués 
dans $\Qpb$ le sont déjà sur $\Qp$.
\end{proof}

\subsection{}

Plaçons nous maintenant dans le cadre des variétés de Shimura de type
P.E.L. non ramifiées propres étudiées dans la première partie. 

\begin{theo}\label{lehteoint}
Soit $\xi\in \Fix(\C)_{\{\g\}}$ un point fixe isolé. Le point $\xi$ appartient
donc à $\Sh_K(\overline{E})$. Supposons que $\xi$ se spécialise sur un
point de la strate $\overline{S}(b)$ où $b\in B(G_{\Qp},\mu_{\Qpb})$. Alors, $\g$
est conjugué dans $G(\Qp)$ à un élément de
 $M(b)(\Qp)$, le centralisateur
du morphisme des pentes.
\end{theo}

\begin{proof}
Soit $(A,\l,\iota)$ le triplet associé à $\xi$ sur une extension
de degré fini du corps reflex. $A$ est une variété abélienne C.M..
Avec les notations précédentes pour
l'uniformisation complexe on peut supposer que $\xi=[x,yK]$ où
$\g.x=x$ et $\g.yK=ygK$. La relation $\g.x=x$ implique que $\g$ induit un
automorphisme de la $\Q$ structure de Hodge polarisée munie de
l'action de $B$ associée à $(A,\l,\iota)$ et que donc, $\g$ induit un
automorphisme du triplet $(\A,\l,\iota)$ sur une extension de
degré fini du corps reflex. Nous noterons $f$ cet automorphisme.
 Rappelons qu'à l'uniformisation complexe est associée un choix
d'isomorphisme de $B$-modules symplectiques :
$$
V\iso H^1(A,\Q)
$$
Il y a alors un diagramme commutatif 
\begin{diagram}[width=2cm]
V & \rTo^\simeq & H^1_{\text{B}} (A,\Q) \\
\dTo^\g & & \dTo_{f^*} \\
 V & \rTo^\simeq & H^1_{\text{B}} (A,\Q)
\end{diagram}
Rappelons que l'on note $E$ le complété en une place $v$ divisant $p$
du corps réflexe. 
Il y a un isomorphisme canonique 
$$
H^1_{\text{B}} (A,\Q) \otimes \Qp \iso H^1_{\text{ét}} (A,\Qp)
$$
où le second membre est une représentation cristalline du groupe de
Galois absolu d'une extension de degré fini de $E$ que nous noterons
$K$.  On a donc un
diagramme commutatif :
\begin{diagram}[width=2cm]
 V\otimes \Qp  & \rTo^\simeq& H^1_{\text{ét}}(A_K,\Qp) \\
\dTo^{\g\otimes 1} & & \dTo_{f^*} \\
 V\otimes \Qp  & \rTo^\simeq& H^1_{\text{ét}}(A_K,\Qp) 
\end{diagram}
Notons $F$ le foncteur de Fontaine entre la catégorie des
représentations cristallines du groupe de Galois absolu de $K$ et celle
des isocristaux filtrés admissibles. $F(H^1_{\text{ét}}(A_K,\Qp))$
s'identifie comme isocristal muni de structures additionnelles à la
cohomologie cristalline de la réduction mod $p$ du triplet
$(A,\l,\iota)$. Dans cette identification, 
$f$ agit sur cet isocristal par un élément de $J_b$. Fixons un
isomorphisme de $L$-modules symplectiques munis d'une action de $B$ :
$$
 F(H^1_{\text{ét}}(A_K,\Qp)) \simeq (V\otimes \Qp)\otimes L 
$$
Il y a donc un diagramme commutatif 
\begin{diagram}[width=2cm]
F(V\otimes \Qp)  & \rTo^\simeq& F( H^1_{\text{ét}}(A_K,\Qp)) &
\rTo^\simeq & (V\otimes\Qp)\otimes L  \\
\dTo^{F(\g\otimes 1)} & & \dTo_{F(f^*)} & & \dTo^g \\
 F(V\otimes \Qp)  & \rTo^\simeq& F(H^1_{\text{ét}}(A_K,\Qp)) &
 \rTo^\simeq  & 
 (V\otimes\Qp)\otimes L  
\end{diagram}
dans lequel $g\in J_b$. 

Trivialisons maintenant le torseur des périodes; c'est à dire étendons
les scalaires à $B_{dR}$. Il y a un isomorphisme canonique 
$$
(V\otimes \Qp)\otimes B_{dR} \iso F(V\otimes \Qp)\otimes_L B_{dR}
$$
d'où au final un diagramme commutatif 
\begin{diagram}[width=2cm]
V_{\Qp}\otimes B_{dR} & \rTo^\simeq & (V_{\Qp}\otimes L)\otimes B_{dR} \\
\dTo^{\g\otimes 1} & & \dTo_{g\otimes 1} \\
V_{\Qp}\otimes B_{dR} & \rTo^\simeq & (V_{\Qp}\otimes L)\otimes B_{dR}
\end{diagram}
Duquel on déduit que $\g$ est stablement conjugué dans $G(L)$ à un
élément de $J_b$.

Rappelons maintenant que dans le cas que nous considérons ($G_{\Qp}$
quasidéployé) la classe de $\s$-conjugaison 
 $b$ provient d'une classe de $M(b)(L)$ encore notée $b$,
et que $J_b=\{ g\in M(b)(L)\; |\; gb\s=b\s g \}$. Rappelons également
que l'on peut supposer $b$ ``decent'' (confère \citeg{RZ}) au sens où 
pour un $s\in \N$, $(b\s)^s=(s.\nu_b)(p) \s^s$. Dans tous ce qui
précède on peut alors remplacer $L$ par $\Q_{p^s}$ l'extension non
ramifiée de degré $s$ de $\Qp$. La relation $b^{-1} g b= g^\s$ dans
$M(b)(\Q_{p^s})$ montre que la classe de conjugaison de $g$ dans
$M(b)(\Q_{p^s})$ est définie sur $\Qp$ et que donc, d'après \citeg{Ko4} ( 
pour les
groupes avec lesquels nous travaillons $M(b)^{der}$ est simplement
connexe), la classe de conjugaison de $g$
dans $M(b)(\overline{\Qp})$ contient un élément de $M(b)(\Qp)$. 

On déduit donc au final que $\g$ est stablement conjugué dans $G(\Qp)$
à un élément de $M(b)(\Qp)$. On conclut alors grâce au lemme \ref{lemstab}.
\end{proof}

\begin{rema}
  Dans la dernière partie de la démonstration, nous n'avons fait
  qu'expliquer l'existence du transfert des classes de conjugaison
  stables de $J_b$ vers celles de sa forme intérieure quasi-déployée
  $M(b)$. Par exemple, lorsque $J_b$ est le groupe des unités d'une
  algèbre à division $D\*$, ce transfert est le transfert des classes de
  conjugaison de $D\*$ vers les classes de conjugaison elliptiques de
  $\GL_n$. 
\end{rema}

\section{Formule de Lefschetz}

\begin{theo}\label{Lefschetz_generique}
Soit $\Sh (G,X)$ une variété de Shimura compacte. 
Soit $f=\otimes_v f_v\in \mathcal{H} ( G(\A_f))$ telle qu'il existe
$v$ telle que $\text{supp} ( f_v) \subset G(\Q_v)_{reg}$. Alors, 
$$
\text{tr}\left ( f; [H^\bullet(\Sh,\LL_\rho)]\right )  
= \sum_{\{\g\}\in \{ G(\Q)\} \atop {\g \text{régulier} \atop {\g
\text{ elliptique dans } G(\R)}}}\; \text{vol} ( G(\Q)_\g \bc G(\A_f)_\g
)\; \text{tr} \rho (\g) \; \O_\g ( f)
$$
\end{theo}
\begin{proof}
On peut supposer que $f$ est de la forme $1_{KgK}$ où $K=\prod_v K_v$
et $K_v g_v K_v\subset G(\Q_v)_{reg}$.
Appliquons la formule des traces de Lefschetz à la correspondance de
Hecke associée. 
 Soit $\g\in G(\Q)$ elliptique
dans $G(\R)$. $\Fix (\C)_{\{ \g\}}\neq \emptyset \impl \g$ est
régulier. En effet, une relation du type $\g yK= ygK$ implique que $\g
y_v K_v =y_v g_v K_v$  qui implique que $y_v^{-1} \g y_v \in K_v g_v
K_v\subset G (\Q_v)_{reg}$. Tous les points fixes sont donc isolés. Et
la somme dans la formule annoncée provient de la décomposition 
$$
\Fix (\C)=\coprod_{{\{\g\}\in \{ G(\Q)\} \atop
 {\text{régulier} \atop {\g
\text{ elliptique dans } G(\R)}}}} \Fix (\C)_{\g}
$$
Reste donc à montrer que 
$$
\sum_{\xi\in \Fix (\C)_{\{\g\}}} \; \text{Lef}_\xi (KgK, \LL_\rho) = 
 \text{vol} ( G(\Q)_\g \bc G(\A_f)_\g
)\; \text{tr} \rho (\g) \; \O_\g ( f)
$$
On sait (corollaire \ref{Lef_comp}) que $\text{Lef}_\xi (KgK,
\LL_\rho)$ est constant sur $\Fix (\C)_{\{\g\}}$ et vaut
$\text{tr}\rho (\g)$. Il faut donc calculer 
$
\# \Fix (\C)_{\{\g\}} 
$. 
$\g$ étant régulier, $\text{Fix} (\g ; X)$ est réduit à un seul
élément. On en déduit que $\Fix (\C)_{\{ \g\}}$ est en bijection avec
$$
G(\Q)_\g \bc \{ yK \; |\; y^{-1}\g y \in KgK \;\}
$$
qui est de cardinal $\text{vol} ( G(\Q)_\g \bc G (\A_f)_\g) \; \O_\g
(f)$. 
\end{proof}

\begin{rema}
  Cette formule se déduit de \cite{Gor_Mac1} en remarquant grâce au
  théorème 5.2 de \cite{Gor_Mac1} que la somme des caractères des
  séries discrètes de $G(\R)$ ayant même caractère infinitésimal que
  $\check{\rho}$ évaluée en un $\g$ régulier est $\text{tr}\rho
  (\g)$. 
\end{rema}

\begin{rema}
  Lorsque le support de $f$ reste dans un compact fixé de $G(\A_f)$,
  les classes $\{\g\}$ intervennant dans la formule ci dessus
  varient dans un ensemble fini ne dépendant que de ce compact. Cela
  peut se voir de deux façons : la première en utilisant la
  proposition 8.2 de \cite{Ko7} qui implique qu'il y a un nombre fini
  de classes de $G(\A)$-conjugaison contenant un élément  de
  $G(\Q)$ et rencontrant un compact de $G(\A)$ fixé (ce qui est le cas
  ici puisque nos classes sont elliptiques à l'infini); la deuxième de
  la même manière que dans la remarque (\ref{finitude_classes_conj}). 
\end{rema}

%%% Local Variables: 
%%% mode: latex
%%% TeX-master: "thesef"
%%% End: 

\chapter[Caractérisation par les traces]
{Caractérisation des éléments de $\Groth (G( \A_f))\times
W_{E_\nu}$ par leurs traces}

\section{Un lemme d'algèbre linéaire}

Soir $K$ un corps et $V$ un $K$-espace vectoriel de dimension
finie. Soient $u\in \End (V)$ et $v\in \Gl (V)$. Le but du lemme qui
suit est de montrer que si l'on connaît $\text{tr}(u v^N)$ pour $N$
grand alors on connaît $\text{tr} (u)$.

Pour cela on s'inspire de la formule classique : si $P(T)=\det ( Id -
T v)$ alors si $F(T)=\sum_{N\geq 0} \text{tr} ( v^{N+1} ) T^N$, $F$
est une fraction rationnelle et $F(T) dT= - dlog P$. Ainsi, puisque $v$
est inversible $$\text{Res}_\infty (F(T) dT)= v_\infty (P^{-1}) = 
\dim (V)=\text{tr} (Id)$$

\begin{lemm}\label{lemalg}
Soit 
$$
F(T)=\sum_{N\geq 0} \text{tr} (u v^{N+1})\; T^N\in K[[T]]
$$
La série formelle 
$F$ est une fraction rationnelle de degré $-1$ et 
$$
\text{Res}_\infty ( F(T) dT)= \text{tr} (u)
$$
\end{lemm}
\begin{proof}
On peut supposer $K$ algébriquement clos de caractéristique $0$. 
 Soit $v=v_s+ v_n$ la
décomposition de $v$ en parties semi-simples et nilpotentes.
Supposons le lemme démontré pour $v$ semi-simple. 
On a alors
\begin{eqnarray*}
F(T)&=&\sum_{N\geq 0} \sum_{k=0}^{N+1} \left ( N+1 \atop k \right )
\text{tr} \left ( u 
v_n^k  v_s^{N+1-k} \right ) T^N \\
&=&  \sum_{k=0}^{\dim V} \sum_{N+1\geq k \atop N\geq 0}  \left ( N+1
  \atop k \right ) 
\text{tr} \left ( u  
v_n^k  v_s^{N+1-k} \right ) T^N \\
&=& \sum_{N\geq 0} \text{tr} ( u v_s^{N+1}) T^N + \sum_{k=1}^{\dim V}
 \underbrace{ \sum_{N\geq k-1}\left ( N+1 \atop k \right ) \text{tr} \left ( u 
v_n^k  v_s^{N+1-k} \right ) T^N}_{\frac{T^{k-1}}{k!} \frac{d^k}{dT^k}
\left ( F_k (T)\right )}
\end{eqnarray*}
où $$F_k (T)=\sum_{N\geq 0} \text{tr} ( u v_n^k v_s^{-k} v_s^{N+1})
T^N$$ 
Le terme de gauche est une fraction rationnelles par application du cas
semi-simple à $u$ et $v_s$, quant au terme de droite cela résulte du
fait que pour tout $1\leq k\leq n$, $F_k(T)$ est une fraction
rationnelle par application du cas semi-simple à $uv_n^k v_s^{-k}$ et
$v_s$. La série formelle $F$ est donc une fraction rationnelle. De
plus, par hypothèse, 
$\deg F_k =-1$ et donc $\dpt{\deg \left ( \frac{d^k}{dT^k} (T F_k (T))
  \right ) \leq -k-1 }$. Cela implique que 
$$
\deg \left ( \frac{T^{k-1}}{k!} \frac{d^k}{dT^k} ( T F_k (T)) \right )
\leq -2 
$$
et est donc de résidu nul en l'infini. La cas semi-simple implique
donc les autres cas. 

Reste à démontrer le lemme dans le cas où $v$ est semi-simple. Soit
donc $v$ semi-simple. Pour tout $w\in \Gl (V)$ 
$$
\text{tr} ( u (w^{-1} v w)^{N+1})= \text{tr} ( u w^{-1} v^{N+1} w)=
\text{tr} ( (wuw^{-1}) v^{N+1}) 
$$
et $\text{tr} (wuw^{-1})=\text{tr} (u)$. On peut donc supposer que
$V=K^n$ et que la matrice de $v$ dans la base canonique de $K^n$ est
la matrice diagonale $diag (\l_1,\dots,\l_n)$ (où $\forall i
\;\l_i\neq 0$). Notons $(a_{ij})_{i,j}$ la matrice de $u$. Alors, 
\begin{eqnarray*}
F(T) &=&  \sum_{N=0}^{+\infty} \sum_{i=1}^n a_{ii} \l_i^{N+1} t^N \\
 &=& \sum_{i=1}^n a_{ii} \frac{\l_i}{1 - \l_i T}
\end{eqnarray*}
or $\dpt{\l_i\neq 0 \impl \text{Res}_\infty \left ( \frac{ \l_i}{1-\l_i T}
dT\right ) =1}$. D'où le résultat.
\end{proof}

\begin{rema}
L'intérêt de ce lemme est que si $P$ est un polynôme alors 
$$
\text{Res}_\infty P(T) dT =0
$$
Ainsi, si l'on modifie la série définissant $F(T)$ dans le lemme
précédent par un nombre fini de termes on peut encore retrouver la
trace de $u$. 
\end{rema}

\section{Séparation des représentations sphériques par les traces
de fonctions à support dans les éléments réguliers} 

\begin{prop}\label{sep_spherique}
Soit $H/\Qp$ un groupe réductif non ramifié tel que $H(\Qp)\drt H_{ad} (\Qp)$ soit surjectif. Soit 
$$
\sum_{\pi\in \Irr(G(\Qp))} a_\pi [\pi] \in \Groth (G(\Qp))\hat{\otimes}\C
$$  
(i.e. $\a_\pi\in \C$ et 
$\forall \; K_p\subset H(\Qp)$ compact ouvert $\#\{\pi \; |\;
a(\pi)\neq 0\;\text{ et } \pi^{K_p}\neq (0)\;\}<\infty$)

Si $\forall f_p\in \mathcal{H}(H(\Qp)), \text{supp}(f_p)\subset
H(\Qp)_{reg}$ on a :
$$
\sum_\pi a_\pi \text{tr}\; \pi (f_p)=0
$$
Alors, $\forall \;\pi$ sphérique $\; a(\pi)=0$.
\end{prop}
\begin{proof} Ce résultat est essentiellement contenu dans 
 \citeg{Lab2} dont nous utiliserons librement les notations. 

Commençons par remarquer que dans \citeg{Lab2}, si le paramètre $t$
est régulier alors la fonction élémentaire $f_t$ est à support dans
les éléments réguliers de $H(\Qp)$. D'après la proposition 5 de
\citeg{Lab2} on a en prenant $f_p=f_t$ avec $t$ régulier :
$$
\sum_{{\pi\text{sous-quotient}} \atop {\text{d'une série principale} \atop
  \text{non ramifiée}}} a_\pi\; \text{tr}\; \pi (f_t)=0
$$
i.e. on peut séparer les sous-quotients des séries principales non
ramifiées des autres éléments de $\Irr (H(\Qp))$.

Remarquant que les éléments réguliers $t$ engendrent le groupe
$A(\Qp)/A(\Zp)$, 
la proposition 7 de \citeg{Lab2} couplée à l'indépendance linéaire des
caractères  (ici l'indépendance linéaire des caractères 
non ramifiés $\l$
à valeurs dans $\C^*$ du groupe des points
 dans $\Qp$ d'un tore déployé maximal ) permet de séparer les $\pi$
 sous-quotients de séries principales non ramifiées associées à des
 $\l,\l'$ tels que$\forall w\in W\;\; \l^w\neq \l'$. 

La proposition 8 (cf. également la remarque qui la suit) permet à $\l$
(modulo l'action de $W$ ) fixé de séparer dans les sous-quotients
irréductibles de la série principale non ramifiée le sous-quotient
sphérique des autres sous-quotients. D'où le résultat.

\end{proof}

\section{Le théorème}

Soit $G$ un groupe de similitudes unitaires sur $\Q$ comme dans le I. 

\begin{theo}\label{Caract_Traces}
Soit $A\in \Groth ( G(\A_f)\times W_{E_\nu})$ tel que pour tout type
supercuspidal $(J,\l)$ de $G(\Qp)$, $\forall g\in \tilde{J}$ le sous-
groupe compact modulo le centre associé, $\forall f^p\in \mathcal{H} (
G( \A_f^p))$ telle qu'il existe $w \neq p$ vérifiant $\text{supp} (
f_w)\subset G(\Q_w)_{reg}$ il existe un $N\in \N$ tel que $\forall
\tau\in W_{E_\nu}$ vérifiant $v(\tau)\geq N$ si 
$$
f=(e_\l * \delta_g)\otimes f^p
$$
on ait 
$$
\text{tr} ( f\times \tau ; A)=0
$$
Alors,
$$
A_{cusp}=0
$$
\end{theo}

\begin{proof}
fixons $\s\in W_{E_\nu}$ tel que $v(\s)=1$. Considérons la
distribution sur $G(\A_f)$
$$
f \mapsto \sum_{N\geq 0} \text{tr} ( f\times \tau \s^{N+1} ; A) T^N
$$
qui est à valeurs dans $\C (T)$ d'après le lemme
\ref{lemalg}. Appliquant l'application linéaire $$
\text{Res}_\infty : \C(T)\ldrt \C
$$
à cette distribution on déduit du lemme \ref{lemalg} et de la remarque
qui le suit que dans l'énoncé
du théorème on peut prendre n'importe quel $\tau\in W_{E_\nu}$
indépendamment de $f^p$. 

Fixons $\pi_0$ une représentation supercupsidale de $G(\Qp)$ et
$(\l,J)$ un type associé. Notons 
$$
A=\sum_{\Pi\in \text{Irr} ( G( \A_f))} \; [\Pi]\otimes [\rho (\Pi)]
$$
où $[\rho (\Pi)]\in \Groth ( W_{E_\nu})$. Fixons $f^p \in \mathcal{H}(
G( \A_f^p))$ vérifiant $\exists w \neq p\; \text{supp} ( f_w)\subset G
(\Q_w)_{reg}$ et considérons la distribution sur $G(\Qp)$: 
$$
f_p \longmapsto \text{tr}( f_p\otimes f^p \times \tau ;A)
=\sum_{\pi\in \text{Irr} ( G (\Qp))} \underbrace{\left( \sum_{\Pi\in \text{Irr} (
G(\A_f)) \atop \Pi_p \simeq \pi} \text{tr} (f^p;\Pi^p)\; \text{tr} (
\tau; \rho (\Pi))\right )}_{a_\pi} \; \text{tr}_\pi ( f_p)
$$
Il résulte de la définition de $\Groth ( G( \A_f)\times W_{E_\nu})$
que pour tout sous-groupe compact ouvert $K_p$ de $G(\Qp)$ 
$$
\#\{ \pi\; |\; \pi^{K_p}\neq (0) \text{ et } a_\pi\neq 0\: \}<+\infty
$$
Si $K_p$ est tel que $e_\l\in \mathcal{H} ( G(\Qp)// K_p)$ alors
$\forall g\in \tilde{J}\; \forall \pi\;\; \pi (e_\l* \delta_g)=\pi
(e_\l) \pi (g)$ et donc, 
$$
\#\{ \pi\; |\; \pi (e_\l*\delta_g)\neq 0 \text{ et } a_\pi\neq 0\: \}<+\infty
$$
Il existe donc un ensemble fini $E$ de caractères non ramifiés de $G(\Qp)$
vérifiant 
$$
\forall \chi,\chi'\in E\;\; \chi\neq \chi' \limpl \pi_0\otimes
\chi\not\simeq \pi_0\otimes \chi'
$$
et 
$$
\sum_{\chi\in E} a_{\pi_0\otimes \chi} \text{tr}_{\pi_0\otimes\chi} ( e_\l * \delta_g)=0 
$$
or, $\text{tr}_{\pi_0\otimes \chi}( e_\l *\delta_g)= \chi (g)
\text{tr}_{\pi_0} ( e_\l*\delta_g)$ et  $\text{tr}_{\pi_0} (
e_\l*\delta_g)\neq 0$ (lemme \ref{lemme_typique}). On en déduit 
$$
\forall g\in \tilde{J} \;\; 
\sum_{\chi\in E} a_{\pi_0\otimes \chi} \;\chi (g)=0
$$
Rappelons (annexe \ref{Annexe_types}) que $\pi_0\otimes \chi \not\simeq \pi_0\otimes \chi'
\ssi \chi_{|\tilde{J}} \neq \chi'_{|\tilde{J}}$ et donc par
indépendance linéaire des caractères 
$$
a_{\pi_0}=0
$$
\\

Fixons maintenant $\Pi_0\in \text{Irr} (G (\A_f))$ vérifiant $\Pi_{0 p}
\simeq \pi_0$. Fixons toujours $\tau \in W_{E_\nu}$. 
 Soit $w$ une place de $\Q$ en laquelle $G(\Q_w)$ est non
ramifié de type adjoint et $\Pi_{0 w}$ sphérique. Fixons $f^{p w}\in
\mathcal{H} ( G (\A_f^{p w}))$. Il résulte de l'égalité $a_{\pi_0}=0$
que l'on a la relation  
$$
\forall f_w\in \mathcal{H} ( G (\Q_w)) \;\;\; \text{supp} ( f_w)\subset G
(\Q_w)_{reg} 
$$
$$
\sum_{\pi\in \text{Irr} ( G (\Q_w))} \left ( 
\sum_{\Pi\in \text{Irr} ( G (\A_f)) \atop \Pi_p\simeq \pi_0}
\text{tr}_{\Pi^{p w}}( f^{pw}) \text{tr} (\tau; \rho (\Pi)) \right )
\text{tr}_\pi (f_w) =0 
$$
On déduit donc  du lemme \ref{sep_spherique} que 
$$
\forall f^{pw}\in \mathcal{H} ( G (\A_f^{pw})) \; \forall \tau\in
W_{E_\nu} \;\; \sum_{\Pi\in\text{Irr} ( G(\A_f)) \atop {\Pi_p\simeq
\Pi_{0 p} \atop{ \Pi_w\simeq \Pi_{0 w}}}} \text{tr}_{\Pi^{pw}} ( f^{pw} )
\text{tr} ( \tau; \rho (\Pi)) =0
$$
Fixons un niveau $K^{pw}$ tel que $(\Pi_0^{pw})^{K^{pw}}\neq
(0)$. Alors, si
$$ X=
\{ \Pi^{pw}\in \text{Irr} ( G (\A_f^{pw})) \; |\;
(\Pi^{pw})^{K^{pw}} \neq (0) \text{ et } \rho ( \Pi_{0,p}\otimes
\Pi_{0,w}\otimes \Pi^{pw})\neq (0)\}
$$
$\# X <+\infty$. Les représentations semi-simples de $W_{E_\nu}$ de
dimension finie sur $\C$ sont équivalentes aux représentations
semi-simples de dimension finie sur $\C$ de la $\C$-algèbre 
$$ R =\underset{i\in \N}{
\limi} \C [ W_{E_\nu}/ \Gal ( \overline{E}_\nu | E_\nu)^{(i)}]
$$
Il résulte alors de l'indépendance linéaire des caractères que les
traces des éléments de  
$$
\mathcal{H} ( G ( \A_f^{pw}) // K^{pw}) \otimes_\C R 
$$
permettent de séparer les 
$$
\Pi^{pw} \otimes \rho (\Pi_{0,p}\otimes \Pi_{0,w}\otimes \Pi^{pw})\;,\;\;\;
\Pi^{pw}\in X
$$
 On en déduit donc que
$$
[\rho ( \Pi_0)]=0
$$
\end{proof}

%%% Local Variables: 
%%% mode: latex
%%% TeX-master: "thesef"
%%% End: 

\chapter{Le théorème}

Soit $\mathcal{D}$ une donnée de type P.E.L. non ramifiée en $p$, et
$\Sh_K$ les variétés de Shimura associées, que nous supposerons
compactes. Nous supposerons pour cela  que $D=\End_B(V)$ est une
algèbre à division. Nous supposerons de plus qu'en toutes les places finies
de $F$, $D$ est soit déployée soit une algèbre à division.
 Nous supposerons également que l'on est
dans le cas (A), c'est à dire dans le  cas unitaire et 
 que le corps C.M. $F$ est de la forme $F^+\mathcal{K}$ où
$\mathcal{K}|\Q$ est une extension quadratique imaginaire.

Rappelons que l'on note 
$$
G(\Qp)=\prod_{i\in I} \Gl_n (F_{v_i}) \times G ( \prod_{j\in J} U (
n; F_{v_j}) )
$$

Le facteur de droite n'étant pas tout à fait un produit, ses
représentations ne se décrivent pas immédiatement à partir des
représentations de chaque facteur.   
Afin de ne pas trop alourdir les notations  nous supposerons donc dans
ce chapitre que $J=\emptyset$ ou bien que $\#J=1$, ce qui est suffisant
pour les applications que nous avons  en vue. 

Si $$A=\sum_{\Pi \in \text{Irr} (G(\A_f))} [\Pi]\otimes [\s(\Pi)] \in
\Groth (G(\A_f)\times W_{E_\nu})$$ on note $$A_{cusp}=\sum_{\Pi \atop
\Pi_p \text{ supercuspidale}} [\Pi]\otimes [\s (\Pi)]
$$

\begin{theo}\label{coho_ind_en_p}
\begin{enumerate}
\item
Supposons $p$ décomposé dans $\mathcal{K}$  (et donc 
 $J=\emptyset$), alors 
$$ 
\left [ \underset{K}{\limi} H^\bullet (\Sh_K ,\mathcal{L}_\rho) \right ]_{cusp}
=  
\left [ \underset{K}{\limi} H^\bullet (\Sh_K^{an} (b_0) ,\mathcal{L}_\rho) \right ]_{cusp}
$$ 
\item
Supposons $p$ inerte dans $\mathcal{K}$ et $\#J=1$.
Supposons démontré l'existence d'une application de changement de base
local en $p$ pour les groupes unitaires non ramifiés 
(ce qui est le cas pour $U(3)$, \ref{BC_local}) et l'existence de
types pour les représentations supercuspidales des groupes unitaires
non ramifiés en $p$ (confère \ref{types_Un}) (c'est le cas pour $U(3)$
lorsque $p\neq 2$, \ref{types_U3}). On a alors l'égalité 
$$
\left [ \underset{K}{\limi} H^\bullet (\Sh_K ,\mathcal{L}_\rho) \right
]_{|W_{E_\nu \mathcal{K}_p},cusp}
=  
\left [ \underset{K}{\limi} H^\bullet (\Sh_K^{an} (b_0) ,\mathcal{L}_\rho) \right ]_{|W_{E_\nu \mathcal{K}_p},cusp}
$$
\end{enumerate}
\end{theo}

\begin{proof}
La démonstration consiste à comparer deux formules des traces, l'une
sur la fibre générique, l'autre sur la fibre spéciale.
\\

Nous donnons la démonstration  dans le cas où $\#J=1$, le cas où
$J=\emptyset$ étant plus simple puisqu'il suffit dans les
démonstrations de remplacer le groupe de similitudes unitaires associé
à la place indexée par $J$ par
$\Qp^\times$. 
\\

Commençons par rappeler quelques notations. Notons $J=\{ j \}$. 
 Fixons une représentation
supercuspidale $\pi_0$ de $G(\Qp)$. Via l'isomorphisme 
$$
G(\Qp)\simeq \prod_{i\in I} \Gl_n (F_{v_i})\times GU (n;F_{v_j})
$$
notons 
$$
\pi_0=\otimes_{i\in I} \pi_{0,i} \otimes \pi_{0,j}
$$
où les représentations $\pi_{0,i}$  et $\pi_{0,j}$ sont supercuspidales. Pour tout $k$
 dans $I\cup J$  soit 
$(\l_k,J_k)$  un type associé à la classe
d'équivalence inertielle de 
$\pi_{0,k}$. Nous noterons $\tilde{J}_k$
le sous-groupe compact modulo 
le centre  associé
(confère l'annexe \ref{Annexe_types})), et  
$\tilde{\l}_k$ l'extension de $\l_k$ à
$\tilde{J}_k$ vérifiant  
$$ \forall i\in I\;\;
\pi_{0,i}\simeq c-Ind^{\Gl_n( F_{v_i})}_{\tilde{J}_i}
\tilde{\l}_i\text{ et } \forall j\in J\;\; \pi_{0,j}\simeq
c-Ind^{ GU(n;F_{v_j})}_{\tilde{J}_j} 
$$
Pour tout $k$ dans $I\cup J$ fixons $g_k\in \tilde{J}_k$, et soit
$e_{\l_k}$ l'idempotent de l'algèbre de Hecke associé défini par 
$$
e_{\l_k}(g)=\left \{ 0\text{ si } g\notin J_k \atop \frac{\dim
\l_k}{vol (J_k)} \chi{\l_k} (g) \text{ si } g\in J_k \right.
$$

Posons 
\begin{eqnarray*} \forall k\in I\cup J\;\;\;
f_k &=& e_{\l_k}* \delta_{g_k}\\
 f_p &=& \otimes_{i\in I} f_i \otimes f_j \in \mathcal{H} (G(\Qp))
\end{eqnarray*}
Soient $f^p \in \mathcal{H} ( G( \A_f^p))$ et $\tau\in W_{E_\nu\mathcal{K}_p}$.

Commençons par établir la formule des traces sur la fibre générique. 
Plus précisément, on cherche une expression pour 
$$
\text{tr} ( f\times \tau ; [H^\bullet (\Sh,\mathcal{L}_\rho)])
$$
Si la représentation $\Pi \in \text{Irr} ( G( \A_f))$ est telle que $\Pi$ intervienne
dans $[ H^\bullet ( Sh,\mathcal{L}_\rho)]$, la non nullité de
$\text{tr} ( f ; \Pi)$ implique que
$\Pi_p$ est supercuspidale et qu'il existe un caractère non ramifié $\chi$
de $G(\Qp)$ tel que 
$\Pi_p\simeq \pi_0 \otimes \chi$.

Il résulte donc de l'annexe A que :
\begin{eqnarray*} 
&\text{tr} & (f\times \tau ; [ H^\bullet ( Sh, \mathcal{L}_\rho) ] ) \\
&=& \sum_{\Pi\in \text{Irr}( G (\A_f)) \atop {\Pi_p \in [\pi_0, G (\Qp)]
\atop \exists \Pi'\in \mathcal{T} (G)_\rho \; \Pi=\Pi'_f} }
m_\Pi\; \text{tr}_\Pi (f) \;\text{tr}\left ( \tau; r_{\mu_{\Qpb}}\circ
\tilde{\s}_\ell ( \text{BC}(\Pi_p) )_{| E_\nu} |.|^{-\frac{\dim Sh}{2}}\right )
\end{eqnarray*}
où $m_\Pi$ est la multiplicité d'une $\Pi'\in \mathcal{T}(G)_\rho$
vérifiant $\Pi\simeq \Pi'_f$ et $BC$ désigne le changement de base
local  (qui existe par hypothèse en la place $v_j$). 

On dispose seulement d'une formule des traces pour la trace d'une
fonction de l'algèbre de Hecke (théorème \ref{Lefschetz_generique}),
mais pas pour $f\times \tau$.
C'est pourquoi nous avons besoin du lemme suivant :

\begin{lemm} Pour tout $\tau\in W_{E_\nu \mathcal{K}_p}$
il existe une fonction $$\tilde{f}_p^{\tau} \in \bigotimes_{i\in I} \left (e_{\l_i} *
\mathcal{H} ( \Gl_n ( F_{v_i}))* e_{\l_i}\right ) \otimes ( e_{\l_j} *
\mathcal{H} ( GU(n; F_{v_j} ) ) * e_{\l_j}) \subset \mathcal{H} ( G(\Qp))$$
telle que $\forall \chi$ caractère non ramifié de $G(\Qp)$ on ait 
$$
\text{tr} ( \tilde{f}_p^{\tau}; \pi_0\otimes \chi)= \text{tr} ( f_p;
\pi_0\otimes \chi) \; \text{tr} \left ( \tau; r_{\mu_{\Qpb}}\circ
\tilde{\s}_\ell (\pi_0\otimes \chi) |.|^{-\frac{\dim Sh}{2}}\right ) 
$$
\end{lemm}
\begin{proof}

Notons $X_{\Qlb} (G(\Qp))$ le groupe des caractères non ramifiés de
$G(\Qp)$ à valeurs dans $\Qlb^\times$. Le groupe $X_{\Qlb} (G(\Qp))$
s'identifie aux points à valeurs dans $\Qlb$ d'une variété algébrique
de type fini sur $\Qlb$ (plus précisément un tore). De plus, soit le
groupe fini 
$$
\GG=\{ \chi\in  X_{\Qlb} (G(\Qp))\; |\; \pi_0\otimes \chi \simeq \pi_0
\}
$$ 
Les fonctions régulières sur le tore quotient par $\GG$ s'identifient à
l'algèbre de Hecke du type $\l=\otimes_i \l_i \otimes \l_j$ de
$G(\Qp)$ (annexe \ref{Fait_type1}).  
 
Considérons maintenant la fonction 
\begin{eqnarray*} F:
X_{\Qlb} (G(\Qp)) & \ldrt & \Qlb \\
\chi & \longmapsto & \text{tr} ( f_p; \pi_0\otimes \chi) \; \text{tr}
( \tau;   r_{\mu_{\Qpb}}\circ
\tilde{\s}_\ell ( \pi_0\otimes\chi )_{| E_\nu} |.|^{-\frac{\dim Sh}{2}})
\end{eqnarray*}

Étant donné que $\tilde{\s}_\ell$ est défini au niveau des classes
d'isomorphisme de représentations, $F$ est invariante par $\GG$ et
descend donc en une fonction sur  $X_{\Qlb} (G(\Qp))/\GG$. Montrons
que cette fonction est régulière. 

Un tel $\chi$ élément de $X_{\Qlb} (G(\Qp))$ est de la forme 
$$
\chi=\otimes_{i\in I} (\chi_i \circ \det) \otimes \chi_j
$$
où $\forall i \; \chi_i$ est un caractère non ramifié de
$F_{v_i}^\times$ et $\chi_j$ de $GU( n; F_{v_j})$. Le caractère $\text{BC} (\chi_j)$ est
un caractère non ramifié de $\Gl_n (F_{v_j})\times \Qp^\times$· Notons
$\text{BC}(\chi_j)=( \chi'_j \circ \det ) \otimes \chi''_j$. L'application
$\chi_j \mapsto \text{BC} (\chi_j)$ est régulière (on peut la décrire
explicitement). 
On a alors, 
$$
( \pi_{0,i}\otimes \chi_i)(f_i)= \pi_{0,i} (f_i)\; \chi_i ( \det g_i)\;
\text{ et } \;( \pi_{0,j}\otimes \chi_j)(f_j)= \pi_{0,j} (f_j)\; \chi_j (  g_j)
$$
duquel on déduit 
$$
\text{tr} ( f_p; \pi_0 \otimes \chi) = \text{tr} ( f_p;
\pi_0). \prod_{i\in I} \chi_i (\det g_i) \;\chi_j (g_j)
$$
De plus, 
$$ \forall i\in I\;\;
\s_\ell ( \pi_{0,i}\otimes \chi\circ \det)=\s_l ( \pi_{0,i})\otimes
\s_\ell (\chi) 
$$
où $\s_\ell (\chi_i): W_{F_{v_i}}\ldrt \Qlb^\times$ est associé à
$\chi_i$ par la théorie du corps de classe local, 
et 
$$
\s_{\ell} ( \text{BC} ( \pi_{0,j}\otimes \chi_j)) = \s_\ell (
\text{BC}(\pi_{0,j}) \otimes ((\chi'_j\circ \det)\otimes \chi''_j))
=\s_\ell (\text{BC} (\pi_{0,j})) \otimes \s_\ell (\chi'_j) \otimes
\s_\ell (\chi''_j) 
$$
Utilisant le fait que pour tout $k$ dans $I\cup J$ l'extension $E_k |\Qp$ est non
ramifiée on en déduit : 
$$
r_{\mu_{\Qpb}}\circ \tilde{\s}_\ell ( \pi_0\otimes \chi) (\tau) =(
r_{\mu_{\Qpb}}\circ \tilde{\s}_\ell (\pi_0) (\tau)) . \prod_{i\in I}
\chi_i (p)^{-v(\tau) \sum_{\tau\in \Hom_{\Qp} (F_{v_i},\Qpb)}
p_{i}}  \chi'_j
(p)^{-v (\tau) p_j}    \chi''_j (p)^{-v (\tau)}
$$
où les $p_i,p_j$ sont des entiers associées à $\mu$. La fonction 
 $F$ est donc une fonction régulière. D'après l'annexe \ref{Fait_type1} il  existe donc $\tilde{f}_p^{\tau}\in
 \mathcal{H} (\l, G(\Qp))$ vérifiant 
$$
\forall \chi\in X_{\Qlb} ( G(\Qp))\;\;\; F( \chi)= \text{tr}_{\pi_0\otimes
\chi} ( \tilde{f}_p^{\tau})
$$
\end{proof}

\begin{rema}
On peut aussi utiliser le théorème de Paley-Wiener scalaire
(\cite{Paley_Wiener})  pour
obtenir l'existence d'une fonction $\tilde{f}_p^{\tau}\in \mathcal{H} ( G(
\Qp))$ vérifiant 
$$
\forall \pi \in \text{Irr} ( G(\Qp))\;\; \text{tr}_\pi ( \tilde{f}_p^{\tau})=
\text{tr}_\pi (f_p)\;  \text{tr}
( \tau;   r_{\mu_{\Qpb}}\circ
\tilde{\s}_\ell ( \pi )_{| E_\nu} |.|^{-\frac{\dim Sh}{2}})
$$
Cependant
nous aurons besoin par la suite de l'annulation des intégrales
orbitales de $\tilde{f}_p^{\tau}$ sur les éléments réguliers non
elliptiques. Pour la fonction $\tilde{f}_p^{\tau}$ du lemme ci dessus ce sera
une conséquence 
facile de l'appartenance de $\tilde{f}_p^{\tau}$ à l'algèbre de Hecke d'un
type supercuspidal. Si l'on veut utiliser le théorème de Paley-Wiener
scalaire et non le lemme ci dessus, cela est une conséquence du
théorème A-(b) de \cite{Kazhdan1}.  
\end{rema}

Posons $\tilde{f}^\tau=\tilde{f}_p^{\tau}\otimes f^p$ qui dépend donc de $\tau$. Il résulte du lemme que
l'on a l'égalité 
$$
\text{tr} ( \tilde{f}^\tau ; [H^\bullet ( \Sh,\LL_\rho)]) = 
\text{tr} ( f\times\tau ; [H^\bullet ( \Sh,\LL_\rho)]) 
$$

Supposons maintenant qu'en une place $w$ de $\Q$, $f_w$ soit à
support dans les éléments réguliers de $G(\Q_w)$ (on suppose $f$
décomposée en un produit $\otimes_{v} f_v$). Il résulte de
l'égalité précédente ainsi que du théorème \ref{Lefschetz_generique} que :
$$
\text{tr}\left ( f\times \tau ; [ H^\bullet (\Sh,\mathcal{L}_\rho )]\right)
= \sum_{\{\g \} \in \{ G(\Q) \} \atop {\g \text{ régulier} \atop \g
\text{ elliptique dans }G(\R)}} vol ( G(\Q)_\g \bc G (\A)_\g)
\text{tr} \rho (\g) \; \O_\g ( \tilde{f}_p^\tau) \; \O_\g ( f^p)
$$
Remarquons maintenant qu'étant donné que $\tilde{f}^\tau_p$ est dans l'algèbre
de Hecke associée à un type supercuspidal, pour tout élément $\g$ de
$G(\Qp)$ semi-simple régulier non elliptique $\O_\g (\tilde{f}_p)=0$.

D'où la formule 
$$
\text{tr}\left (f\times\tau; [ H^\bullet (\Sh,\mathcal{L}_\rho)]\right )
=
 \sum_{\{\g \} \in \{ G(\Q) \} \atop { \g \text{ régulier} \atop \g
\text{ elliptique dans }G(\R) \text{ et } G(\Qp) }} vol ( G(\Q)_\g \bc G (\A)_\g)
\text{tr} \rho (\g) \; \O_\g ( \tilde{f}_p^\tau) \; \O_\g ( f^p)
$$
\\

Passons maintenant à la formule des traces sur la fibre spéciale.

Soit donc $f_p$ comme précédemment fixée définitivement. Soit
$\tilde{K}\subset G( \A_f^p)$ 
un compact. Il résulte du théorème \ref{Lefschetz_speciale} appliqué à
toute la variété de Shimura qu'il existe un entier $N$ 
 tel que $\forall \tau\in W_{E_\nu\mathcal{K}_p}$, $ \; v(\tau)\geq N$, $ \;
\forall f^p\in \mathcal{H} ( G(\A_f^p))\; \text{supp} ( f^p)\subset
\tilde{K}$
$$
\text{tr}\left (f\times \tau; [H^\bullet (\Sh,\mathcal{L}_\rho)]\right
) = \sum_\phi \sum_{\{\g\}\in \{ I^\phi (\Q)\}} \l_{\phi,\{\g\},\tau} \;
\O_\g ( f^p) 
$$
où $\l_{\phi,\{\g\},\tau}=0$ sauf pour un nombre fini de $(\phi,\{\g\})$
(ce nombre fini tendant vers l'infini lorsque $\tilde{K}$
grossi et $v(\tau)$ varient). Nous choisirons de plus $N$ suffisemment grand de tel façon que le
théorème \ref{Lefschetz_speciale}  s'applique également à la strate basique.

On a donc l'égalité  
$$
\text{tr}\left (f\times \tau; [H^\bullet (\Sh^{an}(b_0),\mathcal{L}_\rho)]\right
) = \sum_{\phi, b(\phi)=b_0}  \sum_{\{\g\}\in \{ I^\phi (\Q)\}} \l_{\phi,\{\g\},\tau} \;
\O_\g ( f^p) 
$$ 
Remarquons que si $\phi$ n'est pas basique et $\g \in I^\phi (\Q)$ alors
$\g$ est dans une forme intérieure d'un Levi propre de $G(\Qp)$, $J_b$. 
\\

Nous pouvons maintenant comparer les formules des traces sur la fibre
générique et sur la fibre spéciale. 

Fixons un compact $\tilde{K}$ dans 
$G(\A_f^p)$.
Soit $v_1\neq p$ une place de $\Q$ décomposée dans
$\mathcal{K}$. Supposons le compact $\tilde{K}$ de la forme $\tilde{K}_{v_1}
\times \tilde{K}^{v_1}$. Soit $f_p$ comme précédemment et
$f^{p,v_1}\in \mathcal{H} ( G ( \A_f^{p,v_1}))$ telle que $\text{supp} (
f^{p,v_1})\subset \tilde{K}^{v_1}$, et en une place $f$ est à support
dans les éléments réguliers.
Soit $\tau\in W_{E_{\nu}\mathcal{K}_p}$ tel que $v(\tau)\geq N$ où $N$ est associé
à $\tilde{K}$. 
Pour une fonction $f_{v_1}$, on note  $\tilde{f}_p^\tau$ la fonction associée à $\tau$ et
$\tilde{f}^\tau= \tilde{f}_p \otimes f^p$. Désormais, $\tau$ et
$f^{v_1}$ sont fixés.  
 Alors, 
$$
\forall f_{v_1}\;\; \text{supp} (f_{v_1})\subset \tilde{K}_{v_1} \;\;\;\;
\sum_{\{\g \} \in \{ G(\Q) \} \atop { \g \text{ régulier} \atop \g
\text{ elliptique dans }G(\R) \text{ et } G(\Qp) }} \a_{\{\g \}}
 \O_\g ( f_{v_1})
= \sum_\phi \sum_{\{\g\}\in \{ I^\phi (\Q)\}} \b_{\phi,\{\g\}} \O_\g ( f_{v_1})
$$
où l'on a posé 
$$
\a_{\{\g\}}=  vol ( G(\Q)_\g \bc G (\A)_\g)
\text{tr} \rho (\g) \; \O_\g ( (\tilde{f}^\tau)^{v_1})
$$
et $b_{\phi,\{\g\}}=\l_{\{\g\},\tau} \; \O_\g ( f^{p,v_1})$ qui dépendent
de $f^{v_1}$ et $\tau$. 

Regroupons les termes par classes de conjugaison dans $G(\Q_{v_1})$:
on obtient 
$$
\sum_{\{\g'\}\in \{ G (\Q_{v_1})\}_{reg}} \left ( 
\sum_{\{\g \} \in \{ G(\Q) \} \atop { \g \text{ régulier} \atop { \g
\text{ elliptique dans }G(\R) \text{ et } G(\Qp)\atop {
\g \sim \g'} }}} \a_{\{\g\}} - \sum_{\phi,\{\g\}\in I^\phi (\Q) \atop
\g \sim \g'} \b_{\phi,\{ \g\}} \right ) \;\O_{\g'} ( f_{v_1})=0
$$
où, étant donné que $\text{supp} ( f_{v_1})$ reste dans le compact
$\tilde{K}_{v_1}$, la somme porte sur un nombre fini de classes de
conjugaisons semi-simples régulières dans $G(\Q_{v_1})$ intersectant
le compact $\tilde{K}_{v_1}$.

Étant donné que les classes de conjugaison sont fermées, disjointes et
en nombre fini on peut les séparer par des éléments $f_{v_1}\in
\mathcal{H} ( G (\Q_{v_1}))$ vérifiant $\text{supp} ( f_{v_1})\subset
\tilde{K}_{v_1}$. 

On obtient donc $\forall \g'\in G (\Q_{v_1})$ semi-simple régulier 
$$
\sum_{\{\g \} \in \{ G(\Q) \} \atop { \g \text{ régulier} \atop { \g
\text{ elliptique dans }G(\R) \text{ et } G(\Qp)\atop {
\g\sim \g'} }}} \a_{\{\g\}} =\sum_{\phi,\{\g\}\in I^\phi (\Q) \atop
\g \sim \g'} \b_{\phi,\{ \g\}}
$$
Remarquons maintenant que si $\phi$ est non basique et $\g_1\in I^\phi
(\Q)$ régulier, $\g_1$ ne peut être conjugué dans $G(\Q_{v_1})$ à un élément
$\g_2$ de $G(\Q)$ elliptique dans $G(\Q_p)$. En effet, 
$\g_1$ vu comme élément de $J_b$ se transfert en un élément $\g'_1$ de
$M(b)$ (le centralisateur du morphisme des pentes) la forme intérieure
quasi-déployée de $J_b$ (l'explicitation de l'existence de ce
transfert est faite à la fin du théorème III.\ref{lehteoint}).
 L'élément 
$\g_2$ serait
alors stablement conjugué à un élément $\g'_1$ d'un Levi propre de
$G(\Qp)$ ce qui d'après le lemme \ref{lemstab} 
est incompatible avec le fait que $\g_2$ est elliptique régulier.

On en déduit donc :
$$
\sum_{\phi\text{ non basique}\atop \{\g\}\in \{I^\phi (\Q)\}} \b_{\phi,\{
\g \}} =0
$$
On peut maintenant appliquer de nouveau le théorème
\ref{Lefschetz_speciale} à la strate basique pour obtenir que 
 $\forall f_{v_1},\text{ supp} (f_{v_1})\subset \tilde{K}_{v_1} $
$$
\text{tr} \left ( f\times \tau ; [H^\bullet
(\Sh,\mathcal{L}_\rho)]\right )
= \text{tr} \left ( f\times \tau ; [H^\bullet
(\Sh^{an}(b_0),\mathcal{L}_\rho)]\right )
$$
Nous pouvons mainteant faire varier $f^{v_1}$ et $\tau$ tels que
$\text{supp} (f^{v_1})\subset \tilde{K}^{v_1}$ et $v(\tau)\geq N$.  
On peut alors appliquer le théorème \ref{Caract_Traces} pour conclure.
\end{proof}

\begin{rema}
  Dans le théorème précédent, dans le second cas, lorsque $G$ est un
  groupe de similitudes unitaires en trois variables la restriction de
  $W_{E_\nu}$ à $W_{E_\nu \mathcal{K}_p}$ est inutile. C'est une
  conséquence du théorème \ref{rep_gal_U3}.
\end{rema}

%%% Local Variables: 
%%% mode: latex
%%% TeX-master: "thesef"
%%% End: 

\part{Application {\`a} la cohomologie des espaces de Rapoport-Zink de
type E.L. et  P.E.L.}
\chapter{}

Nous démontrons dans cette partie les principaux résultats de
cette thèse, à savoir les conjectures de Kottwitz 
(\cite{Rapoport2} conjecture 5.1, \cite{ECM} conjecture 5.3 et la
conjecture 5.4 qui est une généralisation de la conjecture de
Kottwitz) pour la partie supercuspidale de la cohomologie des espaces
de Rapoport-Zink basiques suivants :
\begin{itemize}
\item Les espaces de type E.L. non ramifiés pour lesquels le groupe
  $J_b$ est anisotrope modulo le centre ou bien égal à $G$ (théorèmes
  \ref{theo_princ1} et \ref{theo_princ2}) 
 (la  restriction sur $J_b$ 
  provient du fait que cette hypothèse simplifie la comparaison des
  formules des traces)
\item Pour les espaces de type P.E.L. non ramifiés  associés au groupe
  unitaire $U(3)$ sur des extensions de degré impaires de $\Qp$
  , $p\neq 2$, nous les démontrons
  pour les représentations 
  supercuspidales stables. Nous démontrons des résultats plus faibles
  pour les autres représentations supercuspidales. (théorème \ref{theo_princ3})  
\item De même pour $U(n)$ en supposant certains faits connus
  concernant l'analyse harmonique sur les groupes unitaires locaux et
  globaux en $n$ variables.
\end{itemize}
Avec
le cas des espaces de type P.E.L. non ramifiés associés au groupe
$U(3)$ nous donnons ainsi  le
premier exemple de loi de réciprocité locale non abélienne construite
géométriquement et associée à des groupes autres que des formes
intérieures du groupe linéaire.

\begin{rema}
  Si l'on sait démontrer la conjecture de Kottwitz pour les espaces de
  Rapoport-Zink associés à des données simples, on sait la démontrer
  pour des produits de tels espaces, c'est à dire des espaces
  associés à des données semi-simples. C'est une simple application de
  la formule de  Künneth. En particulier, les résultats annoncés  
  restent valables pour des produits des espaces non ramifiés de type
  E.L. et P.E.L. pour lesquels le résultat est démontré. 
\end{rema}

\section{Espaces de type E.L. non ramifiés}
\subsection{Construction de données globales à partir de données
  locales}
\subsubsection{Le résultat de Clozel}\label{Construction_unitaire}

Commençons par rappeler les résultats de la seconde section de
\cite{Clozel1} qui donnent une condition nécessaire et suffisante pour
qu'une famille de formes intérieures locales d'un groupe unitaire
presque-partout localement triviale provienne d'un groupe unitaire
global. Ces conditions sont déduites du formalisme développé par
Kottwitz pour étudier la cohomologie galoisienne des groupes réductifs
sur les corps de nombres.

Soit donc $F$ un corps C.M. de sous corps totalement réel maximal
$F^+$. Supposons nous donnée une famille $(U_v)_v$, où $v$ parcourt les
places de $F^+$, de groupes unitaires en $n$ variables sur $F^+_v$
(nous entendons par là que $U_v$ est une forme intérieure du groupe
$U(n)_{F^+_v}$ où $U(n)/F^+$ est le groupe unitaire quasi-dépolyé en $n$
variables associé à l'extension quadratique $F|F^+$). Soit $\Phi$ un
type C.M. de $F$. Lorsque $v$ parcourt les places infinies de $F^+$,
les groupes $U_v$ sont donc donnés par des couples d'entiers
$(p_\tau,q_\tau)_{\tau\in \Phi}$ tels que si $\tau$ induit $v$ alors
$$
U_v \simeq U(p_\tau,q_\tau)
$$
Si $v$ est une place de $F^+$ décomposée dans $F$, $U_v$ est du type 
$\GL_a(D)$ où $a|n$ et $D/F^+_v$ est une algèbre à division
d'invariant $\dpt{\frac{r}{n/a}}$ avec $r\wedge n/a=1$. Avec les
notations de (2.3) de \cite{Clozel1}, lorsque $n$ est pair, 
l'invariant local associé dans
$\mu_n^D\simeq \Z/n\Z$ est $ra \text{ mod } n$. Dans ce cas, 
la composante en $v$ de l'invariant global
associé est donc $ ra\text{ mod } 2$. Qui est non nul ssi $a$ est
impair (puisque $a$ impaire implique $r$ impair). 

Si $v$ est une place finie de $F^+$ inerte dans $F$, si $n$ est impaire, il
y a un unique groupe unitaire en $n$ variables sur $F^+_v$, le groupe
unitaire quasi-déployé. Le groupe $U_v$ est donc ce groupe.
 Si $n$ est pair, il y en a deux : la forme
quasidéployée et l'autre (suivant que le discriminant de la forme
hermitienne définissant le groupe est une norme de l'extension quadratique ou non). Dans ce cas
là, il y a donc deux possibilités pour $U_v$. 

Les résultats de la section 2 de \cite{Clozel1} s'énoncent ainsi : 

\begin{prop}\label{prop_Clozel1}
  Supposons que pour presque tout $v$, $U_v$ est quasi-déployé. 

Si $n$ est impair, il existe un groupe unitaire $U/F^+$ tel
que $\forall v\;\; U_{F^+_v}\simeq U_v$. 

Si $n$ est pair, notons $A$ le cardinal de l'ensemble des places
finies déployées  de $F^+$ telles que $U_v$ soit de la forme $\GL_a (D)$ avec $a$
impair, notons $B$ le cardinal des places finies inertes de $F^+$ telles que
$U_v$ soit non quasi-déployé. Il existe un groupe unitaire $U/F^+$ tel
que $\forall v\;\; U_{F^+_v}\simeq U_v$ ssi
$$
 \frac{n}{2} [F^+ : \Q] +\sum_{\tau\in \Phi} p_\tau \equiv A + B
 \text{ mod } 2
$$

Dans tous les cas, $U$ est le groupe unitaire associé à $(\mathcal{B},*)$ où
$\mathcal{B}/F$ est une algèbre simple et $*$ une involution de
seconde espèce.
\end{prop}

\subsubsection{Un lemme}

\begin{lemm}\label{tot_reel}
  Soit $F'|\Qp$ une extension de degré fini. Il existe un corps
  de nombres totalement réel $F^0$ tel que $p$ soit inerte dans $F^0$
  et $F^0_p\simeq F'$. 
  
\end{lemm}
\begin{proof}
  Écrivons pour cela $F'=\Q_p(\a)$ avec $P\in\Q_p[X]$ le polynôme
minimal de $\a$. On peut choisir $\a$ tel que $P$ soit unitaire à
coefficients dans $\Z_p$. Notons $d$ le degré de $P$.
 Alors, le lemme de Krasner montre que
pour $Q\in\Z_p[X]$ suffisamment proche de $P$, $Q$ est
irréductible et une racine de $Q$ engendre $F'$.
\\

Par densité de $\Z$
dans $\Z_p$ le lemme se ramène alors  à montrer que pour
tout polynôme $Q\in\Z[X]$ de degré $d$,
 pour tout voisinage $p$-adique de $Q$
dans les polynômes à coefficients entiers
de degré $d$ il existe un polynôme
dont toutes les racines soient réelles.
Le fait que dans dans l'extension construite à partir d'une racine
d'un tel polynôme, $p$ reste inerte est conséquence du fait que le
polynôme de $\Z[X]$ est irréductible dans $\Qp[X]$. 
\\

Soit   $E=\R[X]_d$ le $\R$-espace vectoriel des polynômes de degré
inférieure ou égal à $d$. Pour tout entier $N$, $p^N.\Z[X]_d$ est
un réseau de $E$. Considérons l'isomorphisme de $\R$-espaces vectoriels
\begin{eqnarray*}
\varphi:E &\ldrt & \R^{d+1} \\
R &\longmapsto & (R(0),\dots,R(d))
\end{eqnarray*}
Le groupe 
$\varphi(p^N.\Z[X]_d)$ est un réseau de $\R^{d+1}$. Et l'on veut montrer
que $\varphi(Q+p^N.\Z[X]_d)$ intersecte l'ensemble $\{x_0<0,x_1>0,x_2<0,\dots\}$.
Or on montre aisément que pour tout réseau $\La$ et tout $x\in\R^{d+1}$,
$x+\La$ intersecte un tel ensemble.
\end{proof}

\subsubsection{Construction d'une donnée de type P.E.L. globale à
  partir d'une donnée de Rapoport-Zink locale de type E.L. non ramifiée simple}

Soit $(F',V',b',\mu')$ une donnée de type E.L. non ramifiée simple
(I.\ref{EL_RZ}) sans sa classe de $\s$-conjugaison. 
 Nous notons $n=\dim_{F'} (V')$. 
Soit $F^+$ un corps de nombres totalement réel dans lequel $p$ est
inerte et tel que $F^+_p\simeq F'$ (lemme \ref{tot_reel}). Nous fixons
définitivement un  isomorphisme $F^+_p\iso F'$.  
Soit $\mathcal{K}$ un
corps quadratique imaginaire dans lequel $p$ est décomposé. Notons $F$
le corps C.M. $F^+ \mathcal{K}$. Notons $p=v_1 v_1^c$ la décomposition
de $p$ dans $F$. Fixons un plongement $\nu : \Qb \hookrightarrow \Qpb$
(où, rappelons le, $\Qb\subset \C$). Définissons le type C.M. $\Phi$
de $F$ de la façon suivante :
$$
\Phi=\{ \tau\in \Hom (F , \Qb)\;|\; \nu \circ \tau : F\hookrightarrow
\Qpb \text{ induit } v_1\;\}
$$
Rappelons que le cocaractère $\mu'$ est donné par des couples
d'entiers $(p'_\tau,q'_\tau)_{\tau\in \Hom_{\Qp} (F',\Qpb)}$. Pour 
  $\tau \in \Phi$, $\tau$ induit un plongement $\nu\circ \tau : F
  \hookrightarrow \Qpb$ qui via l'isomorphisme $F_{v_1}\simeq F'$
  (induit par l'isomorphisme $F^+_p\simeq F'$) induit un plongement 
 $\a (\tau)\in \Hom_{\Qp}(F',\Qpb)$. Nous noterons $(p_{\tau},q_\tau)=(p'_{\a
   (\tau)},q'_{\a(\tau}))$. 

Construisons un système $(U_v)_v$ comme dans le
\ref{Construction_unitaire}. Pour toute place $v$ de $F^+$ divisant l'infini,
il existe un unique $\tau\in \Phi$ tel que $\tau$ induise $v$. Posons
alors 
$$
U_v=U(p_\tau,q_\tau)
$$ 
Posons également 
$$
U_p=\Gl_n ( F')
$$
Choisissons une place finie $w$ de $F^+$ différente de $p$,
 décomposée dans $F$, et posons 
$U_w$ comme étant égal au groupe des unités d'une algèbre à division
sur $F^+_w$. Posons pour toute place finie $w'$ de $F^+$ décomposée
dans $F$, 
différente de $p$ et $w$, $U_{w'}=\GL_n (F^+_{w'})$. 
Si $n$ est impair, faisons un choix quelconque (presque partout
quasi-déployé) 
 de groupes unitaires $U_{w''}$ en les places $w''$ de $F^+$ inertes
 dans $F$. Si $n$ est pair, on peut toujours choisir le nombre de
 $U_{w"}$ non quasi-déployés de tel manière que l'équation de la
 proposition \ref{prop_Clozel1} soit satisfaite. 

A ce choix de groupes locaux,
la proposition \ref{prop_Clozel1} associe un couple
$(\mathcal{B},*)$ où $\mathcal{B}$ est une algèbre simple sur $F$ et
$*$ une involution de seconde espèce sur $\mathcal{B}$. 

Étant donné qu'en la place $w$,
$\mathcal{B}$ est une algèbre à division, $\mathcal{B}$ est une
algèbre à division. Posons $$B=\mathcal{B}^{opp}$$ qui est donc muni de
l'involution $*$. 
 Sur $\R$ l'involution $*$
devient conjuguée à une involution positive. D'après le lemme 2.8 de
\cite{Ko1}, l'ensemble des $b\in B^{*=1}\otimes \R$ tels que
le produit symétrique  $(x,y)_b=\text{tr}_{B_\R/\R} (x b y^*)$ 
sur $B_\R\times B_\R$ soit
positif est un cône ouvert non vide. Par un argument de densité, cet
ensemble 
rencontre donc l'ensemble $B^\times\cap B^{*=1}$. 
On en déduit qu'il existe  un
$b\in B^\times$ tel que $b^*=b$ et tel que l'involution $\#$ 
définie sur $B$ par 
$$ \forall x\in B\;\;\;
x^\#=b x^* b^{-1}
$$
est une involution positive de seconde espèce. Fixons un tel $b$ et
notons $\#$ l'involution 
associée.
 Soit $\beta\in F$ tel
que $\b^\#=-\b$. 
 Posons $V=B$ vu comme $B$-module à gauche. Considérons le produit
 $<.,.>$ défini sur $V\times V$ par 
$$
\forall x,y\in V\;\;\; <x,y>=\text{tr}_{B/\Q} ( x b^{-1} \beta y^\# )
$$
C'est un produit symplectique qui fait de $V$ un $(B,\#)$-module
hermitien. De plus, 
\begin{eqnarray*}
\forall a\in \mathcal{B}=\End_B (V)\;\forall x,y\in V\;\;
<a(x),y> & =& < xa, y> \\
 &=& \text{tr}_{B/\Q} ( x a b^{-1} \beta y^\# ) \\
 &= & \text{tr}_{B/\Q} ( x b^{-1} \beta (bab^{-1}) y^\# ) \\
 &=& \text{tr}_{B/\Q} (x b^{-1} \beta (y
 \underbrace{(bab^{-1})^\#}_{a^*} )^\# ) \\
 &=& <x, a^* (y)>
\end{eqnarray*}
On en déduit donc que $\mathcal{B}=\End_B (V)$ muni de l'adjonction
associée à $<.,.>$ est égal au couple $(\mathcal{B},*)$. Soit $G$ le
groupe des similitudes unitaires associé à $(B,\#,V,<.,.>)$ et $G_1$
le groupe unitaire. Le groupe $G_1$ est donc le groupe unitaire associé à
$(\mathcal{B},*)$. Étant donné qu'à l'infini les signatures
$(p_\tau,q_\tau)$ sont fixées, il existe une unique classe de
$G_1(\R)$-conjugaison $h$ 
 de morphisme d'algèbre à involution de $(\C,c)$ dans
$(\mathcal{B},*)$ telle que le produit $<.,h(i).>$ soit défini
positif. Fixons un tel $h:\C\ldrt \End_B (V)$. Et posons 
$$
\mathcal{D}=(B,\#,V,<.,.>,h)
$$
qui est une donnée de type P.E.L. au sens du chapitre I. 

D'après le choix fait pour $\Phi$ et les $(p_\tau,q_\tau)$ en fonction
des $(p'_\tau,q'_\tau)$ on en déduit la proposition suivante :

\begin{prop}\label{existence_donnee}
  Étant donnée une donnée locale de type E.L. non ramifiée simple sur
  une extension de $\Qp$, 
  il existe une donnée globale de type P.E.L. $\mathcal{D}
=(B,\#,V,<.,.>,h)$, un plongement $\nu:\Qb\hookrightarrow \Qpb$ tel
  que  via $\nu$, $\mathcal{D}$ induise (i.e. après équivalence de
  Morita, confère I.\ref{Dec_glo_loc}) la donnée locale. On peut de
  plus supposer que $\End_B(V)$ est une algèbre à division qui est en
  toutes les places finies soit déployée, soit une algèbre à
  division. On peut également supposer le corps C.M. $F$ de la forme
  $F^+\mathcal{K}$ où $\mathcal{K}$ est un corps quadratique
  imaginaire et où $p$ est inerte dans $\F^+$ et décomposé dans $\mathcal{K}$. 
\end{prop}

\subsection{Le théorème principal}

Si $H$ est un groupe réductif p-adique, si $V$ est un $H$-module
lisse, 
on peut définir le sous $H$-module $V_{cusp}$ de $V$, la partie
supercuspidale, en utilisant par exemple le centre de Bernstein :
$$
V_{cusp}= \bigoplus_{e} e.V
$$
où $e$ parcourt les idempotents du centre de Bernstein de $H$ associés
aux classes d'équivalences inertielles des représentations
supercuspidales de $H$. 

\begin{theo}\label{theo_princ1}
  Soit $(V,F,b,\mu)$ une donnée locale  de Rapoport-Zink de type
  E.L. non ramifiée simple basique. Soit $E\subset \Qpb$ le corps
  réflex associé.  
Supposons le groupe réductif $J_b$ anisotrope modulo
  son centre, c'est à dire
 $J_b=D^\times$ pour une algèbre à division
  $D$ sur $F$ d'invariant calculé en fonction de $\mu$. 
Notons $JL$ la correspondance de Jacquet-Langlands. 
Soit $\pi$ une
  représentation irréductible de $J_b$
telle que $JL(\pi)$ soit supercuspidale. Il y a une égalité dans le
  groupe de Grothendieck $\Groth ( G(\Qp)\times W_E)$ 
\begin{eqnarray*}
 \sum_i (-1)^i 
[\underset{K}{\limi} \Hom_{J_b} ( H^i_c ( \M_K\hat{\otimes}\Cp,\Qlb),
\pi)]_{cusp} 
 = 
[\text{JL} (\pi)] \otimes [r_\mu\circ \tilde{\s}_\ell (\text{JL}
(\pi))_{|E}]\; |.|^{\frac{-\sum_\tau {p_\tau q_\tau}}{2}} 
\end{eqnarray*}
où
$\tilde{\s}_\ell$ la correspondance de Langlands locale ``absolue''
$\ell$-adique pour le
groupe linéaire (\ref{Lang_Rev}) et $r_\mu$ la représentation du
L-groupe associé (\ref{le_r_mu}) sur $E$. 
\end{theo}

Rappelons avec les notations du I.\ref{Isoc_GL} que la classe de
conjugaison de  $\mu$ est donnée par
des entiers $(a_i)_{i\in \Z/d\Z}$, et que $J_b$  est
anisotrope  ssi 
$$
\left (\sum_{i\in \Z/d\Z} a_i \right )\wedge n =1
$$
Alors, $J_b= D^\times$ où l'invariant local de $D$ est $\frac{1}{n}
\sum_i a_i\in \Q/\Z$. 

Il peut également arriver, lorsque $n|\sum_i a_i$,  que $J_b$ soit égal
à $G(\Qp)$. Cela ne veut pas dire que les groupe p-divisibles associés
sont étales (ce qui est la cas ssi $\forall i\; a_i=0$),
car le groupe $J_b$ tient compte de l'action de $F$.
 Dans ce cas là, la correspondance de Jacquet-Langlands est
 l'identité et on a le théorème suivant : 

 \begin{theo} \label{theo_princ2}
    Soit $(V,F,b,\mu)$ une donnée locale  de Rapoport-Zink de type
  E.L. non ramifiée simple basique. Soit $E\subset \Qpb$ le corps
  réflex associé.  
Supposons le groupe réductif $J_b$ égal à $G=\GL_n$. Soit $\pi$ une
  représentation irréductible supercuspidale de $J_b$. Il y a une
  égalité dans $\Groth ( G(\Qp)\times W_E)$ 
$$
\sum_i (-1)^i 
[\underset{K}{\limi} \Hom_{J_b} ( H^i_c ( \M_K\hat{\otimes}\Cp,\Qlb),
\pi)]_{cusp} = [\pi]\otimes  [r_\mu\circ \tilde{\s}_\ell (
\pi)_{|E}]\; |.|^{-\sum_\tau \frac{p_\tau q_\tau}{2}}
$$
\end{theo}

Voici la démonstration des deux théorèmes ci dessus :

\begin{proof}
 Soit $\mathcal{D}$ une donnée de type P.E.L. globale induisant la
 donnée locale, comme dans la proposition \ref{existence_donnee}. Soit
 $\Sh$ la variété de Shimura associée.
 Le groupe $G(\Qp)$ est $\GL_n
 (F_v)\times \Qp^\times$ (où $F_v$ est le corps local noté $F$ dans
 les énoncés). Notons $\M(b,\mu)$ l'espace de Rapoport-Zink local
 associé à la donnée locale. L'espace de Rapoport-Zink associé à
 $\mathcal{D}$ en $p$ 
 et à la classe basique $b\in B(G_{\Qp},\mu_{\Qpb})$
 est donc $$\M(\mathcal{D}_{\Qp},b)=\M(b,\mu)\times \Qp^\times
 /\Zp^\times$$ 
D'après les résultats de l'appendice A,  dans le groupe de
 Grothendieck $\Groth (G(\A_f)\times W_{E_\nu})$, si 
$$
[H^\bullet (\Sh,\LL_\rho)] =\mid \ker^1 (\Q,G) \mid \; \sum_{\Pi\in \text{Irr} (G(\A_f))}
[\Pi]\otimes [R_{\ell,\rho,\mu}(\Pi)]
$$
la représentation semi-simple $ [R_{\ell,\rho,\mu}(\Pi')]$ est non nulle ssi
il existe $\Pi\in \mathcal{T}(G)_\rho$ telle que $\Pi'\simeq
\Pi_f$. De plus, pour une telle $\Pi\in \mathcal{T}(G)_\rho$,
contribuant à $[H^\bullet (\Sh,\LL_\rho)]_{cusp}$, i.e. telle que
$\Pi_p$ soit supercuspidale, 
$$
 [R_{\ell,\rho,\mu}(\Pi_f)]= m_\Pi [r_\mu\circ \tilde{\s}_\ell
 (\Pi_p)] \; |.|^{-\frac{\dim Sh}{2}}
$$
où $m_\Pi$ désigne la multiplicité de $\Pi$ dans l'espace des
formes automorphes sur $G$. On obtient donc l'égalité suivante 
$$
[H^\bullet (\Sh,\LL_\rho)]_{cusp}=\mid \ker^1 (\Q,G)\mid \sum_{\Pi\in
  \mathcal{T}(G)_\rho \atop { \Pi_p \text{ supercuspidale}}} [\Pi_f]
\otimes [r_\mu\circ \tilde{\s}_\ell
 (\Pi_p)] \; |.|^{-\frac{\dim Sh}{2}}
$$
 Soit $\phi$ une classe d'isogénie
intervenant dans la strate basique et $I^\phi$ le groupe réductif
associé.
 Rappelons que $I^\phi (\R)$ est la forme compacte modulo le
centre de $G(\R)$, que $I^\phi(\Qp)=J_b\times \Qp^\times$ et que
$I^\phi (\A_f^p)=G(\A_f^p)$.

Le corollaire \ref{coro_Matsushima} couplé au théorème
\ref{coho_ind_en_p} montre que l'on a l'égalité suivante : 
$$
[H^\bullet (\Sh,\LL_\rho)]_{cusp}=\mid \ker^1 (\Q,G)\mid 
\sum_{\Pi\in \mathcal{T}(I^\phi) \atop {\Pi_\infty=\check{\rho} \atop { 
 \Pi_p \text{ supercuspidale}}}}
 [\underset{K}{\limi} \text{Hom}_{J_b\times\Qp^\times} ( H^\bullet_c
(\M_K (\mathcal{D}_{\Qp},b),\Qlb),\Pi_p)]\otimes [\Pi^p] 
$$
Si $\Pi_p=\pi_1\otimes \chi$ où $\pi_1$ est une représentation
irréductible de $J_b$ et $\chi$ un caractère de $\Qp^\times$, il y a
un isomorphisme
$$
\underset{K}{\limi} \text{Hom}_{J_b\times\Qp^\times} ( H^\bullet_c
(\M_K (\mathcal{D}_{\Qp},b),\Qlb),\Pi_p)
=\underset{K}{\limi} \Hom_{J_b} ( H^\bullet_c
(\M_K (b,\mu),\Qlb),\pi_1)\otimes \chi 
$$
comme représentation de $G(\Qp)=\GL_n (F)\times \Qp^\times$. 

Fixons $\rho=1$.
Il existe une représentation automorphe $\Pi$ de $I^\phi$ telle
que $\Pi_p=\pi\otimes 1$ et $\Pi_\infty=1$ (\cite{Clozel2}). 
 Utilisons le théorème \ref{comparaison_formes_int} appliqué aux
   groupes formes intérieurs $I^\phi$ et $G$. On obtient alors le résultat en égalant
   les deux expressions ci dessus pour $[H^\bullet 
   (\Sh,\LL_\rho)]_{cusp}$ et en prenant la partie
   $\Pi_f^p$-isotypique (après avoir divisé par $\mid \ker^1 (\Q,G)\mid
   m_\Pi$, où $m_\Pi$ désigne la multiplicité de $\Pi$ dans l'espace
   des formes automorphes sur $I^\phi$  qui est égale
  à celle de $\pi_\infty\otimes JL(\Pi_p)\otimes \Pi^p_f$ dans l'espace des formes
  automorphes sur $G$ pour une représentation de la série discrète
  $\pi_\infty$ ayant de la cohomologie dans $\rho$, 
d'après \ref{comparaison_formes_int}). 
\end{proof}

\section{espaces de type P.E.L.}
\subsection{Construction de données globales à partir de données
  locales} 

Soit $F'|\Qp$ une extension non ramifiée de degré pair. On veut
trouver un corps de nombres totalement réel $F^+$ et un corps
quadratique imaginaire $\mathcal{K}|\Q$ dans lesquels $p$ est inerte, 
tels que $F^+_p$ et $\mathcal{K}_p$ soient linéairement disjointes
(c'est à dire $p$ reste inerte dans $F^+ \mathcal{K}$) et tels que 
$F'=F^+_p \mathcal{K}_p$. Cela est
possible ssi $[F':\Qp]$ n'est pas un multiple de $4$ (utiliser le
lemme \ref{tot_reel}). De la même manière que pour la proposition
\ref{existence_donnee} on montre alors la proposition suivante :

  \begin{prop}\label{donnee_ccc}
    Étant donnée une donnée locale de type P.E.L. non ramifiée simple
    sur une extension de $\Qp$ de degré qui n'est pas un multiple de
    $4$, il existe une donnée globale de type
    P.E.L. $\mathcal{D}=(B,\#,V,<.,.>,h)$,  un
    plongement $\nu:\Qb\hookrightarrow \Qpb$ tel que via $\nu$,
    $\mathcal{D}$ induise la donnée locale. On peut de plus supposer
    que $\End_B(V)$ est une algèbre à division qui est en toutes les
    places finies soit déployée, soit une algèbre à division. On peut
    également supposer le corps C.M. $F$ de la forme $F^+\mathcal{K}$
    où $\mathcal{K}$ est un corps quadratique imaginaire. 
  \end{prop}

\subsection{Le théorème principal pour $GU(3)$}

\begin{theo}\label{theo_princ3}
  Soit $(F,*,V,<.,.>,b,\mu)$ une donnée locale de type P.E.L. non
  ramifiée simple basique. Supposons $p\neq 2$, $[F:\Qp]/2$
  impair et $\dim_F (V)=3$. Soit $E\subset \Qpb$ le corps réflex associé.
\begin{enumerate}
\item 
 Soit
  $\psi:W_{\Qp}\ldrt \L G$ une classe de conjugaison de L-paramètre
  vérifiant : $\im (\psi )$ n'est pas contenue dans $\L B$ pour un
 sous groupe de  Borel $B$ de $G$. Soit $\Pi (\psi )$ le L-paquet
  supercuspidal de représentations de $G(\Qp)$ associé
  (\ref{paquets_sup}). Rappelons que $J_b= G(\Qp)=GU(3)$. Il y a une
  égalité dans $\Groth (G(\Qp)\times W_{E})$ : 
$$ 
\sum_{\pi\in \Pi (\psi)}\left [\underset{K}{\limi} \Hom_{J_b} \left ( H^\bullet_c
( \M_K\hat{\otimes}\Cp,\Qlb),\pi\right ) \right ]_{cusp} = \sum_{\pi\in \Pi
  (\psi)} [\pi]\otimes [ r_\mu \circ \psi_{|W_E}]\;
|.|^{-\frac{\sum_\tau p_\tau q_\tau}{4}}
$$ 
En particulier, si $\Pi (\psi)$ est stable la conjecture de Kottwitz
est vérifiée. 
\item
Soit $\xi$ un caractère du groupe endoscopique $H$ de $U(3)$
(appendice \ref{Resultats_Rog}) et $St_H (\xi)$ la représentation de Steinberg
associée. Soit $\Pi ( St_H (\xi)) = \{ \pi^2 (\xi), \pi^s (\xi) \}$ le
L-paquet transfert endoscopique à $U(3)$ (\ref{Class_sup}). Soit
$\Pi=\{\pi^2,\pi^s\}$ un L-paquet de représentations de $G(\Qp)$ tel
  que $\pi^2_{|U(3)}=\pi^2 (\xi),\pi^s_{|U(3)}=\pi^s (\xi)$. A $\Pi$
  est associé un L-paramètre $\psi:W_{\Qp}\times \text{SL}_2 (\C)
  \ldrt \L G$. Il y a une égalité dans $\Groth (G(\Qp)\times W_{E})$
  : 
  \begin{eqnarray*}
  &\,&  \left [\underset{K}{\limi} \Hom_{J_b} \left ( H^\bullet_c
    (\M_{K}\hat{\otimes} \Cp ,\Qlb),\pi^s\right )\right ]_{cusp} + \left[
    \underset{K}{\limi} \text{Ext}^{\, \bullet}_{J_b-\text{lisse}} (
    H^\bullet_c (\M_{K}\hat{\otimes} \Cp ,\Qlb),\pi^2)\right ]_{cusp}
    \\
 &=& [\pi^s]\otimes [r_\mu\circ \psi_{|W_E}] \;
    |.|^{-\frac{\sum_\tau p_\tau q_\tau }{4}}
  \end{eqnarray*}
\item 
Soit $\chi$ un caractère de $G(\Qp)$ et $St_G(\chi)$ la représentation
de Steinberg associée. Il y a une égalité dans $\Groth ( G(\Q)\times
W_{E})$ : 
$$
\left [ \text{Ext}^\bullet_{J_b-\text{lisse}} \left (  H^\bullet_c
( \M_K\hat{\otimes}\Cp,\Qlb), St_G (\chi)\right )\right ]_{cusp} = 0 
$$
\end{enumerate}
\end{theo}

\begin{proof}
  Soit $\mathcal{D}$ une donnée globale de type P.E.L. induisant la
  donnée locale comme dans la proposition \ref{donnee_ccc}. Soit $\Sh$
  la variété de Shimura associée. Il y a un isomorphisme d'espaces
  analytiques $\M (\mathcal{D}_{\Qp},b)\simeq \M (b,\mu)$. Soit $\phi$
  une classe d'isogénie intervenant dans la classe basique. 

Dans les trois cas de l'énoncé, pour chacun des L-paquets locaux en
$p$, $\Pi$ (dans le troisième cas $\{ St_G (\chi) \}$ est un L-paquet,
confère page 174 de \cite{Rogawski1}),  
 on peut trouver une représentation irréductible de $G(\A_f^p)$, $\Pi^p$,
telle que 
$$ \forall \pi_p\in \Pi\;\;
\text{Triv}_\infty \otimes \pi_p\otimes \Pi^p \in \mathcal{T}
(I^\phi)_{Triv_\infty} 
$$
En effet, il suffit de choisir une représentation $\pi_p\in \Pi$, de la
globaliser (\cite{Clozel2}) (dans le troisième cas cela est une
conséquence du fait que $St_G (\chi)$ est de carré intégrable). Alors,
les L-paquets globaux de représentations de $G$ étant stables (c'est
une conséquence du fait que $\End_B(V)$ est une algèbre à division,
\cite{Rogawski1} théorème 14.6.3) tous les éléments d'un même L-paquet
apparaissent avec la même multiplicité dans le spectre discret de
$G$ (multiplicité qui est $1$). 

Utilisons le théorème \ref{rep_gal_U3}
Le résultat s'en déduit comme dans la démonstration du théorème
\ref{theo_princ1} en prenant la partie $\Pi^p$ isotypique dans les
deux expressions différentes pour $[H^\bullet (\Sh,\Qlb)]_{cusp}$  en
termes  de représentations automorphes de $G$ et $I^\phi$, et en
utilisant le fait que si deux représentations automorphes de $G$,
resp. $I^\phi$, sont isomorphes hors $p$ elle sont dans le même
L-paquet global et que donc leur composante en $p$ sont dans le même
L-paquet local (\cite{Rogawski1})(et en utilisant également le fait
énoncé ci dessus que 
la multiplicité est constante pour tous les éléments d'un même L-paquet
global). 
\end{proof}

\subsection{Spéculations}

Le théorème précédent ne donne une démonstration de la conjecture de
Kottwitz que dans le cas des représentations supercuspidales $\pi$
telles que $\{ \pi \}$ soit un L-paquet, c'est à dire telles que le
changement de base de $\pi$ soit une représentation supercuspidale de 
$\GL_3 (F)\times K^\times$. La raison de ``l'échec'' de la méthode dans le
cas  des autres représentations supercuspidales de $G(\Q)$ est que
quelque soit la globalisation d'une représentation locale en une représentation
globale de $G$, le L-paquet global obtenu (il faudrait plutôt parler
de A-paquet) 
 est stable et donc les
multiplicités de chaque élément sont égales. On obtient donc dans tous
les cas dans les membres de gauche du théorème \ref{theo_princ3} une
somme sur tous les éléments du L-paquet local qui ne permet pas de
séparer chaque terme. 

Voila une stratégie permettant d'attaquer la conjecture de Kottwitz
dans les autres cas : il faudrait travailler avec des groupes unitaires
$G$ qui ne sont plus associés à une algèbre à division mais à $M_3
(F)$. En effet, dans ce cas là le spectre discret de $G$ est beaucoup plus
riche que dans le cas précédent. Par exemple, avec les notations de
\ref{paquets_sup} si $\rho$ est une représentation supercuspidale
stable de $H(\Qp)$, si $\Pi (\rho)$ désigne le L-paquet associé de
$U(3)$ celui ci est de cardinal $2$. On peut globaliser $\rho$ en
une représentation de l'analogue global du groupe endoscopique local
$H$ (\cite{Rogawski1}) dont le changement de base à $GL_2$ est cuspidale, puis considérer le transfert endoscopique global
de cette représentation à $U(3)$.  Notons $\Pi$ le L-paquet global
obtenu. Avec les notations de \cite{Rogawski1}, $\Pi\in \Pi_e (G)$.
 En $p$, le L-paquet local associé à $\Pi$ est le L-paquet
local $\Pi (\rho)$. De plus, les multiplicités de chaque élément de ce
L-paquet global dans le spectre discret de $G$ sont calculées dans
\cite{Rogawski1} (confère également \cite{Rogawski2}). Ces
multiplicités ne sont pas constantes. On peut ainsi espérer obtenir
suffisamment de relations pour isoler chaque facteur dans les sommes des
membres de gauche du théorème \ref{theo_princ3}.
Cependant, pour l'algèbre à division $M_3 (F)$ les variétés de Shimura
associées ne sont pas en général compactes.  
Pour avoir des variétés de Shimura compactes, il faut imposer 
que le groupe unitaire à l'infini possède un 
facteur simple compacte (c'est à dire $\exists \tau\; p_\tau=0$ ou
$q_\tau=0$). Cela implique que la variété de Shimura est propre sur la
fibre générique. Il resterait à vérifier que c'est également le cas
pour les modèles entiers (un exercice en théorie des modèles de
Néron).

%%% Local Variables: 
%%% mode: latex
%%% TeX-master: "thesef"
%%% End: 

\appendix

\chapter[Représentations automorphes des
groupes
unitaires]
{Résultats concernant les
représentations automorphes des groupes unitaires, leur cohomologie et
les représentations galoisiennes associées}

\section{Hypothèses générales}

Soit $\mathcal{D}=(B,*,V,<.,.>,h)$ une donnée de type P.E.L. globale
comme dans le premier chapitre. 

 Nous supposerons dans cet appendice que 
 que $\End_B(V)$ est une algèbre
 à division $D$ sur le corps C.M. $F$. Les variétés de Shimura
 associées sont donc propres. Nous supposerons de plus que le
 corps $F$ est de la forme $F^+ \mathcal{K}$ où $\mathcal{K}|\Q$ est
 une extension quadratique imaginaire. 
 Nous ferons
 également l'hypothèse qu'en toute place $v$ de $F$, l'algèbre simple
 $D_v$ est soit 
 déployée soit une algèbre à division. 

 Soit
$\Phi$ un type C.M. de $F$.
 Rappelons que  
$$
G_1(\R)\simeq \prod_{\tau\in \Phi} U(p_\tau,q_\tau)
$$
où $p_\tau+q_\tau=n$ et que 
$$
G(\C)\simeq \prod_{\Phi} GL_n (\C) \times \C^\times 
$$
Rappelons que $\rho$ est une représentation algébrique irréductible de
dimension finie de $G_{\C}$ et que nous notons également $\rho$ sa
restriction à $G(\R)$. 

Si $v$ est une place finie de $F$ divisant un premier $p$ de $\Q$
décomposé dans $\mathcal{K}$, le groupe $\GL_n (F_v)$ est un facteur
de $G(\Qp)$. Si $\Pi$ est une représentation automorphe de $G$ on peut
donc définir $\Pi_v$ comme représentation irréductible de $\GL_n (F_v)$.

Le groupe $G$ étant anisotrope modulo son centre, l'espace des formes
automorphes sur 
$G$ de caractère central fixé
est décomposé discrètement; il n'y a pas de spectre continu. Si $\Pi$
est une représentation automorphe 
de $G$ on peut donc définir sa multiplicité comme étant la
multiplicité de la classe d'isomorphisme de $\Pi$ dans le
$(\mathfrak{g}_\infty,K_\infty)\times G(\A_f)$-module admissible des
formes automorphes se transformant selon le caractère central de
$\Pi$. Nous noterons $m_\Pi$ cette multiplicité.

\section[Résultats sur les représentations galoisiennes]
{Résultats de Clozel, Harris et Labesse sur la pureté de la
cohomologie des représentations automorphes} 

Certains des résultats énoncés ci dessous sont démontrés par Clozel
dans \cite{Clozel1} (on pourra également consulter le corollaire VII.2.7
de \cite{Har4}) . Ils sont démontrés dans le cas de  signature
quelconque par Harris et Labesse dans des notes non publiées. 

\begin{defi}
Nous noterons $\mathcal{T}(G)_{\rho}$ l'ensemble des représentations
automorphes $\Pi$ de $G$ vérifiant que $\Pi_\infty$ a de la
cohomologie dans $\rho$. 
\end{defi}

Bien que son énoncé ne le fasse pas apparaître, la démonstration du
résultat qui suit utilise la cohomologie des variétés de Shimura de
type P.E.L., et notamment les
résultats de Kottwitz sur la cohomologie en une place de bonne
réduction (\cite{Ko6}), 
 couplés aux estimations connues sur les valeurs propres des
 composantes locales des représentations automorphes unitaires
 cuspidales de 
$\GL_n$ (c'est à dire les estimations connues concernant la conjecture
de Ramanujan). 
Plus précisément, si $\Pi$ est une représentation automorphe de $G$
intervenant dans la cohomologie de la variété de Shimura associée, 
le lien démontré par Kottwitz entre les valeurs propres
de Frobenius en une place 
de bonne réduction (qui vérifient certaines conditions grâce au
théorème de pureté de Deligne) 
 et les paramètres de Hecke associés permettent de
démontrer que le motif (tout du moins sa réalisation $\ell$-adique, De
Rham...)
découpé par $\Pi_f$ dans la variété de Shimura est pure de poids
$\dim \Sh$.

\begin{theo}\label{purete}
Soit $\Pi\in \mathcal{T}(G)_{\rho}$ telle qu'il existe
un premier $p$ de $\Q$ décomposé dans $\mathcal{K}$ et 
 une place finie
$v$ de $F$ divisant $p$  
vérifiant :  $\Pi_v$ soit dans la série discrète. Alors, 
$$
\dim H^i (\mathfrak{g}_\infty,K_\infty ; \Pi_\infty\otimes \rho) =
\left \{ { 1 \text{  si } i=\sum_\tau p_\tau q_\tau =\dim \Sh \atop 
 0 \text{  sinon} } \right.
$$
De plus, $\Pi_\infty$ appartient au L-paquet de la série discrète de
$G(\R)$ formé des représentations de la série discrète de $G(\R)$
ayant de la cohomologie dans $\rho$. Ce L-paquet est exactement la
fibre du changement de base vers $G(\C)$ en l'unique représentation
tempérée de $G(\C)$ ayant de la cohomologie dans $\rho$ . Son cardinal
est égal à $\dim r_\mu$ où $r_\mu$ est la représentation du L-groupe
de $G$ associée à $\mu$ (\cite{Langlands1})(confère \ref{mutau}). 
\end{theo}

\begin{theo}\label{purete_supercusp}
Sous les hypothèses du théorème précédent, supposons de plus que
$\Pi_v$ soit supercuspidale. Alors, pour toute représentation
$\pi'_\infty$ de $G(\R)$ dans le L-paquet de la série discrète ci
dessus $\Pi'=\pi'_\infty\otimes \Pi_f\in \mathcal{T}(G)_{\rho}$. De
plus, les multiplicités $m_\Pi$, $m_{\Pi'}$ de $\Pi$ et $\Pi'$ sont égales.
\end{theo}

\section[Changement de base quadratique stable]
{Résultats de Clozel, Harris et Labesse sur le changement de
base quadratique stable}\label{CHL}

Nous énonçons ici un théorème de Harris et Labesse démontré dans des
notes non publiées. Le cas  de signature $(1,n-1)$ est
démontré par Clozel dans \cite{Clozel1} (dans le  VII.2 de
\cite{Har4} le résultat de Clozel est expliqué en détails, il y est
également expliqué 
comment passer du cas des groupes unitaires au cas des groupes de
similitudes unitaires). Les morphismes de L-groupes associés aux
fonctorialités de Langlands que démontrent ce théorème sont exposés
dans la démonstration du théorème \ref{gal_coho_sh}.
\\

 Soit  $p$  un  premier de $\Q$ 
 décomposé dans $\mathcal{K}$, $p=w_0 w_0^c$. Supposons l'algèbre à
 division $D=\End_B (V)$ déployée en $p$. 
Si $(v_i)_{i\in I}$ désigne un ensemble de représentants des classes
de places $v$ divisant $p$ modulo l'action de $c$, 
$$
G(\Qp)=\prod_{i\in I} \Gl_n (F_{v_i}) \times \Qp^\times 
$$
et si $\pi=(\otimes_i \pi_i)\otimes \chi$ est une représentation
irréductible de $G(\Qp)$ on peut définir le changement de base 
quadratique de $\pi$ comme étant 
$$
BC(\pi) =( \otimes_{i\in I}  \pi_i\otimes \check{\pi}_i^c )\otimes
 (\chi\otimes \chi^{-1}) 
$$ 
Lorsque $D$ n'est plus déployée en $p$ on peut encore définir un
changement de base quadratique stable en composant le changement de
base qui vient d'être décrit avec la correspondance de
Jacquet-Langlands locale. 

Si maintenant  $p$ un premier de $\Q$ en lequel $G_{\Qp}$ est non ramifié, si
$\pi$ est une représentation sphérique de $G(\Qp)$ on peut alors
définir le changement de base non ramifié de $\pi$ comme
représentation de 
$$
\prod_{v|p} \Gl_n( F_v) \times \prod_{w|p} \mathcal{K}_w^\times 
$$
où $v$ parcourt des places de $F$ et $w$ de $\mathcal{K}$.
Dans le cas où $p$ est décomposé dans $\mathcal{K}$ on retrouve le cas
précédent restreint aux représentations sphériques. Dans le cas où $p$
n'est pas décomposé dans $\mathcal{K}$, si $BC(\pi)=(\otimes_v \pi'_v)
\otimes \chi_w$, on a les relations 
$$
\forall v \;\;\pi_v'\simeq \check{\pi}'_{v^c}  
$$

Le théorème qui suit énonce l'existence d'un changement de base
quadratique stable pour des représentation automorphes de $G$ vers celles
d'un groupe linéaire. Le mot stable signifie que l'on a
composé un changement de base quadratique avec une correspondance de
Jacquet-Langlands.

\begin{theo}\label{BCCC}
Soit $\Pi\in\mathcal{T}(G)_{\rho}$ telle qu'il existe un premier $p$
décomposé dans $\mathcal{K}$ et une place $v$ de $F^+$ divisant $p$
telle que $\Pi_v$ soit supercuspidale. Il existe alors une 
représentation automorphe cuspidale $\Pi'$ de $\Gl_n(\A_F)$ ainsi
qu'un caractère de Hecke $\chi$ de $\A_{\mathcal{K}}^\times$ 
 vérifiant 
\begin{itemize}
\item $\chi=\chi_{\Pi|\A_{\mathcal{K}}^\times}^c$
\item Si $\Pi_p$ est non
ramifiée alors $(\Pi'\otimes\chi)_p$ est le changement de base non ramifié
de $\Pi_p$
\item $\Pi'_\infty$ est de la forme $$\prod_{\tau\in \Phi}
\pi'_{\tau}$$
où $\prod_\tau \pi'_\tau$ est l'unique représentation tempérée de
$G_1(\C)$ ayant de la cohomologie dans $\rho_{|G_1(\C)}$. De plus,  $$
\chi_\infty=\rho_{|\C^\times}^{-1}
$$
où $\C^\times\subset G(\C)$ désigne l'inclusion du facteur de similitude.
\item $\Check{\Pi}' \simeq {\Pi'}^c$ 
\item Si $p$ est décomposé dans
$\mathcal{K}$ alors 
$$
(\Pi'\otimes \chi)_p= BC (\Pi_p)
$$
comme représentation de $(\text{Res}_{F/\Q} \Gl_n \times
\text{Res}_{\mathcal{K}|\Q} \Gm)(\Qp)$.
\item $\chi_{\Pi'|\A_{\mathcal{K}}^\times}=\chi/\chi^c$
\end{itemize}
\end{theo}

Par multiplicité un forte, une telle représentation automorphe $\Pi'$ est unique.

\section{Résultats de Harris et Labesse sur la comparaison entre la
formule des traces pour deux groupes unitaires formes intérieures} 

 Soit $v$ une place finie de $F$ vérifiant $v\neq v^c$. 
Soient $G_1$ et $G_2$ deux groupes unitaires comme ci dessus, formes
intérieurs l'un de l'autre tels que $G_1(\A_f^v)=G_2 (\A_f^v)$. 
 Nous notons $JL$ le transfert local des
représentations supercuspidales de $G_{1v}$ et $G_{2 v}$ vers les
représentations de leur forme intérieur quasidéployée : $\Gl_n
(F_v)$. 

\begin{theo}\label{comparaison_formes_int}
 Soit $\Pi_1\in
\mathcal{T}(G_1)_\rho$ telle que $JL(\Pi_{1v})$ soit supercuspidale. Alors,
pour toute représentation $\pi_\infty$ de $G_2 (\R)$ faisant partie du
L-paquet de la série discrète ayant de la cohomologie dans $\rho$, il
existe une $\Pi_2\in \mathcal{T}(G_2)_\rho$ telle que 
\begin{itemize}
\item $\Pi_{2\infty}=\pi_\infty$
\item $\Pi_2^{v}\simeq \Pi_1^v$
\end{itemize}
De plus, pour une telle $\Pi_2$
 $$JL(\Pi_{2v})=JL(\Pi_{1v})
$$
et les multiplicités de $\Pi_1$ et $\Pi_2$ sont égales.
\end{theo} 

On déduit en particulier de ce théorème, du théorème \ref{purete} et
de la formule de Matushima que si $\Pi_1$ est comme dans
l'énoncé, la multiplicité de $\Pi_{1,f}$ dans la cohomologie de la
variété de Shimura associée est 
$$
\text{multiplicité} (\Pi_1) \times \dim r_\mu
$$
Nous verrons en fait dans le théorème \ref{gal_coho_sh} que le facteur
$\dim r_\mu$ est la dimension d'une représentation galoisienne telle que la
représentation galoisienne décomposée par $\Pi_{1,f}$ dans la
cohomologie 
 est $\text{multiplicité} (\Pi_1)$-fois cette
représentation. Dans le cas de la signature $(1,n-1)$ ce théorème de
divisibilité est démontré par Taylor en utilisant la théorie de
Hodge-Tate p-adique (à l'origine dans une lettre de Taylor à Clozel,
cet argument est 
esquissé dans \cite{Har1} page 102; pour plus de détails, confère
\cite{Har4} et notamment le lemme II.2.2 et la proposition VIII.1.8). Dans le cas général, nous
déduirons ce résultat du cas de signature $(1,n-1)$ (théorème \ref{gal_coho_sh}).

\section{Résultats de Clozel, Harris, Kottwitz et Taylor sur les représentations
galoisiennes}

Voici une partie de l'énoncé du théorème VIII.1.9 de \cite{Har4} qui
généralise le théorème principal de \cite{Clozel1}. 

\begin{theo}\label{theogalo}
Soit $\Pi$ une représentation automorphe cuspidale de $\Gl_n(\A_F)$
satisfaisant les hypothèses suivantes :
\begin{itemize}
\item $\check{\Pi}\simeq \Pi^c$
\item $\Pi_\infty$ a même caractère infinitésimal que
$\check{\rho}_{|G_1(\C)}$
\end{itemize}
Il existe alors une représentation continue
$$
R_\ell (\Pi):\Gal(\bar{F}|F)\ldrt \GL_n (\Qlb)
$$
vérifiant 
\begin{itemize}
\item En toute place $v$ de $F$ ne divisant pas $\ell$ la restriction de
$R_\ell (\Pi)$ à $W_{F_v}$ est $\s_\ell( \Pi^\vee_v)\otimes
|.|^{\frac{1-n}{2}}$ où $\s_\ell$ désigne la 
correspondance de Langlands locale pour $\GL_n(F_v)$. 
\item En toute place $v$ où $R_\ell (\Pi)$ est non ramifiée, les
valeurs propres de Frobenius sont algébriques pures de poids $n-1$
\end{itemize}
\end{theo}

\section{Résultat de Kottwitz sur les représentations galoisiennes
  obtenues dans la cohomologie des variétés de Shimura en une place de
  bonne réduction}\label{theo_kott_sh}

Dans l'article \cite{Ko6} Kottwitz démontre qu'en presque toutes les
places de bonne réduction, la correspondance de Langlands locale non
ramifiée est réalisée dans la cohomologie des variétés de Shimura que
nous considérons. Ce résultat est le corollaire des nombreux travaux
de Kottwitz sur la stabilisation de la formule des traces ainsi que de ses
travaux sur le comptage des points des variétés de Shimura sur les
corps finis 
(\cite{Ko1}). 
 Nous renvoyons le lecteur au théorème 1 de \cite{Ko6}
 (et au nombreux travaux antérieurs de Kottwitz) pour 
l'énoncé. Le théorème \ref{gal_coho_sh} est une généralisation 
de ce théorème à  d'autres places. 

\section[Application au cas P.E.L. de signature quelconque]
{Application aux représentations galoisiennes obtenues dans 
la cohomologie des variétés de Shimura de type P.E.L. de signature quelconque}

Nous montrons que la représentation galoisienne associée à une
représentation automorphe d'un groupe unitaire dans une variété de
Shimura de signature quelconque est bien celle prédite dans
\cite{Langlands1}, du moins en une place décomposée ou bien en
d'autres places sous certaines
hypothèses. 

\subsection{Correspondance de Langlands locale pour $\GL_n$ revisitée}
\label{Lang_Rev}

Soit $F|\Qp$ une extension  de degré fini.
 Soit 
$$
\s_\ell: \mathcal{A}_\ell (n,F) \iso \mathcal{G}_\ell (n,F)
$$
la correspondance de Langlands locale pour $\Gl_n (F)$ $\ell$-adique, au
sens où des deux cotés les représentations sont à valeurs dans des
$\Qlb$-espaces vectoriels et où nous prenons comme définition de
$\mathcal{G}_\ell (n,F)$ : l'ensemble des classes d'isomorphisme de
représentations 
$
\rho: W_F\ldrt \Gl_n(\Qlb)
$
Frobenius semi-simples. L'ensemble $ \mathcal{A}_\ell (n,F)$ est
défini de la même façon que pour des coefficients complexes en
remplaçant formellement $\C$ par $\Qlb$. 

Considérons le groupe $H=\Res_{F/\Qp} (\Gl_{n/F})$, et 
$$ 
\L H=\left ( \prod_{\tau\in I_F} \Gl_n (\Qlb) \right ) \rtimes W_{\Qp}
$$
son $L$-groupe $\ell$-adique où $I_F=\Hom_{\Qp} (F,\Qpb)$ sur lequel 
$\Gal (\Qpb |\Qp)$ agit par composition. 

Le lemme de Shapiro 
$$
H^1(W_{\Qp}, \L H^0) \simeq H^1 (W_F , \Gl_n (\Qlb))
$$
implique alors l'existence d'une bijection 
$$
\pi\longmapsto \tilde{\s}_\ell (\pi)
$$
entre $\mathcal{A}_\ell (n,F)$ et l'ensemble des classes de morphismes
continus de
$W_{\Qp}$ dans $\L H$ induisant l'identité après projection sur $W_{\Qp}$, tels que
l'image du Frobenius de $F$ dans $\L H^0$ soit semi-simple, et où deux
tels morphismes sont dans la même classe ssi ils sont conjugués par un
élément de $\L H^0$. La correspondance $\tilde{\s}$ est la
correspondance de Langlands locale ``absolue'', c'est à dire sur $\Qp$. 

\subsection{Le morphisme $r_\mu$ local} \label{le_r_mu}

Nous explicitons la définition de \cite{Langlands1} dans notre cas.

Soit $H$ comme précédemment et 
$$
\mu:\Gm{/\Qpb} \ldrt H_{/\Qpb}
$$
un cocaractère. Soit $E\subset \Qpb$ le corps de définition de la
classe de conjugaison de $\mu$. 
$$
H_{\Qpb}\simeq \prod_{\tau\in I_F} \Gl_{n\Qpb}
$$
et si via cet isomorphisme $\mu=(\mu_\tau)_{\tau\in I_F}$ où $\forall
\tau \;\mu_\tau\in \Z^n/\S^n \simeq X_* (\Gl_n)/\text{conjugaison}$,
$\Gal (\Qpb | E)$ est le stabilisateur de $(\mu_\tau)_\tau$ dans $\Gal
(\Qpb |\Qp)$. 

\begin{exem}
Si $\tau_0\in I_F$ est fixé et $\forall \tau\neq \tau_0\;\;$
$\mu_\tau=0$, $\mu_{\tau_0}\neq 0$ alors 
$E=F^{\tau_0}$.
\end{exem}

$\mu$ définit un poids dominant $\check{\mu}$ de $\L H^0$ dont la
classe de conjugaison est invariante par $\Gal (\Qpb | E)$. 

Supposons $\check{\mu}$ minuscule donné par des entiers
$(p_\tau)_{\tau\in I_F}$, $0\leq p_\tau \leq n$, et $\mu_\tau
=(\underbrace{1,\dots, 1}_{p_\tau}, 0, \dots, 0)$. 
A $\check{\mu}$ est associé la représentation irréductible de
dimension finie de $\L H^0$ de plus haut poids $
(1,\dots,1)-\check{\mu}$ (et non $\check{\mu}$, à cause de la
convention homologique de Kottwitz qui identifie $\mu$ au caractère
définissant la filtration de Hodge sur le $H_1$ et non le $H^1$).  
Étendons cette représentation en une représentation $r_\mu$ de $\L
H_E$ en posant que l'action de $W_E$ sur l'espace de plus haut poids
 est triviale. 

Plus explicitement :
$$
{r_\mu}_{|\L H^0}= \bigotimes_{\tau\in I_F} \left ( \bigwedge^{p_\tau} St
\right ) ^* 
$$
où $St$ désigne la représentation standard. 
Pour un élément  $\s\in W_E$ 
$$
 (1 \rtimes \s). ( \otimes_\tau v_\tau) = \otimes_\tau v_{\s^{-1}\tau}
$$
(cette expression est bien définie car si $\s\in W_E,
p_{\s^{-1}\tau}=p_\tau$).  

Si $\pi\in \mathcal{A}_\ell (n,F)$ on peut donc définir 
$$
r_\mu \circ \tilde{\s}_\ell (\pi)_{|E}  \in \mathcal{G}_\ell (\dim r_\mu , E)
$$
\\

Si plus généralement 
$ H=\prod_{i\in I} \Res_{F_i/\Qp} ( \Gl_{n F_i} )$
pour des corps locaux $p$-adiques $F_i$
, si $\mu=\prod_{i} \mu_i$ est un cocaractère de $H_{\Qpb}$, si $E_i$ est le
corps réflex de la composante $\mu_i$ alors, le corps réflex de $\mu$,
encore noté $E$, vérifie
$E\subset\prod_i E_i =: M$. De plus, la représentation $r_\mu$ du
L-groupe de $H$ est
encore définie et 
$$
r_{\mu_{|M}}=\bigotimes_{i\in I} r_{\mu_i |M}
$$
Ainsi, si $\pi=\otimes_i \pi_i$ est un représentation irréductible de
$H(\Qp)$, la représentation 
 $\tilde{\s}_\ell (\pi) : W_{\Qp} \ldrt
\L H$ est définie et l'on a l'égalité 
$$
r_\mu \circ \tilde{\s} (\pi)_{|M} =\bigotimes_{i\in I} \left (r_{\mu_i} \circ
\tilde{\s}_\ell (\pi_i)_{|E_i} \right)_{|M}
$$

\subsection{Le théorème} \label{mutau}

 Nous
 supposerons maintenant que le type C.M.  $\Phi$ est induit à partir d'un 
type C.M. de $\mathcal{K}$, c'est à dire qu'il existe un $\tau_0\in
 \Hom (\mathcal{K},\C)$ tel que $\Phi=\{\tau\in \Hom (F,\C)\;|\;
 \tau_{|\mathcal{K}}=\tau_0 \}$.  
Soit $\mu$ le cocaractère de
$G_{\Qb}$ associé à la donnée de Shimura. 
Rappelons que le cocaractère  $\mu$ est donné par un uplet
 $(p_\tau,q_\tau)_{\tau\in \Phi}$ 
où $p_\tau+q_\tau=n$. Nous reprenons les notations du I. 
\\

Rappelons la formule de Matsushima qui décompose la cohomologie de De
Rham de la
variété de Shimura en représentations irréductibles de $G(\A_f)$ : 
$$
\forall i\;\;\;\; \underset{K}{\limi} H^i (\Sh_K,\LL_\rho) =
\bigoplus_{\Pi\in \mathcal{T} (G)_\rho} \dim H^i
(\mathfrak{g}_\infty, K_\infty; \Pi_\infty \otimes \rho). \;\; \Pi_f 
$$
L'action de $G(\A_f)\times \Gal ( \Qb | E)$ sur cette cohomologie
permet quant à elle de décomposer la cohomologie étale $\ell$-adique en 
$$
\forall i\;\;\;\; \underset{K}{\limi} H^i (\Sh_K,\LL_\rho) = 
\bigoplus_{\Pi' \in \text{Irr} (G(\A_f))} \Pi' \otimes
R_{\ell,\rho,\mu}^i ( \Pi') 
$$
où $R_{\ell,\rho,\mu}^i ( \Pi_f)$ est une représentation continue de
$\Gal ( \Qb | E)$. 
De la formule de Matsushima on déduit que 
$$ \exists i\;
R_{\ell,\rho,\mu}^i ( \Pi')\neq 0 \limpl \exists \Pi\in \mathcal{T}
(G)_\rho \;\; \Pi_f\simeq \Pi'
$$
Du théorème \ref{purete_supercusp}, on
déduit que si $\Pi\in \mathcal{T} (G)_\rho$ et 
si en une place $v$ décomposée de $F$, $\Pi_v$ est supercuspidale 
$$
 R_{\ell,\rho,\mu}^i ( \Pi_f) =0 \text{ si } i\neq \dim \Sh
$$
et 
$$
 \dim R_{\ell,\rho,\mu}^{\dim \Sh} (\Pi_f)= m_\Pi \; \dim r_\mu
$$
où $m_\Pi$ désigne la multiplicité de $\Pi$. 
Nous poserons $[R_{\ell,\rho,\mu}^i ( \Pi_f)]$
comme étant la représentation galoisienne semi-simplifiée de la représentation précédente. 

\begin{theo}\label{gal_coho_sh}
Soit $\Pi\in\mathcal{T} (G)_\rho$, soit $[R_{\ell,\rho,\mu} (\Pi_f) ]$ la
représentation semi-simple de $\Gal(\Qb | E)$ associée à $\Pi_f$ en
degré moitié dans la 
cohomologie de la variété de Shimura de type P.E.L. associée à $\mu$.
Supposons qu'en une place finie $v$ de $F$ divisant un premier de $\Q$  décomposé
dans $\mathcal{K}$, $\Pi_v$ soit supercuspidale.
Soit $p$ un premier de $\Q$ décomposé dans $\mathcal{K}$. Notons 
$G(\Qp)=\prod_{i\in I} \Gl_n (F_{v_i})\times \Qp^\times$ et
$\dpt{\Pi_p=\otimes_v \Pi_v \otimes \chi}$.  
 Rappelons qu'un plongement $\nu$ de $\Qb$ dans $\Qpb$ est fixé,
  et qu'on note $\mu_{\Qpb}=\prod_{i\in I} \mu_i$ le cocaractère de Shimura local
associés. Il y alors une égalité 
$$
[R_{\ell,\rho,\mu} (\Pi_f) ]_{|W_{E_\nu}} = [ r_{\mu_{\Qpb}} \circ
\tilde{\s}_\ell (\Pi_p)_{|E_\nu}]^{m_\Pi}\otimes |.|^{-\frac{\dim Sh}{2}} 
$$
comme représentation semi-simple de $W_{E_\nu}$.

De plus, si pour $i$ dans $I$, $E_i$ désigne le corps réflex de $\mu_i$, $E_\nu\subset
 \prod_i E_i =: M$ et 
$$
[R_{\ell,\rho,\mu} (\Pi_f) ]_{|W_M} =
 m_\Pi \bigotimes_{i\in I} [ r_{\mu_i |M} \circ \tilde{\s}_\ell
(\Pi_{v_i})_{|M} ] 
\otimes \s_\ell (\chi)^{-1}_{|M}\; |.|^{-\frac{\sum_\tau p_\tau q_\tau}{2}}
$$
\end{theo}

\begin{proof}
La deuxième formule est une conséquence immédiate de la première.

Commençons par quelques préliminaires de L-groupes ``$\ell$-adiques''
(étant donné que nous ne disposons que de systèmes compatibles de
représentations galoisiennes et pas de représentation du groupe de
Weil global, nous ne pouvons travailler directement avec des L-groupes
complexes). 
Soit le L-groupe $\ell$-adique de $G$ 
$$
\L G=\left ( \prod_{\tau\in \Phi} \GL_n (\Qlb) \times \Qlb^\times
\right )\rtimes
\Gal (\Qb | \Q) 
$$
où $\forall \s\in \Gal (\Qb |\Q)\;
\s.((g_\tau)_\tau,z)=((h_\tau)_\tau,z')$ avec 
$$
\forall \tau\in\Phi\;\; h_\tau= \left \{ g_{\s^{-1}\tau} \text{  si }
\s^{-1}\tau \in \Phi \atop \,w ^t(g_{c\s^{-1}\tau})^{-1}w^{-1}\text{  si }
\s^{-1}\tau\notin \Phi \right.
$$
où 
$$
w= (( -1)^{i-1} \delta_{i, n-j+1})_{i,j}
$$
et
$$
\; z' = z \prod_{\tau\in \Phi \atop \s^{-1} \tau\notin \Phi } \det (
g_\tau) 
$$
L'élément $w$ a été choisi de manière à ce qu'il envoie
l'image par l'automorphisme $g\mapsto \,^t g^{-1}$ de 
l'épinglage standard de $\Gl_n$ sur l'épinglage standard. 
Soit maintenant $H=( \Res_{F/\Q} \Gl_n)\times
\Res_{\mathcal{K}/\Q}\Gm$. Le groupe réductif $\Res_{\mathcal{K}/\Q}
G_\mathcal{K}$  
étant une forme intérieure de $H$ il y a un isomorphisme de L-groupes
$$
\L \left ( \Res_{\mathcal{K}/\Q} G_\mathcal{K} \right ) \iso \L H
\simeq \left ( \prod_{\tau\in I_F} \Gl_n (\Qlb) \times \C^\times\times\C^\times\right )
\rtimes \Gal (\Qb |\Q)
$$
auquel est associé la correspondance de Jacquet-Langlands. 
Il y a également un morphisme de L-groupes 
\begin{eqnarray*}
\L G & \ldrt & \L  \left ( \Res_{\mathcal{K}/\Q} G_\mathcal{K} \right
) \\
((g_\tau)_\tau,z)\rtimes \s & \longmapsto & ((g_\tau)_\tau, w\,^t
(g_\tau^{-1})_{c\tau} w^{-1},z,z\prod_{\tau\in\Phi} \det g_\tau )\rtimes \s 
\end{eqnarray*}
 auquel est associé le changement de base quadratique. Le composé des deux
morphismes ci dessus est associé au changement de base stable de la
section \ref{CHL}. 

Rappelons que  $E$ désigne le corps réflex de $\mu$. Ce corps ne
contient pas forcément le corps réflex de $(F,\Phi)$. 
Comme dans la section précédente dans le cas local, au cocaractère
$\mu$ est associé une 
représentation $r_\mu$ de $\L G_E$ vérifiant
$$
r_{\mu|\L G^0}=\bigotimes_{\tau\in \Phi} \left ( \bigwedge^{p_\tau} St
\right )^* \otimes (St)^{-1}
$$ 
et 
$$
\forall \s\in \Gal(\Qb | E)\;\;\; r_\mu (\s)\left (\bigotimes_{\tau\in \Phi}
v_\tau\otimes z \right )= \bigotimes_{\tau\in \Phi \atop \s^{-1}\tau \in
  \Phi} v_{\s^{-1}\tau} \bigotimes_{\tau\in \Phi \atop
  \s^{-1}\tau\notin \Phi} \,^{t}\left (\bigwedge^{p_\tau} w^{-1} \right
) ( v_{c\s^{-1}\tau}) \otimes z 
$$

\begin{rema}
  Le morphisme $r_\mu$ est bien défini : on doit vérifier l'égalité
  pour tout $g$ et $\s$ 
$$
r_\mu (\s) r_\mu (g) r_\mu (\s)^{-1} = r_\mu ( g^{\s} )
$$
Cela ne pose pas de problème au niveau des composantes du produit
tensoriel 
 indexées par des $\tau\in \Phi$ vérifiant $\s^{-1} \tau \in
 \Phi$. Cependant, si $\s^{-1}\tau\notin \Phi$, il faut remarquer
 qu'étant donné que $\s$ stabilise la classe de conjugaison $\mu$, 
 $q_{c\s^{-1}\tau} = p_{\tau}$, et que si $g$ est un endomorphisme d'un
 espace vectoriel de dimension $n$ alors pour tout $i$,
$\det (g)^{-1}.
 \bigwedge^{i} g = \,^t \left ( \bigwedge^{n-i} g \right )$.
 \end{rema}

Rappelons que $F^\Phi$ désigne le corps réflex de $(F,\Phi)$.
Dans notre cas $F^\Phi=\tau_0(\mathcal{K})$ pour un plongement
$\tau_0$ de $\mathcal{K}$ dans $\C$. Nous
allons nous intéresser au corps composé $E F^\Phi$. 
On peut définir une représentation $r_\mu'$ de $\L H_{EF^\Phi}$ de
la façon suivante :
$$
\L H \simeq \left (\prod_{\tau\in \Phi} \Gl_n (\Qlb) \times
\prod_{\tau\in c \Phi} \Gl_n (\Qlb) \times \Qlb^\times \times
\Qlb^\times \right ) \rtimes \Gal(\Qb |\Q)
$$
Posons 
$$
r'_{\mu | \L H^0}=\bigotimes_{\tau\in \Phi} \left ( \bigwedge^{p_\tau} St
\right )^* \bigotimes_{\tau\in c \Phi} 1 \otimes (St)^{-1}\otimes 1
$$
et $\forall \s\in \Gal (\Qb | E F^\Phi)$ 
$$
r'_\mu (\s) (\bigotimes_{\tau\in \Phi} v_\tau \otimes \l)=\bigotimes_{\tau\in
\Phi} v_{\s^{-1}\tau} \otimes \l
$$
Et remarquons maintenant que le diagramme suivant commute 
\begin{diagram}
\L G_{EF^\Phi} & \rTo & \L (\Res_{\mathcal{K}/\Q} G_{\mathcal{K}})_{EF^\Phi} &
\rTo^\sim & \L H_{EF^\Phi} \\
 & \rdTo^{{r_\mu}_{| E F^\Phi} } & & \ldTo_{r'_\mu} \\
& & \Gl_{dim (r_\mu)} (\Qlb) 
\end{diagram}
où la première application horizontale est l'application associée au
changement de base et la seconde à la correspondance de
Jacquet-Langlands décrites auparavant. 
\\

Notons maintenant $R_1=[R_{\ell,\rho,\mu} (\Pi_f)]$
la représentation de $\Gal ( \Qb | E)$ 
associée à $\Pi_f$ dans la cohomologie de la variété
de Shimura. Soit
$(\Pi',\chi)$ le changement de base stable de $\Pi$ associé à $\Pi$
grâce au théorème \ref{BCCC}. Soit $R'$ la représentation
$\ell$-adique de $\Gal (\overline{F} |F)$ associée à $\Pi'$ par le
théorème \ref{theogalo}. Il résulte du théorème \ref{BCCC} que $\chi$
est un caractère de Hecke algébrique. On peut donc également lui associer un
caractère $\ell$-adique $\chi':\Gal(\overline{\mathcal{K}} |\mathcal{K})\ldrt
\Qlb^\times$.
Faisons maintenant la remarque fondamentale suivante qui résulte de
l'égalité $q_\tau=n- p_\tau$ : 
$$
a:=\sum_{\tau\in \Phi}p_\tau (n-1)-\sum_{\tau\in \Phi} p_\tau q_\tau \equiv 0 \text{ mod } 2
$$
 A l'entier $\frac{a}{2}$ est associé 
 le caractère cyclotomique $\ell$-adique tordu $\frac{a}{2}$-fois (qui est
 le caractère Galoisien $\ell$-adique associé au caractère de Hecke
 algébrique $ ||.||^{a/2}$); nous noterons $(-)\mapsto
 (-)(\frac{a}{2})$ la torsion par ce caractère. 
 Le couple
$(R',\chi')$ fournit (par le lemme de Shapiro) un morphisme 
$$
R'_2:\Gal(\Qb |\Q) \ldrt \L H
$$
Posons 
$$
R_2=\left ( r'_\mu\circ R'_{2|E F^\Phi }\right )^{m_\Pi} (\frac{a}{2})
$$
comme représentation de $\Gal(\Qb | EF^\Phi)$. 

Montrons que $R_{2|\Gal ( \Qb | EF^\Phi)}=R_1$ ce qui conclura la démonstration  d'après les
propriétés de $R'$ 
énoncées  dans le théorème  \ref{theogalo}, la forme du changement de
base local en un $p$ décomposé, la définition de $r'_\mu$, le choix
de $a$ et qu'étant donné que $p$ est décomposé dans $\mathcal{K}$, $(E
F^\Phi)_{\nu}= E_\nu$.

Pour montrer que $R_{2|\Gal ( \Qb | EF^\Phi)}=R_1$ il suffit d'après le
  théorème de densité 
de Tchebotarev de montrer que pour presque toutes les places $v$ de $EF^\Phi$,
$ R_{2|W_{(EF^\Phi)_v}}=R_{1 |W_{( EF^\Phi)_v}}$. Restreignons nous aux places $v$
divisant des premiers $p_0$ de $\Q$ tels que $\Pi_{p_0}$ soit non
ramifiée. L'égalité ci dessus résulte alors de la commutativité du
diagramme global ci dessus restreint à $W_{(EF^\Phi)_v}$ et du théorème de Kottwitz
calculant la représentation galoisiennes non ramifiée en une place de
bonne réduction.
\end{proof}

\begin{coro}
Soit $v_0$ une place de $F$ fixée.
Soit $\nu$ un plongement de $\Qb$ dans $\Qpb$ tel que si $\tau\in
\Phi$ vérifie $\nu\circ \tau : F\hookrightarrow \Qpb$ n'induise pas
$v_0$ alors $p_\tau=0$ (autrement dit, la donnée de Shimura locale
associée est étale hors $v_0$). Alors, $E_\nu=E(\mu_{v_0})$ et 
$$
[R_{\ell,\rho,\mu} (\Pi_f) ]_{|W_{E_\nu}} = m_\Pi [r_{\mu_{v_0}}\circ
\tilde{\s}_\ell (\Pi_{v_0})_{|E(\mu_{v_0})}] \otimes \s_\ell (\chi)_{|E_\nu}^{-1} \;
|.|^{\frac{-\sum_\tau p_\tau q_\tau}{2}} 
$$
\end{coro}

\begin{rema}
Les représentations $\ell$-adiques des groupes de Weil locaux considérées ci dessus
sont semi-simples et pas seulement Frobenius semi-simples : on a
semi-simplifié la monodromie de ces représentations. Ainsi, par
exemple si $\rho$ est une représentation galoisienne locale 
$$
[\rho \otimes sp (r)]=r [\rho]
$$
Cependant le théorème ci dessus devrait être encore vrai sans
semi-simplification de la monodromie, c'est à dire il devrait y avoir
un isomorphisme comme représentations du groupe de Weil-Deligne. 
\end{rema}

Nous allons maintenant nous intéresser à la restriction de la
représentation galoisienne obtenue dans la cohomologie de la variété
de Shimura en une place non décomposée. Soit donc $p$ un premier de
$\Q$ inerte dans $\mathcal{K}$. Avec les notations du I, il y a une
décomposition 
$$
G(\Qp)=\prod_{i\in I} \GL_n (F_{v_i}) \times G( \prod_{j\in J} U(n , F_{v_j}))
$$
où l'on suppose que l'ensemble $J$ paramétrant les places de $F$
égales à leur complexe conjuguée est non vide. Si $\Pi$ est une
représentation automorphe de $G$ satisfaisant les hypothèses du
théorème \ref{BCCC}, le changement de base stable global $BC (\Pi)$
consiste en une représentation automorphe $\Pi'$ de $\GL_{n/ F}$ et un
caractère de Hecke $\chi$ de $\A_{\mathcal{K}}^\times$. La composante en
$p$ de $\Pi'\otimes \chi$ est donc une représentation du groupe
p-adique 
$$
\prod_{i\in I} \left (\GL_n ( F_{v_i})\times \GL_n (F_{v_i^c}) \right )\times
\prod_{j\in J} \GL_n ( F_{v_j}) \times \mathcal{K}_p^\times 
$$
de la forme 
$$
(\Pi'\otimes \chi)_p = \bigotimes_{i\in I} \left ( \pi_i\otimes
  \check{\pi}_i \right )  \otimes \bigotimes_{j\in J} \pi_j \otimes
\chi_p 
$$
où $\forall j\in J \; \pi_j\simeq \check{\pi}_j$. Nous noterons alors
$$
\text{BC} ( \Pi )_p^0 = \bigotimes_{i\in I} \pi_i \otimes
\bigotimes_{j\in J} \pi_j \otimes \chi_p 
$$
comme représentation de $\prod_{i\in I} \GL_n (F_{v_i}) \times
\prod_{j\in J} \Gl_n ( F_{v_j}) \times \mathcal{K}_p^\times$ qui est
le groupe des $\Qp$-points du groupe algébrique $\prod_{i\in I} \Res_{F_{v_i}/\Qp}
  \GL_n \times\prod_{j\in J} \Res_{F_{v_j}/\Qp} \GL_n \times
  \Res_{\mathcal{K}_p /\Qp} \Gm$. On peut en
particulier définir la représentation $\tilde{\s}_\ell (\text{BC} ( \Pi
)_p^0 )$ de $\Gal ( \Qpb |\Qp)$ dans le L-groupe de ce groupe
réductif (confère \ref{Lang_Rev}). Dans la démonstration du théorème
précédent, l'égalité $R_{2|\Gal(\Qb | E F^\Phi )}=R_1$ implique alors le
théorème suivant : 

\begin{theo}
  Plaçons nous dans le cadre du théorème précédent en supposant cette
  fois ci que $p$ est inerte dans $\mathcal{K}$. Il y a une égalité
  pour toute $\Pi\in \mathcal{T} (G)_\rho$ : 
$$
[R_{\ell,\rho,\mu} (\Pi_f)]_{| W_{(E_\nu \mathcal{K}_p)}} =[
r_{\mu_{\Qpb}} \circ \tilde{\s}_\ell ( \text{BC} (\Pi)^0_p )_{|E_\nu
  \mathcal{K}_p}]^{m_\Pi} \otimes |.|^{-\frac{\dim \Sh}{2}}
$$
et si $\text{BC} (\Pi)^0_p =\otimes_{i\in I} \pi_i \otimes_{j\in J}
\pi_j \otimes \chi$, si $M$ désigne le produit des corps réflex
locaux, on a l'égalité 
$$
[R_{\ell,\rho,\mu} (\Pi_f)]_{| W_{(M \mathcal{K}_p)}}]
= m_\Pi \bigotimes_{i\in I} [ (r_{\mu_i})_{|M \mathcal{K}_p} \circ
\tilde{\s}_\ell ( \pi_i)_{|M \mathcal{K}_p}\bigotimes_{j\in J}  (r_{\mu_j})_{|M \mathcal{K}_p} \circ
\tilde{\s}_\ell ( \pi_j )_{| M \mathcal{K}_p} \otimes \s_\ell
(\chi)_{|M\mathcal{K}_p} 
$$
\end{theo}

\begin{rema}
Dans le cas d'un groupe unitaire en $3$-variables les travaux de
Rogawski (\cite{Rogawski1} et appendice \ref{BC_local}) montrent que
dans le théorème précédent $BC (\Pi)_p^0$ ne dépend que de $\Pi_p$ et
que donc la représentation galoisienne découpée par $\Pi_f$ restreinte
à un groupe de décomposition en $p$ ne dépend
que de $\Pi_f$. 
\end{rema}

Nous allons préciser ce théorème dans le cas de d'un groupe unitaire
en trois variables. 

\begin{theo}\label{rep_gal_U3}
Supposons $n=3$ et $G(\Qp)=GU(3)$ le groupe non ramifié en trois
variables. Soit $\Pi\in\mathcal{T}(G)_\rho$. Supposons  $\Pi_p$
supercuspidale ou bien  une représentation de Steinberg de
$G(\Qp)$. Soit $\psi:W_{\Qp}\ldrt \L GU(3)$ le L-paramètre local
associé (confère l'appendice \ref{Resultats_Rog}; dans le cas de la représentation
notée $\pi^{s} (\xi)$ dans l'appendice  \ref{Resultats_Rog} ou de la représentation
de Steinberg, le L-paramètre n'est pas trivial sur le facteur $SL_2
(\C)$;  dans ce cas on note $\psi$ la restriction à $W_{\Qp}$ de ce
L-paramètre). Supposons qu'en une place finie de $\Q$ décomposée dans $\mathcal{K}$
, $\Pi$ est supercuspidale. Il y a a alors une égalité 
$$
[R_{\ell,\rho,\mu}(\Pi_f)]_{|W_{E_\nu}}= [r_\mu\circ \psi_{|W_{E_\nu}}]
\; |.|^{-\frac{\dim \Sh}{2}}
$$
\end{theo}

\begin{proof}
Notons $(\Pi',\chi)$ le changement de base stable de $\Pi$
(\ref{CHL}) et $R'$ la représentation $\ell$-adique de $\Gal
(\overline{F}|F)$ associée à $\Pi'$ par le théorème
\ref{theogalo}. Étant donné que $\check{\Pi}'\simeq {\Pi'}^c$,
$\check{R}'\simeq {R'}^c$ (utilisant qu'en une place finie de $F$,
$R'$ réalise la correspondance de Langlands locale on en déduit que
cet isomorphisme est vrai en toute les places finies
de $F$ et qu'il est donc vérifié par le théorème de densité de
Tchebotarev). La représentation
$\Pi$ étant supercuspidale en une place finie décomposée, $\Pi'$ est
supercuspidale en une place finie, et donc la restriction de $R'$ à un
groupe de décomposition en cette place est irréductible (il s'agit d'une
propriété de la correspondance de Langlands locale). La représentation
$R'$ est donc irréductible. Comme dans le lemme 15.1.2 de
\cite{Rogawski1}, il existe un signe $c(R')\in\{ \pm 1 \}$ obstruction à
ce que $R'$ se relève en un L-paramètre $\Gal(\Qb | F^+)\ldrt \L
G_1$. Montrons que $c(R')=1$. Le caractère $\det (R')$ est le
caractère $\ell$-adique associé au caractère central de $\Pi'$
(puisque cela est vérifié en toutes les places finies d'après une
propriété standard de la correspondance de Langlands locale).  Comme
dans le lemme 15.1.2 de \cite{Rogawski1}, $c(R')=1$ ssi ce caractère
central est trivial sur $(F^+)^\times \bc \A_{F^+}^\times$, ce qui est le
cas d'après la dernière propriété du théorème \ref{CHL}. Il existe
donc une unique extension de $R'$, 
$$
\varphi':\Gal (\Qb | F^+) \ldrt \L G_1
$$
où $G_1$ est vu comme groupe sur $F^+$. D'où un L-paramètre 
$$
\varphi'': \Gal (\Qb |\Q) \ldrt \L G_1
$$
où cette fois ci $G_1$ est un groupe sur $\Q$. Le caractère algébrique
$\chi$ fournit un caractère $\ell$-adique $\chi'$ qui couplé à
$\varphi''$ nous donne un L-paramètre 
$$
\varphi : \Gal (\Qb |\Q) \ldrt \L G
$$
L'application de changement de base locale de Rogawski est injective
(\cite{Rogawski1} théorème 13.2.1). On montre de même que si $GU(3)$
désigne le groupe de similitudes unitaires $p$-adique non
ramifiée $3$ variables, si $\varphi:W_{\Qp}\ldrt \L GU(3)$ est un L-paramètre,
$\varphi_{|W_K}$ détermine $\varphi$ de façon unique, où $K$ est
l'extension quadratique non ramifiée de $\Qp$. On en déduit en particulier que
$$
\varphi_{|W_{\Qp}}= \psi
$$
et qu'en toutes les places $v$ où $\Pi$ est non ramifiée,
$\varphi_{|W_{\Q_v}}$ est associé à $\Pi_v$ via la correspondance de
Langlands locale non ramifiée. 

Considérons maintenant $r_\mu\circ \varphi$. Comme dans la démonstration
précédente, d'après le théorème de Kottwitz (\ref{theo_kott_sh}),
$[r_\mu\circ \varphi]$ est égal à $[R_{\ell,\mu,\rho} (\Pi_f)]\;
|.|^{-\frac{\dim \Sh}{2}}$ en
presque toutes les places où $\Pi$ est non ramifiée. D'où le résultat
d'après le théorème de densité de Tchebotarev. 
\end{proof}

\chapter[Théorie des types]{Résultats concernant les types associés aux
représentations des groupes p-adiques} \label{Annexe_types}

\section{Quelques propriétés générales des types associés aux
représentations supercuspidales des groupes $p$-adiques (\cite{BushKutz})}

Soit $F$ un corps local et $G$ les points à valeurs dans $F $ d'un
groupe réductif sur $F$. Soit $\pi$ une représentation supercuspidale
de $G$. Notons $\mathfrak{s}=[\pi,G]$ sa classe d'équivalence
inertielle. 

Nous ferons l'hypothèse suivante :

\begin{enonce}{Hypothèse} \label{hypo_types}
 Il existe un $\mathfrak{s}$ type $(J,\l)$, un 
sous-groupe compact ouvert modulo le centre de $G$, $\tilde{J}$ tel que
$J=\tilde{J}\cap G^1$, et une extension $\tilde{\l}$ de $\l$ à
$\tilde{J}$ tels que 
$$
\pi\simeq c-Ind^G_{\tilde{J}} \; \tilde{\l}
$$
\end{enonce}

Soient maintenant $\chi_1,\chi_2$ deux caractères non ramifiés de $G$.

\begin{Fait} On a l'équivalence 
$$
\pi\otimes\chi_1\simeq \pi\otimes\chi_2 \lssi \chi_{1|\tilde{J}}=
\chi_{2|\tilde{J}}
$$
\end{Fait}

En effet, l'implication de la droite vers la gauche résulte de
l'isomorphisme 
$
(c-Ind^G_{\tilde{J}} \pi)\otimes \chi\simeq c-Ind^G_{\tilde{J}} ( \pi\otimes\chi_{|\tilde{J}})
$.
Dans l'autre sens, si $\pi\otimes\chi_1\simeq \pi\otimes\chi_2$, alors 
$(\pi\otimes \chi_1)_{|\tilde{J}} \simeq (\pi\otimes
\chi_2)_{|\tilde{J}}$ or, pour $i=1,2$ il résulte de la proposition
5.6 de \cite{BushKutz} que $(\pi\otimes\chi_i)_{|\tilde{J}}=
  \tilde{\l}\otimes\chi_{i|\tilde{J}}$. $\tilde{\l}$ étant irréductible
  de dimension finie, on en déduit en utilisant le caractère de
  $\tilde{\l}$ que $\chi_{1|\tilde{J}}=
\chi_{2|\tilde{J}}$.
\\

Notons $X_\C (G)$ le groupe des caractères non ramifiés de $G$. Le
groupe $X_\C
(G)$ s'identifie aux points d'un tore algébrique sur $\C$.

\begin{Fait}\label{Fait_type1} 
Soit $\GG=\{ \chi\in X_\C (G)\; |\; \pi\otimes\chi\simeq \pi \}$, un
sous-groupe fini de $X_\C (G)$. L'algèbre de Hecke du type,
$\mathcal{H} (\l,G)$ est isomorphe à l'algèbre des fonctions
régulières sur la composante connexe du centre de Bernstein indexée
par $\mathfrak{s}$. L'algèbre $\mathcal{H} (\l,G)$ s'identifie donc aux
fonctions régulières sur le tore complexe 
$X_\C(G)/\GG$. 
\end{Fait}

Notons $e_\l\in \mathcal{H} (G)$ l'idempotent associé à
$\l$ et  $\mathcal{H} (\l,G)=e_\l*\mathcal{H} (G)* e_\l$.  Pour toute
représentation $\pi'$ dans la classe d'équivalence 
$\mathfrak{s}$, l'espace $\pi' (e_\l). V_{\pi'}$ est un 
 $\mathcal{H} (\l,G)$-module sur lequel  $\mathcal{H} (\l,G)$ opère
 par un caractère qui correspond au point associé à $\pi'$. 

L'application suivante est alors un isomorphisme 
\begin{eqnarray*}
 \mathcal{H} (\l,G) & \iso & \C [ X_\C (G)/\GG ] \\
f & \mapsto & \left [ [\chi ]\mapsto \text{tr}_{\pi\otimes\chi}
(f)\right ]
\end{eqnarray*}
\\

Nous aurons besoin du lemme suivant :

\begin{lemm}\label{lemme_typique} On a l'inégalité 
$$
\forall g\in\tilde{J}\;\;\;\text{tr}_\pi (e_\l * \delta_g)\neq 0
$$
\end{lemm}
\begin{proof}
Commençons par remarquer que 
$$
\text{tr}_\pi ( e_\l*\d_g)=\text{tr}_\pi ( e_\l*\d_g*e_\l)\text{ où }
e_\l*\d_g*e_\l \in \mathcal{H} (\l,G)
$$
On vérifie de plus que 
$$
\forall \chi\in X_\C (G) \;\;\; \text{tr}_{\pi\otimes\chi} (e_\l
*\d_g) = \chi (g) \;\text{tr}_{\pi} ( e_\l*\d_g)
$$
et donc, $\text{tr}_\pi (e_\l*\d_g)=0 \lssi $ la fonction régulière
associée à $e_\l*\d_g*e_\l$ est nulle $\lssi$ $e_\l*\d_g*e_\l=0$ (grâce
à l'isomorphisme ci dessus).

Or, 
$$
(e_\l*\d_g*e_\l)*\d_{g^{-1}}= e_\l * (\d_g*e_\l *\d_{g^{-1}}) = e_\l *e_{\l^g}
$$
Étant donné que $\l$ s'étend en $\tilde{\l}$ et que $g\in \tilde{J}$,
$\l^g\simeq \l$ et donc $e_{\l^g}=e_\l$. D'où
$$
(e_\l*\d_g*e_\l)*\d_{g^{-1}} = e_\l *e_\l =e_\l\neq 0 \limpl e_\l* \d_g
* e_\l\neq 0
$$
\end{proof}

\section{Types pour $\Gl_n$}

Rappelons l'un des théorèmes principaux de \cite{BushKutzGln} 

\begin{theo}[\cite{BushKutzGln}] Pour $G=\Gl_n$, l'hypothèse \ref{hypo_types}
de la section précédente est vérifiée. 
\end{theo}

Nous remarquerons de plus que tous les groupes $\tilde{J}$ obtenus
sont contenus dans des normalisateurs de groupes parahoriques dans
$\Gl_n$ (cela est valable pour n'importe quel sous-groupe compact modulo le
centre de $\Gl_n$ et résulte de la classification des sous-groupes
compacts maximaux de $\text{PGL}_n$). 

\section{Types pour les groupes unitaires non ramifiés}
\subsection{Le groupe unitaire en trois variables}\label{types_U3}

Dans \cite{Moy1} le théorème suivant est démontré :

\begin{theo}
  Supposons $p$ différent de $2$. Toutes les représentations
  supercuspidales du groupe
  unitaire quasi-déployé $U(3)$ sur un corps local p-adique possèdent
  un type vérifiant l'hypothèse \ref{hypo_types}. 
\end{theo}

Nous utiliserons des types pour les groupes de similitudes unitaires
non ramifiés. Bien que le théorème précédent soit démontré dans le cas
des groupes unitaires, on peut l'étendre aux groupes de similitude
unitaires en utilisant la section \ref{pass_u_gu} 

\subsection{Le groupe unitaire général}\label{types_Un}

Dans l'article \cite{Kim1}, Ju-Lee Kim construit des types 
vérifiant l'hypothèse \ref{hypo_types} pour des
représentations supercuspidales des groupes unitaires non
ramifiés. Elle conjecture qu'elle obtient ainsi toutes les
représentation supercuspidales de ces groupes. Ses résultats
s'étendent aux groupes de similitudes unitaires (bien que cela n'est
pas été publié).

\chapter[U(3)]{Résultats de Rogawski sur $U(3)$}\label{Resultats_Rog}

Nous rassemblons dans cette annexe les résultats de \cite{Rogawski1}
que nous utilisons. La lecture des articles \cite{Rogawski2} et
\cite{Rogawski3} est également utile. 

\section{}
\subsection{Notations}

Soit $F|F_0$ une extension quadratique de corps locaux
$p$-adiques. Soit $G=U(3)$ le groupe unitaire en trois variables
quasidéployé vu comme groupe sur $F_0$. 

Soit $H=U(2)\times U(1)$ l'unique groupe endoscopique elliptique de $G$. Il est
donc muni d'un morphisme de L-groupes $\L H \ldrt \L G$ décrit dans le
chapitre 4.6 et 4.7 de \cite{Rogawski1}.

Soit $C=U(1) \times U(1) \times U(1)$ l'unique groupe endoscopique
elliptique de
$H$. 

Les groupes $C$ et $H$ sont les deux groupes endoscopiques 
associés à  $G$ (plus généralement, les
groupes endoscopiques associés aux groupes unitaires $U(n)$ sont
décrits dans la partie 1.2 de \cite{Rogawski2}). 

\subsection{Transfert endoscopique local}
\subsubsection{De $C$ à $H$}

Dans le chapitre 11 de \cite{Rogawski1}, l'existence du transfert
endoscopique des représentations admissibles irréductibles de $C$ vers
celles de $H$ est énoncé (les démonstrations étant identiques à celles
de \cite{Lab3}). Il s'agit d'une application injective 
\begin{eqnarray*}
  \Pi (C) & \ldrt & \Pi (H) \\
\theta & \longmapsto & \rho (\theta )
\end{eqnarray*}
où $\Pi (C)$ désigne les représentation irréductibles de $C(F_0)$
(i.e. les caractères $\theta=\theta_1\otimes \theta_2 \otimes
\theta_3$ de $U(1)\times U(1) \times U(1)$) et $\Pi (H)$ les L-paquets
de représentations admissibles irréductibles irréductibles de
$H(\Qp)$. 

\subsubsection{De $H$ à $G$}

Dans le chapitre 13 de \cite{Rogawski1}, l'existence du transfert
local 
\begin{eqnarray*}
  \Pi (H) & \ldrt & \Pi (G) \\
\rho & \longmapsto & \Pi (\rho) 
\end{eqnarray*}
est démontré, où $\Pi(G)$ désigne l'ensemble des L-paquets associés à
$G(F_0)$. 

\subsection{Changement de base local}\label{BC_local}

Dans le chapitre 13.2 de \cite{Rogawski1}, Rogawski démontre
l'existence d'une application de changement de base quadratique 
$$
BC: \Pi (G) \ldrt \Pi (\GL_{n /F})=\mathcal{A} (3,F)
$$
compatible avec le changement de base global de l'annexe A. 
Son image est constituée de représentations de $\GL (3,F)$, $\pi$
vérifiant $\check{\pi}\simeq \pi^c$ et $\omega_{\pi | F_0^\times}=1$
où $c$ désigne l'automorphisme non trivial de $F|F_0$ et $\omega_\pi$
le caractère central de $\pi$. 

\subsection{Classification des L-paquets supercuspidaux de $G$}

Dans la section 12-2 de \cite{Rogawski1}, il est donné une description
détaillée des L-paquets de représentations de $G(F_0)$. Nous rappelons
cette classification dans le cas des représentations supercuspidales
qui est le seul cas que nous considérons. 

\subsubsection{L-paquets supercuspidaux pour H}[proposition 11.1.1 de
\cite{Rogawski1}] 

\begin{defi}
Un L-paquet de représentations de $H$ contient une représentation
supercuspidale ssi ce L-paquet est formé de représentations
supercuspidales (ce fait est faux en général). Nous appellerons un tel
L-paquet un L-paquet supercuspidal. 
\end{defi}

Les L-paquets supercuspidaux de représentations de $H$  se scindent en
deux ensembles : 
\begin{itemize}
\item Les L-paquets de cardinal $1$ ou encore L-paquets stables.
 Ce sont les $\{ \rho\}$ où $\rho$
  est une représentation supercuspidale de $H$ telle que le changement
  de base quadratique $BC(\rho)$ soit une représentation supercuspidale
  de $\GL_2 (F)\times F^\times$.
\item Si $\theta$ est un caractère de $C$ de la forme
  $\theta_1\otimes\theta_2\otimes\theta_1$ où $\theta_1\neq \theta_2$,
  le L-paquet transfert endoscopique $\rho (\theta)$ est supercuspidal
  de cardinal $2$. Lorsque $\theta$ parcourt les caractères de $C$ du
  type précédent on obtient ainsi tous les L-paquets supercuspidaux de
  $G$ de cardinal plus grand que $1$. Un L-paquet supercuspidal $\Pi$ est de la
  forme $\rho(\theta)$ ssi $BC(\Pi)$ n'est par supercuspidale
  ($BC(\Pi)$ est alors un quotient d'une série principale de $\GL_2
  (F)\times F^\times$ calculée en fonction de $\theta$). 
\end{itemize}

\subsubsection{L-paquets supercuspidaux pour G} \label{paquets_sup}

\begin{defi}
Un L-paquet de représentations de $G$ est appelé un L-paquet
supercuspidal s'il est constitué de représentations supercuspidales.
\end{defi}

Les L-paquets supercuspidaux de représentations de $G$ se scindent en trois sous ensembles
:

\begin{itemize}
\item Les L-paquets de cardinal $1$ ou encore L-paquets stables. Ce
  sont les $\Pi=\{\pi\}$ où 
  $\pi$ est une représentation supercuspidale telle que $BC(\pi)$
  soit une représentation supercuspidale de $\GL_3 (E)$.
\item Les L-paquets de cardinal $2$. Ils sont construits de la façon
  suivante : soit $\rho$ un L-paquet supercuspidal de cardinal $1$ de
  $H$. Alors, le transfert endoscopique $\Pi (\rho)$ est un L-paquet
  supercuspidal de cardinal $2$ de $G$. Un L-paquet supercuspidal
  $\Pi$ est de cette forme ssi $BC(\Pi)$ est une représentation de
  $\GL_3 (F)$ dont le support supercuspidal est une représentation du
  sous groupe de Levi $\GL_2 (F)\times \GL_1 (F)$.
\item Les L-paquets de cardinal $4$. Ils sont construits ainsi : 
Soit $\rho=\rho (\theta)$ un L-paquet supercuspidal de cardinal $2$ de
$H$ tel que si $\theta=\theta_1\otimes \theta_2\otimes \theta_3$ alors
$\theta_1, \theta_2, \theta_3$ sont deux à deux distincts. 
 Le L-paquet obtenu par transfert de $H$ à $G$ de $\rho$, noté
$\Pi (\rho)$ est alors de cardinal $4$. Un L-paquet supercupsidal
$\Pi$ de $G$ est de ce type là ssi $BC (\Pi)$ est un sous-quotient
d'une série principale de $\GL_3 (F)$. 
\end{itemize}

\subsection{Correspondance de Langlands locale pour les
L-paquets  supercuspidaux de $G$}

Dans le chapitre 15 de \cite{Rogawski1}, il est expliqué comment 
construire une correspondance de Langlands locale pour les
L-paquets supercuspidaux de $G$ à partir de la correspondance
de Langlands locale pour $\GL_3 (F)$, et la classification précédente.

Rappelons cette correspondance :
\begin{itemize}
\item L'application $BC:\Pi (G)\ldrt \mathcal{A} (3,F)$ induit une
  bijection entre les L-paquets supercuspidaux de cardinal $1$ et les
  représentations supercuspidales $\pi'$ de $\GL_3 (F)$ vérifiant
  $\check{\pi}'\simeq {\pi'}^c$ (où $c$ est l'automorphisme non trivial
  de $F|F_0$) et $\omega_{\pi'|F_0^\times}=1$. 
 Si $\{\pi\}\in \Pi(G)$ est un tel L-paquet, si
  $\varphi:W_F\ldrt \GL_3 (\C)$ est le paramètre de Langlands associé à
  $BC(\pi)$ (une représentation irréductible de dimension $3$), $\varphi$
  défini un morphisme de L-groupe $W_F \ldrt \L G_F=\GL_3 (\C)\times
  W_F$ qui possède une unique extension en un L-morphisme
  $W_{F_0}\ldrt \L G$. Cela établit une bijection entre les L-paquets
  supercuspidaux de cardinal $1$ de $G$ et les classes de
  $\L G^0$-conjugaison $\varphi$ de L-morphismes de $W_{F_0}$ dans $\L G$
  telles que $\varphi_{|W_{F_0}}$ soit irréductible. 
\item On établit de même une correspondance de Langlands locale pour
  les L-paquets supercuspidaux de $H$ de cardinal $1$. Alors, si
  $\Pi=\Pi (\rho)$ est un L-paquet supercuspidal de cardinal $2$ de
  $G$, si $\varphi:W_{F_0}\ldrt \L H$ est le paramètre de Langlands
  associé à $\rho$, le paramètre associé à $\Pi$ est celui composé 
$W_{F_0}\ldrt \L H \ldrt \L G$. Cela établit une bijection entre les
L-paquets supercuspidaux de $G$ de cardinal $2$ et les classes de $\L
G^0$-conjugaison de L-morphismes $\varphi:W_{F_0}\ldrt \L G$
telles que $\im (\varphi)$ ne soit pas contenue dans $\L B$, où $B$ est
un sous-groupe de Borel de $G$, et telles que
$\varphi_{|W_F}$ soit somme de deux représentations irréductibles.
\item La correspondance de Langlands locale pour $C$ se déduit de la
  théorie du corps de classe local. Si $\Pi=\Pi (\rho (\theta))$ est
  un $L$-paquet supercuspidal de cardinal $4$ de $G$, où $\theta$ est
  un caractère de $C$, si $\varphi:W_{F_0}\ldrt \L C$ est le paramètre
  associé à $\theta$, on associe à $\Pi$ le paramètre composé 
$W_{F_0}\ldrt \L C \ldrt \L H \ldrt \L G$. Cela définit une bijection entre
les L-paquets supercuspidaux de $G$ de cardinal $4$ et les classes de
$\L G^0$-conjugaison de L-morphismes $\varphi: W_{F_0}\ldrt \L G$
telles que $\im (\varphi)$ ne soit pas contenue dans $\L B$, où $B$ est
un sous-groupe Borel de $G$, et
 telles que
$\varphi_{|W_{F_0}}$ soit somme de trois caractères.
\end{itemize}

Il y donc une bijection entre les
L-paquets supercuspidaux 
 de $G$ et les classes de $\L G^0$-conjugaison
de L-morphismes $\varphi:W_{F_0}\ldrt \L G$ telles que $\im (\varphi)$ ne
soit pas contenue dans $\L B$ pour un sous-groupe de Borel $B$ de
$G$. L'application de changement de base se lit au niveau des
L-paramètres comme l'application de restriction $\varphi\longmapsto
\varphi_{|W_{F}}$. 

\subsection{Classification des représentations supercuspidales de $G$}
\label{Class_sup}

Comme nous l'avons vu, les représentations supercuspidales de $H$ sont
toutes membres d'un L-paquet supercuspidal. Cela est faux pour $G$. Il
existe une famille de 
 représentations supercuspidales de $G$, $(\pi^s (\xi))_{\xi}$,
toutes membres de L-paquets de
cardinal $2$  qui ne sont pas supercuspidaux.

Ces représentations supercuspidales se construisent ainsi : soit $\xi$
un caractère de $H$. A $\xi$ est associé la représentation de
Steinberg $St_H (\xi)$ de $H$ définie comme constituante d'une série
principale réductible de $H$ (l'autre constituant étant $\xi$, confère
\cite{Rogawski1} 12.1). L'ensemble $\{St_H (\xi) \}$ est un L-paquet
stable de carré intégrable de représentations de $H$. 
Considérons le L-paquet de représentations de $G$, $\Pi (St_H (\xi))$,
transfert endoscopique de $\{St_H (\xi) \}$. C'est un L-paquet de
cardinal $2$ noté 
$$
 \Pi (St_H (\xi))=\{ \pi^2 (\xi),\pi^s (\xi) \}
$$
où $\pi^s (\xi)$ est une représentation supercuspidale et $\pi^2
(\xi)$ une représentation de carré intégrable. Alors,

\begin{enonce}{Fait}
Les représentations supercuspidales de $G$ sont soit membres d'un
L-paquet supercuspidal soit de la forme $\pi^s (\xi)$ pour une unique
représentation de dimension $1$ de $H$, $\xi$. 
\end{enonce}

\subsection{Passage de $U(3)$ à $GU(3)$}\label{pass_u_gu}

  Soit $G=GU(3)$ un groupe de similitudes unitaires non ramifié en $3$
  variables défini
  sur $\Qp$ et $G_1\subset G$ le groupe unitaire noyau du facteur de
  similitude. Soit $K|\Qp$ l'extension quadratique non ramifiée. Il y
  a une inclusion $K^\times\hookrightarrow  Z_G$. 

  \begin{lemm}
    On a l'égalité suivante 
$$
G(\Qp)=G_1 (\Qp) K^\times
$$
  \end{lemm}
  \begin{proof}
    Étant donné que $G$ est non ramifié il y a une décomposition
    d'Iwasawa $G(\Qp)=C.A (\Qp). N(\Qp)$, où $C$ est compact, $A$ est
    un tore déployé maximal et $N$ unipotent. Donc, $v_p ( c
    (G(\Qp)))=v_p ( c (A(\Qp)))$. Avec les notations de l'appendice D, 
$A(\Qp)=\{\text{diag}(t,x,t^{-1} x^2 )\;|\; t\in \Qp^\times \; x\in
    \Qp^\times \;\}$ et $c=x^2$. Donc, $v_p(c ( G(\Qp)))=2 \Z$. 

De plus, $N_{K/\Qp} ( K^\times )=p^{2\Z} \Zp^\times $ et donc
$c(G(\Qp))\subset N_{K/\Qp} ( K^\times )$. D'où le résultat puisque
$c_{|K^\times} =N_{K/\Qp}$.   
  \end{proof}

  \begin{coro}
    Il y a une bijection entre les représentations irréductibles de
    $G(\Qp)$ et les couples $(\pi,\chi)$ où $\pi$ est une
    représentation irréductible du groupe unitaire $G_1 (\Qp)$ et
    $\chi$ un caractère de $K^\times$ tel que $\omega_{\pi | K^\times}=\chi$.
  \end{coro}

  \begin{coro}
    La description des représentations supercuspidales de $G(\Qp)$ est
    identique à celle de $G_1 (\Qp)$. En particulier, il y a une
    paramétrisation des paquets supercuspidaux de $G(\Qp)$ par des
    L-paramètres $\psi: W_{\Qp}\ldrt \L G$ tels que $\im (\psi)$ n'est
    pas contenu dans $\L B$ pour un sous groupe de Borel $B$ de $G$. 
  \end{coro}

\chapter[Isocristaux munis de structures additionnelles]
{Isocristaux munis de structures additionnelles : le cas
 des groupes unitaires}\label{app_cl_sig}

Nous traduisons dans cet appendices les résultats de \cite{Ko2} et
\cite{Ko3} dans le cas des groupes de similitudes unitaires quasidéployés.
Nous utiliserons librement les notations de ces deux articles. 

\section{Groupes unitaires en un nombre pair de variables}
\label{Isoc_AU}

Soit $F|\Qp$ non ramifiée de degré pair $[F:\Qp]=2d$. Nous notons
comme d'habitude $\s$ le Frobenius de $F|\Qp$ et $\GG$ le groupe de
Galois $\Gal (F|\Qp)$. 
Soit $*=\s^d$
l'involution de degré $2$ et $F_0$ le corps fixe de $*$. Notons $G$ le
groupe des similitudes unitaires quasidéployé en $n$ variables, $n$
pair, 
$$
G=\{g\in \Gl_n(F) \;|\; {}^t g^* J g=c(g) J, c(g)
\in\Qp\* \}
$$
où $J=\left (
  \begin{array}{cc}
0 & I_{n/2} \\
I_{n/2} & 0 \\
  \end{array}
\right )
$

Un tore maximal de $G$ défini sur $\Qp$ est 
$$
T=\left \{\left ( \begin{array}{cccccc}
t_1 \\
& \ddots \\
& & t_{n/2} \\
& & & c (t_1^{-1})^* \\
& & & & \ddots \\
& & & & & c (t_n^{-1})^* 
\end{array} \right ) \;\;|\;\; t_i\in F\*, c\in\Qp\* 
\right \} \simeq (F\*)^{n/2}\times \Qp\*
$$
Un tore déployé maximal est $A=(\Qp\*)^{n/2}\times \Qp\*\subset
(F\*)^{n/2}\times \Qp\*= T$.

$G$ est déployé sur $F$,
$$
G_{F}\simeq \prod_{0\leq i < d-1} \Gl_n(F) \times F\*
$$
Le L-groupe connexe de $G$ est donné par 
$$
\widehat{G}=\prod_{0\leq i <d } \Gl_n(\C) \times \C\*
$$
sur lequel $\s=\bar{1}\in \Z/2d\Z$ opère via 
$$
\s.((g_0,\dots,g_{d-1}),c)=(( w {}^t g_{d-1}^{-1} w
,g_0,\dots,g_{d-2}),c\det g_{d-1})
$$
où $w= ( (-1)^{i-1} \delta_{i, n-j+1})_{i,j}$. 
$$
X_*(T)=\{ (x_{ij})_{i\in \Z/2d\Z, 1\leq j \leq n} \; | \;
x_{i,j}+x_{i+d,j+\frac{n}{2}}=c\in\Z \text{ constant } \}
$$
$\GG=\Z/2d\Z$ agissant par translations de coordonnée $i$ sur les
$(x_{ij})$

Le groupe de Weyl absolu est 
$$
W=\S_n^d
$$
et 
$$
X_*(T)/W=\{ (x_{ij})\;|\; \forall 0\leq i<d \;\;\forall j\leq j'\;\;
x_{ij}\geq x_{ij'} \}
$$
Le groupe de Weyl relatif sur $\Qp$ est
$W_{\Qp}=(\Z/2\Z)^{\frac{n}{2}}\rtimes\S_{n/2}$
et 
$$
X_*(A)/W_{\Qp} = \{ (y_i)_{1\leq i \leq \frac{n}{2}} \; |\; \forall
  i\leq i' \; y_i\geq y_{i'} \text{ et } \forall i \; y_i\geq 0 \}\times \Z  
$$
$G^{der}$ est simplement connexe et le cocentre $D=G/G^{der}$ est
isomorphe à 
$$
\{(x,t)\in F\*\times \Qp\* \; |\; N_{F/F_0}(x)=t^n\} \text{ via } 
g\mapsto ({\det}_F(g),c(g))
$$
$D\simeq \{\ker N_{F/F_0}:F\*\drt F_0\*\}\times\Qp\*$ puisque $n$ est pair et 
$N_{F/F_0}=t^n\ssi N_{F/F_0}(xt^{-n/2})=1$. On en déduit aisément que,
$$
X^*(Z(\widehat{G})^\Gamma)=X_*(D)_\Gamma =B(D)\simeq \Z/2\Z\times \Z 
$$
Nous normaliserons cet isomorphisme de la façon suivante : 
avec les notations précédentes 
$$
Z (\widehat{G})=\prod_{0\leq i <d} \C^\times \times \C^\times 
$$
et on vérifie aussitôt avec l'action donnée de $\s$ sur $\widehat{G}$
que l'on a :
$$
Z( \widehat{G})^\GG = \mu_2 \times \C^\times 
$$
où $\mu_2 \hookrightarrow (\C^\times)^d$ est le plongement diagonal et
l'isomorphisme est fixé en posant $\mu_2 =\Z/2\Z$. 
\\

Nous allons maintenant nous intéresser aux classes de $\s$ conjugaison
$b\in B(G)$ vérifiant $v_p(c(g))=1$ ce qui équivaut à dire que via
l'application $\kappa:B(G)\drt \Z/2\Z\times \Z$ la seconde composante
de $\kappa(b)$ est $1$. 

De telle classes de $\s$ conjugaison classifient des isocristaux
$(N,\Phi)$ munis d'une action de $F$ et d'une polarisation
$<.,.>:N\times N \drt \Qp(1)$ telle que si 
$$
N=\bigoplus_{i\in \Z/d\Z} N(i)
$$
est la décomposition $N(i)=\{n\in N| \forall x\in F \; x.n=\s^i(x) n
\}$, 
$\Phi:N(i)\drt N(i+1)$, alors 
$$
\forall i\neq j+d \; <N(i),N(j)>=0 \text{ et } <.,.> \text{ est
  parfait sur }
N(i)\times N(i+d) 
$$
\vspace{1mm}

{\bf Classes basiques : }

D'après Kottwitz, il y a un isomorphisme 
$$
\kappa:
B(G)_b\xrig{\; \sim \;} \Z/2\Z \times \Z 
$$
où la seconde composante est $b\mapsto v_p ( c(g))$. 
Il y a donc une bijection entre les classes basiques telles que
$v_p(c(g))=1$ et $\Z/2\Z$. Pour $b\in B(G)_b$, notons $\kappa (b)_1\in
\Z/2\Z$
la première composante de $\kappa (b)$. 

L'unique pente du $\Qp$ isocristal sous-jacent associée à une telle classe
est $\frac{1}{2}$.

Le groupe $J_b$ associé est une forme intérieure de $G$ qui est un
groupe de similitudes unitaires qui dépend de $\kappa (b)_1$. On a
en fait le lemme suivant valable pour n'importe quel $G$ et dont la vérification est facile :
\begin{lemm}
Soit $\kappa : B(G)_b \iso X^* ( Z (\widehat{G})^\GG)$. Via
l'identification $\widehat{G_{ad}}\simeq (\widehat{G})_{sc}$
considérons l'isomorphisme de Kottwitz 
$ H^1 ( F, G_{ad}) \iso X^*( Z ((\widehat{G})_{sc})^\GG)$.
Soit l'application $B(G)_b\ldrt H^1 (F,G_{ad})$ qui à une classe
basique $b$ associe la classe de cohomologie d'un cocyle définissant la
forme intérieure $J_b$. Le diagramme suivant est alors commutatif 
\begin{diagram}
B(G)_b & \rTo^\kappa & 
X^* ( Z (\widehat{G})^\GG) \\
\dTo & & \dTo \\
H^1 (F,G_{ad}) & \rTo^\sim  &  X^*( Z( (\widehat{G})_{sc})^\GG) 
\end{diagram}
où l'application verticale de droite est la composition d'un caractère
de $ Z (\widehat{G})^\GG $ avec l'application naturelle $
 Z ((\widehat{G})_{sc})^\GG \drt  Z (\widehat{G})^\GG$ induite par
l'application Galois équivariante $ Z((\widehat{G})_{sc}) \drt  Z
(\widehat{G})$, elle même déduite du morphisme $(\widehat{G})_{sc}
\drt \widehat{G}$.
\end{lemm}

On en déduit donc le corollaire suivant dans la cadre de notre groupe
unitaire :
\begin{coro}
$$ \forall b\in B(G)_b\;\;
J_b = \left \{ G \text{ si } \kappa (b)_1 =\bar{0} \atop G' \text{ si
} \kappa (b)_1 =\bar{1}  \right.
$$
où $G'$ désigne le groupe des similitudes unitaires en $n$ variables
non quasidéployé.  
\end{coro}

{\bf Classes quelconques : } 
Si $b\in B(G)$, $v_p(c(g))=1$, et si les pentes du $F$-isocristal 
sous-jacent (avec les notations précédentes cet isocristal est
$(N(0),\Phi^{2d})$ ) sont : 
$$
\l_1=\frac{d_1}{h_1} >\dots > \l_r=d > d-\l_{1} >\dots >d-\l_{r-1} 
$$
de multiplicités $m_1,\dots , m_{r-1},m_r,m_{r-1},\dots,m_1$. Avec les
notations précédentes pour $X_* (A) /W_{\Qp}$ on a  
$$
Newt(b)=(\underbrace{\frac{\l_1}{2d},\dots,\frac{\l_1}{2d}}_{m_1 h_1},\dots
,\underbrace{\frac{1}{2},\dots,\frac{1}{2}}_{m_r/2}) \times (1) \in X_*(A)_\Q/W_{\Qp}
$$
Si $m_r\neq 0$
le centralisateur du morphisme des pentes est le sous-groupe de Lévy 
\begin{eqnarray*}
M &=& \Res_{F/\Qp}(\Gl_{m_1 h_1})\times \dots\times \Res_{F/\Qp}
(\Gl_{m_{r-1}h_{r-1}})
\times GU(F;2m_r) \\
&=& \left \{ \left ( 
\begin{array}{ccccccc}
A_1 \\
& \ddots \\
& & A_{r-1} \\
& & & B \\
& & & & c {}^t A_1^{-*} \\
& & & & & \ddots \\
& & & & & & c {}^t A_{r-1}^{-*}
\end{array} \right ) \; |\; A_i\in\Gl,\; B\in GU\right \}
\end{eqnarray*}
où $GU(F;2m_r)$ est le groupe des similitudes unitaires quasidéployé en $2m_r$
variables. 
Si $m_r=0$,
$$
M=\Res_{F/\Qp}(\Gl_{m_1 h_1})\times \dots\times \Res_{F/\Qp}
(\Gl_{m_{r-1}h_{r-1}})\times \Qp\*
$$

\begin{center}
\setlength{\unitlength}{3947sp}%
\begingroup\makeatletter\ifx\SetFigFont\undefined%
\gdef\SetFigFont#1#2#3#4#5{%
  \reset@font\fontsize{#1}{#2pt}%
  \fontfamily{#3}\fontseries{#4}\fontshape{#5}%
  \selectfont}%
\fi\endgroup%
\begin{picture}(3024,1824)(5089,-5173)
\thinlines
\put(5101,-4861){\line( 1, 0){3000}}
\put(5101,-4861){\line( 0, 1){1500}}
\put(5101,-4861){\line( 1, 0){1500}}
\put(6601,-4861){\line( 1, 1){1500}}
\put(5101,-4861){\line( 4, 1){1200}}
\put(6301,-4561){\line( 2, 1){600}}
\put(6901,-4261){\line( 4, 3){1200}}
\multiput(6301,-4561)(0.00000,-8.95522){68}{\makebox(1.6667,11.6667){\SetFigFont{5}{6}{\rmdefault}{\mddefault}{\updefault}.}}
\multiput(6901,-4261)(0.00000,-9.00000){101}{\makebox(1.6667,11.6667){\SetFigFont{5}{6}{\rmdefault}{\mddefault}{\updefault}.}}
\multiput(6301,-5161)(109.09091,0.00000){6}{\line( 1, 0){ 54.545}}
\put(6901,-5161){\vector( 1, 0){0}}
\put(6301,-5161){\vector(-1, 0){0}}
\multiput(5101,-5161)(114.28571,0.00000){11}{\line( 1, 0){ 57.143}}
\put(6301,-5161){\vector( 1, 0){0}}
\put(5101,-5161){\vector(-1, 0){0}}
\multiput(6901,-5161)(114.28571,0.00000){11}{\line( 1, 0){ 57.143}}
\put(8101,-5161){\vector( 1, 0){0}}
\put(6901,-5161){\vector(-1, 0){0}}
\multiput(5101,-4861)(0.00000,-9.09091){34}{\makebox(1.6667,11.6667){\SetFigFont{5}{6}{\rmdefault}{\mddefault}{\updefault}.}}
\multiput(8101,-3361)(0.00000,-9.00000){201}{\makebox(1.6667,11.6667){\SetFigFont{5}{6}{\rmdefault}{\mddefault}{\updefault}.}}
\put(5701,-5086){\makebox(0,0)[b]{\smash{\SetFigFont{12}{14.4}{\rmdefault}{\mddefault}{\updefault}
(1)
}}}
\put(6601,-5086){\makebox(0,0)[b]{\smash{\SetFigFont{12}{14.4}{\rmdefault}{\mddefault}{\updefault}
(2)
}}}
\put(7501,-5086){\makebox(0,0)[b]{\smash{\SetFigFont{12}{14.4}{\rmdefault}{\mddefault}{\updefault}
(3)
}}}
\end{picture}
\begin{itemize}
\item (2) est autoduale de pente 1/2 et provient d'un groupe unitaire
\item (1) et (3) sont en dualité de pentes $\l, 1-\l$ et proviennent
  d'un $\Gl$  
\end{itemize}
\end{center}

$b$ est donc déterminé par les pentes $\l_1,\dots,\l_r$ les
multiplicités et par l'élément de $\Z/2\Z$ associé à la classe basique
de $GU(F;2m_r)$. On a
alors 
$$
J_b=\Gl_{m_1}(D_{\l_1})\times\dots\times \Gl_{m_{r-1}}(D_{\l_{r-1}})
\times H_b
$$
où $H_b$ est un groupe de similitudes unitaires en un nombre pair de
variables déployé ou non suivant l'invariant dans $\Z/2\Z$ (et où l'on
pose $H_b= \Qp^\times$ si $m_r =0$).
\\

{\bf Polygone de Hodge : } Soit $\mu\in X_*(T)$ tel que $\mu$ soit
donné par $\mu=(x_{ij})_{i\in\Z/2d\Z,1\leq j\leq n}$ où
$$\forall 1\leq
i<d \;\;\forall j\leq a_i\;\; x_{ij}=1,\; \forall j>a_i \;\; x_{ij}=0
$$
où les $a_0,\dots,a_{d-1}$ sont des entiers compris entre $1$ et $n$,
et où $c\circ \mu=1$ i.e. $x_{i+d,j+\frac{n}{2}}=1-x_{i,j}$.

$\mu_1=\widehat{\mu}_{Z(\widehat{G})^\Gamma}\in \Z/2\Z\times \Z$ est
donné par 
$$
\mu_1=(\sum_{0\leq i <d-1} a_i \text{ mod } 2, 1 ) 
$$
et $\mu_2\in X_*(A)_\Q=\Q^{n/2}$ est donné par 
$$
\mu_2=\frac{1}{d} \sum_{0\leq i <d } l_i
$$
où
\begin{itemize}
\item
  $\dpt{
l_i=(\underbrace{1,\dots,1}_{a_i},\frac{1}{2},\dots,\frac{1}{2})}$
  si $a_i\leq \frac{n}{2}$
\item $\dpt{
l_i=(\underbrace{\frac{1}{2},\dots,\frac{1}{2}}_{n/2-a_i},1,\dots,1)
}$ si $a_i>\frac{n}{2}$
\end{itemize}

Et donc finalement : 
$$
 B(G,\mu)   =\left \{
 b\; | \;v_p(c(g))=1,\,
 \kappa(b)_1=\sum_{0\leq i <d} a_i \text{ mod }
  2\text{  et } 
 \left  (  
 \underbrace{ \l_1 ,\dots,\l_1}_{m_1 h_1} ,\dots,\underbrace{\frac{1}{2},
\dots,\frac{1}{2}
}_{m_r/2} \right )
\leq \frac{1}{d} \sum_{0\leq i< d} l_i \right \}
$$
au sens où le membre de droite doit d'abord être ramené dans la
chambre de Weyl positive de $X_*(A)_\Q$ avant d'être comparé.

L'unique classe basique de $B(G,\mu)$ est la classe basique $b$ telle
que $v_p(c(g))=1$ et $\kappa (b)_1\equiv\sum_i a_i \in \Z/2\Z$.
\\

\begin{rema} Dans la description de $B(G,\mu)$ les points terminaux
des polygones n'interviennent pas. En effet,
 la condition que les point terminaux des polygones de
Hodge et de Newton soient égaux est automatique puisqu'un polygone
symétrique  a pour dimension le double de sa
hauteur.
\end{rema}

\section{Groupes unitaires en un nombre impair de variables}

Soit $F|\Qp$ comme précédemment et 
$$
G=\{ g\in \Gl_n (F)  \; |\; {}^tg^* J g=c(g) J,\;  c(g)\in\Qp\* \}
$$
où $J=\left ( \begin{array}{ccc}
0 & 0 & I_{\frac{n-1}{2}} \\
0 & 1 & 0 \\
I_{\frac{n-1}{2}} & 0 & 0 \\
\end{array} \right )
$

Un tore maximal défini sur $\Qp$ est 
$$
T=diag(t_1,\dots,t_{\frac{n-1}{2}},c (t_1^{-1})^*,\dots,c (t_n^{-1})^*)
$$
où $t_i,x\in F\*, xx^*=c\in\Qp\*$.
$T\simeq (F\*)^{\frac{n-1}{2}}\times \{x\in F\* \; | \;
N_{F/F_0}(x)\in\Qp\*\}$.

Un tore déployé maximal  est
$A=(\Qp\*)^{\frac{n-1}{2}}\times\Qp\*\subset T$. 

$\widehat{G}$ et ${}^L G$ sont de la même forme que précédemment.

Notons 
$$
X_*(T)=\{ (x_{ij})_{i\in \Z/2d\Z,1\leq j\leq n}\; |\; \forall j\leq
\frac{n-1}{2}\;\; x_{i,j}+x_{i+d,j+\frac{n+1}{2}}=c,
x_{i,\frac{n+1}{2}}+x_{i+d,\frac{n+1}{2}}=c \}
$$
$W=\S_n^d$ et 
$$
X_*(T)/W = \{ (x_{ij}) \; |\; \forall 0\leq i <d-1\;\; \forall j\leq
j'\; 
\;
x_{ij}\geq x_{ij'} \}
$$
Le groups de Weyl relatif s'identifie quant à lui  à
$W_{\Qp}\simeq (\Z/2\Z)^{\frac{n-1}{2}}\rtimes\S_{\frac{n-1}{2}}
$ et 
$$
X_*(A)/W_{Qp} =\{ (y_j)_{1\leq j\leq \frac{n+1}{2}}\; | \; 
\forall 1\leq j\leq j' <\frac{n+1}{2}\;\; y_j\geq y_{j'}, \; y_1\geq 0\}
$$
\vspace{1mm}

Le cocentre $D$ est isomorphe comme précédemment à 
$$
\{(x,t)\in F\*\times \Qp\*\; |\; N_{F/F_0}(x)=t^n \}
$$
qui est isomorphe à 
$$
\{x\in F\*\;|\; N_{F/F_0}(x)\in \Qp\* \}\times \Qp\*
$$
On en déduit aisément que 
$$
X^*(Z(\widehat{G})^\Gamma)\simeq \Z
$$
et que $\forall b\in B(G)\; \kappa(b)=v_p(c(g))$. 
\\

{\bf Classes basiques : } 
Soit l'isomorphisme de Kottwitz 
$$
\kappa :
B(G)_b\xrig{\;\sim\;}\Z
$$
On en déduit  qu'il existe une unique classe basique $b$ telle que
$v_p(c(g))=1$. Et on a alors, puisque $G$ ne possède pas d'autre forme
intérieure que lui même,  
$$
J_b=G
$$

{\bf Classes quelconques : } La discussion pour les pentes, $M$ et
$J_b$ est la même que précédemment à ceci près que la classe basique
du groupe unitaire est unique. On a 
$$
Newt(b)=\left (
  \underbrace{\frac{\l_1}{2d},\dots,\frac{\l_1}{2d}}_{m_1 h_1},\dots,
\underbrace{\frac{1}{2},\dots,\frac{1}{2}}_{\frac{m_r-1}{2}},\frac{1}{2}\right )
$$
\\

{\bf Polygone de Hodge : } Supposons donné $\mu$ par
$\mu=(\mu_{ij})_{i,j}$ tel que 
$$\forall 1\leq
i<d \;\forall j\leq a_i\; x_{ij}=1,\; \forall j>a_i \; x_{ij}=0
$$
et où $c\circ \mu=1$ i.e. $$x_{i,j}+x_{i+d,j+\frac{n+1}{2}}=1,
x_{i,\frac{n+1}{2}}+x_{i+d,\frac{n+1}{2}}=1$$

$\mu_1=\widehat{\mu}_{|Z(\widehat{G})^\Gamma}$ est donné par
$\mu_1=1$ et $\mu_2\in X_*(A)_\Q = \Q^{\frac{n+1}{2}}$ est donné par 
$$
\mu_2=\frac{1}{d} \sum_{0\leq i <d} l_i 
$$
où \begin{itemize}
\item
  $l_i=(\underbrace{1,\dots,1}_{a_i},\frac{1}{2},\dots,\frac{1}{2},\frac{1}{2} )$
si $a_i\leq \frac{n-1}{2}$ 
\item $l_i=(1,\dots,1,\frac{1}{2})$ si $a_i=\frac{n+1}{2}$
\item
  $l_i=(\underbrace{\frac{1}{2},\dots,\frac{1}{2}}_{a_i-\frac{n+1}{2}},1,\dots,1,\frac{1}{2})$ si $a_i>\frac{n+1}{2}$ 
\end{itemize}
Posons $l'_i=l_i$ moins le $\frac{1}{2}$ à droite, alors 
$$
 B(G,\mu) = \left \{ b\;| \; v_p(c(g))=1\text{ et } 
\left (
  \underbrace{\frac{\l_1}{2d},\dots,\frac{\l_1}{2d}}_{m_1 h_1},\dots,
\underbrace{\frac{1}{2},\dots,\frac{1}{2}}_{\frac{m_r-1}{2}}\right )
 \leq \frac{1}{d}\sum_{0\leq i< d} l'_i \right \}
$$
\\

\chapter{Espaces adiques et espaces analytiques}\label{an_ad}
\section{Différents espaces}

Soit $k$ un corps valué complet muni d'une valuation de rang $1$. 

Rappelons qu'il y a  un foncteur pleinement fidèle (\citeg{Berk1}
 1.6.1)
$$
X\longmapsto s(X)
$$ 
de la catégorie des espaces analytiques de Berkovich Hausdorff (i.e. $|X|$ est
séparé) dans la catégorie des variétés rigides quasiséparées.

Si $\AA$ est une algèbre affinoïde, $s(\mathcal{M}(\AA))=Spm(\AA)$. 

Rappelons également (\citeg{Hu1} 1.1.11) qu'il y a une équivalence de
catégorie 
$$
Y\longmapsto r(Y)
$$
entre la catégorie des variétés rigides quasiséparées $Y$ sur $k$ et la
catégorie des espaces adiques quasiséparés localement de type fini sur
$\Spa(k,k^0)$. Nous noterons 
$$
X\longmapsto X^{rig}
$$
la composée $X\mapsto r(s(X))$ qui est pleinement fidèle.

Rappelons la proposition suivante 
\begin{Fait}(\citeg{Hu1} 8.3.3)
  L'image essentielle de $X\mapsto X^{rig}$ est formée de la catégorie
  des espaces adiques quasiséparés localement de type fini tendus
  (``taut'') sur $\Spa(k,k^0)$.
\end{Fait}
La condition d'être tendu portant uniquement sur l'espace topologie de
l'espace adique (\citeg{Hu1} 5.1.2).
\\

Supposons maintenant la valuation de $k$ discrète. 

A $\X$ un schéma formel localement formellement de type fini sur
$\spf(k^0)$ est associé un $k$-espace analytique paracompact
$\X^{an}$. Il est muni d'un morphisme de spécialisation 
$$
sp:\X^{an}\ldrt \X
$$
Si $Y\subset \X_s$ est un sous-schéma localement fermé de la fibre
spéciale de $\X$ on a un isomorphisme canonique 
$$
(\X_{/Y}^\wedge)^{an}\iso sp^{-1}(Y)
$$
où $sp^{-1}(Y)$ est un domaine analytique localement fermé dans
$\X^{an}$. 

Si $Y$ est ouvert, $sp^{-1}(Y)$ est fermé ($Y=D(f)\impl
sp^{-1}(Y)=\{|f|=1\}$.

Si $Y$ est fermé, $sp^{-1}(Y)$ est ouvert ($Y=V(f)\impl
sp^{-1}(Y)=\{|f|<1\}$. 
\\

A $\X$ est également associé fonctoriellement
un espace adique $t(\X)$. $t$ étant défini localement par 
$t(\spf(A))=\Spa( A,A)$. On peut alors considérer l'espace adique 
$$
d(\X)=t(\X)_a
$$
l'ouvert de $t(\X)$ formé des points analytiques i.e. des valuations
ayant un support non ouvert. Si $I$ est un idéal de définition de $A$, 
$$
d(\spf(A))=\Spa(A,A)\setminus V(I)
$$
On a alors un morphisme de spécialisation d'espaces annelés 
$$
(d(\X),\O_{d(\X)})\xrig{\l} (\X,\O_{\X})
$$
l'application au niveau des espaces topologiques étant celle qui
associe à un point $x\in d(\X)$ le support d'une 
 spécialisation de $x$ dans
$t(\X)$ formée d'un point non analytique.

Si $\X/\spf(k^0)$ est adique c'est la bonne définition de la fibre
générique de $\X$ au sens des espaces adiques.

Dans le cas général
$$
\widetilde{d}(\X)=t(\X)\setminus V(\pi) \subset d(\X)
$$
est un espace adique ouvert dans $d(\X)$ muni du morphisme de
spécialisation obtenu par restriction de $\l$ : 
$$
\widetilde{\l}:\widetilde{d}(\X)\ldrt \X
$$
C'est la bonne définition de la fibre générique au sens de Berthelot
et on a des isomorphismes canoniques commutant aux morphismes de
spécialisation 
$$
(\X^{an})^{rig}\simeq \widetilde{d}(\X)
$$
Nous noterons donc $\widetilde{d}(\X)=\X^{rig}$. 

Si maintenant $Y\subset \X_s$ et un sous-schéma fermé, 
$\widetilde{\l}^{-1}(Y)\subset \X^{rig}$ est un ensemble fermé constructible
et on a un isomorphisme canonique 
$$
(\X_{/Y}^\wedge)^{rig}\simeq (\widetilde{\l}^{-1}(Y))^0
$$
(l'intérieur de $\widetilde{\l}^{-1}(Y)$).
Ainsi $(\X_{/Y}^\wedge)^{rig}$ s'identifie à un ouvert de $\X^{rig}$
qui est l'ouvert $sp^{-1}(Y)^{rig}$.

Si $\mathcal{Y}\subset \X$ est ouvert, $\mathcal{Y}^{rig}\subset
\X^{rig}$ est l'ouvert $\widetilde{\l}^{-1}(Y)$.
\\

\begin{exem}
 Si $\X=\spf(k^0<X_1,\dots,X_n>)$, $\X^{an}=\BB^n$ la
boule fermée analytique, $\X^{rig}=\BB^n$ la boule fermée adique
(qui est ouverte dans $\A^n$) 
. Si
$Y=V(0)$, $\X_{/Y}^\wedge=\spf(k^0\llbracket X_1,\dots,X_n\rrbracket)$, 
$sp^{-1}(Y)=\tilde{\BB}^n$ la boule ouverte analytique
 obtenue par recollement
d'une union croissante de boules fermées (ce n'est pas l'intérieur de
la boule fermée analytique pour la topologie de Berkovich ! ). 

$\widetilde{\l}^{-1}(Y)=\{|X_i|<1\}\subset \BB^n$ un fermé (ce n'est
pas un espace adique, c'est un espace pseudo-adique) et
$(\X_{/Y}^\wedge)^{rig}$ est la boule ouverte adique obtenue comme
union croissante des boules ouvertes de rayon strictement inférieur à
$1$. $\widetilde{\l}^{-1}(Y)\setminus \widetilde{\l}^{-1}(Y)^0$ est
constitué de valuations de rang $\geq 2$ qui ne sont pas dans l'espace
de Berkovich.
\\

On remarquera que la stratification associée à un fermé $Y\subset
\X_s$ de la fibre spéciale est plus compliquée dans le cas de l'espace
adique :
\begin{itemize}
\item La stratification de $\X^{an}$ associée est
$
sp^{-1}(Y) \subset \X^{an}$ où $sp^{-1}(Y)$ est un domaine
analytique ouvert et $\X^{an}\setminus sp^{-1}(Y)$ est un domaine
analytique fermé. Les deux strates sont donc des espaces analytiques
\item La stratification de $\X^{rig}$ associée est 
$$\widetilde{\l}^{-1}(Y)^0\subset 
\widetilde{\l}^{-1}(Y)\subset \X^{an}$$ où $\widetilde{\l}^{-1}(Y)$
est un espace pseudo-adique (i.e. un germe d'espace adique) et le
complémentaire est un espace adique ouvert dans $\X^{rig}$.
\end{itemize}
\end{exem}

\section{Différents morphismes étales}

Si $j:\BB^n\hookrightarrow \Aff^n$ désigne l'inclusion de la boule fermée
dans l'espace affine, $j$ est une immersion ouverte au sens des
espaces adiques et est donc étale. L'espace adique $\BB^n$ est lisse
sur $\Spa(k,k^0)$. Cependant, $j$ n'est pas étale au sens des espaces
analytiques car $\partial(\BB^n /\Aff^n)\neq\emptyset$. Ainsi, $\BB^n$
n'est pas lisse comme espace analytique sur $k$.
\\

Rappelons :
\begin{defi}\citeg{Berk2} 
  Soit $f:X\drt Y$ un morphisme d'espace analytique. $f$ est
  quasi-étale si tout point $x\in X$ possède un voisinage qui est un
  domaine analytique fermé étale au dessus de $Y$.
\end{defi}

\begin{Fait}(\citeg{Berk2} 2.3, resp. \citeg{Hu1} 3.5.1) 
  Si $\X\drt\mathcal{Y}$ est un morphisme de schémas formels,
  $\X^{an}\drt \mathcal{Y}^{an}$ est quasi-étale  et $\X^{rig}\drt
  \mathcal{Y}^{rig}$ est étale.
\end{Fait}
Par exemple, si $\X\hookrightarrow \mathcal{Y}$ est une immersion
ouverte, donc étale, 
elle induit une inclusion d'un domaine analytique fermé
$\X^{an}$ dans $\mathcal{Y}^{an}$ qui n'est pas étale mais
quasi-étale. Par contre elle induit une immersion ouverte d'espaces adiques.
\\

Le lien entre les deux notions est le suivant :
\begin{Fait} (\citeg{Hu1} 8.3.4)
  Si $f:X\drt Y$ est quasi-étale, $f$ est étale ssi
  $f^{rig}:X^{rig}\drt Y^{rig}$ est partiellement propre.
\end{Fait}
Rappelons qu'un morphisme d'espaces adiques est partiellement propre
s'il vérifie un critère valuatif de propreté (et est localement de
type fini) et qu'un morphisme est propre ssi il est partiellement propre et
quasicompact. 
\\

Ainsi, un domaine analytique ouvert induit une immersion ouverte
d'espace adique qui est partiellement propre et à un ouvert analytique
d'un espace analytique de bord vide sur $k$ est du point de vue adique partiellement
propre (c'est en particulier le cas d'un ouvert analytique de l'espace
analytique associé à un schéma)
\\

Les deux exemples typique d'espace partiellement propre sont les suivants :
\begin{Fait}
  Si $X/k$ est un schéma localement de type fini,
  $X^{rig}/\Spa(k,k^0)$ est partiellement propre.
\end{Fait}

\begin{Fait}
Si $\X/\spf (k^0)$ est un schéma formel localement formellement de
type fini séparé, 
$$
\X_s \text{ propre } \limpl \X^{rig} \text{ partiellement propre}
$$
Lorsque $\X$ est adique sur $k^0$, il s'agit même d'une équivalence.
\end{Fait}

\begin{exem}
Le tube au dessus d'un fermé propre de la fibre spéciale d'un schéma
formel est partiellement propre.
\end{exem}

Du point de vue des espaces adiques de Huber 

\begin{Fait}
  Soit $f:X\drt Y$ un morphisme d'espaces analytiques. $f^{rig}$ est
  partiellement propre ssi $\partial (X/Y)=\emptyset$. 
\end{Fait}

Si $X=\mathcal{M}(\AA)\xrig{\;\; f\;\;} Y=\mathcal{M}(\mathcal{B})$, $f$ est
partiellement propre s'il existe une immersion fermée 
\begin{diagram}
  X & \rInto & \BB^n(0,\e) \times Y \\
& \rdTo_f   & \dTo \\
 & & Y
\end{diagram}
pour $0<\e< 1$ telle que $X$ soit Zariski fermé dans $\BB^n
(0,1)\times Y$.

\chapter{Cohomologie $\ell$-adique des espaces analytiques}
\label{coho_l} 

Dans cette section nous donnons des démonstration des résultats non
publiés de \citeg{Berk4} en suivant \citeg{Hu2}.

\section{Propriétés générales}

Soit $k$ un corps valué complet et $X$ un $k$-espace analytique
Hausdorff (i.e. $|X|$ est séparé). Soit $\La$ un anneau de valuation
discrète complet d'idéal maximal $\mm$ et de caractéristique
résiduelle $\ell\neq p$ où $p=char(k^0/k^{00})$.
\\

\begin{defi}
  Un $\La_\bullet$ faisceau ou encore un faisceau $\La$-adique sur $X_{\et}$ 
 est un système projectif
  $(\F_n)_{n\in \N}$ de faisceaux en $\La$-modules 
 sur $X_{\et}$ vérifiant $\forall n\;\;
  \mm^n\F_n=0$. Nous noterons $\Lpf_{/X^{\et}}$ la catégorie associée.
\\
On définit de même $\La_\bullet -mod$ la catégorie des $\La_\bullet$-modules
\\
On notera également $\Lf_{/X^{\et}}$ la catégorie des
faisceaux de $\La$-modules.
\end{defi}
\vspace{4mm}
$\Lpf_{/X^{\et}}$ est  une catégorie abélienne
$\La$ linéaire possédant suffisamment d'injectifs.
\\

On a un foncteur exact à gauche 
\begin{eqnarray*}
  \pi_*:\Lpf\Xet &\ldrt & \Lf\Xet \\
 (\F_n)_n &\longmapsto & \limp \F_n \\
\end{eqnarray*}

Si $(\F_n)_n \in \Lpf\Xet $ et $X$ n'est pas compact il n'est pas
raisonnable de poser 
$\GG_c(X,(\F_n)_n)$ comme étant égal à $\limp \GG_c (X,\F_n)$ car si 
$(s_n)_n\in \limp \GG_c(X,\F_n)$ il se peut que l'adhérence de 
 $\cup_n supp(s_n)$ ne
soit pas compact.

Nous adopterons donc la définition suivante :
\begin{defi}
  \begin{eqnarray*}
    \GG_!:\Lf\Xet & \ldrt & \Ab \\
  \F &\longmapsto & \GG_c(X,\F) 
  \end{eqnarray*}
 $$\GG_c=\GG_!\circ \pi_*$$ 
On a donc $\GG_c(X,(\F_n)_n)=\{(s_n)_n
\in \limp \GG(X,\F_n) \; |\;\overline{ 
\cup_n supp(s_n)}$ est compact $\}$. 

On posera $$H^p_c(X,(\F_n)_n)=R^p(\GG_c(X,-))((\F_n)_n)$$  
\end{defi}

\begin{rema}
  Si $X$ est compact $\GG_c(X,(\F_n)_n)=\underset{n}{\limp}
\GG(X,\F_n)$. 
\end{rema}

\begin{lemm}\label{flasque}
  Soit $\F$ un faisceau abélien sur $X_{\et}$. Si $\forall x\in X\;\;
  \F_x$ est un $\Gal(\mathcal{K}(x)^s|\mathcal{K}(x))$-module acyclique
  et si ${\F_{|}}_{|X|}$ est un faisceau mou, $\F$
  est $\GG_c$ acyclique.
\end{lemm}
\begin{proof} Il suffit de considérer le morphisme de site $f:X_{\et}\drt
|X|$. La première hypothèses montre que $\forall p>0\;\; R^p f_* (\F)
=0$ car ses fibres en tous les points de $|X|$ sont nulles
(cf. \citeg{Berk1}, 4.2.4)
. La suite
spectrale de Leray associée dégénère donc en
$H^p_c(X_{\et},\F)\simeq H^p_c(|X|,\F_{|X|})$ mais  $|X|$ est localement
paracompact donc ${\F_{|}}_{|X|}$ étant mou il est $\GG_c(|X|,-)$ 
acyclique.
\end{proof}

L'intérêt du lemme qui suit est qu'il permet de calculer la
cohomologie à support  compact d'un faisceau
$\La$-adique comme l'hypercohomologie à support compact d'un 
complexe de faisceaux de $\La$-modules.

\begin{lemm}\label{Rpp}
  $\R\GG_c=\R\GG_!\circ \R \pi_*$
\end{lemm}

\begin{proof}
 En suivant la preuve du i) du lemme 2.3
 de  \citeg{Hu2}  on voit qu'il
suffit de montrer que si $\forall n\;$
$\I_n$ est un injectif de $\Lmf\Xet$ alors $\prod_{k=0}^\infty \I_k$
est $\GG_!$ acyclique.
Pour cela nous allons appliquer le lemme \ref{flasque}.
 
$$\forall x\in X\;\;\;\;
(\prod_k \I_k)_x\simeq \left ( \prod (\I_k)_x\right )^{disc} 
$$
comme $\Gal(\mathcal{K}(x)^s | \mathcal{K}(x))$-module galoisien
discret, où si $M$ est un module galoisien $M^{disc}$ désigne sa partie
discrète. Or, si $(M_k)_{k\in K}$ sont des $\Gal(\mathcal{K}(x)^s| \mathcal{K}(x))$-modules
discrets 
$$
H^q \left
(\Gal(\mathcal{K}(x)^s |\mathcal{K}(x) ),\left ( \prod_k M_k \right
)^{disc}
\right )\simeq \prod_k H^q( \Gal(\mathcal{K}(x)^s |\mathcal{K}(x) ),M_k)
$$
$\I_k$ étant injectif comme faisceau de $\La/\mm^n$-modules, 
$\forall x \;\I_{k,x}$ est un $\Gal (\mathcal{K}(x)^s|\mathcal{K}(x))$-module galoisien 
acyclique (\citeg{Berk1}, 4.2.5)
. On en déduit donc que $\left (\prod_k \I_k \right )_x$
en est également un.

De plus,  ${(\prod_k \I_k)_{|}}_{|X|}$ est mou car  il est
flasque et $|X|$ est localement paracompact.
\end{proof}
\vspace{3mm}

Dégageons quelques propriétés simples de la cohomologie $\La$-adique :

\begin{lemm}
  Soit $\U(X)$ un ensemble d'ouverts de $X$ tels que $\forall U\in
  \U(X), \overline{U}$ est compact, $\forall U_1,U_2\in \U(X)
  \;\exists U_3\in\U(X)\;\; U_1\cup U_2\subset U_3$ et $\bigcup_{U\in
  \U(X)} U=X$. On a alors,
$$
\GG_c(X,(\F_n)_n)=\underset{U\in \U(X)}{\;\limi\;} \underset{n}{
\limp} \GG_c(U,\F_n)
$$
\end{lemm}

\begin{proof}
 Cela résulte de ce que tout compact de $X$ est contenu
dans un élément de $\U(X)$ et réciproquement.
\end{proof}

Nous choisirons généralement $\U(X)=\{$ ouverts $U$ de $X$, $\overline{U}$ est
compact $\}$ ou bien $\U(X)=\{$ ouverts distingués de $X\}$ où
rappelons 
\begin{defi}
  \citeg{Berk4} 
Un ouvert $U$ de $X$ est distingué s'il peut s'écrire 
$U=V_1\setminus V_2$ où $V_1$ et $V_2$ sont des domaines analytiques compacts.
\end{defi}
\vspace{5mm}

\begin{lemm}
  \begin{enumerate}
\item Si $i:Y\hookrightarrow X$ est une immersion fermée d'espaces
  analytiques
ou une immersion d'un domaine analytique fermé, il y a un
isomorphisme canonique 
$$
H^p_c(Y,(\F_n)_n)\iso H^p_c(X,(i_*\F_n)_n)
$$ 
\item Si $j:U\hookrightarrow X$ est une immersion ouverte, pour tout
  $(\F_n)_n\in \Lpf\Xet$ 
 il y  a une application naturelle 
$$
H^p_c(U,(\F_{n|U})_n)\ldrt H^p_c(X,(\F_n)_n)
$$
et $\forall (\F_n)_n\in \Lpf_{/U_{\text{èt}}}$ une application
naturelle 
$$
H^p_c(U,(\F_n)_n)\ldrt H^p_c(X,(j_!\F_n)_n)
$$
\item Si $\mathbb{V}$ est un recouvrement ouvert de $X$ tel que
  $\forall V_1,V_2\in \mathbb{V}\;\;\exists V_3\in \mathbb{V}\;\;
  V_1\cup V_2\subset V_3$ 
$$
H^p_c(X,(\F_n)_n)=\underset{V\in \mathbb{V}}{\limi} H^p_c(V,(\F_{n|V})_n)
$$
\end{enumerate}
\end{lemm}

\begin{proof}
 1) Soit le foncteur
\begin{eqnarray*}
  i_*:\Lpf_{/Y^{\text{ét}}} & \ldrt & \Lpf\Xet \\
 (\F_n)_n & \longmapsto & (i_*\F_n)_n 
\end{eqnarray*}
il est exact et possède un adjoint à gauche $(\F_n)_n\longmapsto
(i^*\F_n)_n$. Il envoie donc les injectifs sur des injectifs. De plus, 
$|Y|$ étant fermé $\GG_c(Y,-)=\GG_c(X,-)\circ i_*$. D'où le 1).
\\

2) Soit le foncteur  
\begin{eqnarray*}
  j^*:\Lpf\Xet &\ldrt & \Lpf_{/U^{\text{èt}}} \\
(\F_n)_n &\longmapsto & (j^*\F_n)_n 
\end{eqnarray*}
qui est exact et possède un adjoint à gauche $(\F_n)_n \mapsto
(j_!\F_n)_n$. On a donc $R^+(\GG_c(U,-)\circ j^*)=R^+\GG_c(U,-)\circ
j^*$. On a de plus une application naturelle 
$$
\GG_c(U,-)\circ j^*\ldrt \GG_c(X,-) 
$$
ce qui conclu la démonstration de la première partie.

Si maintenant $(\F_n)_n\in \Lpf_{/U^{\text{èt}}}$ il y a une
application canonique d'adjonction 
$$
(\F_n)_n\ldrt (j^*j_!\F_n)_n
$$
L'application cherchée est alors la composée 
$$
H^p_c(U,(\F_n)_n)\drt H^p_c(U,(j^*j_!\F_n)_n)\drt H^p_c(X,(j_!\F_n)_n)
$$
où la dernière application est celle définie auparavant.
\\

Le 3) est clair car $\underset{V\in\mathbb{V}}{\limi}\circ
\underset{U\in\U(X)}{\limi} =\underset{U\in \U(X)}{\limi}$ et les
limites étant filtrées elles sont exactes.
\end{proof}
7
\begin{rema} En général l'application 
$$
H^p_c(U,(\F_n)_n)\ldrt H^p_c(U,(j_!\F_n)_n)
$$
n'est pas bijective. Cela est du au fait qu'en général $j$ ne possède
pas de bonnes propriétés de finitude. Déjà dans le cas algébrique,
pour $X\xrig{ f } Y$ un morphisme de type fini, pour montrer que $(R^+
f_!)((\F_n)_n)=(R^+ f_!\F_n)_n$ pour $(\F_n)_n$ constructible et
$R^+\GG_c\circ R^+f_!=R^+\GG_c$ il faut utiliser les propriétés de finitude
de $f_!$: $f_!($constructible$)$ est constructible.

Par exemple, en général $H^p_c(\breve{\F}^{ad},\Ql)\neq
H^p_c(\breve{\F},j_!\Ql)$ (cf. \citeg{Hu2}, exemple 2.7 )
\end{rema}

Comme dans le lemme 2.3 de \citeg{Hu2} on a une factorisation 
$R^+\GG_c=\underset{\U(X)}{\limi}\circ \text{Ind}(R^+\rho)\circ R^+\s$
de laquelle on déduit :
\vspace{1mm}

\begin{prop}
  $\forall (\F_n)_n\in \Lpf\Xet$ il y a des suites exactes : $\forall
  p\in \N$ 
$$
0\ldrt \underset{U\in \U(X)}{\limi} \underset{n}{\limp}^{\!\!\!\! 1}\;\;
H^{p-1}_c(U,\F_n) 
\ldrt H^p_c(X,(\F_n)_n) \ldrt \underset{U\in
  \U(X)}{\limi}\underset{n}{\limp} H^p_c(X,\F_n) \ldrt 0
$$
\end{prop}

\begin{coro} \label{corlim}
  \begin{enumerate} \item[] 
  \item[a)] Si $\forall n\; \forall p \;\;\; \forall U\in \U(X)\;\;\;
    H^p_c(U,\F_n)$ est un $\La/\mm^n$-module de type fini 
$$
\forall p\in \N \;\;\; H^p_c(X,(\F_n)_n)\iso \underset{U\in
  \U(X)}{\limi}
\underset{n\in \N}{\limp} H^p_c(U,\F_n)
$$
\item[b)] $\forall p>2 \dim (X)+\text{cd}_\ell (k)+1 \;\; \forall
  (\F_n)_n\in \Lpf\Xet \;\; H^p_c(X,(\F_n)_n)=0 $ 
et sous l'hypothèse du a) cela est vrai pour tout $p>2 \dim (X)
+\text{cd}_\ell (k)$
  \end{enumerate}
\end{coro}
\vspace{1mm}

\begin{exem} Si $X=\widetilde{\BB}^n_k$ est la boule ouverte de
dimension $n$, $\U(X)=\{\widetilde{\BB}^n_k(0,\e)\; |\; 0<\e<1\; \}$
les hypothèses du a) sont vérifiées pour $\Zl=(\Z/\ell^n \Z)_n$  et 
\begin{eqnarray*}
  H^p_c(\widetilde{\BB}^n_k,\Zl) &=& \underset{\e\drt 1}{\limi}
  \underset{
n\in \N}{\limp} H^p_c(\widetilde{\BB}^n_k(0,\e),\Z/\ell^n\Z) \\
&=&  \left \lbrace \begin{array}{l}
  0  \text{ si } p\neq 2n \\
  \Zl \text{ si } p=2n \end{array}\right.
\end{eqnarray*}
\end{exem}

\section{Cohomologie $\ell$-adique des espaces analytiques
  quasi-algébriques}

Les hypothèses faites ici : $X$ quasi algébrique et $\F$ localement
constant ne sont là que pour palier à l'existence d'une théorie des
faisceaux constructibles pour les espaces analytiques comme c'est le
cas pour les espaces adiques (\citeg{Hu1},\citeg{Hu2}).
\\

\begin{defi}\label{def_quasi_alg}
  (\citeg{Berk2},\citeg{Berk4}) $X$ est quasi-algébrique si $\forall x\in X$
  $\;\exists V_1,\dots,V_n$ des domaines affinoïdes tels que $x\in
  \bigcap_{i=1}^n V_i, \bigcup_{i=1}^n V_i$ est un voisinage de $x$ et
  tout $V_i$ est quasi-étale au dessus d'un domaine affinoïde de
  l'espace analytique associé à un schéma de type fini sur $k$. 
\end{defi}

\begin{rema}
  D'après \citeg{Berk4} cela implique que l'on peut trouver des $V_i$
  qui sont eux mêmes des domaines affinoïdes dans l'espace analytique
  associé à un schéma de type fini.
\end{rema}

\begin{rema}
  Tout espace analytique lisse sur $k$ est quasi-algébrique.
\end{rema}
\vspace{4mm}

\begin{defi}
  $(\F_n)_n\in\Lpf\Xet$ est localement constant si 
$$
\forall m\geq n \;\;\; \F_m\otimes \La/\mm^n\iso \F_n
$$
et $\forall n\;\; \F_n$ est localement constant fini 
sur $X_{\text{ét}}$ comme $\La/\mm^n$-module.
\end{defi}
\vspace{4mm}

\begin{lemm}
  Soit $X$ quasi-algébrique et $U$ un ouvert distingué de $X$ ou un
  domaine analytique compact. $\forall \F\in \Lmf\Xet$ localement
  constant fini, $\forall p\in\N\;\;\; H^p_c(U,\F)$ est un
  $\La/\mm^n$-module de type fini. 
\end{lemm}

\begin{proof}
 Si $U$ est un domaine compact c'est le corollaire 5.6 de
\citeg{Berk2}.

Si $U=V_1\setminus V_2$ est un ouvert où $V_1$ et $V_2$ sont des
domaines compacts cela résulte de la suite exacte longue de
cohomologie
\begin{diagram}
 & \rTo & H^p_c(U,\F) & \rTo & H^p( V_1,\F) & \rTo & H^p (V_1\cap
  V_2 ,\F) & \rTo  &
\end{diagram}
et de l'assertion pour un domaine compact.
\end{proof}

\begin{coro}
Si $X$ est quasi-algébrique et $(\F_n)_n$ est localement constant 
$$
H^p_c(X,(\F_n)_n)=\underset{U\text{
    distingué}}{\limi}\;\;\underset{n\in \N}{\limp}\;\;\; H^p_c(U,\F_n)
$$
En particulier si $X$ est compact 
$$
H^p(X,(\F_n)_n)=\underset{n\in\N}{\limp} H^p(X,\F_n)
$$
\end{coro}
\begin{proof}
 Cela résulte du lemme précédent et du corollaire
\ref{corlim}.
\end{proof}

\begin{prop} Supposons $cd_\ell (k)<+\infty$. 
  Soit $U$ un ouvert distingué ou un domaine compact d'un espace
  quasi-algébrique X, $\F$ localement constant sur $X$. 
 $\forall p\in\N\;\;\; 
(H^p_c(U,\F_n))_n$ est AR $\mm$-adique et $\underset{n\in\N}{\limp}
  H^p_c(U,\F_n)$ est un $\La$-module de type fini.
\end{prop}
\begin{proof}
 $H^p_c(U,\F_1)$ étant de type fini, la première partie de
la proposition entraîne la seconde.

Pour tout $n$, soit $C^\bullet (\F_n)$ la résolution de Godement de
$\F_n$ associée aux points de $|X|$ (si $x\in|X|$ il y a un morphisme
de sites $i_x:\spec(\mathcal{K}(x)^s)_{\et}\drt \spec
(\mathcal{K}(x))_{\et} \drt X_{\et}$ et $C^0(\mathcal{G})=\prod_x
{i_x}_* i_x^*\mathcal{G}$ ).

 En utilisant le lemme
\ref{flasque} on vérifie que c'est une résolution $\GG_c(U,-)$
acyclique de $\F$. Cette résolution est scindée sur chaque
fibre. Etant donné que $\forall n\;\; \F_n$ est localement constant, 
 pour tout $n\;\;C^\bullet(\F_n)$ est $\La/\mm^n$ plat et 
$$\forall m\geq
 n\;\;\;\;
C^\bullet(\F_m)\otimes \La/\mm^n\iso C^\bullet(\F_n)
$$

Si $m\geq n$ considérons la résolution périodique $D^\bullet$ de
$\La/\mm^n$ dans la catégorie des $\La/\mm^m$-modules (et
faisceaux) 
\begin{diagram}
  \cdots & \rTo  & \La/\mm^m & \rTo^{\varpi_\La^n} & \La/\mm^m &
\rTo^{\varpi_\La^{m-n}} & \La/\mm^m  & \rTo^{\varpi_\La^n} &\La/\mm^m  & \rTo 0 
\end{diagram}

Pour tout $p$ appliquons la formule de projection (\citeg{Berk1} 5.3.9) :
\begin{eqnarray*}
  C^p(\F_m)\overset{\mathbb{L}}{\otimes}_{\La/\mm^m} D^\bullet &= &
C^p(\F_m)\otimes \La/\mm^n \simeq C^p(\F_n) \\
\limpl
R^\bullet\GG_c(C^p(\F_m)\overset{\mathbb{L}}{\otimes}_{\La/\mm^m}
D^\bullet )&\simeq & \GG_c(U,C^p(\F_n)) 
\end{eqnarray*}
et 
$$
R^\bullet\GG_c(U,C^p(\F_m)) \overset{\mathbb{L}}{\otimes}_{\La/\mm^m}
D^\bullet
=\GG_c(U,C^p(\F_m))\overset{\mathbb{L}}{\otimes} D^\bullet
$$
\begin{diagram}
 =  \cdots & \rTo  & \GG_c(U,C^p(\F_m)) & \rTo^{\varpi_\La^n} &\GG_c(U,C^p(\F_m))  &
\rTo^{\varpi_\La^{m-n}} &\GG_c(U,C^p(\F_m))    &
\rTo^{\varpi_\La^n} & \GG_c(U,C^p(\F_m))  & \rTo 0 
\end{diagram}

On en déduit donc que 
$$\GG_c(U,C^p(\F_m))/\pi_\La^n\GG_c(U,C^p(\F_m)) \iso \GG_c
(U,C^p(\F_n))
$$
et que $\GG_c(U,C^p(\F_m))$ est un $\La/\mm^m$-module plat (par
récurrence sur $m$).
\\

$X$ étant quasi-algébrique et $\overline{U}$ compact, $\dim
(U)<+\infty$ et donc pour $d>2\dim(U)+cd_\ell(k)$ $\;\; \tau_{\leq d}
\GG_c(U,C^\bullet(\F_n))$ calcule $H^\bullet_c(U,\F_n)$ qui est un
$\La/\mm^n$-module de type fini. On peut alors appliquer le lemme
12-14 de \citeg{FK} pour conclure.
\end{proof}
\vspace{3mm}
\begin{defi}
  Soit $L$ le corps des fractions de $\La$. 
 Un faisceau $L$-adique est un
  faisceau $\La$-adique modulo $\mm$ torsion. 
  On notera pour $\F\in \Lpf$  $\;\;\F\otimes L$ le
  faisceau $ L$-adique associé et on posera $H^p_c(X,\F\otimes
  L)=H^p_c(X,\F)\otimes L$.

Un faisceau $L$-adique localement constant est un faisceau $L$-adique
de la forme $\F\otimes L$ avec $\F$ localement constant.
\end{defi}
\vspace{2mm}
\begin{coro}\label{Cont}
  Soit $X$ quasi-algébrique muni d'une action continue (\citeg{Berk2},
  6)
  d'un groupe $G$ tel que $G$ possède un pro-p sous-groupe
  ouvert. Soit $\F$ un faisceau $L$-adique  localement
  constant muni d'une action de $G$ compatible à celle sur $X$. Alors,
$H^p_c(X,\F)$ est un $G$-module lisse. 
\end{coro}
\begin{proof}
 Si $\F=\G\otimes L$, $\G$ localement constant, 
 $$
H^p_c(X,\F)=\underset{U\text{ distingué}}{\limi}\; \left
    [ \underset{n}{\limp} H^p_c(U,\G) \right ]\otimes L $$ 

Soit donc $s\in H^p_c(X,\F)$. Il existe $U$ distingué tel que $$s\in
\im \left ( H^p_c(U,\G)\otimes L\drt H^p_c(X,\G)\otimes L\right
)$$

Si $U=V_1\setminus V_2$ où $V_1$ et $V_2$ sont des domaines compacts,
par continuité de l'action $\exists K\subset G$ un sous-groupe ouvert
tel que $\forall g\in K\;\; g.V_1=V_1,g.V_2=V_2$ et donc $g.U=U$. 
D'après \citeg{Berk2} 7.7,  $\forall n\;\; H^p_c(U,\G_n)$ est un
$K$-module  lisse. Etant un $\La/\mm^n$-module de type fini c'est un
$K$-module  discret. $K$ agit donc continuement sur le $\La$-module de type
fini $\underset{n}{\limp}\; H^p_c(U,\F_n)$. Si $K'\subset K$
 est un sous-pro-p-groupe ouvert, un sous-groupe ouvert $K''$ de $K'$
 s'envoit sur un pro-$\ell$ groupe et agit donc trivialement. Donc,
 $\forall g\in K''\;\; g.s=s$.
\end{proof}

On montre de même :
\begin{coro}
  Sous les mêmes hypothèses $\Gal(\kb|k)$ agit continûment sur
  $H^p_c(X\otimes \widehat{\kb},\F)$ au sens où  $H^p_c(X\otimes
  \widehat{\kb},\F)$ est une union croissante de $\Gal(\kb|k)$-modules
  de dimension finie continus.
\end{coro}

%%% Local Variables: 
%%% mode: latex
%%% TeX-master: "thesef"
%%% End: 

\backmatter

\bibliographystyle{smfplain}
\bibliography{biblio}
\end{document}